\input epsf
\input amssym
\font\goth=eufm10

\def\g{{\hbox{\goth g}}}
\def\d{{\hbox{\goth d}}}
\def\z{{\hbox{\goth z}}}
\def\a{{\hbox{\goth a}}}
\def\b{{\hbox{\goth b}}}
\def\t{{\hbox{\goth t}}}

\def\s{{\hbox{\goth s}}}
\def\m{{\hbox{\goth m}}}

\def\h{{\hbox{\goth h}}}
\def\sl{{\hbox{\goth sl}}}
\def\su{{\hbox{\goth su}}}

\bf \centerline{KAC--MOODY AND VIRASORO ALGEBRAS}
\vskip .1in
\centerline{Antony Wassermann, Michaelmas 1998.}
\vskip .3in
\rm This course develops the representation theory of affine
Kac--Moody algebras and the Virasoro algebra. These
infinite--dimensional Lie algebras play an important r\^ole in string
theory and conformal field theory; they are the Lie algebras of the
loop groups and diffeomorphism group of the circle. We adopt a
unitary viewpoint and use supersymmetry as the main technique, as
suggested by the supersymmetric coset constructions of
Goddard--Kent--Olive and Kazama--Suzuki. Even in the case of
finite--dimensional simple Lie algebras, this approach is fruitful.
Not only does it give a streamlined route to the classical 
Weyl character formula (taken from unpublished notes of Peter Goddard), 
but it also has a natural geometric interpretation 
in terms of Dirac operators and index theory. 
\vskip .1in
\bf \centerline{CONTENTS}
\vskip .05in
\noindent \rm Chapter I. Clifford algebras, fermions and the spin group.
\vskip .05in
\noindent 1. Tensor, symmetric and exterior algebras.

\noindent 2. Inner products and tensors.

\noindent 3. The double commutant theorem.

\noindent 4. Fermions and Clifford algebras.

\noindent 5. Quantisation and the spin group.

\noindent 6. Matrix groups and their Lie algebras.

\noindent 7. The odd--dimensional case.

\noindent 8. The spin representations of Spin(V). 
\vskip .05in
\noindent Chapter II. Compact matrix groups and simple Lie algebras.
\vskip .05in
\noindent Part I. Elementary structure theory.
\vskip .05in
\noindent 1. Compact Lie algebras.

\noindent 2. Examples of simple compact Lie algebras.

\noindent 3. Maximal tori.

\noindent 4. Representations of SU(2) and sl(2).

\noindent 5. The root system.

\noindent 6. The Weyl group as a reflection group.

\noindent 7. Geometric apprach to Weyl chambers and simple roots.

\noindent 8. Weyl's uniqueness theorem.

\noindent 9. Classification of compact simple Lie algebras.
 
\vskip .05in
\noindent Part~II. Representation theory.
\vskip .05in
\noindent 10. Root and weight lattices.

\noindent 11. Poincar\'e--Birkhoff--Witt theorem.

\noindent 12. Highest weight vectors.

\noindent 13. Eigenvalues of the Casimir operator.

\noindent 14. Lie algebraic construction of irreducible representations.

\noindent 15. Projective representations and covering groups.

\noindent 16. The Dirac operator and supersymmetry relations.

\noindent 17. The square of the Dirac operator.

\noindent 18. Weyl's character and denominator formulas.

\noindent 19. Remarks on connections and Dirac operators.

\noindent 20. Remarks on Dirac induction and Bott's principle.
\vskip .05in
\noindent Chapter~III. Representations of affine Kac--Moody algebras.
\vskip .05in
\noindent 1. Loop algebras and the Witt algebra.

\noindent 2. Positive energy representations and Kac--Moody algebras.

\noindent 3. Complete reducibility.

\noindent 4. Classification of positive energy representations.

\noindent 5. Sugawara's formula for $L_0$.

\noindent 6. Sugawara's construction of the Virasoro algebra.

\noindent 7. Weights, roots and the quantum Casimir operator.

\noindent 8. The affine Weyl group.

\noindent 9. Construction of irreducible representations.

\noindent 10. Garland's ``no--ghost'' theorem on unitarity.

\noindent 11. The character of a positive energy representation.

\noindent 12. Bosons and fermions on the circle.

\noindent 13. The Kazama--Suzuki supercharge operator.

\noindent 14. The square of the supercharge operator.

\noindent 15. Kac's character and denominator formulas.
\vskip .05in
\noindent Chapter~IV. Representations of the Virasoro algebra.
\vskip .05in
\noindent 1. Positive energy representations of the Virasoro algebra.

\noindent 2. The Goddard--Kent--Olive construction.

\noindent 3. Character of the multiplicity space.

\noindent 4. The Kac determinant formula.

\noindent 5. The Friedan--Qiu--Shenker unitarity criterion for h.

\noindent 6. The multiplicity one theorem.

\noindent 7. The Feigin--Fuchs character formula 
for the discrete series.

\noindent 8. The Friedan--Qiu--Shenker unitarity criterion for c.

\vfill\eject

\bf \centerline{CHAPTER~I. CLIFFORD ALGEBRAS, FERMIONS AND THE SPIN
GROUP.}
\vskip .1in

\it We develop here the theory of fermions in finite--dimensions. This
provides the first example of the principle of quantisation: if an
algebraic object acts irreducibly on a Hilbert space and a group of
automorphisms of the object preserves the equivalence calss of the
representation, then the automorphism group is implemented by a unique
projective representation on the Hilbert space. It is convenient to
develop the linear algebra for bosons and fermions in parallel (it is
even possible to develop a simultaneous super--theory).  Roughly
speaking bosons are operators that satisfy the commutation 
relation $aa^*-a^*a=I$ and fermions are operators that satisfy the
anticommutation relation $aa^*+a^*a=I$. In this chapter we use
fermionic quantisation to construct the spin group, a double covering of the
special orthogonal group. In Chapter~III, we carry out an analogous
treatment for bosons. Bosonic quantisation leads to a construction of
the metaplectic group, a double covering of the symplectic group; as
applications we will prove the modular transformation properties of
theta functions and quadratic reciprocity. Fermionic and bosonic
quantisation is also very important in infinite dimensions. The
infinite--dimensional theory proceeds very much as in the
finite--dimensional case, except that not all automorphisms can be
quantised (they must satisfy a `Hilbert--Schmidt' quantisation
criterion). In addition there is a remarkable equivalence between the
bosonic and fermionic theories in infinite dimensions which has been
used to explain the KdV and KP hierarchies in soliton theory. 

\vskip .1in

\noindent \bf 1. TENSOR, SYMMETRIC AND EXTERIOR ALGEBRAS. \rm 
\vskip .1in
\noindent \bf Tensor products. \rm If $V$ and $W$ are
finite--dimensional vector spaces over ${\Bbb R}$ or ${\Bbb C}$, we
defining their tensor product $V\otimes W$ by taking bases $(v_i)$ and
$(w_j)$ in $V$ and $W$ and then decreeing $V\otimes W$ to be the
vector spec with basis $v_i\otimes w_i$. In general we set
$(\sum a_i v_i)\otimes (\sum b_j w_j) =\sum a_i b_j \, v_i\otimes
w_j$, so that $v\otimes w$ is defined for any $v\in V$, $w\in W$. This
definition is up to isomorphism independent of the choice of basis.
Clearly ${\rm dim}(V\otimes W) ={\rm dim}(V){\rm dim}(W)$. Iterating
we get a similar definition of a $k$--fold tensor product $V_\otimes
\cdots \otimes V_k$. By definition these is a natural one--one correspondence
between the vector space of multilinear maps $V_1\times \cdots \times
V_k\rightarrow U$ and ${\rm Hom}(V_1\otimes \cdots \otimes V_k,U)$;
this could equally well be used as the universal property
characterising the tensor product.

The tensor product has various obvious functorial properties. Thus
for example $V_1\otimes V_2\cong V_2\otimes V_1$,
 $(V_1\otimes V_2)^*=V_1^*\otimes V_2^*$, $V_2\otimes V_2^*
\cong {\rm Hom}(V_2,V_1)$, $V_1\otimes V_2 \cong {\rm
Hom}(V_2^*,V_1)$, ${\rm Hom}(V_1\otimes V_2,V_3) \cong {\rm
Hom}(V_1,V_2^*\otimes V_3)$. Moreover if $f_i:U_i\rightarrow V_i$ are
linear maps, then we have $f_1\otimes f_2:U_1\otimes U_2\rightarrow
V_1\otimes V_2$ sending $u_1\otimes u_2$ to $f_1(u_1)\otimes
f_2(u_2)$. 
\vskip .1in
\noindent \bf The tensor algebra. \rm Let $T^n(V) =V^{\otimes n}
=V\otimes \cdots \otimes V$ ($n$ times) and $T(V)=\bigoplus
V^{\otimes n} =\bigoplus T^n(V)$, the tensor algebra. Multiplication
$T^a(V)\rightarrow T^b(V)\rightarrow T^{a+b}(V)$ is defined by
concatenation, so that $(v_1\otimes \cdots \otimes v_a)\times 
(w_1\otimes \cdots \otimes w_b) = v_1\otimes \cdots \otimes v_a
\otimes w_1 \otimes \cdots \otimes w_b$. This makes $T(V)$ into a
non--commutative associative algebra.
\vskip .1in
\noindent \bf Action of $S_n$ on $V^{\otimes n}$. \rm The symmetric
group $S_n$ acts on $V^{\otimes n}$ by permuting the tensor factors. 
Thus $\sigma(v_1\otimes \cdots \otimes v_n) = v_{\sigma 1} \otimes
\cdots \otimes v_{\sigma n}$ for $\sigma\in S_n$. Define $\varepsilon
:S_n\rightarrow \{\pm 1\}$ to be the sign homomorphism, assigning $+1$
to an even permutation and $-1$ to an odd permutation. Let 
$$S\omega ={1\over n!} \sum_{\sigma \in S_n} \sigma\omega, \qquad
A\omega = ={1\over n!} \sum_{\sigma \in S_n} \varepsilon(\sigma)\sigma\omega$$
be the symmetrising and antisymmetrising operators on $V^{\otimes n}$. 
\vskip .1in
\noindent Symmetric and exterior algebras. \rm Let $S^n(V) =\{\omega\in
V^{\otimes n}: \sigma \omega =\omega \forall \sigma\in
S_n\}=SV^{\otimes n}$ and
$\Lambda^k(V) =\{\omega\in
V^{\otimes n}: \sigma \omega =\varepsilon(\sigma)\omega \forall
\sigma\in S_n\}=AV^{\otimes n}$. $S(V)=\bigoplus S^n(V)$ and $\Lambda(V)=\bigoplus
\Lambda^n(V)$ are called the symmetric and exterior algebras. Their
multiplication is defined on homogenous elements by $\omega_1 \cdot
\omega_2= S(\omega_1\otimes \omega_2)$ or $\omega_1\wedge \omega_2 =
A(\omega_1\otimes \omega_2)$ and extended bilinearly to the whole of
$S(V)$ or $T(V)$. It is easy to check that $a\cdot (b\cdot c) =
S(a\otimes b\otimes c) = (a\cdot b)\cdot c$ and that
$a\wedge (b\wedge c) = A(a\otimes b\otimes c) = (a\wedge b)\wedge c$,
so that $S(V)$ and $\Lambda(V)$ become associative algebras.
\vskip .1in
\noindent \bf Lemma. \it $S(V)$ is a commutative ring and $\Lambda(V)$
is a graded commutative ring.
\vskip .05in
\def\d{{\partial}}
\noindent \bf Proof. \rm The first result follows straight from the
definitions and is just part of the fact that $S(V)$ coincides with
the algebra of polynomial functions on $V^*$ (see below). The algebra
$\Lambda (V)$ is ${\Bbb Z}_2$--graded into even or odd elements,
according to degree of homogeneous elements. We set $\partial a=0$ or
$1$ according as $a$ is even or odd. Graded commutativity is just the
statement that 
$a\wedge b =(-1) \d a\d b b\wedge a$, which is immediate from the
definitions. 
\vskip .1in
\noindent \bf Concrete realisations of $S(V)$ and $\Lambda(V)$. \rm We map
$S(V)$ into polynomial functions on $V^*$. Note that
$S^k(V)^*=S^k(V^*)$. We need
\vskip .1in
\noindent \bf Polarisation Lemma. \it The tensors $v^{\otimes m}$ with
$v\in V$ span $S^mV$
\vskip .05in
\noindent \bf Proof \rm Note that if $X$ is a subspace and
$f(\lambda_1,\dots,\lambda_m)$ is a polynomial function of
$\lambda_1,\dots,\lambda_m$ with values in $X$, then 
${\partial^{|\alpha|}\over \partial \lambda^\alpha}f$ also lies in $X$
for any multinomial $\alpha$, since $X$ is finite--dimensional, so
closed. Take $v_1,\dots,v_m\in W$ and consider $f(\lambda)=(\sum
\lambda_i v_i)^{\otimes m}$. Up to a constant non--zero factor
${\partial^m f\over \partial \lambda_1\cdots\partial \lambda_m} $ is
the symmetrisation of $v_1\otimes\cdots \otimes v_m$. This shows that
the symmetrisation of any elementary tensor (and hence any tensor)
lies in the subspace $X$ of $S^m V\subset V^{\otimes m}$ spanned by
the tensors $v^{\otimes m}$. 
\vskip .1in
In particular, $S^k V^*$ is spanned by tensors $x^{\otimes n}$ with
$x\in V^*$. Hence the map $f\mapsto f(x^{\otimes m})=f(x)$ defines an
injection of $S^k(V)$ into the polynomials of degree $k$ on $V^*$. The
map is clearly surjective, so we may identify $f\in S^k(V)$ with the
polynomial $f(x)$. It is easy to see that under this identification
$f\cdot g (x) = f(x)g(x)$, so that as a commutative algebra $S(V)$ can
be identified with the algebra of polynomial functions on $V^*$.

Note that if $v_1,\dots,v_n$ is a basis of $V$, then a basis of
$\Lambda^k(V)$ is given by $v_{i_1}\wedge v_{i_2} \wedge \cdots \wedge
v_{i_k}$. Thus ${\rm dim}\Lambda^k(V) ={n\choose k}$ and ${\rm
dim}\Lambda(V) =2^n$. In particular $\Lambda^m(V)=0$ for $m>n$ and
$\Lambda^n(V)$ is one--dimensional. We can also identify $\Lambda^kV$
with alternating multilinear functionals on $V^*\times \cdots \times
V^*$. 
If $f$ and $g$ are homogeneous of degree $a$ and $b$
respectively, then exterior multiplication is given by the formula 
$$f\wedge g(x_1,\dots, x_{a+b}) = {1\over (a+b)!}\sum_{\sigma\in
S_{a+b}}\varepsilon(\sigma) f_1(x_{\sigma 1}, 
\dots , x_{\sigma a}) g(x_{\sigma(a+1)},\dots , x_{\sigma(a+b)}).$$
(Actually the sum can be reduced to a sum over the coset space
$S_{a+b}/S_a\times S_b$ since $\sigma(f\otimes g) =\varepsilon(\sigma)
f\otimes g$ for $\sigma \in S_a\times S_b$.)
\vskip .1in
Finally note that $V\rightarrow S(V)$ and $V\rightarrow \Lambda(V)$
are {\it functors} from the additive category of vector spaces to the
multiplicative tensor category of vector spaces. This will not be
important for us, although it is the key to quantisation in quantum
field theory. As Nelson said, first quantisation is a mystery while
second quantisation is a functor. This functoriality appears in the
isomorphism $S(V\oplus W)=S(V)\otimes S(W)$ and $\Lambda(V\oplus W) =
\Lambda (V)\otimes \Lambda(W)$ between (graded) commutative algebras. 
We need $(a\otimes b)(c\otimes d) =(-1)^{\d b \d c} ac\otimes bd$ to
define the tensor product of graded algebras. The basic rule in
discussing graded objects is that if we move a symbol of degree $\d_1$
past a symbol of degree $\d_2$, then a sign $(-1)^{\d_1\d_2}$ must be
introduced. The functor $S$ corresponds to bosons which
satisfy the canonical commutation relations while the functor
$\Lambda$ corresponds to fermions which satsify the canonical
anticommutation relations. The basic idea of supersuymmetry is that
the bosonic and fermionic theory can be developed in parallel at each
stage, so that any concept introduced in one theory has its natural
counterpart in the other. 
\vskip .1in

\noindent \bf 2. INNER PRODUCTS AND TENSORS. \rm If $U$ and $V$
are real or complex inner product spaces, we can define an inner
product on $U\otimes V$ by taking any positive multiple of the
inner product $(u_1\otimes v_1,u_2\otimes v_2)=(u_1,u_2) (v_1,v_2)$.
In particular we define the inner product on $T^k(V)=V^{\otimes k}$ by
$(a_1\otimes \cdots \otimes a_k,b_1\otimes \cdots \otimes a_k) = k!
\prod (a_i,b_i)$. (The factor of $k!$ is essential here to guarantee 
$(\exp(a),\exp(b)) =\exp{(a,b)}$ for $a,b\in V$.) This inner product
extends to $T(V)$ by declaring the $T^k(V)$ to be mutually orthogonal.
Note that, since $S(V),\Lambda(V)\subset T(V)$, there are naturally
induced inner products on $S(V)$ and $T(V)$. The definition
immediately give the following
explicit formulas in the functional realisations above.
\vskip .1in
\noindent \bf Lemma. \it (a) In $\Lambda(V)$, we have $(a_1\wedge
\cdots \wedge a_m, b_1\wedge 
\cdots \wedge b_n) =\delta_{nm} \det (a_i,b_j)$.

\noindent (b) In $S(V)$, we have $(x^m, y^n)= \delta_{mn} n! (x,y)^n$.
\vskip .1in
\rm We will see that regarded as polynomial functions on $V^*={\Bbb C}^n$, the
inner product in $S(V)$ agrees with the inner product
$(f,g) = \pi^{-n} \int_{{\Bbb C}^n} f(z)\overline{g(z)} \,
e^{-|z|^2}$, so that $S(V)$ can be identified with so--called
holmorphic Fock space (see below). Part (a) of the lemma shows that if
$(e_i)$ is an orthonormal basis of $V$, then $e_{i_1}\wedge
e_{i_2}\wedge \cdots \wedge e_{i_k}$ ($i_1< \cdots <i_k$) is an
orthonormal basis for $\Lambda^k(V)$. 

Now both on $\Lambda(V)$ and $S(V)$ we have the operation of
multiplication by $v\in V$. We now work out their adjoints.
\vskip .1in
\bf \noindent Theorem (adjoint derivations). \it 
(a) The adjoint $e(v)^*$ of $e(v)$ is the graded derivation
$d_v( v_1\wedge \cdots \wedge v_k )=\sum (-1)^{i+1} (v_i,v) v_1
\wedge \cdots \wedge v_{i-1} \wedge v_{i+1} \wedge \cdots \wedge v_k$
with $d_v(1)=0$.

\noindent (b) The adjoint of multiplication by $v$ is the derivation 
$\partial_v$ given by $\partial_v(x_1\cdots x_n) =\sum (x_i,v)
\prod_{j\ne i} x_j$ for $x_j\in V$ 
with $\partial_v (1)=0$.

\vskip .05in
\noindent \bf Proof. \rm (a) We have
$$\eqalign{
(e(w_1)^*& v_1\wedge \cdots \wedge v_{n+1}, w_2\wedge \cdots \wedge
w_{n+1})= (v_1\wedge \cdots \wedge v_{n+1} , w_1\wedge \cdots \wedge
w_{n+1}) \cr
 & = \det(v_i,w_j) \cr
& =\sum (-1)^{i+1} (v_i,w_1) (v_1\wedge \cdots \wedge v_{i-1}\wedge
v_{i+1}\wedge \cdots \wedge v_{n+1}, w_1\wedge \cdots \wedge
w_{n+1}). \cr}$$
expanding the determinant by the first column. This proves the formula
for $e(v)^*$. This is usually called ``contraction'' with $v$
or ``interior multiplication''. It is routine to check from the
definition of $d_v$ that, if $\omega_1$ and $\omega_2$ are
homogeneous, then $d_v(\omega_1\wedge \omega_2) =d_v(\omega_1)\wedge
\omega_2 +(-1)^{\d \omega_1} \omega_1 \wedge d_v\omega_2$. This means
that $d_v$ is a graded derivation, the signs being compatible with our
previous convention since $d_v$ is odd. Note that $d_v$ is uniquely
determined once we declare that it is a graded derivation, $d_v(1)=0$
and $d_v w=(w,v)$ for $w\in V$.

\noindent (b) This can be checked directly using the inner product as
in (a). When $V$ is a complex inner product space, it is also obvious
in the functional realisation in terms of 
polynomials on ${\Bbb C}^n$ with the above inner product, for there
clearly $z_i$ has adjoint $\partial/\partial z_i$. 

\vskip .1in
\bf \noindent Theorem (real and complex wave representation). \it (a) Let $V$
be an inner product space. Then if $a,b\in V$, the operators
$e(a),e(b)$ on $\Lambda(V)$ satisfy the canonical anticommutation relations
$e(a)e(b)+e(b)e(a)=0$, 
$e(a)^*e(b)^*+e(b)^* e(a)^*=0$ and $e(a) e(b)^*+ e(b)^*e(a)=(a,b)$. 

\noindent (b) Let $V$ be an inner product space. Then, if
$v,w\in V$, the operators
$z$ and $\partial_w$ on $S(V)$ satisfy the canonical commutation
relations $zw-wz=0$, $\partial_z\partial_w -\partial_w\partial_z=0$
and $\partial_w z - z\partial_w=(z,w)$.
\vskip .05in
\noindent \bf Proof. \rm (a) Clearly $e(a)$ and $e(b)$ anticommute, so
taking adjoints so too do $e(a)^*$ and $e(b)^*$. Now
$$(e(a)e(b)^*+e(b)^* e(a))\omega= a\wedge e(b)^*\omega + e(b)^*(a\wedge
\omega) =a\wedge e(b)^*\omega + (a,b)\omega - a\wedge (eb)^*\omega
=(a,b)\omega.$$

\noindent (b) Clearly $z$ and $w$ commute. hence so do their adjoints
$\partial_z$ and $\partial_w$. Now
$$(\partial_w z - z\partial_w ) p = z\partial_w p +(z,w)p -z\partial_w
p =(z,w)p.$$
This proves the last commutation relation.

\vskip .1in

\bf \noindent Theorem (irreducibility of wave
representation). \it (a) If $V$ is 
an complex inner product space, the
operators $e(v)$ and $e(v)^*$ act irreducibly on $\Lambda(V)$.

\noindent (b) If $V$ is a complex inner product space, the operators
$v$ and $\partial_v$ act irreducibly on $S(V)$.
\vskip .05in
\noindent \bf Proof. \rm (a) Let $U\ne (0)$ be an invaraiant subspace
and take $\omega\ne 0$ in $U$. Then $\omega =\sum a_I v_{i_1}\wedge \cdots
\wedge v_{i_k}$ with respect to some orthonormal basis $(v_i)$. Pick a
non--zero term of maximal degree, $a_I v_{i_1}\wedge \cdots
\wedge v_{i_k}$. Then $e(v_{i_k})^* \cdots e(v_{i_1})^*\omega = a_I $,
so that $1\in U$. Since all of $\Lambda(V)$ can be obtained by
applying $e(v)$'s to $1$, we see that $U=\Lambda(V)$.

\noindent (b) We have to show that the operators $z_i$ and
$\partial/\partial z_j$ act irreducibly on the polynomial algebra
${\Bbb C}[z_1,\dots,z_n]$. Let $U$ be an invariant subspace and take
$p(z)\ne 0$ in $U$. Then $p(z)=\sum a_\alpha z^\alpha$. Pick a
non--zero term of 
maximal degree $a_\alpha z^\alpha$. Then $\partial_\alpha p(z) =
\alpha! $, so that $1\in U$. Since all polynomials can be obtained by
multiplying $1$ by $z_i$'s, we see that $U={\Bbb C}[z_1,\dots,z_n]$.

\def\H{{\cal H}}

\vskip .1in
\noindent \bf 3. THE DOUBLE COMMUTANT THEOREM. \rm  Let $V$ be a
finite--dimensional inner product space over ${\Bbb C}$ and let
$A\subseteq {\rm End}\, V$ be a *--subalgebra of ${\rm End}\,V$. This
means that $I\in A$ and $A$ is a linear subspace closed under
multiplication and the adjoint operation $T\mapsto T^*$. For any
subset ${\cal S}\subseteq {\rm End}\,V$, we define the {\it commutant}
of ${\cal S}$ by
$${\cal S}^\prime={\rm End}_{\cal S}(V)=\{T\in {\rm End}\, V:
Tx=xT\quad\hbox{for all $x\in {\cal S}$}\}.$$
\vskip .1in
\noindent \bf Schur's Lemma. \it (i) $A$ acts irreducibly
on $V$ (i.e.~has no invariant subspaces) iff $A^\prime={\Bbb C}$. 

\noindent (ii) If $A$ acts on two irreducible subspaces
$V_i$ and $T\in {\rm Hom}_A(V_1,V_2)$ (i.e.~commutes with $A$), then
$T=0$ or is an isomorphism. 
\vskip .05in
\noindent \bf Proof. \rm (Spectral Theorem.) (i) Say $A$ does not
act irreducibly and $U\subset V$ be a proper subspace invariant under $A$
(i.e.~$U$ is an $A$--submodule). Then, if $P$ is the orthogonal
projection onto $U$, we have $P\in A^\prime$. So $A^\prime \ne {\Bbb
C}$. 

Conversely if $T\in A^\prime$, then, since $A^\prime$ is a *--algebra,
both ${\rm Re}\, T = T+T^*/2$ and ${\rm Im}\, T= T-T^*/2i$ lie in
$A^\prime$. By the spectral theorem for self--adjoint matrices, so
does any projection onto an eigenspace (i.e.~a spectral projection).
So if $T\notin {\Bbb C}$, we have produced a projection $P\in
A^\prime$ with $P\ne 0, I$. The corresponding subspace is invariant.

\noindent (ii) If $v_1$ and $V_2$ are irreducible and $T$ is an
intertwiner, then so is $T^*$ (simply take adjoints of the
intertwining relation and replace $a$ by $a^*$). But then $TT^*$ and
$T^*T$ are also intertwiners, i.e.~$T^*T\in \pi_1(A)^\prime$ and
$TT^*\in \pi_2(A)^\prime$. They must be scalars by (i), so either both
zero or both the same multiple of the identity.

\vskip .1in
\noindent \bf Double commutant theorem. \it If $A\subset {\rm End}(V)$
is a *--algebra, then $A^{\prime\prime}=A$. 
\vskip .05in
\noindent \bf Proof. \rm (1) \it If $U$ is a subspace of $V$ invariant
under $A$, then so is $U^\perp$. In particular $V$ is a direct sum of
irreducible $A$--submodules.
\vskip .05in
\noindent \bf Proof. \rm Say $\xi\in U^\perp$ and $a\in A$. Let
$\eta\in U$. Then $\langle a\xi,\eta\rangle =\langle \xi,a^*\eta\rangle =0$ since $a^*\eta\in U$
and $\xi\perp U$. So $a\xi\perp U$, i.e.~$a\xi\in U^\perp$. So
$V=U\oplus U^\perp$ with $U$ and $U^\perp$ $A$--modules. We continue
this game if $U$ or $U^\perp$ fail to be irreducible. 

\vskip .05in
\noindent (2)\it If $S\in A^{\prime\prime}$ and $v\in V$, there is a
$T\in A$ such that $Tv=Sv$.
\vskip .05in
\noindent \bf Proof. \rm In fact let $W=Av\subseteq V$. This is an
$A$--submodule of $V$. The orthogonal projection onto $W$ gives a
projection $E\in {\rm End}\,V$ ($E^2=E=E^*$) which commutes with $A$,
from (1). So $E\in A^\prime$. But $S\in A^{\prime\prime}$, so $SE=ES$.
This means that $S$ leaves $W$ and $W^\prime$ invariant. (Note that
$I-E$ is the orthogonal projection onto $W^\perp$.) But $v\in W$. So
$Sv\in W=Av$. So $Sv=Tv$ for some $T\in A$.
\vskip .05in
\noindent (3) Let $V^\prime=V\oplus\cdots \oplus V$ ($m$ times) with
$A$ acting diagonally, $a(\xi_1,\dots,\xi_n)=(a\xi_1,\dots,a\xi_m)$.
This means we can identify $A$ with a *--subalgebra of ${\rm
End}\,V^\prime$ (for the initiated, $V\oplus\cdots\oplus V= V\otimes
{\Bbb C}^m$). It's easy to check that $\pi(A)^\prime=A^\prime\otimes
M_m({\Bbb C})=M_m(A^\prime)$, if we write elements of ${\rm End}\,
V^\prime$ as $m\times m$ matrices with entries in ${\rm End}\, V$. We
go on to check that
$$\pi(A)^{\prime\prime}=(\pi(A)^\prime)^\prime=\pi(A^\prime)
=\left\{\pmatrix{x&&&&\cr
                  &x&&&\cr
                  &&\cdot &&\cr
                   &&&\cdot &\cr
                    &&&&x\cr}\right\},$$ 
where here (and above) $\pi$ denotes the embedding ${\rm End}\,
V\rightarrow {\rm End}\,V^\prime$ taking operators to diagonal
operators. Take $m={\rm dim}\,V$. Set $v=\pmatrix{e_1\cr
\cdot\cr\cdot\cr\cdot\cr e_m\cr}$ where $e_1,\dots,e_m$ is a basis of
$V$. By step (2), we have $\pi(A)v=\pi(A)^{\prime\prime}v$. But
$\pi(A)^{\prime\prime}=\pi(A^{\prime\prime})$ from the above. So given
$S\in A^{\prime\prime}$ we can find $T\in A$ such that
$\pi(S)v=\pi(T)v$. Hence 
$$\pmatrix{S&&&&\cr
                  &S&&&\cr
                  &&\cdot &&\cr
                   &&&\cdot &\cr
                    &&&&S\cr}
\pmatrix{e_1\cr
\cdot\cr\cdot\cr\cdot\cr e_m\cr}
=\pmatrix{T&&&&\cr
                  &T&&&\cr
                  &&\cdot &&\cr
                   &&&\cdot &\cr
                    &&&&T\cr}
\pmatrix{e_1\cr
\cdot\cr\cdot\cr\cdot\cr e_m\cr}.$$
Thus $Se_i=Te_i$ for all $i$ and hence $S=T$.
\vskip .1in
\noindent \bf Corollary~1. \it A *--algebra $A$ acts irreducibly iff
$A={\rm End}(V)$.
\vskip .1in
\noindent \bf Corollary~2. \it All *--representations of
${\rm End}(V)$ are on direct sums of copies of $V$.
\vskip .05in
\noindent \bf Proof~1. \rm Since all representations are sums of
irreducibles, it suffices to show that $V$ is the only irreducible
representation of ${\rm End}(V)$. But if $W$ is another inequivalent
irreducible, the commutant on $V\oplus W$ must be ${\rm End}(V)\oplus
{\rm End}(W)$ by Schur's lemma and the double commutant theorem. But
the image of ${\rm End}(V)$ must coincide with its double commutant, a
contradiction.
\vskip .05in
\noindent \bf Proof~2. \rm Choose matrix units in $A={\rm End}(V)$ and
let $W$ be an $A$--module. Set $W_1=e_{11}W$ and let $(w_j)$ be 
basis of $W_1$. Consider the map
$T:\oplus V\otimes W_1 \rightarrow W$, $\oplus \mu_{ij} e_i \otimes w_j
\mapsto \sum \mu_{ij} e_{i1}w_j$. $T$ is surjective since
$I=\sum e_{j1}e_{11}e_{1j}$, so that $AW_0=W$. $T$ is also injective,
for $\sum \mu_{ij} e_{i1}w_j=0$ forces $\mu_{ij}w_j=0$ for each $i$
(premultiply by $e_{1i}$) and hence $\mu_{ij}\equiv 0$. By
construction it commutes with the actions of $A$. It is even unitary
if $(w_j)$ is chosen orthonormal. So $W$ is a direct sum of copies of
$V$.

\vskip .1in
\noindent \bf Remarks. \rm In the exercises this theorem is used to
determine the structure of finite--dimensional *--subalgebras of
${\rm End}(V)$, $V$ a complex inner product space. There is also an
infinite--dimensional Hilbert space version of the 
double commutatant theorem due to John von Neumann which is the
starting point of the modern theory of operator algebras.
\vskip .1in
\noindent \bf Corollary~3 (Schur--Weyl duality). \it if $A$ is the
*--algebra of linear combinations of 
$g^{\otimes m}$'s as $g$ ranges over $GL(V)$ and $B$ is the
*--algebra of linear combinations of the $\sigma$'s as $\sigma$ ranges
over $S_m$, we have $A=B^\prime$ and $B=A^\prime$, so that $A$ and $B$
are each other's commutants.
\vskip .05in
\rm \noindent Proof. \rm Since $A$ and $B$ are *--algebras, by the
double commutant theorem, 
$A=A^{\prime\prime}$ and $B=B^{\prime\prime}$. So to prove $A^\prime
=B$, it is equivalent to check that $B^\prime =A$. The algebra $A$ is
a finite--dimensional subspace, so closed. Now any non--invertible
matrix is the limit of invertible matrices: for $x+\varepsilon I$ for
all $\varepsilon$ sufficiently small. So $A$
contains all tensors $w\otimes\cdots \otimes w$ even if $w$ is not
invertible. So $C$ coincides with the fixed points of $S_m$ in ${\rm
End}\, V^\otimes m$, i.e.~the commutant of $S_m$. (Note that
conjugation by $\sigma$ gives the permutation action of $S_m$ on
$W^{\otimes m}$.) Thus $B^\prime=A$ as required. 
\vskip .1in

\noindent \bf 4. FERMIONS AND CLIFFORD ALGEBRAS. \rm
\vskip .1in
\noindent Real Clifford algebras. \rm Let $V$ be a real
$2n$--dimensional inner product space. Operators $c(v)$ on a real or complex
inner product space $W$ are said to satisfy the real Clifford algebra
relations iff $v\mapsto c(v)$ is ${\Bbb R}$--linear, $c(v)^*=c(v)$ and
$c(a)c(b) + c(b) c(a)=2(a,b) I$. 
\vskip .1in
\noindent \bf Lemma~1. \it If the operators $c(v)$ satisfy the real
Clifford algebra relations, then the real *--algebra $A$ they generate is
spanned by products
$c(v_{i_1}) c(v_{i_2}) \cdots c(v_{i_k})$ with $i_1<i_2<\cdots <i_k$
and $(v_i)$ a basis of $V$. Moreover ${\rm dim}(A)\le 2^{{\rm
dim}(V)}$.
\vskip .05in
\noindent \bf Proof. \rm Clearly the algebra generated by the
$c(v_i)$'s is a *--algebnra since each $c(v)$ is self--adjoint. Thus
it suffices to prove that $A_0={\rm lin}_{\Bbb C}(c(v_{i_1})
c(v_{i_2}) \cdots c(v_{i_k}))$ is closed under multiplication by
$c(v_i)$. This, however, is obvious from the Clifford relations. Hence
$A=A_0$. Clearly ${\rm dim}(A_0) \le 2^{{\rm dim}(V)}$.
\vskip .1in
\bf \noindent Lemma~2. \it Let $c(v) =e(v) +e(v)^*$ acting on
$W=\Lambda_{\Bbb R} V$. 

\noindent (a) The $c(v)$'s satisfy the real Clifford algebra relations. 

\noindent (b) The vector $\Omega=1\in \Lambda^0(V)$ is cyclic for a
the real*--algebra $A$ generated by the $c(v)$'s, i.e.~$A\omega
=\Lambda(V)$. 

\noindent (c) The operators $c(v_{i_1}) c(v_{i_2}) \cdots c(v_{i_k})$
with $i_1<i_2<\cdots <i_k$ are linearly independent and the vector
$\Omega=1$ is separating for $A$, i.e. 
$a\Omega=0$ implies $a=0$. 

\noindent (d) $c(v_{i_1}) \cdots c(v_{i_k}) \omega = v_{i_1} \wedge
v_{i_2} \wedge \cdots \wedge v_{i_k}\wedge \omega +$ lower order terms
modulo two.

\vskip .05in
\bf Proof. \rm (a) By the canonical anticommutation relations,
$$c(a)c(b) + c(b) c(a) = (e(a) + e(a)^*)(e(b)+e(b)^*) + (e(b) +
e(b)^*) (e(a) + e(a)^*) = 2(a,b) I.$$
(b) Let $W_k={\rm lin}\{c(x_1)c(x_2) \cdots c(x_j)\Omega: j\le k\}$
for $k\ge 0$. We prove by induction that $W_k=\oplus_{j=0}^k
\Lambda^j(V)$. For $k=0$, this is trivial. For $k>0$, $c(x_1)
x_2\wedge \cdots \wedge x_k= x_1\wedge
\cdots \wedge x_k $ plus a term in $\Lambda^{k-2}(V)$. Thus, by
induction,
$x_1\wedge \cdots \wedge x_k $ lies in $c(x_1)W_{k-1} + W_{k-2}
\subset W_k$, as required.

\noindent (c) Since $2^{{\rm dim}(V)} \le {\dim}\Lambda(V) ={\rm
dim}\,A\Omega=\le {\rm dim}(A) \le 2^{\rm dim}(V) $, this is
obvious from (b) and Lemma~1. 

\noindent (d) This follows easily by induction on $k$. 
\vskip .1in
We define the real Clifford algebra ${\rm Cliff}(V)$ to be the real
*--algebra generated by the $c(v)$'s on 
$\Lambda(V)$. We show that ${\rm Cliff}(V)$ has a similar universal 
property to the group algebra ${\Bbb C}[G]$. This is defined as the
algebra of operators on $\ell^2(G)$ generated by left translations.
Any finite--dimensional unitary representation of $G$ gives rise to a
*--representation of ${\Bbb C}[G]$ and conversely, so that ${\Bbb
C}[G]$ is the universal algebra for representations of $G$. We claim
that any given Clifford algebra relations $C(v)$ on $W$, there is a
unique *--representation of ${\rm Cliff}(V)$ sending $c(v)$ to $C(v)$.
Uniqueness is clear, since the $c(v)$'s generate ${\rm Cliff}(V)$; to
prove existence, we take a basis $(v_i)$ of $V$ and send the basis
element $c(v_{i_1}) \cdots c(v_{i_k})$ of ${\rm Cliff}(V)$ to $C(v_{i_1})
\cdots C(v_{i_k})$. This is clearly a homomorphism of *--algebras. 
If $W=\Lambda(V)$, a real inner product space, we have a natural
complexification $W_{\Bbb C}=W\otimes_{\Bbb R}{\Bbb C}$. This is just
obtained by taking an orthonormal basis for $V$ and hence $\Lambda W$
and extending the scalars and inner product in the obvious way. The
algebra $A={\rm Cliff}(V)$ and its complexification ${\rm Cliff}_{\Bbb
C}(V)=A_{\Bbb C} = A\oplus i A$ acts on $W_{\Bbb C}$. $A_{\Bbb C}$ is
a complex *--algebra and $\Omega$ is again cyclic and separating for
$A_{\Bbb C}$. This means $A_{\Bbb C}$ cannot act irreducibly; for if it
did, $A_{\Bbb C}={\rm End}(W_{\Bbb C})$ and $\Omega$ is not separating
for ${\rm End}(W_{\Bbb C})$. Note tthat $A\rightarrow A\Omega$ gives
an isomorphism between ${\rm Cliff}(V)$ and $\Lambda(V)$ as linear
spaces. This allows us to speak about the degree of an element of
${\rm Cliff}(V)$. Note the following immediate consequence of Lemma~2 (d).
\vskip .1in
\noindent \bf Corollary. \it If $\omega_1\in \Lambda^a(V)$ and
$\omega_2\in \Lambda^b(V)$, then $\omega\cdot \omega_2 =
\omega_1\wedge \omega_2 +$ lower degree terms modulo two. 

\vskip .1in

\rm We now show how introducing a complex
structure on $V$ allows us to produce an irreducible representation of
the real Clifford algebra relations. By definition a complex structure
on $V$ is a map $J\in {\rm End}(V)$ such that $J^2=-I$ and $J$ is
orthogonal. Since ${\rm dim}(V)=2n$ is even, such maps always exist.
We can then define a complex inner product space $V_J$ from $V$ by
taking $J$ to be multiplication by $i$ and taking the complex inner
product on $V$ as 
$(v,w)_{\Bbb C} = (v,w)_{\Bbb R} -i(Jv,w)_{\Bbb R}$, where
$(v,w)_{\Bbb R}$ denotes the original real inner product on $V$.
\vskip .1in
\noindent \bf Lemma. \rm $V_J$ is a complex inner product space with
$(v,v)_{\Bbb R} =(v,v)_{\Bbb C}$. 
\vskip .05in
\bf \noindent Proof. \rm Clearly $(v,w)$ is ${\Bbb R}$--bilinear.
Moreover 
$\overline{(v,w)_{\Bbb C}} = (v,w)_{\Bbb R}+i(Jv,w)_{\Bbb R} =
(w,v)_{\Bbb R} -i(Jw,v)_{\Bbb R}=(w,v)_{\Bbb C}$. Since $(Jv,w)
=i(v,w)$, it follows that $(v,w)_{\Bbb C}$ is ${\Bbb C}$--linear in 
$v$ and conjugate linear in $w$. Now $(Jv,v)_{\Bbb R} =-(v,Jv)_{\Bbb
R}=(Jv,v)_{\Bbb R}$, so that $(Jv,v)_{\Bbb R}=0$. Hence
$(v,v)_{\Bbb C}=(v,v)_{\Bbb R}$ and $V_J$ is a complex inner product
space.

\vskip .1in
\noindent \bf Theorem. \it The formula $C(v)=e(v)+e(v)^*$ gives a
faithful(=injective) irreducible representation of ${\rm Cliff}_{\Bbb
C}(V)$ on $S=\Lambda(V_J)$, called the ``spin module''. In particular ${\rm
Cliff}_{\Bbb C}(V) \cong {\rm End}(S)$. 
\vskip .05in
\noindent \bf Proof. \rm Clearly $v\mapsto C(v)$ is ${\Bbb
R}$--linear, $C(v)^*=C(v)$ and  $C(v)C(w) + C(w) C(v) =$
$$(e(v) +e(v)^*)(e(w)+e(w)^*) + (e(w) +e(w)^*)
(e(v) +e(v)^*) = 2{\rm Re} (v,w)_{\Bbb C} I=2 (v,w)_{\Bbb R} I.$$
Hence $C(v)$ satifies the real Clifford algebra relations and
therefore we get *--homomorphism of ${\rm Cliff}(V)$ into ${\rm
End}(\Lambda(V_J))$. Now the relation $C(v)=e(v) +e(v)^*$ implies
$C(Jv) =e(v) +e(Jv)^*=e(v)-ie(v)^*$. Hence $e(v)= {1\over 2}(C(v)
-iC(Jv))$ and $e(v)^*={1\over 2}(C(v)+iC(Jv))$. But the $e(v)$'s and
$e(v)^*$'s act irreducibly on $\Lambda(V_J)$ (it is the complex wave
representation), so the $C(v)$'s must also act irreducibly. Therefore
the $C(v)$'s generate ${\rm End}(S)$. Thus the 
image of ${\rm Cliff}_{\Bbb C}(V_J)$ has ${\Bbb C}$--dimension ${\rm
dim}(S)^2= 2^{2{\rm dim}_{\Bbb C}(V_J)} = 2^{{\rm dim}_{\Bbb R}(V)}$.
But this is the ${\Bbb C}$--dimension of ${\rm Cliff}_{\Bbb C} (V)$,
so the representation of ${\rm Cliff}_{\Bbb C} (V)$ is faithful and
surjective. Hence ${\rm Cliff}_{\Bbb C} (V)\cong {\rm End}(S)$.
Moreover the representation must a fortiori be faithful on the real
subalgebra ${\rm Cliff}(V)$. 

\vskip .1in

\bf \noindent 5. QUANTISATION AND THE SPIN GROUP. \rm 
\vskip .1in
\noindent \bf Bogoliubov automorphisms of ${\rm Cliff}(V)$. \rm
Consider the compact group $SO(V)$. 
\vskip .1in
\noindent \bf Lemma. \it $SO(V)$ is connected. 
\vskip .05in
\noindent \bf Proof. \rm Any matrix in $SO(V)$ is conjugate to a
block diagonal matrix with $2\times 2$ diagonal blocks
$D_i=\pmatrix{\cos x_i & \sin x_i\cr -\sin x_i & \cos x_i\cr}$, so can
be connected by a continuous path to $I$ by the path of matrices with
blocks $D_i=\pmatrix{\cos tx_i & \sin tx_i\cr -\sin tx_i & \cos
tx_i\cr}$.
\vskip .1in
If $g\in SO(V)$, $v\mapsto c(gv)$ also satisfies the real Clifford
algebra relations, so induces an automorphism of ${\rm Cliff}(V)$. In
fact $SO(V)$ acts orthogonally on $\Lambda(V)$ via $g(x_1\wedge \cdots
\wedge x_k)=gx_1\wedge \cdots \wedge gx_k$, so that
$ge(v)g^{-1}=e(gv)$ and hence $gc(v)g^{-1} =c(gv)$ since
$c(v)=e(v)+e(v)^*$. Thus $SO(V)$ normalises ${\rm Cliff}(V)$ on
$\Lambda(V)$. We write $\alpha_g$ for the automorphism of ${\rm
Cliff}(V)$ and ${\rm Cliff}_{\Bbb C}(V)$ induced by ${\rm Ad}\, g$,
$a\mapsto g ag^{-1}$. In particular $g_0=-I$ acts and gives a period
two automorphism $\gamma=\alpha_{-I}$ of ${\rm Cliff}(V)$ satisfying 
$\alpha c(v) =-c(v)$. This automorphism gives rise to a ${\Bbb
Z}_2$--grading on $A={\rm Cliff}(V)$, because we can take the $\pm 1$
eigenspaces $A_\pm$ of $\gamma$. Clearly $A_+A_+\subset A_+$,
$A_+A_-\subset A_-$, $A_-A_+\subset A_-$ and $A_-A_-\subset A_+$.
Under the identification $A\equiv \Lambda(V)$, $A_+=\Lambda^{\rm
even}(V)$ and $A_-=\Lambda^{\rm odd}(V)$. 

Now if $v\mapsto C(v)$ is the irreducible representation of the
Clifford algebra 
relations on the spin module $S$, $v\mapsto C(gv)$ will given another
irreducible representation on $S$. By uniqueness we can find $U_g\in
U(S)$ such that $C(gv) =U_g C(v) U_g^*$ for all $v$. 
Note that $g\in SO(V)$ commutes with the complex structure $J$ iff $g\in
SU(V_J)$. In this case $g$ is canonically implemented on
$S=\Lambda(V_J)$ by $g(v_1\wedge \cdots \wedge v_k) = gv_1\wedge
\cdots \wedge gv_k$. In particular $g_0=-I$ commutes with all $J$'s,
so is canonically implemented on each $S$: $g_0$ acts as $\pm
1$ on $S^\pm$.

The choice of $U_g$ is not unique. If $U^\prime_g$
is another possible choice, then $U_g^*U^\prime_g$ must commute with
all $C(v)$'s and hence must be a scalar matrix by Schur's lemma. Thus
$U_g$ is uniquely determined up to a phase in ${\Bbb T}$, so that
$U_g$ really gives 
a homomorphism of $SO(V)$ into $U(S)/{\Bbb T}=PU(S)$, the projective
unitary group. This is what is meant by quantisation. The prequantised
action on $V$ can be implemented on Fock space $S$ by a unitary; the
phase represents the anomoly that usually arises when we quantise. As
we shall see, we really get a 2--valued representation of $SO(V)$ or
equivalently a representation of a double cover, called ${\rm
Spin}(V)$, which we now construct. Observe first that $U_g C(v)
U_g^*=C(gv)$, so that $U_g$ normalises the real 
subalgebra $A={\rm Cliff}(V)$ of ${\rm End}(S)$. 
\vskip .1in
\noindent \bf Theorem (Noether--Skolem). \it $g\in {\rm End}(S)$
normalises $A$ iff $g\in A^*\cdot {\Bbb C}^*$, where $A^*$ denotes the
invertible elements in $A$. 
\vskip .05in
\noindent \bf Proof. \rm We know that ${\rm End}(S)=A\oplus iA$, a
direct sum of real vector spaces. Let $g=a+ib$ with $a,b\in A$ and set
$\alpha(a)=gag^{-1}$. Then $(a+ib)x =\alpha(x) (a+ib)$. Hence
$ax=\alpha(x)a$ and $bx=\alpha(x) b$. Consider the polynomial $p(t)={\rm
det}(a+t b)$. Since $p(i)\ne 0$, we can find $t\in {\Bbb R}$ such that
$p(t)\ne 0$. Let $h=a+tb\in A$ and let $h^{-1}=u+iv$. Then
$h(u+iv)=I$, so that $hv=0$ and hence $v=0$. Thus $h^{-1}\in A$. Since
$hx=(a+tb)x=\alpha(x)(a+tb)=\alpha(x)h$, it follows that $z=h^{-1}g$
commutes with $A$, so lies in ${\Bbb C}^*$. Hence $g=hz$ as claimed.
\vskip .05in
\noindent \bf Corollary. \it For each $g\in SO(V)$, there is a unitary
element $u_g\in A^*$ uniquely determined up to a sign such that $u_g
c(v) u_g^*=c(gv)$. 
\vskip .05in
\noindent \bf Proof. \rm Suppose $u_g=\lambda U_g$. Then
$u_gu_g^*=|\lambda|^2 = u_g^*u_g$. Scaling $u_g$, we may therefore
arrange that $u_g$ is unitary. Since $A^*\cap {\Bbb C}^*={\Bbb
R}^*$, $u_g$ is uniquely determined up to sign.
\vskip .1in
\noindent \bf The spin group. \rm Let ${\rm Spin}(V)=\{\pm u_g:g\in
SO(V)\}\subset {\rm Cliff}(V)$, the spin group. 
\vskip .1in
\noindent \bf Lemma. \it ${\rm Spin}(V)$ consists of unitaries 
$u\in{\rm Cliff}(V)$ normalising $c(V)$ such that the orthogonal
transformation $g$ defined by $c(gv)=uc(v)u^*$ lies in $SO(V)$. In
particular ${\rm Spin}(V)$ is a closed subgroup of the unitary group
of $A$, so compact.
\vskip .05in
\noindent \bf Proof. \rm Clearly any element of ${\rm Spin}(V)$
satisfies these conditions. The converse is obvious from the
corollary by uniqueness. 
\vskip .1in
The map ${\rm Spin}(V)\rightarrow SO(V)$ is a surjective contiuous
homorphism, by construction. Its kernel is ${\pm I}$, so that ${\rm
Spin}(V)$ is a double cover of $SO(V)$.
\vskip .1in
\noindent \bf Theorem. \it (a) ${\rm Spin}(V)$ is connected.

\noindent (b) ${\rm Spin}(V)\subset {\rm Cliff}^+(V)$.
\vskip .05in
\noindent \bf Proof. \rm (a) Let $f:{\rm Spin}(V)\rightarrow {\Bbb Z}$
be a continuous function; we must show it is constant. If we show that
$f(-g)=f(g)$ for all $g$, then $f$ will drop to a continuous map of
$SO(V)$ into ${\Bbb Z}$ and hence be constant, by the connectivity of
$SO(V)$. But $x(t)=\cos \pi t +c(e_1) c(e_2)\sin \pi t$ ($t\in [0,1]$)
is a continuous path in ${\rm Spin}(V)$ from $I$ to $-I$. (To see this
either use the representation $C(e_1)=\pmatrix{1 & 0\cr 0 &
-1\cr}$, $C(e_2)=\pmatrix{0 & 1\cr 1 & 0\cr}$ to write
$x(t)=\pmatrix{ \cos \pi t & \sin \pi t \cr -\sin \pi t & \cos \pi
t\cr}$; or note that $x(t)=\exp \pi c(e_1)c(e_2)t$ with $c(e_1)c(e_2)$
skew--adjoint.) Hence $t\mapsto f(gx(t))$ is continuous so
constant. Hence $f(g)=f(-g)$. 

\noindent (b) Let $u_0$ be the element of ${\rm Cliff}(V)$
implementing the grading automorphism $\gamma$. Thus $u_0 c(v)
u_0^*=-c(v)$. But $\alpha_g(u_0)$ also implements $\gamma$, so that
$\alpha_g(u_0)=\lambda(g) u_0$ with $\lambda(g)=\pm 1$. Thus
$\lambda(g)$ is a continuous homorphism $SO(V)\rightarrow \{\pm 1\}$.
Since $SO(V)$ is connected, $\lambda(g)\equiv 1$. Since
$\alpha_g(u_0)=u_gu_0u_g^*$, this implies that $u_g$ commutes with
$u_0$. But then $\gamma(u_g)=u_0u_gu_0^*=u_g$, so that $u_g\in {\rm
Cliff}^+(V)$.
\vskip .1in
\noindent \bf Remark. \rm There is also an infinitessimal version of
the action of ${\rm Spin}(V)$ in terms of bilinear combinations of
fermions (see the exercises).
\vskip .1in
\def\gl{{\rm gl}}
\noindent \bf 6. MATRIX GROUPS AND THEIR LIE ALGEBRAS. \rm We start by
proving von Neumann's theorem on closed subgroups of $GL(V)$. We
define $\gl(V)={\rm End}\,V$ with the usual operator norm.
\vskip .1in
\noindent \bf Lemma (Lie's formulas). \it If $a,b\in {\rm End}(V)$
then $(\exp(a/n) \exp(b/n))^n\rightarrow \exp(a+b)$ and \break
$(\exp(a/n)\exp(b/n)\exp(-a/n)\exp(-b/n))^{n^2}\rightarrow \exp[a,b]$.
\vskip .05in
\noindent \bf Proof. \rm Recall that $\exp(a) =\sum a^n/n!$ for all $a$ and
$\log(1+x) = \sum (-1)^{n+1}x^n/n$ for $\|x\|<1$. For $\|a\|$
sufficiently small, we have $\log\exp a=a$ and for $x$ sufficiently
small $\exp\log (1+x)=1+x$. Then
$$\log([\exp(a/n)\exp(b/n)]^n) =n \log(1 +(a+b)/n +O(1/n^2)) = a+b
+O(1/n)\rightarrow a+b,$$
and 
$$\log([\exp(a/n)\exp(b/n) \exp(-a/n)\exp(-b/n)]^{n^2})=
n^2\log(1+[a,b]/n^2 + O(1/n^3)) = [a,b] +O(1/n) \rightarrow [a,b].$$

\vskip .1in
\noindent \bf Theorem (von Neumann). \it Let $G$ be a closed
subgroup of $GL(V)$ and let 
$${\rm Lie}(G)=\{X\in {\rm End}\, V \, |\,
\exp(tX)\in G\,\hbox{for all $t$}\}.$$
Then ${\rm Lie}(G)$ is a linear subspace of ${\rm End}\, V$
closed under 
the Lie bracket $[a,b]=ab-ba$ and
$\exp({\rm Lie}(G))$ is a neighbourhood of $1$ in $G$. In fact if $U$ is a
sufficiently small open neighbourhood of $0$ then $\exp(U)$ is an open
neighbourhood of $1$ in $G$ and $\exp$ gives a homeomorphism between
$U$ and $\exp(U)$. 
\vskip .05in
\bf \noindent Proof. \rm  Lie's formulas applied to $tX$ and $tY$
immediately show that ${\rm Lie}(G)$ is a 
subspace closed under the bracket $[X,Y]=XY-YX$. 

It remains to show that $\exp({\rm Lie}(G))$ is a neighbourhood of $1$ in $G$.
Let ${\rm Lie}(G)^\perp$ be a vector subspace complementing ${\rm Lie}(G)$
$\gl(V)$, so that 
$\gl(V) = {\rm Lie}(G)\oplus {\rm Lie}(G)^\perp$. By the
inverse function theorem,
$X\oplus Y\mapsto \exp(X)\exp(Y)$ gives a homeomorphism between a
neighbourhood of $0$ in ${\rm End}(V)$ and $1$ in $GL(V)$ (its
derivative is $I$). If $\exp({\rm Lie}(G))$ is
not a neighbourhood of $1$ in $G$, then we can find $g_n\in G$ with
$g_n\rightarrow 1$ but $g_n\notin \exp({\rm Lie}(G))$. Write
$g_n=\exp(X_n)\exp(Y_n)$ with $X_n\in {\rm Lie}(G)$, $Y_n\in {\rm
Lie}(G)^\perp$. By 
assumption $Y_n\ne 0$ for all $n$. But since $\exp(X_n)$ and $g_n$ are
in $G$, it follows that $\exp(Y_n)\in G$ for all $n$. Since
$g_n\rightarrow 1$, we must have $Y_n\rightarrow 0$. 
By compactness, we may assume by passing to a subsequence if necessary
that $Y_n/\|Y_n\|\rightarrow Y\in{\rm Lie}(G)^\perp$ with $\|Y\|=1$. Since
$\|Y_n\|\rightarrow 0$, we can choose integers $m_n$ such that $m_n
\|Y_n\| \rightarrow 
t$. Then $\exp(m_n Y_n)=\exp(Y_n)^{m_n}\in G$ has limit $\exp(tY)$.
Since $G$ is closed, $\exp(tY) \in G$ for $t>0$ and hence for all $t$
on taking inverses. So by definition $Y$ lies in ${\rm Lie}(G)$, a
contradiction. 

\vskip .1in
This result says that matrix groups are Lie groups. If $G$ is a matrix
group, we denote its Lie algebra by ${\rm Lie}(G)$. We shall be
interested in matrix groups that are closed subgroups of $O(n)$.
Since $O(n)\subset U(n)$, they are also closed subgroups of $U(n)$.
\vskip .1in
\noindent \bf Corollary. \it Let $G$ and $H$ be matrix groups and
$\pi:G\rightarrow H$ a continuous homomorphism. Then there is a unique
Lie algebra homomorphism $\pi:{\rm Lie}(G)\rightarrow {\rm Lie}(H)$
such that $\pi(\exp(X)) =\exp\pi(X)$ for $X\in {\rm Lie}(G)$.
\vskip .05in
\noindent \bf Proof. \rm Uniqueness follows because we may replace
$X$ by $tX$ and take the coefficient of $t$. Conversely note that
$\pi\exp(tX)$ is a one parameter subgroup in $H$. Now $H$ is a closed
subgroup of $U(n)$; since commuting unitaries can be simultaneously
diagonalised, it follows that $\pi\exp (tX)= \exp t A$ for some matrix
skew--adjoint matrix $A$. But then by definition $A$ lies in ${\rm
Lie}(H)$. We define $\pi(X)=A$. From Lie's formulas, the map $X\mapsto
\pi(X)$ is a Lie algebra homomorphism. 
\vskip .1in
There is also an infinitesimal version of
the action of ${\rm Spin}(V)$ in terms of bilinear combinations of
fermions.
\vskip .1in
\def\gl{{\rm gl}}

\noindent \bf Proposition. \it (a) ${\rm Lie}(SO(V))=
\{A\in {\rm End}(V):A^t=-A\}$. 

\noindent (b) ${\rm Lie}({\rm Spin}(V))=\{x\in {\rm
Cliff}^+(V)|x^*=-x, \,[x,c(V)] \subset c(V)\}= {\rm lin}_{\Bbb
R}\{c(a)c(b)-c(b)c(a):a,b\in V\}$. A basis is given by $c(e_i)c(e_j)$
with $i<j$. 

\noindent (c) If $\pi:{\rm Spin}(V)\rightarrow SO(V)$ is the double
cover, $\pi^{-1}(A)={1\over 4} \sum_{i\ne j} a_{ij} c(e_i) c(e_j)$.
In other worsds, if $A\in so(V)$, then $\pi^{-1}(A)={1\over 4} \sum_{ij}
c(A\cdot e_i) c(e_j)$. 

\noindent (d) If $X\in {\rm Lie}({\rm Spin}(V))$, then $[X,c(v)]=
c(\pi(X	)v)$. 
\vskip .05in
\noindent \bf Proof. \rm (a) is obvious. To prove (b) and (c), note
that if $e^{yt}$ lies in ${\rm Spin}(V)$, then $y$ is even, $y^*=-y$
and $e^{y t} c(v) e^{-y t}=c(e^{At}v)$ for some 
$A\in {\rm Lie}(SO(V))$. Taking the coefficient of $t$, we get
$[y,c(v)]=c(Av)$. Let $(Ae_j,e_i)=a_{ij}$, so that $a_{ij}$ is
antisymmetric and real, and let $x={1\over 2}\sum_{i,j} a_{ij}
c(e_i)c(e_j)$. Then
$$\eqalign{[x,c(v)] &={1\over 4} \sum_{i,j} a_{ij} [c(e_i)c(e_j),c(v)]
={1\over 4} \sum_{i,j} a_{ij} [-\{c(e_i),c(v)\} c(e_j)
+c(e_i)\{c(e_j),c(v)\}]\cr 
&={1\over 4} \sum_{i, j} a_{ij}[(v,e_i) c(e_j) -(v,e_j)c(e_i)]
=c(Av),\cr}$$
using the graded commutator $\{c(u),c(v)\}=c(u)c(v)+c(v)c(u) =2(u,v)I$
and the rules for computing graded commutators.
Thus $[y-x,c(v)]=0$ for all $v\in V$ and therefore $y-x$ must be a
real scalar. Since $(y-x)^*=-(y-x)$, we deduce that $y=x$ as
required. The map between Lie algebras is ${1\over 2} \sum a_{ij}
c(e_i)c(e_j) \rightarrow (a_{ij})$ by uniqueness.
Finally, if $A\in {\rm Lie}({\rm Spin}(V))$, then $e^{At} c(v) e^{-At}
= c(\pi(e^{At})v)=c(e^{\pi(A)t}v)$. Taking coefficients of $t$, we get
$[A,c(v)]=c(\pi(A)v)$, so (d) follows.
\vskip .2in
\noindent \bf 7. THE ODD--DIMENSIONAL CASE. \rm The structure of Clifford
algebras and spin groups for odd--dimensional inner product spaces can
easily be deduced from the even--dimensional case. Let $V$ be a real inner
product space with odd dimension. We may write $V=V_0\oplus {\Bbb
R}e_0$, where $V_0$ is even--dimensional and $e_0$ is a unit
vector. Let $e_1,\dots,e_{m}$ be an  
basis of $V_0$, with $m=2n$. Suppose that the operators $c(v)$
satisfy the real Clifford algebra
relations $c(v)c(w) + c(w) c(v)=2(v,w)$, $c(v)^*=c(v)$. 
\vskip .1in
\noindent \bf Lemma. \it The element $z=c(e_0)c(e_1)c(e_2) \cdots
c(e_m)$ commutes with all $c(v)$'s and satisfies $z^2=(-1)^n$.
\vskip .1in
\noindent \bf Proof. \rm It is immediate that $c(e_i)z=zc(e_i)$ for
all $i$, since there are an odd number of $e_j$'s. Hence $c(v)z=zc(v)$
for all $v$. A simple induction argument shows that
$c(e_0)  \cdots c(e_k) = (-1)^{k(k-1)/2} c(e_k) c(e_{k-1})
\cdots c(e_0)$. Thus $c(e_0)c(e_1)\cdots c(e_m)=(-1)^n c(e_m)\cdots
c(e_0)$, so that $z^2=(-1)^n$. 
\vskip .1in
Suppose that the $c(v)$'s act irreducibly on $W$, a complex inner
prouct space. By Schur's lemma, $z$ must be a scalar, so $\pm
(i)^n$. Thus 
$$c(e_0)=\mp i^n c(e_1)\cdots c(e_m).\eqno{(1)}$$ It follows that
$c(e_1),\cdots, c(e_m)$ already act irreducibly on $W$. Thus ${\rm
dim}(W)= 2^n$ is the standard irreducible representation of ${\rm
Cliff}_{\Bbb C}(V_0)$. Conversely let $W$ be an irreducible
representation of ${\rm Cliff}(V_0)$ and define $c(e_0)$ by (1). It is
easy to check that the $c(e_i)$'s satisfy the Clifford relations and
hence we get irreducible representations $W_\pm$ of ${\rm
Cliff}(V)$. The representations are inequivalent because
$c(e_0)\cdots c(e_m) = z$ with $z=\pm i^n$. The maximum dimension of
${\rm Cliff}_{\Bbb C}(V)$ is $2^{m+1}$ and by the double commutant
theorem we have a surjection onto ${\rm End}(W_+)\oplus {\rm
End}(W_-)$. Since this space also has dimension $2^{m+1}$, this map is
an isomorphism: 
$${\rm Cliff}_{\Bbb C}(V)\cong {\rm End}(W_+)\oplus {\rm End}(W_-).$$
\vskip .1in
\noindent \bf Lemma. \it The inclusion ${\rm Cliff}_{\Bbb
C}^+(V)\subset {\rm Cliff}_{\Bbb C}(V)$ induces isomorphisms 
${\rm Cliff}^+_{\Bbb C}(V) \cong {\rm End}(W_\pm)$. The spaces $W_\pm$
give equivalent irreducible representations of 
${\rm Cliff}^+_{\Bbb C}(V)$. 
\vskip .1in
\noindent \bf Proof. \rm Taking any reordering of the $e_i$'s, the
previous lemma shows that $$c(e_i) =\pm (i)^n c(e_0)\cdots \widehat{
c(e_i)} \cdots c(e_m)\in \pi({\rm Cliff}_{\Bbb C}^+(V)).$$ 
It follows that ${\rm Cliff}_{\Bbb C}^+(V)$ acts irreducibly on
$W_\pm$. Thus we 
have a surjection of ${\rm Cliff}_{\Bbb C}^+(V)$ onto ${\rm
End}(W_\pm)$. Since both space have dimension $2^m$, this is an
isomorphism, so the first assertion follows. The second follows
because a matrix algebra has a unique irreducible representation.
\vskip .1in
The theory of the spin group for odd--dimensional spaces now proceeds
exactly as in the even--dimensional case, but based on ${\rm
Cliff}^+(V)$ rather than ${\rm Cliff}(V)$. The group $SO(V)$ acts by
automorphisms on ${\rm Cliff}^+(V)$. On the other hand ${\rm
Cliff}^+(V)$ has a unique irreducible representation $W$. Therefore
for each $g\in SO(V)$ there is a unitary $U_g\in U(W)$, unique up to a
scalar multiple, such that $c(gv)c(gw)=U_g c(v)c(w) U_g^*$ for all
$v,w\in V$. Since every $c(e_i)$ can be expressed as a product of $n$
elements $c(v)c(w)$, this is equivelent to the condition that $U_g
c(v) U_g^* = c(gv)$. Let $A=\pi({\rm Cliff}^+(V))\subset {\rm
End}(W)$. Then 
$A+i A={\rm End}(W)$ and $A\cap iA=(0)$. The Noether--Skolem argument
then implies that we can find $u_g\in A^*$ with $u_gu_g^*=I$ such that
$u_g c(v) u_g^*=c(gv)$. Moreover $u_g$ is unique up to a sign. Let
${\rm Spin}(V)=\{\pm u_g: g\in SO(V)\}$. The proofs of the following
results are as before.
\vskip .1in

\noindent \bf Lemma. \it ${\rm Spin}(V)$ consists of unitaries 
$u\in{\rm Cliff}(V)$ normalising $c(V)$ such that the orthogonal
transformation $g$ defined by $c(gv)=uc(v)u^*$ lies in $SO(V)$. In
particular ${\rm Spin}(V)$ is a closed subgroup of the unitary group
of $A$, so compact.
\vskip .1in
\rm The map ${\rm Spin}(V)\rightarrow SO(V)$ is a surjective contiuous
homorphism, by construction. Its kernel is ${\pm I}$, so that ${\rm
Spin}(V)$ is a double cover of $SO(V)$.
\vskip .1in
\noindent \bf Theorem. \it (a) ${\rm Spin}(V)$ is connected.

\noindent (b) ${\rm Spin}(V)\subset {\rm Cliff}^+(V)$.
\vskip .1in

\noindent \rm There is also an infinitesimal version of
the action of ${\rm Spin}(V)$ in terms of bilinear combinations of
fermions.
\vskip .1in

\noindent \bf Proposition. \it (a) ${\rm Lie}(SO(V))=
\{A\in {\rm End}(V):A^t=-A\}$. 

\noindent (b) ${\rm Lie}({\rm Spin}(V))=\{x\in {\rm
Cliff}^+(V)|x^*=-x, \,[x,c(V)] \subset c(V)\}= {\rm lin}_{\Bbb
R}\{c(a)c(b)-c(b)c(a):a,b\in V\}$. A basis is given by $c(e_i)c(e_j)$
with $i<j$. 

\noindent (c) If $\pi:{\rm Spin}(V)\rightarrow SO(V)$ is the double
cover, $\pi^{-1}(A)={1\over 4} \sum c(A\cdot e_i) c(e_i)$.

\noindent (d) If $A\in {\rm Lie}({\rm Spin}(V))$, then $[\pi^{-1}(A),c(v)]=
c(\pi(A)v)$.
\vskip .2in
\noindent \bf 8. THE SPIN REPRESENTATIONS OF ${\rm Spin}(V)$. \rm We
treat the even and odd dimensional cases separately. 
\vskip .1in
\noindent \bf Proposition (irreducibility of spin
representations). \it If ${\rm dim}V$ is even, the spin 
representations $S^\pm$ of ${\rm Spin}(V)$ are irreducible. If ${\rm
dim}(V)$ is odd, the spin representations $S$ of ${\rm Spin}(V)$ is
irreducible. 
\vskip .1in
\noindent \bf Proof. \rm (1) If ${\rm dim}(V)$ is odd, we have seen
that the $c(v)c(w)$'s act irreducibly on $S$. Since the algebra these
generate is the double commutant of ${\rm Spin}(V)$, ${\rm Spin}(V)$
acts irreducibly.

\noindent (2) Since $A={\rm Cliff}^+(V)$ is generated by the
$c(v)c(w)$'s as a unital algebra, it is generated by ${\rm spin}(V)$
as a unital algebra. We claim that ${\rm Cliff}^+(V)$ acts irreducibly
on $S^\pm$ and that these are inequivalent representations. In fact
${\rm Cliff}(V)$ acts irreducibly on $W^+\oplus W^-$. There are two
ways to see this. (a) We may introduce $\gamma = c(v_1) \cdots
c(v_{2m})\in A$, where $(v_i)$ is an orthonormal basis of $V$. Then
$\gamma c(v) = -c(v) \gamma$ for $v\in V$. Thus $\gamma$ commutes 
with $A$ but is not a scalar. On the other hand $\gamma^2$ is central
in ${\rm Cliff}(V)$ and unitary, so $\gamma^2 =\pm I$. Thus $\gamma$
must take distinct values on $W^+$ and $W^-$, proving their
inequivalence. (b) Take a unit vector $v\in V$ and set $g=c(v)$. Then
$g^2=I$ and $gAg^{-1}=A$. If the two representations $W^\pm$ are
equivaelent on $A$, then $A$ would be isomorphic to ${\rm
End}(W^\pm)$ and hence have $1/4$ times the dimension of ${\rm
Cliff}(V)$. But ${\rm Cliff}^+(V)=A$ and ${\rm Cliff}^-(V)=Ag$, so 
that ${\rm dim}\, {\rm Cliff}(V)=2\cdot {\rm dim}\,A$, a
contradiction. Finally  ${\rm Cliff}(V)$ is generated by $c(v)c(w)$'s
and therefore 
by the image ${\rm spin}(V)$. It follows that the representations
$W^\pm$ are irreducible and inequivalent on ${\rm spin}(V)$ and hence
${\rm Spin}(V)$. 
\vskip .1in
Every matrix $g\in SO(N)$ is conjugate to a matrix with $2\times 2$
blocks $\pmatrix{\cos\theta_j & \sin\theta_j \cr -\sin\theta_j
&\cos\theta_j\cr}$ down the diagonal where $j=1,\dots,[N/2]$. Note that
there is an additional $1$ on the diagonal if $N$ is odd. Thus the
complex eigenvalues of $g$ are $e^{\pm i \theta_j}$. In any irreducible
projective representation $\pi(g)$ of $SO(N)$, the eigenvalues of
a generic block diagonal element $\pi(g)$ are called the {\it weights}
of the representation (see Chapter~II for a more precise definition).

\vskip .1in
\noindent \bf Lemma (weights of spin representation). \it (1) If ${\rm
dim}(V)$ is even, the weights of the spin representation $W^\pm$ are
$\exp{i {1\over 2} \sum \pm \theta_k}$ where the number of plus signs is even
for $W^+$ and odd for $W^-$.

\noindent (2) If ${\rm dim}(V)$ is odd, the weights of the spin
representation $W$ are $\exp{ {1\over2} \sum \pm i\theta_k}$. 
\vskip .1in
\bf \noindent Proof. \rm Let $v_1,\dots ,v_m$ be an orthonormal basis
for $V_J$ and 
set $v_I=v_{i_1}\wedge \cdots \wedge v_{i_k}$ for
$i_1<\cdots<i_k$. Then
$c(v_k)=e(v_k) + e(v_k)^*$ and $c(iv_k) = i(e(v_k)-e(v_k)^*)$. The
generators of Lie algebra of the maximal torus are given by
$T_j={1\over 4} (c(v_j)c(iv_j) - c(iv_j)c(v_j))={1\over 2} c(v_j)
c(iv_j)$. Thus $T_j v_I = i/2 \, v_I$ if $j\in I$ and $T_j v_I =-i/2\,
v_I$ if $j\notin I$. The corresponding self--adjoint operators 
$S_j$ satisfy $T_j=iS_j$ so that $S_jv_I=\pm {1\over 2} v_I$. Since
$W^+=\Lambda^{\rm even}(V)$ is spanned by $v_I$'s with $I$ even and 
$W^-=\Lambda^{\rm odd}(V)$ is spanned by $v_I$'s with $I$ odd, the
result follows. 

\noindent (2) Note that the maximal torus of $SO(V_0)$ coincides with
the maximal torus of $SO(V)$. On the other hand the spin
representation of $SO(V)$ equals $W^+\oplus W^-$, where $W^\pm$ are
the spin representations of $SO(V_0)$. So the result follows
immediately from (1). 
\vskip .1in

\noindent \bf Lemma (grading operator). \it If $(e_i)$ is any
orthonormal basis of $V$, 
then $c(e_1)\cdots c(e_n)$ equals $\pm u_0$, the operator
implementing the grading. Moreover $u_0^2=(-1)^{{1\over 2}
{\rm dim}(V)} I$. The grading operator on $S$ is given by $\lambda
u_0$ where $\lambda =(i)^{{1\over 2} {\rm dim}(V)}$.
\vskip .05in
\noindent \bf Proof.\rm If ${\rm dim}(V)=2m$, then the elements
$a_i=c(e_{2i-1})c(e_{2i})$ commute and satisfy $a_i^2 =-1$. 
Hence $g=c(e_1)\cdots c(e_n)$ satisfies $g^2=(-1)^{m}I$. Moreover $g
c(e_i) =- c(e_i) 
g$. Hence $g=\pm u_0$ and $u_0^2=g^2=(-1)^{m}I$. Now a multiple
$\lambda u_0$ of $u_0$ 
acts as $\pm 1$ on $S^\pm$. Since $u_0^2=(-1)^m$, we get
$\lambda^2=(-1)^m$. 

\vskip .05in
\noindent \bf Corollary. \it ${\rm Spin}(V)=\{u\in {\rm Cliff}^+(V):
uu^*=u^*u=I,\, uc(V)u^*=c(V)\}$.
\vskip .05in
\noindent \bf Proof. \rm Suppose that $u\in {\rm Cliff}^+(V)$ is unitary
and that the orthogonal transformation $g$ with $uc(v)u^*=c(gv)$ has
determinant $-1$. Define $h\in O(V)$ by $he_1=-e_1$ and $he_i=e_i$ for
$i>1$. Then $x=g^{-1}h\in SO(V)$, so corresponds to $v\in {\rm
Spin}(V)$. Hence $h=gx$ corresponds to $w=uv\in {\rm Cliff}^+(V)$. But
$wc(e_1)\cdots c(e_n) w^*=c(he_1) \cdots c(he_n) =-c(e_1)\cdots
c(e_n)$, so that $\gamma(w)=-w$, a contradiction.

\vskip .1in
\noindent \bf Caveat. \rm If ${\rm dim}(V)$
is even and $J$ is a complex structure on $V$, then $U(V_J)\subset
SO(V)$ is the subgroup of $SO(V)$ commuting 
with $J$. Thus $U(V)$ acts canonically on $W_\pm=\Lambda^\pm V_J$ fixing
$\Lambda^0 V_J={\Bbb C}\Omega$. Denote this representation  by
$\pi$. Note however that the representation 
of $U(V)$ obtained by restricting the spin representations $W_\pm$ of
${\rm SO}(V)$ is given by ${\rm det}(g)^{-1/2} \pi(g)$.

\vfill\eject

\noindent \bf \centerline{CHAPTER~II. COMPACT MATRIX GROUPS AND SIMPLE
LIE ALGEBRAS} 
\vskip .2in
\noindent \bf PART~1. ELEMENTARY STRUCTURE THEORY. 
\vskip .1in
\noindent \bf 1. COMPACT LIE ALGEBRAS. \rm Let $G$
be closed subgroup 
of $O(V)$ or $U(V)$, where $V$ is a real or complex inner product
space. Since $G$ is closed and $O(V)$ and $U(V)$ compact, $G$ must
also be compact. We call $G$ a compact matrix group. Recall that we
have defined the Lie
algebra of $G$ as 
$${\rm Lie}(G)=\g=\{X\in {\rm End}\, V \, |\,
\exp(tX)\in G\,\hbox{for all $t$}\}.$$
By von Neumann's theorem, $\g$ is a linear subspace closed under the
Lie bracket and the exponential is locally a homeomorphism between
neighbourhoods of $0\in \g$ and $1\in G$. Thus $\g$ is indeed a Lie
algebra! It is 
routine to check that the bracket defines the Jacobi identity
$$[[X,Y],Z] +[[Y,Z],X] + [[Z,X],Y]=0$$
for $X,Y,Z\in \g$. An abstract course in Lie algebras might take the
antisymmetry of the bracket ($[X,Y]=-[Y,X]$ ) and the Jacobi identity
as the defining axioms of a Lie algebra. Since every
finite--dimensional Lie algebra can be realised as a Lie subalgebra of
matrices, we prefer a more concrete approach. We also proved that
every continuous homorphism between matrix groups  
$\pi:G\rightarrow H$ gave rise to a unique
Lie algebra homomorphism $\pi:{\rm Lie}(G)\rightarrow {\rm Lie}(H)$
such that $\pi(\exp(X)) =\exp\pi(X)$ for $X\in {\rm Lie}(G)$. The
converse of this statement is also true when suitably interpreted. We
shall return to this point later.

Note that if $G$ is a closed subgroup of
$U(V)$ for some $V$ and hence $\g$ carries a real inner product,
namely ${\rm Re}\, {\rm Tr}(XY^*)=(X,Y)$. This inner product is
invariant under $G$ and therefore under $\g$, i.e.~$([X,Y],Z)
+(Y,[X,Z])=0$. We define an involution on $\g$ by $X^*=-X$. This
extends to a conjugate linear involution on $\g_{\Bbb C}=\g +i\g$,
$(X+iY)^*= X^*-i Y^*=-X +iY$. The inner prouct extends to a complex
inner product on $\g_{\Bbb C}$ such that ${\rm ad}(X^*)={\rm
ad}(X)^*$. 

Now suppose that $\g$ is a Lie algebra with an invariant
inner product $(X,Y)$. We call $\g$ a compact Lie algebra. Clearly if
$\b$ is an ideal in $\g$ then so is 
$\b^\perp$, and $\g=\b\oplus \b^\perp$. 
\vskip .1in
\noindent \bf Lemma. \it Every compact Lie algebra is the direct sum
of an Abelian algebra (its centre) and simple Lie algebras (of compact
type). 
\vskip .1in
\noindent \bf Proof. \rm $\g$ acts as a *--representation on $\g$ via
${\rm ad}$. We can therefore decompose $\g$ into a direct sum of
irreducibles. Grouping together the copies of the trivial
representation into $\g_0$, we have
$$\g = \g_0 \oplus \g_1 \oplus \cdots \oplus \g_n,$$
with $\g_0=\z$, the centre of $\g$. Clearly $[\g,\g_i]\subset \g_i$
and $[\g_j,\g_i]=(0)$ for $i\ne j$. By definition of $\g_0$, we must
therefore have $[\g_i,\g_i]\ne (0)$. Any ${\rm ad}\g_i$--submodule of
$\g_i$ is clearly ${\rm ad}\g$--invariant, so that $\g_i$ is simple
and non--Abelian. Indeed $[\g_i,\g_i]=\g_i$, since it is a non--zero
ideal; this implies $[\g,\g]=\bigoplus_{i>0} \g_i =\z^\perp$.  We call
$\g$ semisimple if its centre is trivial. 
Clearly if $\g$ is semisimple (and compact), then $[\g,\g]=\g$. (Weyl
proved that any complex semisimple Lie algebra is the complexification
of a compact Lie algebra; see below.) 
\vskip .1in
\noindent \bf Lemma. \it If $\g$ is semisimple and compact iff the
Killing form
$B(X,Y)={\rm Tr}({\rm ad}(X){rm ad}(Y))$ is negative definite.
\vskip .1in
\noindent \bf Proof. \rm If $B$ is negative definite, $(X,Y)=-B(X,Y)$
gives an invariant inner product on $\g$. If $X\in \z$, then $(X,Y)=0$
for all $Y$, so $X=0$ and hence $\g$ is semisimple.
Conversely if $B$ is semisimple and compact, let $(e_i)$ be an
orthonormal basis of $\g$ for the invariant inner product $(X,Y)$.
$$B(X,X)={\rm Tr}({rm ad}\, X)^2=\sum (({\rm ad}\, X)^2 e_i,e_i)
=-\sum \|{rm ad}(X)e_i\|^2.$$
Thus $B(X,X)\le 0$ with equality iff ${\rm ad}(X)e_i=0$ for all
$i$ iff ${\rm ad}(X)=0$ iff $X=0$. 

\vskip .1in

\noindent \bf Lemma. \it If $\g$ is compact and semisimple, every
derivation of $\g$ is inner. 
\vskip .1in
\noindent \bf Proof. \rm Let $\d\in so(\g)$ be the Lie algebra of
derivations preserving the inner product. Clearly $\g\equiv {\rm
ad}(\g)\subseteq \d$. Write $\d={\rm ad}(\g)\oplus {\rm
ad}(\g)^\perp$, orthogonal ideals. If $D\in {\rm ad}(\g)^\perp$, then
$[D,{\rm  ad}(\g)]=0$. Thus ${\rm ad}(D\g)=0$, so that $D\g=0$ and
hence $D=0$. Thus $\d={\rm ad}(\g)$. 

\vskip .1in

\vskip .1in
\noindent \bf 2. EXAMPLES OF SIMPLE COMPACT LIE ALEGBRAS. 
\vskip .1in
\noindent \bf Classical compact groups. \rm The classical compact simple Lie
groups are $SU(N)$, $SO(2N+1)$, $U({\Bbb H}^n)$, $SO(2N)$. It is easy
to compute their Lie algebras and verify that they are simple.

\vskip .1in
\noindent \bf Compact Lie algebras constructed from lattices. \rm We
shall now give a method of constructing a compact Lie algebra from a
lattice, essentially due to Tits. Let $\Lambda$ be a lattics in the
real inner product space $V$ such 
that $\Lambda$ is integral, i.e.~$(\alpha,\beta)\in {\Bbb Z}$ for
$\alpha,\beta\in \Lambda$, and $\Lambda$ is even,
i.e.~$(\alpha,\alpha)\in 2{\Bbb Z}$ for $\alpha\in \Lambda$. The
lattice has a natural bicharacter
$B(\alpha,\beta)=(-1)^{(\alpha,\beta)}$ with values in ${\Bbb
Z}_2=\{\pm 1\}$, i.e.~a bilinear form 
$\Lambda\times \Lambda\rightarrow {\Bbb Z}_2$. An
$\varepsilon$--factor is a bilinear map $\varepsilon:\Lambda\times
\Lambda\rightarrow {\Bbb Z}_2$ such that
$B(\alpha,\beta)=\varepsilon(\alpha,\beta)\varepsilon(\beta,\alpha)$
and $\varepsilon(\alpha,\alpha)=(-1)^{\|\alpha\|^2/2}$. 
\vskip .1in
\noindent \bf Lemma. \it Every $\varepsilon$--factor has the form
$\varepsilon(\alpha,\beta)=(-1)^{b(\alpha,\beta)}$ where
$b:\Lambda\times \Lambda\rightarrow{\Bbb Z}$ is a bilinear map such
that $b(\alpha,\alpha)\equiv \|\alpha\|^2/2$ modulo $2$. 

\vskip .1in
\noindent \bf Proof. \rm Let $\alpha_1,\dots,\alpha_m$ be a ${\Bbb
Z}$--basis of 
$\Lambda$. Define $b:\Lambda\times \Lambda\rightarrow {\Bbb Z}$
bilinear by $b(\alpha_i,\alpha_i)=\|\alpha_i\|^2/2$,
$b(\alpha_i,\alpha_j)=(\alpha_i,\alpha_j)$ if $i<j$ and
$b(\alpha_i,\alpha_j)=0$ if $i>j$. Clearly
$\varepsilon(\alpha,\beta)=(-1)^{b(\alpha,\beta)}$ is an
$\varepsilon$--factor. Conversely given an $\varepsilon$--factor,
choose $b(\alpha_i,\alpha_j)\in {\Bbb Z}$ such that
$\varepsilon(\alpha_i,\alpha_j)=(-1)^{b(\alpha_i,\alpha_j)}$.
Extending $b$ bilinearly to $\Lambda\times \Lambda$, we evidently have
$\varepsilon = (-1)^b$. 
\vskip .1in 
Using the inner product each $\alpha\in \Lambda$ defines a
real linear form $\alpha(H)=(H,alpha)$ on $V$ which extends by complex
linearity to $\h=V\oplus iV$. Let $\Phi=\{\alpha\in
\Lambda:\|\alpha\|^2=2\}$ and let 
$$\g=\h
\oplus_{\alpha\in \Phi} {\Bbb C}E_\alpha.$$ 
Define non--trivial
brackets by $[H,E_\alpha]=\alpha(H)E_\alpha$, $[E_\alpha,E_{-\alpha}]=
-\alpha$ and
$[E_\alpha,E_\beta]=\varepsilon(\alpha,\beta) E_{\alpha+\beta}$ if
$\alpha+\beta$ is a root. 
\vskip .1in
\noindent \bf Proposition. \it The above brackets make $\g$ into a
complex Lie algebra. 
\vskip .1in
\noindent \bf Proof. \rm It is easy to verify that the bracket
satisfies $[X,Y]=-[Y,X]$ by taking $X$ and $Y$ to be basis
elements. We therefore have only to check that the Jacobi identity 
$$[[X,Y],Z] + [[Y,Z],X] + [[Z,X],Y]=0$$ 
is satisfied when $X$ and $Y$ are basis elements. This identity is
obvious if at least two of $X$, $Y$, $Z$ lie in $\h$. If $Z$ lies in
$\h$ and $X=E_\alpha$, $Y=E_\beta$ then the left hand side is
$[X,Y]\{(\alpha+\beta)(Z)-\alpha(Z)-\beta(Z)\}=0$. So we may suppose
that $X$, $Y$, $Z$ all are $E_\alpha$'s.
Let $\g_\alpha={\Bbb C}E_\alpha$. Note that by definition
$[\h,\g_\alpha]\subseteq \g_\alpha$ and $[\g_\alpha,\g_\beta]\subseteq
\g_{\alpha+\beta}$. Note also that if $\alpha,\beta,\alpha+\beta\in
\Phi$, then $(\alpha,\beta)=-1$ and so
$\varepsilon(\alpha,\beta)=-\varepsilon(\beta,\alpha)$. 

Suppose then that  $X=E_\alpha$, $Y=E_\beta$ and $Z=E_\gamma$. If
$\alpha+\beta+\gamma\notin \Phi\cup\{0\}$, then each term in the 
Jacobi identity must be zero. So either $\alpha+\beta+\gamma=0$ or
$\alpha+\beta+\gamma\in \Phi$. If $\alpha+\beta+\gamma=0$, then
$$\varepsilon(\alpha,\beta)
=\varepsilon(\beta,\gamma)=\varepsilon(\gamma,\alpha).$$
By symmetry only one of these equalities needs to be proved; the first
holds because $\gamma =-\alpha-\beta$ and $\varepsilon(\alpha,\alpha)=-1$,
$\varepsilon(\beta,\alpha)=- \varepsilon(\alpha,\beta)$.
Hence the left hand side of the Jacobi inequality is proportional to
$$[E_{\alpha+\beta},E_\gamma] + [E_{\beta+\gamma},E_\alpha] +
[E_{\gamma+\alpha} ,E_\beta]=[E_{-\gamma},E_\gamma] +
[E_{-\alpha},E_\alpha] +[E_{-\beta},E_\beta]=\gamma +\alpha +
\beta=0.$$
Now suppose that $\delta=\alpha+\beta+\gamma\in \Phi$. Expanding
$\|\delta\|^2=2$, we get
$$(\alpha,\beta)+(\beta,\gamma) + (\gamma,\alpha)=-2.$$
Let $a=(\beta,\gamma)$, $b=(\gamma,\alpha)$ and
$c=(\alpha,\beta)$. Any permutation of $\alpha,\beta,\gamma$ results
in a distinct permuation of $a,b,c$. We have $a+b+c=-2$ and $-2\le
a,b,c\le 2$. We may therefore assume that $c\ge 0$. If $c=0$, we must
have (after permuting) $a=-1$ and $b=-1$ or $a=0$ and $b=-2$. If
$c=1$, we must have (after permuting) $a=-2$ and $b=-1$. If $c=-2$,
we must have $a=-2=b$. Thus there are four possibilities:
\vskip .05in
\noindent (1) $(a,b,c)=(-1,-1,0)$. 
Thus $(\alpha,\beta)=0$, $\beta+\gamma,\alpha+\gamma\in
\Phi$. Thus $\alpha+\beta\notin \Phi\cup \{0\}$. The left hand side of
the Jacobi identity becomes
$$ 0 + \varepsilon(\beta,\gamma) [E_{\beta+\gamma},E_\alpha] -
\varepsilon(\gamma,\alpha) [E_{\gamma+\alpha},E_{\beta}] 
=[\varepsilon(\beta,\gamma)
\varepsilon(\beta+\gamma,\alpha)-\varepsilon(\gamma,\alpha)
\varepsilon(\beta,\gamma+\alpha)]E_\delta=0.$$

\noindent (2) $(a,b,c)=(-1,-2,1)$. Thus
$(\alpha,\beta)=1$ so that $\alpha+\beta\notin \Phi\cup\{0\}$,
$\gamma=-\alpha$ and $\beta+\gamma\in \Phi$. The left hand side of the
Jacobi identity is
$$\varepsilon(\beta,\gamma) [E_{\beta+\gamma}, E_\alpha]
+(\alpha,\beta) E_\beta
=-\varepsilon(\beta,\gamma) \varepsilon(\alpha,\beta+\gamma) E_\beta -
E\beta= (-(-1)^{(\beta,\gamma)} +1) E_\beta=0.$$

\noindent (3) $(a,b,c)=(0,-2,0)$. In this case
$\alpha,\gamma\perp\beta$ and $\alpha=-\gamma$. Thus $\beta\pm
\alpha\notin\Delta\cup\{0\}$. So two terms vanish and the remaining
term vanishes because $[\alpha,E_\beta]=(\alpha,\beta)E_\beta=0$.

\vskip .05in
\noindent (4) $(a,b,c)=(-2,-2,2)$. In this case
$\beta=\alpha=-\gamma$. The left hand side of the Jacobi identity is
trivially zero by skew symmetry of the bracket. This completes the proof.
\vskip .1in
\noindent \bf Proposition. \it Define a complex inner product on $\g$
by extending the real inner product on $V$ to a complex inner product
on $\h=V+iV$ and then decreeing the $E_\alpha$'s to be orthonormal and
orthogonal to $\h$. Define a conjugate--linear map $X\mapsto X^*$ on $\g$
by $\alpha^*=-\alpha$ and $E_\alpha^*=-E_{\-alpha}$. Then
$[X,Y]^*=[Y^*,X^*]$ for $X,Y\in \g$, $\overline{(X,Y)}=(X^*,Y^*)$ and
${\rm ad}(X^*)={\rm 
ad}(X)^*$ for $X\in \g$.
\vskip .1in
\noindent \bf Proof. \rm This is a routine verification. 
\vskip .1in
\noindent \bf Corollary. \rm Let $\g_0=\{X\in \g: X^*=-X\}$. Then
$\g_0$ is compact Lie algebra with invariant real inner product
$(X,Y)$.
\vskip .1in
\noindent \bf Proof. \rm Clearly $\g_0$ is closed under bracket and
real scalar multiplication. Since $\overline{(X,Y)}=(X,Y)$ for $X,Y\in
\g_0$, it follows that the inner product is real on $\g_0$. It is
invariant since ${\rm ad}(X)^*=-{\rm ad}(X)$ for $X\in \g_0$.
\vskip .1in
In particular it follows that $\g$ is the direct sum of its centre and
a set of simple algebras. We now determine the centre and each of the
simple summands.
\vskip .1in
\noindent \bf Proposition. \it The centre $\z$ of $\g$ is contained in $\h$
and equals $\Phi^\perp\subset \h$. Thus $\g$ has no centre iff $\Phi$
spans $\h$. 
\vskip .1in
\noindent \bf Proof. \rm The adjoint action of $\h$ on $\g$ is
diagonal: the eigenspace decomposition is $\h \oplus \bigoplus
\g_\alpha$, with $\h$ the $0$--eigenspace. Thus if $X$ is central, it
must lie in the $0$--eigenspace of $\h$, i.e.~$\h$. But then we need
$0=[X,E_\alpha]=\alpha(X)E_\alpha$ for all $\alpha\in \Phi$. This
happens iff $X\perp \Phi$ as required.
\vskip .1in
Define $\alpha,\beta\in \Phi$ to be adjacent if $(\alpha,\beta)\ne
0$. Define $\alpha,\beta\in \Phi$ to be connected if there is a chain
of adjacent elements of $\Phi$ linking $\alpha$ and $\beta$. This
celarly determines an equivalence relation on $\Phi$. Let
$V_0=\Phi^\perp$. Let these equivalence classes be
$\Phi_1,\dots,\Phi_m$, let $V_i$ be the real--linear span of $\Phi_i$
and let $\Lambda_i$ be the ${\Bbb Z}$--linear span of $\Phi_i$.
Thus $V=V_0\oplus V_1\oplus \cdots \oplus V_m$ is an orthogonal direct
sum. 
\vskip .1in
\noindent \bf Proposition. \rm  Let $\g_i=V_i\oplus i V_i \oplus
\bigoplus_{\alpha\in \Phi_i} {\Bbb C} E_\alpha$. Then $\g_i$ is a
simple non--Abelian Lie algebra and is an ideal in $\g$, Moreover
$\g=\z \oplus \g_1\oplus \cdots \oplus \g_m$. These ideals are
invariant under $*$ and mutually orthogonal.
\vskip .1in
\noindent \bf Proof. \rm Clearly each $\g_i$ is the Lie algebra
constructed from the lattice $\Lambda_i$ in $V_i$. Thus to prove the
first part we must show, if $\Phi$ spans $V$ and any two elements of
$\Phi$ are connected, that $\g$ is simple. Let $\a$ be an ideal in
$\g$. Note the following:

\vskip .05in
\item{(1)} If $E_\alpha\in \a$, then $\alpha\in \a$ (since
$[E_\alpha,E_{-alpha}]=-\alpha$).

\item{(2)} If $\alpha\in \a$, then $E_\alpha\in \a$ (since
$[\h,E_\alpha]={\Bbb C}E_\alpha$).
\vskip .05in
We claim that $E_\alpha\in \a$ for some $\alpha$. Suppose not. Since
$\a$ is invariant under ${\rm ad}(\h)$ so can be decomposed
into eigenspaces. If no $E_\alpha$ lies in $\a$, $\a$ is the zero
eigenspace, so that $\a\in \h$. But then $\alpha(\a)\ne 0$ for some
$\alpha$, so that $[\a,E_\alpha]={\Bbb C}E_\alpha$. Hence $E_\alpha\in
\a$ a contradiction. 

Since $E_\alpha\in \a$, so is $\alpha$. We claim that if $E_\beta$ is
in $\a$ and $\beta$ is adjacent to $\gamma$, then $E_\gamma$ is in
$\a$. In fact $E_\beta\in \a$, so $\beta\in \a$. But
$[\beta,E_\gamma]= (\beta,\gamma) E_\gamma$, with $(\beta,\gamma)\ne
0$. Thus $E_\gamma$ lies in $\a$. Since all elements of $\Phi$ are
connected to $\alpha$, it follows that $E_\beta\in \a$ for all
$\beta\in \Phi$ and hence $\beta\in \a$ for all $\beta\in \Phi$. Thus
$\a=\g$ and $\g$ is therefore simple.

\vskip .1in
\noindent \bf Corollary. \it $\g$ is simple iff $\Phi$ spans $V$ and
any two elements of $\Phi$ are connected.

\vskip .1in
This construction gives all the simple algebras of type $A,D,E$
(the so--called simply laced algebras). The remaining simple algebras
arise as fixed point algebras of lattice automorphisms of these
algebras: any automorphism of the
lattice $\Lambda$,
preserving the inner product and the $\varepsilon$--factor,
canonically induces an automorphism of the Lie algebra constructed
above. (See the exercises.)
\vskip .1in
\noindent \bf 3. MAXIMAL TORI. \rm Let $G$ be a compact matrix
group (not necessarily semisimple). A torus $T$ in $G$ is a closed
connected Abelian subgroup. Thus 
if $\t$ is the Lie algebra of $T$, we have $T=\exp(\t)$ and thus
$T\cong \t/\Lambda$ where $\Lambda={\rm ker}(\exp)$ is a lattice in
$\t$. By Kronecker's theorem, $T$
is generated toplogically by a single element $t\in T$ (called a
toplogical generator). [In fact if $x\in T ={\Bbb R}^n/{\Bbb Z}^n $
satisfies ${\bf e}_m(x)=e^{2\pi i m\cdot x}\ne 1$ for all non--zero $m\in {\Bbb
Z}^n$, then $x$ is a topological generator. Indeed let $H$ be the
closed subgroup generated by $x$. Then $H/H^0$ is finitely generated by
$x$, so cyclic. Taking the appropriate power of $x$, we do not change the
hypotheses but now $H$ is connected. If $H\ne T$, then $\h\subset
\t={\Bbb R}^n$ and the kernel of $\exp$ is $\Gamma=\h\cap {\Bbb Z}^n$, a
lattice in $\h$. Since ${\Bbb Z}^n/\Gamma$ is finitely generated and
has a free part, it has a (non--trivial) homomorphism onto ${\Bbb Z}$. Hence
there is a homomorphism $f$ of ${\Bbb Z}^n$ onto ${\Bbb Z}$ with $\Gamma$
in its kernel. Necessarily $f(x)=x\cdot m$ for some $m\in {\Bbb
Z}^n$. But then ${\bf e}_m=1$ on $H$, a contradiction.]

We say that $T$ is a maximal torus in $G$ if it
not properly contained in any other torus of $G$. Note that $T_1\subset T_2$
iff $\t_1 \subset \t_2$, so maximal tori always exist. 
\vskip .1in
\noindent \bf Lemma. \it $T$ is a maixmal torus iff $\t$ is a maximal
Abelian subalgebra in $\g$.
\vskip .1in
\noindent \bf Proof. \rm If $\t$ is maximal Abelian, $T$ cannot
properly be contained in another torus. If $\t$ is not maximal
Abelian, then $\t\subset \t_1$ with $\t_1$ Abelian. Then
$T^\prime=\overline{\exp \t_1}$ is connected, closed and Abelian, so a
torus, with $\t^\prime\supseteq \t_1\supset \t$, so that
$T^\prime\supset T$, so $T$ is not maximal Abelian.
\vskip .1in
\noindent \bf Theorem. \it If $T$ is a maximal torus in $G$ with Lie
algebra $\h$, then
$\g=\bigcup_{g\in G}ghg^{-1}$.
\vskip .1in
\noindent \bf Proof. \rm Take $X\in \g$. Choose $Y\in\h$ such that
$\exp Y$ is a topological generator of $T$. Thus the centraliser of
$Y$ in $\g$ is $\h$. Next choose $g\in G$ so
that $\|gX g^{-1} -Y\|^2$ is minimised, since $G$ is
compact. Replacing $X$ by $gXg^{-1}$, we may assume this minimum
occurs for $g=1$. Looking at a small variation $\exp(A) X\exp(-A)$, we
must have $([X,A],Y)-(Y,[X,A])=0$ for all $A$. Hence $(A,[X,Y])=0$ for
all $A$, so that $[X,Y]=0$. Hence $X\in \h$, as required. 
\vskip .1in
\noindent \bf Theorem. \it If $T$ is a maximal torus in $G$, then
$G=\bigcup_{g\in G}gTg^{-1}$.
\vskip .1in
\noindent \bf Differential geometric remark. \rm
If we knew that $G=\exp(\g)$, this would 
follow immediately from the previous theorem. Surjectivity of the
exponential map can be proved by a geometric argument (the Hopf--Rinow
theorem).

\vskip .1in
\noindent \bf Proof. \rm We show that $B=\bigcup_{g\in G} gTg^{-1}$ is
open and closed in $G$. Since $G$ is connected, we must have $G=B$.
Now clearly $B$ is closed as the continuous image in $G$ of the compact set
$G\times T$ under the map $(g,t)\mapsto gtg^{-1}$. So we need only
show it is open and for this it is enough to show that each $t\in T$
is an interior point. Let $A=C_G(t)^0$. We consider two extreme cases:
$A=G$ and $A=T$. In the first case $t=\exp X$ is central and so, if
$Y\in \g$, 
$\exp(X+Y)$ lies in $\bigcup g\exp( \h
)g^{-1}$, since $X\in \h$ and $\exp(X)$ is central. Thus $B$ contains
an open neighbourhood of $t$. In the second case, consider the map
$f;h\oplus \h^\perp\rightarrow B$ given by
$f(X,Y)=\exp(Y)t\exp(X)\exp(-Y)=t\exp(t^{-1}Yt)\exp (X) \exp(Y)$. The
derivative of this map at $(0,0)$ is $f^\prime(X,Y)=X \oplus (Y
-t^{-1} Yt)$. Since $C_\g(t)=\h$, the map $Y\mapsto Y -t^{-1} Y t$ is
an automorphism of $\h^\perp$. Thus $f^\prime_{(0,0)}$ is an isomorphism
and $f$ is locally a diffeomorphism. This provides an open
neighbourhood of $t$ in $B$.

To handle the general case, we combine these two ideas. Let $t\in T$
and $A=C_G(s)^0$. Thus $T\subseteq A \subseteq G$. Note that $t$ is
central in the maximal torus $T$ of $A$. Thus $t$ is an interior point
of $\bigcup_{g\in A} gTg^{-1}$. Thus if $X\in \a$ is sufficiently
small, $t\exp(X)$ lies in $B$. Now take $Y\in \a^\perp$ and consider
the map $f(X,Y) =\exp(Y)t\exp(X)\exp(-Y)=t
\exp(t^{-1}Yt)\exp(X)\exp(-Y)$. Again $f_{(0,0)}(X,Y)=X\oplus (Y
-t^{-1} Yt)$, which is an isomorphism since $\a=C_\g(t)$. Thus $f$ is
locally a diffeomorphism at $(0,0)$ and therefore the image of an open
ball around $(0,0)$ provides an open neighbourhood of $t$ in $B$.
\vskip .1in

\noindent \bf Corollary. \it Every element of $G$ lies in a maximal
torus. In particular $G=\exp\, \g$. 
\vskip .1in
\noindent \bf Corollary. \it Any two maximal tori are conjugate.
\vskip .1in
\noindent \bf Proof. \rm Since a torus is topologically cyclic, one
must be contained in a conjugate of the other. Since a conjugate of a
maximal torus is a maximal torus, the result follows.

\vskip .1in
\noindent \bf Lemma (centralisers). \it (1) $x\in C_G(x)^o$ for
all $x\in G$.

\noindent (2) $C_G(X)$ is connected for all $X\in \g$.
(More generally
if $\a$ is an Abelian subalgebra of $\g$, then $C_G(\a)$ is connected.)
\vskip .1in
\noindent \bf Proof. \rm (1) Let $T$ be a maximal torus containing
$x$. Then $x\in T\subseteq C_G(x)^o$. 

\noindent (2) Let $A=\overline{\exp(\a)}$, a torus. Suppose
$x\in C_G(\a)$. Then $A\subset C_G(x)^o=H$, so $A$ is contained in a
maximal torus $T$ in $H$. But $x$ is central in $H$, so $x\in
T$. Hence $x\in T\subseteq C_G(\a)^o$.
\vskip .1in
\noindent \bf Remark. \rm Note that the more general statement in (2)
could also be proved inductively using the single element statemenr by
successively passing to centralisers in centralisers using a basis of
$\a$. 

\vskip .1in
\noindent \bf Corollary. \rm (a) $x\in T$ is contained in exactly one
maximal torus iff $C_G(x)^o$ is a maximal torus iff $C_\g(x)$ is
maximal Abelian.

\noindent (b) $C_\g(X)=\h$ iff $C_G(X)=T$.

\noindent (c) A maximal torus is maximal Abelian (but not conversely).

\vskip .1in
\bf \noindent Proof. \rm (a) Since $x$ is central in $C_G(x)^o$ and
$C_G(x)^o$ is the union of all maximal tori containing $x$
(it is the union of its maximal tori and they all contain $x$), the
result follows.  

\noindent (b) This follows because the Lie algebra of $C_G(X)$ is
$C_\g(X)$ and $C_G(X)$ is connected.

\noindent (c) Say $x$ commutes with $T$. Then $T\subseteq C_G(x)^o$ must
contain $x$, since $x$ is central in $C_G(x)^o$.

\vskip .1in

\bf \noindent Corollary. \it The Weyl group $W=N(T)/T$ is finite. 
\vskip .1in
\noindent \bf Proof. \rm By the previous corollary, the continuous map
$W\subset {\rm Aut}(T)=PGL_m({\Bbb Z})$ is injective. Since $W$ is
compact, its image is compact and discrete, so finite.
\vskip .1in
\bf \noindent Corollary. \it $t_1,t_2\in T$ are conjugate in $G$ iff
they are conjugate under the Weyl group $N(T)/T$. (Thus the space of
conjugacy classes $G/{\rm Ad}\, G$ is homeomorphic to $T/N(T)$. )

\vskip .1in
\noindent \bf Proof. \rm Let $H=C_G(t_2)^o$ and suppose that
$t_2=gt_1g^{-1}$. Thus $T, gTg^{-1}\subset H$. Since $T$  and
$gTg^{-1}$ are maximal
tori in $H$, we can find $h\in H$ such that $T=hgTg^{-1}h^{-1}$. 
Let $x=hg$. Then $t_2=gtg^{-1}=hgtg^{-1}h^{-1}=xtx^{-1}$ and $x\in
N(T)$. The last statement follows because the map $T/N(T)\rightarrow G/{\rm
Ad}\, G$ is a continuous bijection between compact spaces.
\vskip .1in

\noindent \bf 4. REPRESENTATIONS OF SU(2) AND $\sl(2)$. \rm Let
$G=SU(2)$, the group of all complex matrices $\pmatrix{\alpha & \beta \cr
-\overline{\beta} & \overline{\alpha}\cr}$ with
$|\alpha|^2+|\beta|^2=1$. Thus $G$ is a compact connected group,
homeomorphic to $S^3$ [so simply connected]. Since $g\in G$ iff ${\rm
det}(g)=1$ and $gg^*=I$, the Lie algebra $\su(2)$ of $G$ is given by
matrices $X$ such that ${\rm tr}(X)=0$ and $X+X^*=0$,
i.e.~skew--adjoint matrices with trace zero. We take as a real basis
of $\su(2)$, $X,Y,T$ with 
$T=\pmatrix{i/2 & 0\cr 0 & -i/2}$, $X=\pmatrix{0 & 1/2\cr -1/2 &
0\cr}$ and $Y=\pmatrix{0 & i/2\cr i/2 & 0\cr}$. These basis elements
are orthogonal with respect to the real inner product $(X,Y)={\rm
tr}(XY^*)$ and satisfy the following relations 
$[X,Y]=T$, $[T,X]=Y$, $[Y,T]=X$. Define the complexification of
$\su(2)$ in $M_2({\Bbb C})$ as $\sl(2)=\su(2)+i\su(2)$. Clearly
$\sl(2)=\{X\in M_2({\Bbb C}): {\rm tr}(X)=0\}$, a
3--dimensional complex Lie algebra. It is the Lie algebra of the
closed matrix group $SL(2,{\Bbb C})$. The natural complex basis (over
${\Bbb C})$ of $\sl(2)$ is  $E=\pmatrix{ 0 & 1\cr 0 & 0\cr}$, $F=\pmatrix{ 0
& 0\cr 1 & 0\cr}$ and 
$H=\pmatrix{1 & 0\cr 0 & -1\cr}$. These are related to the real basis by
$E=X-iY$, $F=-X-iY$ and $H=-2iT$ (so that $2T=iH$). Thus 
$H^*=H$, $E^*=F$ and 
$$[E,F]=H,\quad [H,E]=2E,\quad [H,F]=-2F.\eqno{(*)}$$
The important point about the complex basis is that $E$ and $F$ become
raising and lowering operators in any finite--dimensional
representation. (We will encounter this phenomenon again when we
consider representations of the Heisenberg algebra and pass from the
real basis $P=x$, $Q=id/dx$ to complex basis $P\pm iQ$. The complex
operators are again raising and lowering operators, usually called
creation and annihilation operators.)

If $V={\Bbb C}^2$, the groups $SU(2)$ and $SL(2,{\Bbb C})$ act on
the tensor power $V^{\otimes n}$ by $T(g)=g^{\otimes n}$. The
corresponding action of the Lie 
algebras is $X\mapsto T(X)=X\otimes I\otimes \cdots \otimes I +
\cdots I \otimes I\otimes \cdots \otimes X$ to be compatible with the
exponential map, since $\exp T(X) =\exp X\otimes \exp X\otimes \cdots
\otimes \exp X$. The action of $SU(2)$ on $V^{\otimes n}$ is unitary,
so completely reducible. The actions of $SL(2,{\Bbb C})$, $\sl_2$ and
$\su_2$ are also by operators invariant under taking adjoints, so are
completely reducible. One of our aims is to classify the irreducible
representations of $SU(2)$ that arise in $V^{\otimes n}$. Weyl's
unitarian trick shows that we only need look at representations of $\sl(2)$.
\vskip .1in
\noindent \bf Lemma. \it $SU(2)$, $SL(2,{\Bbb C}$, $\su_2$ and $\sl_2$
have the same centraliser algebra on $V^{\otimes n}$ and hence the
same invariant subspaces.
\vskip .05in
\noindent \bf Proof. \rm Any $g\in SL(2,{\Bbb C})$ has a polar
decomposition $g=up$ where 
$p=(g^*g)^{1/2}$ and $u=gp^{-1}\in SU(2)$.  The unitary $u$ can be
written as $\exp(x)$ with $x\in \su_2$ and $p$ can be written as
$\exp(iy)$ with $y\in \su_2$. Note that  $T$
commutes with $\exp(tA)$ for all $t\in {\Bbb R}$ iff $T$ commutes $A$
(take the coefficient of $t$ in $T\exp(tA)=\exp(tA)T$). 
Thus $T\in {\rm End}\,V^{\otimes n}$ commutes with $SL(2,{\Bbb C})$
iff $T$ commutes with $T(A)$ for every $A\in \sl_2$ iff $T$ commutes
with $T(A)$ for every $A\in \su_2$ (by complex 
linearity) iff $T$ commutes with $SU(2)$. 
\vskip .1in
We now tackle the problem of classifying finite--dimensional
irreducible representations of $\sl(2)$. Thus we have operators $E$,
$F$ and $H$ on $V$ satisfying $(*)$. We shall temporarily abandon the adjoint
conditions, retaining only the property that {\it $H$ is
diagonalisable}. It is easy to check the following commutation relations.
\vskip .1in
\noindent \bf Lemma. 
\noindent (a) $[E^n,F]=nE^{n-1} (H+n-1)=n(H-n+1)E^{n-1}$.

\noindent (b) $[F^n,E]=-nF^{n-1}(H-n+1)=-n(H+n-1)F^{n-1}$.

\noindent (c) $[H,E^n]=2n E^n$ and $[H,F^n]=-2n F^n$.

\vskip .1in
\rm An eigenvector of $H$ is called
a {\it weight vector} and the eigenspaces {\it weight spaces}. Thus if
$Hv=\lambda 
v$, $v$ is a weight vector with weight $\lambda$. Note that
$HEv=(\lambda+2)Ev$ and $HFv=(\lambda-2)Fv$. Thus $E$ increases the
weight by $2$ and $F$ decreases the weight by $2$. For this reason $E$
and $F$ are called raising and lowering operators.
\vskip .1in 

\noindent \bf Lemma. \it Let $V$ be an $\sl(2)$--module and let $v\in
V$ satisfy $Hv=\lambda v$ and $Ev=0$. Let $v_j=(j!)^{-1} F_j v$. Then
$Hv_j=(\lambda-2j)v$ and $Ev_j =(\lambda-j+1) v_{j-1}$.
\vskip .1in
\noindent \bf Proof. \rm Immediate from previous lemma. 
\vskip .1in
\noindent \bf Theorem. \it The irreducible finite--dimensional
representations of $\sl(2)$ 
are classified by their highest weight, a non--negative integer $d$.
The representation $V_d$ has dimension $d+1$ and has a unique highest
weight vector $v$ (up to a scalar multiple). If $v_0=v$ and $v_j =
(j!)^{-1} F^j v$ for $j=0,\dots,d$, then the $v_j$'s form a basis of
$V_d$ and 
$$H\cdot v_j=(d-2j)v_j,\quad F\cdot v_j=(j+1)v_{j+1}, \quad E\cdot
v_j=(d+1-j) v_j. \eqno{(**)}$$
\vskip .1in
\noindent \bf Proof. \rm  Let $v=v_0$ be a vector
of highest weight. Thus $Hv=\lambda v$ and $Ev=0$. The vector
$v_k=(k!)^{-1}F^kv$ 
has weight $\lambda -2k$, all distinct, so by finite--dimensionality
$F^{d+1}v=0$  
for some smallest $d\ge 0$. Since $F^{d+1}v_0=0$, we must have $v_{d+1}=0$.
But by the lemma, $Ev_{d+1}=(\lambda-d)v_d$. Since $v_d\ne 0$, we get
$\lambda=d$. On the other hand it is easy to verify directly that
$(**)$ defines a representation of $\sl(2)$ on ${\Bbb
C}^{d+1}=\bigoplus {\Bbb C}v_i$. It is irreducible,
because if $U$ is an invariant subspace, it must be a sum of
eigenspaces of $H$ and hence contain some eigenvector. Applying raising
and lowering operators we see that all basis vectors lie in $U$.
\vskip .1in
\noindent \bf Adjoint conditions. \rm If one puts in the
self--adjointness conditions $E^*=F$ and $H=H^*$, one can give a
``no--ghost'' argument for $\lambda=d$:
\vskip .1in
\noindent \bf Lemma. \it Let $E,F,H$ be operators on an inner product
space $V$ with $E^*=F$, $H^*=H$ satisfying $(*)$. If $v\in V$
satsifies $Ev=0$ and $Hv=\lambda v$, then $\lambda$ must be a
non--negative integer. 
\vskip .1in
\noindent \bf Proof. \rm By induction on $k$, we have
$[E,F^{k+1}]=(k+1)F^k (H-kI)$ for $k\ge 0$. Hence
$$(F^{k+1}v,F^{k+1} v)=(F^*F^{k+1}v, F^k v)=(EF^{k+1}v,
F^kv)=(k+1)(\lambda-k) (F^kv,F^kv).$$
For these norms to be non--negative for all $k\ge 0$, $\lambda$ has to
be a non--negative integer. 
\vskip .1in
\noindent \bf Character of a representation. \rm  The representation
$V_d$ coincides with $S^d V$ since they have the 
same highest weight and dimension. (It can also be seen directly that
$S^d V$ is irreducible, because this is the linear action on
two variable polynomials of degree $d$.) In particular every
irreducible finite--dimensional representation of $\sl(2)$ comes from
a representation of $SU(2)$ (and even $SL(2,{\Bbb C})$). The character of a
representation $\pi$ of $G$ is given by $\chi(g)={\rm tr}(\pi(g))$. It
is invariant under conjugation, since the trace is. On the other hand
every element of $SU(2)$ is conjugate to a diagonal matrix
$\pmatrix{\zeta & 0 \cr 0 & \overline{\zeta}\cr}$, so it is enough to
know the character on the diagonal element.
From the theorem the character of $V_d$ on 
$\pmatrix{\zeta & 0 \cr 0 & \overline{\zeta}\cr}$ is
$\chi_d(\zeta)=(\zeta^d -\zeta^{-d})/\zeta-\zeta^{-1}$. It follows
that every completely reducible representation is completely specified
by its character. In particular this applies to all representations
arising as subrepresentations of $V^{\otimes n}$ and hence tensor
products of $V_d$'s. By multiplying and expanding the characters we
get the celebrated 
Clebsch--Gordan rules: $$V_r\bigotimes V_s \cong V_{|r-s|} \bigoplus
V_{|r-s|+2} \bigoplus \cdots \bigoplus V_{r+s}.$$
\vskip .1in
\noindent \bf The Casimir Operator. \rm \noindent The Casimir element
is defined as $C=H^2 + 
2(EF+FE)=H^2+2H+4FE$ on any representation. (Note that $C=-2(X^2+Y^2+T^2)$.)
The commutation relations imply that $C$ commutes with $H,E,F$ and
hence with $\sl(2)$. 
By Schur's lemma, the Casimir is therefore a constant on $V_d$ and the
constant can be computed by applying $C$ to $v_0$. We get $Cv=(d^2
+2d)v_0$. Thus the Casimir distinguishes irreducible representations.

\vskip .1in
\noindent \bf Complete Reducibility Theorem. \it Every
finite--dimensional representation $W$ of $\sl(2)$ is completely reducible.
\vskip .05in
\noindent \bf Remark. \rm If we considered
representations on inner 
product spaces
satsifying $E^*=F$, $H^*=H$, this would be immediate; this applies to
all the examples 
above as well as the famous Hodge theory action on hermitian exterior
algebras (see Wells, for example). The theorem implies that every
finite--dimensional representation extends to $SU(2)$ and $SL(2,{\Bbb
C}$ and particular has an invariant inner product.
\vskip .05in
\noindent \bf Proof (van der Waerden--Casimir). \rm We may assume $C$
has only one eigenvalue on 
$W$. Let $W_1$ be a subspace of $W$ that is a direct sum of
irreducibles of maximal possible dimension. If $W_1\ne W$, find an
irreducible subspace $\overline{V}$ of $W/W_1$. Because of the
eigenvalue assumption on $C$, $\overline{V}$ and the irreducible
summands of $W_1$ are all isomorphic to $V_d$ some $d$. Suppose
$\overline{V}=V/W_1$. Then $E^{d+1}\overline{V}=(0)$ so that
$E^{d+1}V\subset W_1$. On the other hand $E^{d+1} W_1=(0)$. Hence $E$
is nilpotent on $V$, say $E^{k+1}=0$, $E^k\ne 0$ on $V$. Thus $k\ge d$. Since
$[E^{k+1},F]=(k+1)(H-k)E^n$, $k$ must be an eigenvalue of $H$ so that
$k\le d$. Hence $k=d$. Now choose $v\in V$ such that
$\overline{v}$ is a highest weight vector in $\overline{V}$. Let
$u=F^d v$ and $v^\prime =E^d u$. Thus $\overline{v^\prime}$ is a
non--zero multiple of $\overline{v}$. On the other hand $Ev^\prime=0$ and
we have just seen that $v^\prime=E^du$ is an eigenvector of $H$. Thus
$v^\prime$ generates a copy of 
$V_d$ not contained in $W_1$, a contradiction. Hence $W$ is completely
reducible.

\vskip .1in

\noindent \bf 5. THE ROOT SYSTEM. \rm Let $G$ be a compact simple (or
semisimple) matrix group with maximal torus $T$. Let $\g$ and $\h$ be
the corresponding Lie algebras. We may write $\g=\h\oplus \m$ where
$\m=\h^\perp$. Since the inner product is ${\rm Ad}$--invariant, $\m$
is invariant under ${\rm Ad}(T)$. It is a real inner product space on
which $T$ acts without fixed vectors (because $\h$ is maimal
Abelian). Hence $\m$ is even--dimensional the orthogonal matrices
${\rm Ad}(t)$ can simultaneously be put in canonical form. Thus there
is an orthonormal basis of $\m$ in which ${\rm Ad}(e^T)$ is block
diagonal with blocks $\pmatrix{\cos\alpha_i(T) & \sin \alpha_i(T)\cr
-\sin\alpha_i(T) & \cos\alpha_i(T)\cr}$ down the diagonal. 
[Alternatively
${\rm Ad}(T)$ is a torus in $SO(\m)$ so contained in a maximal torus.]
The eigenvalues of this matrix are $e^{\pm i \alpha_i(T)}$ and the linear
function $\alpha(T)=\pm \alpha_i(T)$ are called roots. Thus each root
$\alpha$ lies in $\h^*$. 
\vskip .1in
For each root $\alpha$, we can find
orthogonal unit vectors $x$, $y$ such that 
$[t,x]=\alpha(t) y$, $[t,y]=-\alpha(t)x$. Consider $[x,y]$. Then
$[t,[x,y]]=[[t,x],y]+[x,[t,y]]=0$ for $t\in \h$. Hence $[x,y]$ lies in
$\h$. To calculate which element of $\h$ it is, consider
$([x,y],t)=-(y,[x,t])=(y,y)\alpha(t)$. Thus $[x,y]=T_\alpha$, where 
$T_\alpha\in \h$ is the element corresponding to $\alpha\in \h^*$,
i.e.~$(t,T_\alpha)=\alpha(t)$. Clearly
$\alpha(T_\alpha)=\|\alpha\|^2$. Define $X=x/\|\alpha\|$,
$Y=y/\|\alpha\|$ and $T=T_\alpha/\|\alpha\|^2$. Thus
$[X,Y]=T$, $[T,X]=Y$ and $[T,Y]=-X$. We can then form the elements
$E$, $F$ and $H$  in the complexification $\g_{\Bbb C}$ as
above. Clearly $[t,E]=i\alpha(t)E$ and 
$[t,F]=-i\alpha(t)F$ for all $t\in \h$ and the elements $E$ and $F$
are orthogonal. These lie in the subspaces $\g_\alpha$ and
$\g_{-\alpha}$. We call $E,F,H$ the copy of $\sl(2)$ in $\g_{\Bbb
C}$. Thus $H_\alpha=-2iT_\alpha/\|\alpha\|^2$ and
$E_\alpha^*=F_\alpha$. Since $\su(2)$ has a unique invariant norm
with $\|H\|^2=2$ and $\|H_\alpha\|^2=4/\|\alpha\|^2$, we get
$\|E_\alpha\|^2 =\|F_\alpha\|^2=\|H_\alpha\|^2/2=2/\|\alpha\|^2$. 
\vskip .1in
\noindent \bf Remark. \rm This construction also makes sense at the
level of groups. If $G\subset U(V)$ and $X,Y,T \in \g$ satisfy
$[X,Y]=T$, $[T,X]=Y$, $[Y,T]=X$ and span a Lie subalgebra $\s\subset
\g$. Since these operators are 
skew--adjoint, $V$ breaks up as a direct sum of irreducible
representations of $\s$. By the $SU(2)$ theory there is a
representation $\pi$ of
$SU(2)$ on $V$ such that its generators $X_1,Y_1,T_1$ are sent to
$X,Y,T$ under $\pi$. Since $\pi(\exp X)=\exp \pi(X)$, it follows that
$\pi(SU(2))\subseteq G$. The image is a closed connected subgroup of
$G$, so a matrix group in its own right. The kernel $Z$ of $\pi$ is a
closed subgroup of $SU(2)$ so a matrix group. Since $\pi$ is
injective on the Lie algebra of $SU(2)$, $Z$ is discrete so
finite. It is normal in $SU(2)$ so $SU(2)$ acts by conjugation on
$Z$. But $Z$ is discrete and $SU(2)$ connected. Hence $Z$ is central,
so that $Z=\{1\}$ or $\{\pm 1\}$. We call this the copy of $SU(2)$ or
$SO(3)$ in $\g$ corresponding to $\s$. For this reason we can loosely
talk about the ``copy of $SU(2)$'' in $G$ correpsonding to a given
root $\alpha$.

\vskip .1in

The above arguments could also have been carried out
directly in the 
complexification $\g_{\Bbb C}$; in fact we get a useful generalisation
which cannot be seen so clearly just working in $\g$.
\vskip .1in
\noindent \bf Lemma. \it Suppose that $E\in \g_\alpha$ and $F\in
\g_{-\alpha}$. Then $[E,F]=-i(F,E^*)T_\alpha=-i(E,F^*)T_\alpha$.
\vskip .1in
\noindent \bf Proof. \rm Clearly $[E,F]$ commutes with $\h$, so lies
in $\h_{\Bbb C}$. We have $([E,F],t)=(F,[E^*,t]) =(F,i\alpha(t)E^*)= 
-i\alpha(t) (F,E^*)=(-i(F,E^*)T_\alpha,t)$. 
\vskip .1in
\bf \noindent Corollary. \it The root $\alpha$ occurs with multiplicity
one in $\g$, so that ${\rm dim}\g_\alpha =1$. 
\vskip .1in
\noindent \bf Proof. \rm If not, we can find further
$E^\prime,F^prime$ orthogonal to $E,F$ but with the same relations
with $H$, But then $[E,F^\prime]=0$, since $(F,E^*)=(F,F^\prime)=0$
and $[H,F^\prime]=-2F^\prime$. This contradicts the $s\ell(2)$ lemma
applied to $V=\g_{\Bbb C}$, with the adjoint representation of
$E,F,H$, and $v=F^\prime$. 
\vskip .1in

\noindent \bf Lemma. \it If $\alpha\ne \pm \beta$, then
$[\g_\alpha,g_\beta]=\g_{\alpha+\beta}$ if $\alpha+\beta$ is a root
and $(0)$ otherwise.

\vskip .1in
\noindent \bf Proof. \rm Take the copy of $s\ell(2)$, $E,F,H$,
corresponding to the root $\alpha$ and let $V=\oplus_{m\in {\Bbb Z}}
\g_{\beta +m\alpha}$. Then $V$ is invariant under $E,F,H$, so can be
written as a direct of irreducible submodules. On the other hand each
weight space of $V$ is at most one dimensional, so that $V$ must
actually be irreducible. The result follows
immediately, because ${\rm ad}(E)$ is a raising operator so is an
isomorphism between weight spaces of $V$. 
\vskip .1in
\noindent \bf Corollary of proof. \it
$n(\alpha,\beta)=2(\alpha,\beta)/(\alpha,\alpha)$ is an integer.
The roots of the form $\beta +m\alpha$ are exactly those with
$m\in [-p,q]$ where $-p\le 0\le q$ and $p-q=n(\alpha,\beta)$. In
particular $\beta-n(\alpha,\beta)\alpha$ is always a root. 
\vskip .1in
\noindent \bf Proof. \rm The first assertion follows because
$n(\alpha,\beta)=\beta(H)$ and ${\rm ad} \, H$ has only integer
eigenvalues. The corresponding irreducible representation of $E,F,H$
has lowest $H$--eigenvalue $n(\alpha,\beta)-2p$ and highest
$H$--eigenvalue $n(\alpha,beta)+2q$. These must be negatives of each
other, so that $n(\alpha,\beta)=p-q$. The last assertion follows
because $q-p\in [q,-p]$.

\vskip .1in
\bf \noindent Proposition. \it The root system $\Phi\in \h^*=V$ has the
following properties.
\vskip .05in
\item{R1} $\Phi$ spans $V$.
\item{R2} If $\alpha\in \Phi$, then $\sigma_\alpha\Phi=\Phi$ where
$\sigma_\alpha$ is the reflection
$\sigma_\alpha(v)=v-(v,\alpha^\vee)/\alpha$ with
$\alpha^\vee=2\alpha/(\alpha,\alpha)$. 
\item{R3} $2(\alpha,\beta)/(\alpha,\alpha)\in {\Bbb Z}$ for all
$\alpha,\beta\in \Phi$.
\item{R4} If $\alpha\in \Phi$, then the only roots proportional to
$\alpha$ are $\pm \alpha$. 
\vskip .1in
\noindent \bf Proof. \rm (R1) If not, we could find a non--zero $t\in
\h$ such 
that $\alpha(t)=0$ for all $\alpha$. But then $t$ would be central, a
contradiction. (R2) and (R3) were proved in the preceding lemmas.
(R4) If $\beta=s \alpha$ is a root, then $2s=n(\alpha,\beta)\in {\Bbb
Z}$. So that $2s\in {\Bbb Z}$. Since $\alpha=s^{-1}\beta$, we
similarly have $2s^{-1}\in \Bbb Z$. Thus we may assume without loss of
generality that $s=2$ and that $3\alpha$ is not a root. 
Then $V=\g_{2\alpha} \oplus \g_\alpha \oplus {\Bbb C}H \oplus
\g_{-\alpha} \oplus \g_{-alpha}$ would give a 5--dimensional
representation of $E,F,H$, necessarily irreducible. However $\g_\alpha
\oplus {\Bbb C}H \oplus \g_{-\alpha}$ is a subrepresentation, a
contradiction. 
\vskip .1in
\noindent \bf Lemma. \it If $\alpha\in \h^*$ is a root and $T_\alpha\in \h$
the corresponding element of $\h$, so that $(t,T_\alpha)=\alpha(t)$,
then reflection in the ray ${\Bbb R}T_\alpha$ lies in $W=N(T)/T$. 
\vskip .1in
\noindent \bf Proof. \rm Take $X,Y,T$ such that $[X,Y]=T$, $[T,X]=Y$
and $[Y,T]=X$. Thus $[t,X]=\alpha(t)Y$ and $[t,Y]=-\alpha(t)X$ for
$t\in \h$. Let $g_s=\exp(sX)$. Then 
$$g_stg_s^{-1} ={\rm Ad}(e^{sX})\cdot t =e^{s{\rm ad}(X)} t.$$
If $\alpha(t)=0$, we therefore have $g_stg_s^{-1}=t$. On the other hand
$g_sTg_s^{-1}= \cos(s) T +\sin(s) X$. Taking $g=g_\pi$, we get
$g\in N(T)$, with $gtg^{-1}=t$ if $t\perp T$ and $gTg^{-1}=-T$. Hence
the image of $g$ is the reflection in the ray ${\Bbb R} T$. 

\vskip .1in
\noindent \bf Example. \rm For $G=SU(N)$, we have $\g_{\Bbb C}=\sl(N)=\{X\in
M_N({\Bbb C})$. We can identify $\h$ with diagonal matrices $\{ix:
x\in {\Bbb R}^N,  \sum
x_i=0\}$. The roots vectors are the matrix units $e_{pq}$ with $p\ne q$. Since
$e^{ix} e_{pq} e^{-ix} =e^{i(x_p-x_q)} e_{pq}$, the corresponding root
is $\alpha_{pq}(x)=x_p-x_q$. 
\vskip .1in
\noindent \bf 6. THE WEYL GROUP AS A REFLECTION GROUP. \rm Our aim is
to obtain a description in terms of roots for a fundamental 
domain of the Weyl group. For each root $\alpha$, let $\H_\alpha$ be
the hyperplane $\{X:\alpha(X)=0\}$ and let $\H_\alpha^+$ be the
closed half--space $\{X:\alpha(X)\ge 0\}$. 
In $\h$ we define the Weyl chambers to be the connected components of
$$\h^\prime=\h\backslash \bigcup \H_\alpha=\{X:\alpha(X)\ne 0\,
(\alpha\in \Phi)\}.$$
Clearly these are open convex cones. The boundary of each
chamber $C$ is contained in $\bigcup \H_\alpha$. 
Let $W_0$ be the normal subgroup of $W$
generated by reflections in the $\H_\alpha$'s.

\vskip .1in
\noindent \bf Theorem. \it $W=W_0$ and $W$ permutes the Weyl chambers
simply transitively.
\vskip .05in

\noindent\bf Proof. \rm (1) \it $W$ permutes the Weyl chambers. \rm
This is clear because $W$
permutes $\Phi$ and therefore leaves $\h^\prime$ invariant. 

\noindent (2) \it $W_0$ permutes the the Weyl chambers transitively.
\rm Let $C_1$ and $C_2$ be two Weyl chambers. Fix $x\in C_1$
and consider the boundary sphere $S$ of a small ball in $C$ around
$x$. The chamber $C_2$ projects onto an open subset of $S$. Each
intersection of distinct hyperplanes $\H_\alpha\cap \H_\beta$ ($\alpha\ne \pm
\beta$) is a subspace of codimension 2 so projects onto a sphere of
codimension $1$ in $S$. There are only finitely many such spheres so
there is a point $y$ in $C_2$ such that the line segment joining $x$
and $y$ misses each double intersection and therefore has only simple
(or empty) intersections with each hyperplane $\H_\alpha$. Clearly the
composition of the reflection in each of the successive hyperplanes
encountered will carry $C_1$ onto $C_2$.

\noindent (3) \it $W=W_0$ and $W$ is simply transitive. \rm To see
this, let 
$W_C=\{\sigma\in W:\sigma 
C=C\}$ be the stabiliser of 
$C$. By (2) $W=W_0\cdot W_C$. Now take $x\in C$ and set $X=|W_C|^{-1}
\sum_{\sigma\in W_C} \sigma x$. Thus $X\in C$ is fixed by $W_C$. Since
$\alpha(X)\ne 0$ for all $\alpha\in \Phi$, $T=C_\g(X)=C_G(X)$. Hence
$W_C\subset C_G(X)/T=\{1\}$, so that $W_C=\{1\}$ and $W=W_0$. 

\vskip .1in
\bf \noindent Corollary. \rm $\overline{C}$ is a fundamental
domain for the Weyl group $W$.
\vskip .1in
\noindent \bf Proof. \rm Let $C_1$ and $C_2$ be Weyl chambers
with $X\in C_1$, $Y\in 
C_2$. If the line segment $[X,Y]$ crosses a hyperplane $\h_\alpha$,
then $\|X-Y\|> \|X-\sigma_\alpha(Y)\|$. Now minimise the distence
$X-\sigma(Y)\|$ over $\sigma\in W$. Any minimum $\sigma Y$ cannot be
separated from $X$ by any walls, so that $\sigma Y\in C_1$. 
Since there is a unique $\sigma$ such that $\sigma Y\in C_1$, it follows
that if $X,Y\in C$, then $\|X-\sigma Y\|\ge \|X- Y\|$ for all
$\sigma\in W$. This result also holds by continuity for $X,Y\in
\overline{C}$; a similar continuity argument shows that $\h=W\cdot
\overline{C}$. Now suppose 
that $X,\sigma X \in \overline{C}$. Let $Y=\sigma X$ and
$\tau=\sigma^{-1}$. Then $0\ge \|X-\tau Y\|\ge \|X-Y\|$. Hence
$X=Y$. Thus the $W$--orbit of any point intersects $\overline{C}$ in
just one point, so that $\overline{C}$ is a fundamental domain.

\vskip .1in
Note that if $C_1, C_2$ are two Weyl
chambers, then the number of hyperplanes intersecting the line
segment joining $x_1\in C_1$ and $x_2\in C_2$ is independent of
the choice of $x_i$. If $\Phi_i^+=\{\alpha: \alpha(x_i)>0\}$, it is
the number of roots in $\Phi_1^+$ lying in $-\Phi_2^+$; this is
because a sign change occurs whenever $x_1$ and $x_2$ lie on opposite
sides of $\H_\alpha$. Denote this
number by $n(C_1,C_2)$. If we fix a Weyl chamber $C$, we define
$n(\sigma) = n(C,\sigma C)$ for $\sigma\in W$. If $x\in C$, we have
$\Phi^+=\{\alpha: \alpha(x)>0\}$, so that
$n(\sigma)=|\{\alpha:\alpha(x)>0,\,\alpha(\sigma x)>0\}|=|\{\alpha>0:
\sigma^{-1}\alpha<0\}|$. 
\vskip .1in
\noindent \bf 9. GEOMETRIC APPROACH TO WEYL CHAMBERS AND SIMPLE
ROOTS. \rm Let $C$ be a fixed Weyl chamber. We call $\H_\alpha$ a
wall of the Weyl chamber $C$ if $\H_\alpha\cap \overline{C}$ has
non--empty interior in $\H_\alpha$. 
We define the simple roots corresponding to $C$ to be those such that
$\alpha(x)>0$ on $C$ and $\H_\alpha$ is a wall of $C$. The
corresponding reflections in the walls are called simple roots.
\vskip .1in
\noindent \bf Lemma. \it If $\sigma\in W$, then $\sigma$ is a product
of $n(\sigma)$ simple reflections. In fact if $\ell(sigma)$ is the
minimal number of simple reflections required for such a product,
$\ell(\sigma) = n(\sigma)$. In particular $W$ is generated by simple
reflections. 
\vskip .1in
\noindent \bf Proof. \rm We prove the result by induction on
$n(\sigma)$, the result being trivial for $n=0$.
Take a generic line segment joining $x\in C$
to $y\in \sigma^{-1} C$, crossing the hyperplanes $\H_{\beta_1}$,
\dots $\H_{\beta_\ell}$ transversely. Then $\sigma^{-1}C
=\sigma_{\beta_\ell} \cdots \sigma_{\beta_1} C$. By simple
transitivity,
$\sigma^{-1} =\sigma_{\beta_\ell} \cdots \sigma_{\beta_1}$. Thus
$\sigma =\sigma_{\beta_1} \cdots \sigma_{\beta_\ell}$. Thus
$n(C,\sigma C)=\ell$. Let
$\tau=\sigma_{\beta_2} \cdots \sigma_{\beta_\ell}$. Then $n(C,\tau
C)=\ell-1$, because of the properties of the line segment
$[x,y]$. So by induction $\tau$ is the product of $\ell-1$ simple
reflections. Since $\sigma=\sigma_{\beta_1}\tau$ and $\beta_1$ is
simple, we see that $\sigma$ is the product of $\ell$ simple
reflections as required. Thus $\ell(\sigma)\le n(\sigma)$.

We now prove that
$n(\sigma)=\ell(\sigma)$. Note that if $x\in C$ and $\h_\beta$ is a
wall of $C$, then $x$ and $\sigma_\beta x $ are only separated by the
hyperplane $\h_\beta$. Transporting structure by $\tau\in W$, we see
that if $\tau x\in \tau C$, $\tau x$ and 
$\sigma_{\tau \beta} x$ are only separated by the hyperplane
$\tau(\h_\beta)$. But $\sigma_{\tau\beta}=\tau
\sigma_\beta \tau^{-1}$, so that $\tau x$ and $\tau \sigma_\beta x$
are only separated by $\tau\h_{\beta}$. If we write
$\sigma=\sigma_1\dots \sigma_\ell$, a product of simple reflections,
it follows that there is only one 
hyperplane separating $x$ and $\sigma_1 x$, one separating $\sigma_1
x$ and $\sigma_1\sigma_2 x$, and so on. Thus there is a piecewise
linear path from $x$ to $\sigma x$ crossing only $\ell$
hyperplanes. This can only cross more hyperplanes than the
straightline joining $x$ and $\sigma x$, so that $\ell(\sigma)\ge
n(\sigma)$. Hence $\ell(\sigma)=n(\sigma)$.

\vskip .1in

\vskip .1in
\noindent \bf Proposition. \it Each Weyl chamber $C$ is the
intersection of the open half spaces corresponding to its walls,
i.e.~$C=\bigcap_{\alpha\in \Delta} \H_\alpha$. Moreover
$\overline{C}= \bigcap_{\alpha\in \Delta} \overline{\H_\alpha}$.
\vskip .1in
\noindent \bf Proof. \rm Suppose that we have $C=\bigcap_{\alpha\in
\Delta^\prime} \H^+\alpha$ where $\Delta^\prime\supset \Delta$. Then
we claim that $\overline{C}= \bigcap_{\alpha\in \Delta^\prime}
\overline{\H^+\alpha}$. Indeed the inclusion 
$\overline{C}\subseteq \bigcap_{\alpha\in \Delta^\prime}
\overline{\H^+\alpha}$ is clear. If on the other hand $\alpha(x)\ge $
for all $\alpha\in \Delta^\prime$ and $z\in C$, then $x_n=x+z/n\in C$
and $x_n\rightarrow x$. So equality holds. 

Now take $\Delta^\prime\supseteq \Delta$ minimal with
$C=\bigcap_{\alpha\in \Delta^\prime} \h^+_\alpha$. If
$\Delta^\prime\ne \Delta$, take $\beta\in \Delta^\prime\backslash
\Delta$. Let
$\Delta^{\prime\prime}=\Delta^\prime\backslash\{\beta\}$. 
We claim that $\H_\beta$ does not intersect 
$\bigcap_{\alpha\in \Delta^{\prime\prime}} \h^+_\alpha$. If not,
suppose they meet in $x$. Thus $\alpha(x)>0$ for $\alpha\in
\Delta^{\prime\prime}$ while $\beta(x)=0$. Thus $x\in
\bigcap_{\alpha\in \Delta^\prime} \overline{\h^+_\alpha}$ by our first
observations. But then $x\in \overline{C}\cap \H_\beta$. Since
$\alpha(x)>0$ for all $\alpha\ne \beta$, this will also be true in a
pneighbourhood of $x$ in $\H_\beta$. So $\H_\beta$ would have to be a
wall, a contradiction since by assumption $\beta\notin \Delta$. 

Since $ C\subset \bigcap_{\alpha\in \Delta^{\prime\prime}} \h^+_\alpha$
and the latter does not meet $\H_\beta$, they both must lie in
$\H_\beta^+$. But then $\bigcap_{\alpha\in \Delta^{\prime\prime}}
\h^+_\alpha\subset \H_\beta^+$, so that $C=\bigcap_{\alpha\in
\Delta^{\prime\prime}} \h^+_\alpha$. This contradicts the minimality of
$\Delta^\prime$. Hence $\Delta^\prime=\Delta$, so that
$C=\bigcap_{\alpha\in \Delta} \h^+_\alpha$. 
\vskip .1in

\noindent Theorem. \it The simple roots form a basis of $\h^*$. Every
positive root is a non--negative integral combination of simple roots.
The Weyl group orbit of any root contains a simple root; equivalently
every hyperplane $\h_\alpha$ is the wall of some Weyl chamber.

\vskip .1in
\noindent \bf Proof. \rm  (1) \it $\Delta$ spans $\h^*$. \rm 
If $x\in \h$, then $\alpha(x)>0$ iff
$\alpha_i(x)>0$. Hence $\alpha(x)\ge 0$ for all $\alpha>0$ iff
$\alpha_i(x)\ge 0$ for all $\alpha_i$. Similarly $\alpha(x)\le 0$ for
all $\alpha\le 0$ iff $ \alpha_i(x)\le 0$ for all $\alpha_i$. Hence
$\alpha(x)=0$ for all $\alpha$ iff $\alpha_i(x)=0$ for all
$\alpha_i$. Thus the $\alpha_i$'s span $\h^*$. 

\noindent (2) \it Each simple root $\alpha$ is non--redundant, i.e.~cannot be
written as $\alpha=\mu \beta +\nu\gamma$ with $\beta,\gamma\in \Phi^+$
non--proportional and $\mu,\nu\ge 0$. 
\rm If
$\alpha=\mu\beta+\nu\gamma$, then $y\in 
\overline{C}\cap \H_\alpha$ implies $\alpha(x)=0$ and
$\beta(x),\gamma(x)\ge 0$. Hence $\beta(x)=\gamma(x)=0$, so that 
$\overline{C}\cap \H_\alpha$ has at least codimension $1$ in $\H_\alpha$
so cannot have non--empty interior. Thus $\H_\alpha$ cannot be a wall.

\noindent (3) \it If $\alpha,\beta$ are simple, $(\alpha,\beta)\le 0$.
\rm
We show that $(\alpha,\beta)\le 0$ for non--redundant
roots. In fact we know $\sigma_\alpha-\beta=\gamma$ is a root where
$\gamma=\sigma_\alpha\beta=\beta-2\|\alpha\|^{-2}(\alpha,\beta)\beta$. If 
$(\alpha,\beta)>0$, then since $\gamma$ or $-\gamma$ is positive,
either $\beta$ or $\alpha$ would be non--redundant. Hence
$(\alpha,\beta)\le 0$. 

\noindent (4) \it The non--redundant roots form a basis of $\h^*$. 
\rm
 If $\Delta$ were not linearly independent, the
existence of a linear relation would yield a subset
$\Delta_0\subset\Delta$ and non--negative reals $c_\alpha$ such
that $\gamma=\sum_{\alpha\in \Delta_0} c_\alpha \alpha = \sum_{\beta\notin
\Delta_0} c_\beta \beta$. Then $\gamma(x)>0$ for $x\in C$ since not
all $c$'s are 
zero. But $\|\gamma\|^2=\sum c_\alpha c_\beta (\alpha,\beta)\le 0$, so that
$\gamma=0$, a contradiction. 

\noindent (5) \it A positive root is simple iff non--redundant.
\rm
The simple roots are spanning and contained in the
linearly independent set of non--redundant roots, the set of simple
roots must coincide with the set of non--redundant roots.

\noindent (6) \it The simple roots form a basis and every positive root is
a non--negative combination  of simple roots. \rm
The non--redundant roots are linearly independent. Since
the simple roots are non--redundant by (2) and span $\h^*$ by (1),
$\Delta$ forms a basis. Let $X_i$ be the dual basis in $\h$, so that
$\alpha_i(X_j)=\delta_{ij}$. If $\alpha\ge 0$, we know that
$\alpha_i(X)\ge 0$ for all 
$i$ implies $\alpha(X)\ge 0$. Hence $\alpha(X_i)\ge 0$. But
$\alpha=\sum \alpha(X_i)\alpha_i$.

\noindent (7) \it Any $\H_\alpha$ is the wall of some Weyl chamber.
\rm Take $x\in \H_\alpha$ with $x\notin \H_\beta$ for
$\beta\ne \pm \alpha$.  Take $y$ in a small ball around $x$ with
$\alpha(y)>0$. Let $C$ be the Weyl chamber containing $y$. Then
$\H_\alpha$ is a wall of $C$ because $\overline{C}$ intersects
$\H_\alpha$ in a neighbourhood of $x$.

\noindent (8) \it $\Phi= W\cdot \Delta$.
\rm By (7) every $\alpha\in \Phi$ is a simple root for some
Weyl chamber $C^\prime$. But $C^\prime=\sigma C$ for some $\sigma \in
W$, so that $\sigma ^{-1}\alpha\in \Delta$. Hence $\alpha\in
W\cdot\Delta$. 

\noindent (9) \it Every positive root is a non--negative integer
combination of simple roots. 
\rm If $\alpha\in \Phi$, we may have
$\alpha=\sigma\alpha_i$ for $\sigma\in W$ and $\alpha_i\in
\Delta$. Since $\sigma$ is a product of simple reflections, it follows
that $\alpha$ is an integer combination of simple roots. By (6) the
coefficients must be either all non--negative or non--positive. 
\vskip .1in
\noindent \bf Corollary. \it For each $x\in \h$, $W\cdot x \cap
\overline{C}$ is a single point. Thus $\overline{C}$ is a fundamental
domain for $W$. Moreover if $x \in \overline{C}$ then $W_x$ is
generated by the simple reflections fixing $x$, i.e.~by the
reflections in the walls of $C$ containing $x$. 
\vskip .1in
\noindent \bf Proof. \rm The result is obvious for $x\in
\h^\prime$. Otherwise take $x_n\in \h^\prime$ with $x_n\rightarrow
x$. Since there are only finitely many Weyl chambers, we may assume
that $x_n\in \sigma C$ for a fixed $\sigma\in W$. Hence $x\in \sigma
\overline{C}$. 

Now say $x\in \partial C$ and $\sigma x \in \overline{C}$ for
$\sigma\ne 1$. We shall assume the result by induction on
$n(\sigma)=n(C,\sigma C)$, the result being trivial when $\ell=0$.
 Since $x\in \partial C$, we have $x\in \H_\beta$ for
some $\beta$ simple. Since $\H_\beta$ is a wall, we can find an
interior point $y\in \H_\beta\cap \overline{C}$. Take $a\in C$ near
$y$ and $b\in \sigma C$ such the line segment $[a,b]$ is crosses the
hyperplanes $\H_\beta=\H_{\beta_1}, \cdots \H_{\beta_\ell}$
transversely. Let $\sigma=\sigma_{\beta_\ell}\cdots \sigma_{\beta_1}$
and $\tau=\sigma_{\beta_\ell}\cdots \sigma_{\beta_2}$, so that
$\sigma=\tau\sigma_\beta$. Let $C_1=\sigma_\beta C$. Then
$n(C,\sigma C)=\ell$ and $n(C_1,\tau C_1)=\ell-1$. Since $\H_\beta$,
$x=\sigma_\beta x$. Thus $x\in \partial C_1$ and $\tau x\in \partial C_1$. By
induction $\tau x =x$. Since $\sigma=\tau\sigma_\beta$, it follows
that $\sigma x =x$, as required. The second assertion follows by
induction  on $\ell(\sigma)$, since $\ell(\tau)=\ell(\sigma)-1$. 

\vskip .1in

\noindent \bf 8. WEYL'S UNIQUENESS THEOREM. \rm It turns out that
every simple complex Lie algebra is the complexification of a compact
simple Lie algebra, unique up to isomorphism. One proof of this
suggested by Cartan minimises the $\ell^2$ norm of the structure
constants over all choices of orthonormal bases with respect to the
Killing form. Weyl's original proof relied on choosing a basis with
real structure constants, similar to the bases in the lattice
construction. Our aim here is to show that a compact simple Lie
algebra is uniquely determined by its root system.
\vskip .1in
\noindent \bf Theorem~A. \it If $\g_1$ and $\g_2$ are compact
semisimple Lie algebras with isomorphic complexifications, then $\g_1$
and $\g_2$ are isomorphic.
\vskip .1in

\noindent \bf Proof. \rm We may assume that $\tilde{g}$ is the common
complexification. let $J_1$ and $J_2$ be the conjugations
corresponding to $\g_1$ and $\g_2$. Let $B(x,y)={\rm Tr}({\rm ad}(X)
{\rm ad} (Y))$ be the complex Killing form on $\tilde{g}$. Evidently
$B$ restricts to the real Killing forms on both $\g_1$ and
$\g_2$. Take the complex inner product $(X,Y)=-B(X,J_1Y)$ on
$\tilde{\g}$. It is real on $\g_1$. Let $T=J_2J_1$. Then $T\in {\rm
Aut}_{\Bbb C} (\tilde{g})$ and $T$ is self--adjoint,
since
$$(TX,Y)=-B(TX, J_1Y)=-B(X,T^{-1}J_1Y)=-B(X,J_1 TY)=(X,TY).$$
Thus $S=T^2\in {\rm Aut}(\tilde{\g})$. We may identify $\tilde{\g}$
with its image $\g={\rm ad}(\tilde{\g})$ in ${\rm End}(\tilde{g})$ with
the operator bracket. Since
$S \g S^{-1}=\g$, it follows that $S^t \g S^{-t}=\g$ for all $t\in
{\Bbb R}$. On the other hand it is easily checked that
$J_i S^t =S^{-t} J_i$ for $i=1,2$. Let $J^\prime_2= S^t J_2 S^{-t}$.
Then $J_1J_2^\prime=J_1 S^t J_2 S^{-t} =J_1 J_2 S^{-2t}T^{-1}
S^{-2t}$ while $J^\prime_2 J_1=S^t J_2 S^{-t} J_1=S^{2t}T$. These are
equal when $S^{4t}=T^{-2}=S^{-1}$, i.e.~when $t=-1/4$. Thus
$\theta=S^{-1/4}$ gives an automorphism of $\tilde{\g}$ such that
$J_1$ and $J_2^\prime =\theta J_2\theta^{-1}$ commute. 

We claim that $\theta(\g_2)=\g_1$. Let $\g_2^\prime=\theta(\g_2)$. Its
conjugation is now $J_2^\prime$ which commutes with $J_1$. Thus
$\g_1=\{X\in \tilde{\g}: J_1 X= -X\}$ and $\g_2=\{X\in \tilde{\g}:
J_2^\prime X=-X\}$. Now $\g_1=\g_1^+\oplus \g_1^-$, where
$J_2=\pm 1$ on $\g_1^\pm$. Thus $\g_2^\prime=\g_1^+\oplus i \g_1^-$. On the
other hand $-B(X,Y)$ has to be positive definite on $\g_2$ (since
$\g_2$ is compact). Since it is positive definite on $\g_1$ and hence
$\g_1^-$, it is negative definite on $i\g_1^-$. Therefore $\g_1^-=(0)$
and so $\g_2^\prime=\g_1$. Thus $\theta(\g_2)=\g_1$ and $\theta$ is an
isomorphism of $\g_1$ onto $\g_2$, as required.
\vskip .1in

\noindent \bf Theorem~B. \it Let $\g$ and $\g^\prime$
be complexifications of compact simple Lie algebras and
$f:\h\rightarrow \h^\prime$ an 
isometric isomorphism between their maximal abelian subalgebras carrying one
root system onto another. Then $f$ extends uniquely to an (complex)
isomorphism of $\g_{\Bbb C}$ onto $\g^\prime_{\Bbb C}$ carrying $E_i$ onto
$E_i^\prime$.
\vskip .1in
\noindent \bf Proof (A.~Winter). \rm Uniqueness follows because
$f(F_i)$ must be sent onto 
a multiple of $F_i^\prime$. Since $f(H_i)=H_i^\prime$ and
$f(E_i)=E_i^\prime$, the relations $[E_i,F_i]=H_i$ and
$[E_i^\prime,F_i^\prime]=H_i^\prime$ force $f(F_i)=F_i^\prime$. Since
the $E_i$'s and $F_i$'s generate $\g$, this uniquely determines $f$. 

To prove the existence of the isomorphism, let $\overline{\g}$ be the
subalgebra of  $\g\oplus \g^\prime$ generated by the elements
$\overline{H}_i=H_i\oplus H_i^\prime$, $\overline{E_i}=E_i\oplus
E_i^\prime$ and $\overline{F}_i =F_i\oplus F_i^\prime$. 
The algebra $\overline{g}$ has projections $\pi$ and $\pi^\prime$ onto
$\g$ and $\g^\prime$. Clearly ${\rm ker}(\pi)\subset (0)\oplus \g^\prime$ and
${\rm ker}(\pi^\prime)\subset \g\oplus (0)$. Being invariant under
${\rm ad}(\overline{E}_i), {\rm ad}(\overline{F}_i)$ and ${\rm ad}(\overline{H}_i)$, it follows
that ${\rm ker}(\pi^\prime)$ is invariant under ${\rm ad}(E_i)$, ${\rm
ad}(F_i)$ and ${\rm ad}(H_i)$ and hence is an ideal in $\g\oplus
(0)$. Similarly ${\rm ker}(\pi)$ is an ideal in $(0)\oplus
\g^\prime$. Since $\g$ and $\g^\prime$ are simple, either (and hence
both) of these kernels is non--trivial iff $\overline{g}=\g\oplus
\g^\prime$. 

Suppose therefore that $\overline{g}=\g\oplus \g^\prime$. Let $\theta$
be the highest roots for $\h$ and $\h^\prime$ with 
corresponding vectors $E_\theta$ and $E^\prime_{\theta}$. Let
$v= E_\theta\oplus E^\prime_{\theta}\in \g\oplus
\g^\prime=\overline{\g}$. Let $V$ 
be the $\overline{\g}$--submodule generated by $v$. Since $v$ is an
eigenvector for the $\overline{H}_i$ and annihilated by the
$\overline{E}_i$'s, it is clear that $V$ is just the space obtained by
applying monomials in the $\overline{F}_i$'s to $v$. As an $\h\oplus
\h^\prime$--module, the weight $(\beta,\beta)$ occurs with mutiplicity
one in $V$, since the $\overline{F}_i$'s are lowering operators. On
the other hand the $\g\oplus \g^\prime$ by $E_\theta\oplus
E^\prime_{\theta}$ is just $\g\oplus \g^\prime$ by simplicity. This
contradiction proves that the kernels are non--trivial and hence that
$\g_{\Bbb C}$ and $\g^\prime_{\Bbb C}$ are isomorphic as complex Lie
algebras. 
\vskip .1in 
\bf \noindent Theorem~C (Weyl). \it A compact simple Lie algebra is
determined up to isomorphism by its root system.
\vskip .1in
\noindent \bf Proof. \rm Immediate from Theorems A and B.
\vskip .1in
\noindent \bf 9. CLASSIFICATION OF COMPACT SIMPLE LIE ALGEBRAS. \rm 
\vskip .1in
\noindent \bf Irreducibility. \rm A root system is said to be
irreducible if it cannot be written as the disjoint union of two
mutually orthogonal proper subsets. 
\vskip .1in
\noindent \bf Lemma. \it A root sytem is irreducible iff its Weyl
group acts irreducibly.
\vskip .1in
\noindent \bf Proof. \rm Suppose that $\Phi\subset V$ is the root system.
Let $V_1$ be a non--zero $W$--invariant subset of $W$. Since
$x-\sigma_\alpha x= (x,\alpha^\vee)\alpha$, either $\alpha\perp x$ for
all $x\in V_1$ or $\alpha\in V_1$. Thus $\Phi_1=\Phi\cap V_1$ and
$\Phi_2= \Phi\cap V_1^\perp$ are orthogonal and have disjoint union
$\Phi$. So if $\Phi$ is irreducible, $W$ acts irreducibly. Conversely
if $\Phi=\Phi_1\cup \Phi_2$ is an orthogonal splitting, any reflection
$\sigma_\alpha$ fixes pointwise the component in which $\alpha$ does
not lie and hence carries the other component into itself. Thus
$\Phi_1$ and $\Phi_2$ span orthogonal invariant subspaces, so $W$ does
not act irreducibly. 
\vskip .1in
\noindent \bf Cartan matrix. \rm Let
$\Phi$ be a root system with simple roots
$\alpha_1,\dots,\alpha_n$. We define the Cartan matrix $N=(n_{ij})$ by
$n_{ij} =2(\alpha_i,\alpha_j)/(\alpha_i,\alpha_i)$. Note that
$n_{ii}=2$ and $n_{ij}\le 0$ if $i\ne j$; moreover $n_{ij}\ne 0$ iff
$n_{ji}\ne 0$. 
\vskip .1in
\noindent \bf Lemma. \it An irreducible root system is uniquely
determined by its Cartan matrix.
\vskip .1in
\noindent \bf Proof. \rm Let $(\alpha_i)$ and $(\alpha_i^\prime)$ be
systems of simple roots in $V,V^\prime$ such that
$n_{ij}=n^\prime_{ij}$. Define $T:V\rightarrow V^\prime$ by
$T(\alpha_i)=\alpha_i^\prime$. Then $\sigma_{\alpha_j}\alpha_i
=\alpha_i -n_{ji}\alpha_j$. Thus $T\sigma_{\alpha_j}T^{-1}
=\sigma_{\alpha^\prime_j}$. hence $TWT^{-1} =W^\prime$. Since
$\Phi=W\cdot \Delta$ and $\Phi^\prime =W^\prime \cdot \Delta^\prime$,
we get $T\Phi=\Phi^\prime$. Since $W$ and $W^\prime$ act irreducibly,
there is an essentially unique invariant inner product on $V$ and
$V^\prime$. Thus $T$ is a scalar multiple of an isometry so the root
systems $\Phi$ and $\Phi^\prime$ are equivalent.
\vskip .1in
\noindent Dynkin diagram. \rm Let $m_{ij}=2\delta_{ij} -n_{ij}$
($1\le i,j\le n$). Thus $M=(m_{ij})$ is the incidence matrix of a
directed graph, called the Dynkin diagram of the root system. Clearly
the Dynkin diagram completely determines the Cartan matrix and hence
the root sytem. Cearly $\Phi$ is irreducible iff the Dynkin diagram is
connected. 
\vskip .1in
\noindent \bf The highest root. \rm Let $\theta$ be the highest weight
of the adjoint representation on $\g$. Since $\theta$ is the highest
root and
$\sigma_{\alpha_i}\theta=\theta-(\theta,\alpha^\vee_i)\alpha_i$ is
also a root, we must have $(\theta,\alpha_i)\ge 0$. 
Since $\theta$ is a positive root, we may
write $\theta=\sum_{i=1}^n d_i\alpha_i$ with $d_i\ge 0$. Since
$\alpha_i$ is also a weight of $\g$, $\theta-\alpha_i\ge 0$. Thus
$d_i\ge 1$ for all $i$. Note that since $\g$ acts irreducibly on $\g$,
there must be a root $\theta-\alpha_j$ for some $j$ since some
lowering operator must act non--trivially on $\g_\theta$. Thus
$(\alpha_j,\theta)<0$. 

\vskip .1in
\bf \noindent Extended Dynkin diagram. \rm The extended Dynkin diagram
arises naturally in the study of affine Lie algebras, but can easily
be defined without reference to them. It probably provides the
simplest method to classifying Dynkin diagrams. Define
$\alpha_0=-\theta$ extend the definition of
$n_{ij}=2(\alpha_i,\alpha_j)/(\alpha_i,\alpha_i)$ to include $i$ or
$j=0$. We still have $n_{ii}=2$ and $n_{ij}\le 0$. As before define
$m_{ij}= 2\delta_{ij}-n_{ij}$. This is the incidence matrix of a
durected graph called the extended Dynkin diagram. The Dynkin diagram
is obtained by deleting the node $0$ from the extended Dynkin diagram,
so it too is connected. As before
$m_{ii}=0$ and $m_{ij}\ne 0$ iff $m_{ji}\ne 0$. Since $\theta=\sum d_i
\alpha_i$, we have $\sum d_i \alpha_i =0$ if we set $d_0=1$. Thus
we obtain the important equation $\sum m_{ij} d_j = 2d_j$.
This equation implies that the extended Dynkin diagram is a directed
graph with spectral radius $2$. It is easy to classify such graphs.
\vskip .1in

\noindent \bf Graphs of spectral radius two. \rm By a graph we shall mean a
directed graph where nodes $i$ and $j$ are joined by $m_{ij}$
links. We required $m_{ii}=0$ for all $i$ (no loops) and $m_{ij}\ne 0$
iff $m_{ji}\ne 0$. The matrix $M=(m_{ij})$ is called the incidence
matrix of the graph. We shall suppose that the graph is connected. By
Perron--Frobenius theory, the eigenvalue of $M$ of largest modulus is
positive and of multiplicity one; it is the
unique eigenvalue corresponding to an eigenvector with strictly
positive entries. We denote this eigenvalue by $r(M)$ (it is the
spectral radius of $M$.) If we take a connected subgraph,
Perron--Frobenius theory implies that its spectral radius will be
strictly smaller. We use these ideas to classify all connected graphs
of spectral radius $2$.
\vskip .1in
\bf \noindent Theorem. \it Figure~1 gives a complete list of
connected graphs with spectral radius $2$.
\def\ld{\hbox{---}}
$$\matrix{
&&1 &&1 &&&&&&&&&&\cr
A_1^{(1)}&&\circ & \Longleftrightarrow & \circ &&&&&&&&&&\cr
&&&&&&&&&&&&&&\cr
&& 1 &&1 &&&&1&&1&&&&\cr
& & \circ &\ld& \circ & \ld & \cdots & \ld& \circ & \ld
&\circ&&&&\cr
A_n^{(1)}&&|&&&&&&&&|&&&&\cr
& & \circ &\ld& \circ & \ld & \cdots & \ld& \circ & \ld
&\circ&&&&\cr
&& 1 &&1 &&&&1&&1&&&&\cr
&&&&&&&&&&&&&&\cr
&&&&1&&&&&&&&&&\cr
&&&&\circ&&&&&&&&&&\cr
&&&&|&&&&&&&&&&\cr
B_n^{(1)}& & \circ &\ld& \circ & \ld & \circ &\ld & \cdots & \ld& \circ & \Rightarrow
&\circ&&\cr
&&1&&2&&2&&&&2&&2&&\cr
&&&&&&&&&&&&&&\cr
C_n^{(1)}& & \circ &\Rightarrow& \circ & \ld & \cdots & \ld& \circ & \Leftarrow
&\circ&&&&\cr
&&1&&2&&&&2&&1&&&&\cr
&&&&&&&&&&&&&&\cr
&&&&1 &&&&&&1&&&&\cr
&&&&\circ&&&&&&\circ&&&&\cr
&&&&|&&&&&&|&&&&\cr
D_n^{(1)}&& \circ &\ld& \circ & \ld & \circ& \ld&\cdots & \ld& \circ & \ld
&\circ&&\cr
&&1&&2&&2 &&&&2&&1&&\cr}$$
\vfill\eject
$$\matrix{&&&&&&2 && 1&&&&&&\cr
&&&&&&\circ&\ld&\circ&&&&&&\cr
&&&&&&| &&&&&&&&&\cr
E_6^{(1)}&& \circ & \ld & \circ & \ld & \circ & \ld & \circ & \ld
& \circ&&&&\cr
&&1&&2&&3&&2&&1&&&&\cr
&&&&&&&&&&&&&&\cr
&&&&&&&& 2   &&&&&&&&\cr
&&&&&&&&\circ&&&&&&&&\cr
&&&&&&&&|    &&&&&&&&\cr
E_7^{(1)}&&\circ&\ld& \circ & \ld & \circ & \ld & \circ & \ld & \circ & \ld
& \circ&\ld&\circ\cr
&&1&&2&&3&&4&&3&&2&&1\cr
&&&&&&&&&&&&&&\cr
&&&&&&3    &&&&&&&&&&\cr
&&&&&&\circ&&&&&&&&&&\cr
&&&&&&| &&&&&&&&&\cr
E_8^{(1)}&& \circ & \ld & \circ & \ld & \circ & \ld & \circ & \ld
& \circ&\ld&\circ&\ld&\circ&\ld&\circ\cr
&&2 && 4&& 6&&5&&4&&3&&2&&1\cr
&&&&&&&&&&&&&&\cr
F_4^{(1)} && \circ&\ld &\circ& \ld &\circ &\Rightarrow & \circ& \ld &\circ&&&&\cr
&&1&&2&&3&&4&&2&&&&\cr
&&&&&&&&&&&&&&\cr
G_2^{(1)} && \circ&\ld& \circ & \Rrightarrow & \circ &&&&&&&&\cr
&&1&&2&&3&&&&&&&&\cr
&&&&&&&&&&&&&&&\cr
&&1&&2&&&&&&&&&&\cr
A_2{(2)}&&\circ & \Rightarrow^4 &\circ&&&&&&&&&&\cr
&&&&&&&&&&&&&&\cr
&& 1&&2 &&&& 2&& 2&&&&\cr
A_{2n}^{(2)}&&\circ&\Rightarrow&\circ&\ld &\cdots &\ld&\circ
& \Rightarrow &\circ &&&&\cr
&&&&&&&&&&&&&&\cr
&&&&1 &&&&&&&&&&\cr
&&&&\circ &&&&&&&&&&\cr
&&&&|&&&&&&&&&&\cr
A_{2n-1}^{(2)} &&\circ& \ld&\circ &\ld &\circ&\ld&\cdots&\ld&\circ&\Leftarrow &\circ &&\cr
&&1&&2&&2&&&&2&&1&&\cr
&&&&&&&&&&&&&&\cr
D_{n}^{(2)}&&\circ&\Leftarrow & \circ & \ld &\cdots&\ld&\circ&\Rightarrow &\circ&&&&\cr
&&1&&1&&&&1&&1&&&&\cr 
&&&&&&&&&&&&&&\cr
E_6^{(2)}&&\circ & \ld &\circ&\Rightarrow &\circ
&\ld&\circ&\ld& \circ&&&&\cr
&&1&&2&&3&&2&&1&&&&\cr
&&&&&&&&&&&&&&\cr
D_4^{(3)} && \circ&\ld& \circ & \Lleftarrow & \circ &&&&&&&&\cr
&&1&&2&&1&&&&&&&&\cr}$$

\centerline{Figure 1}

\vskip .1in
\noindent \bf Proof. \rm In these graphs a simple bond means that
$m_{ij}=1=m_{ji}$. Otherwise the multiplicity of a bond is indicated
by the number above it. The numbers in the circles give (a multiple
of) the Perron--Frobenius eigenvector corresponding to the eigenvalue $2$.
It is immediately verified by inspection that all the above graphs
have norm $2$, so we just have to show that the list is exhaustive. We
shall consistently use the fact that a graph of spectral radius $2$
cannot have a proper subgraph with spectral radius $2$. 

If the graph contains a cycle with three or more nodes, it must
contain and hence equal a subgraph $A_n^{(1)}$.
Thus we may assume that there are no such cycles. If it contains a
node of connected to four or more other nodes, it must contain and
hence equal a subgraph $D_4^{(1)}$. Thus we may assume that each node
is connected to at most three other nodes. 
Suppose next that the graph has a bond of multiplicity greater than or
equal to $4$. If so it has a
subgraph of type $A^{(2)}_2$, which it must equal. Now suppose the
graph has a bond of multiplicity $3$. If the other bond between these
nodes had 
multiplicity greater than or equal to two, then $A_1^{(1)}$ would
be a proper subgraph and the graph would have spectral radius greater
than two.  So it must have multiplicity one, so that the graph would have to
be contained and hence equal either $G_2^{(1)}$ or $D_4^{(2)}$.
So we may assume all bonds have multiplicity $1$ or
$2$. Suppose it has at least $2$ bonds of multiplicity $2$. If they
are between the same nodes, $A_1^{(1)}$ is a subgraph so the whole
graph. If they are between different nodes, then the graph must
contain and hence equal one of
$C_n^{(1)}$, $A_{2n}^{(2)}$ or $D_{n+1}^{(2)}$. So we may assume that
there is only one bond of multiplicity $2$. Suppose that there is a
node of valency three. It must be connected to one of the nodes in the
multiplicity two bond. Thus the graph contains and hence equals one of
the graphs $B_n^{(1)}$ or $A_{2n-1}^{(2)}$. So we may assume it has no
nodes of valency three, so that graphs is just one long
string. Neither of the nodes in the multiplicity $2$ bond 
can be an endpoint of the string, for the graph would be a proper
subgraph of $B_n^{(1)}$ or $A_{2n-1}^{(2)}$ and hence have spectral
radius strictly less than $2$. Thus each is connected to a further
point; if there no additional points, the graph would be a proper subgraph of
$F_4^{(1)}$ (or equally well $E_6^{(2)}$), and hence have spectral
radius strictly less than $2$. So there must be a fifth point, 
so the graph contains and hence equals $F_4^{(1)}$ or $E_6^{(2)}$.
Thus we may assume the graph has only bonds of multiplicity one and
only trivalent vertices. If it has two trivalent vertices, it must
contain and hence equal $D_n^{(1)}$. If it had no trivalent vertices,
it would be a proper subgraph of $A_n^{(1)}$, which would make its
spectral radius less than $2$. So we may assume it has exactly one
trivalent vertex. If all the branches from that vertex have length
greater than or equal to two, it contains and hence equals
$E_6^{(1)}$. So some branch has length one. If two had length one,
then the graph would be a proper subgraph of $D_n^{(1)}$, so have
spectral radius less than $2$. Hence one of the branches has length
$1$ and the two other length greater than or equal to $2$. If they
both have length greater than or equal to $3$, then the graph contains
and hence equals $E_7^{(1)}$. So one of branches must have length
equal to $2$. If the remaining branch has length greater than or equal
to $5$, the graph contains and hence equals $E_8^{(1)}$. If it had
length less than $5$, it would be contained in $E_8^{(1)}$, so would
have spectral radius less than $2$. This completes the proof (which
could profitably arranged in a flow chart). 
\vskip .1in
\noindent \bf Classification of Dynkin diagrams. \rm Figure~2 gives
a complete list of possible Dynkin diagrams.
\vskip .1in
\noindent \bf Proof. \rm These are the only diagrams that arise from
the list of the preceding theorem by removing one node so as not too
disconnect the graph. The lattice construction and its fixed point
refinement yield root systems corresponding to each Dynkin diagram.
\vfill\eject
$$\matrix{A_n& & \circ &\ld& \circ & \ld & \cdots & \ld& \circ & \ld
&\circ&&&&\cr
&&&&&&&&&&&&&&\cr
&&&&&&&&&&&&&&\cr
&&&&&&&&&&&&&&\cr
B_n& & \circ &\ld& \circ & \ld & \cdots & \ld& \circ & \Rightarrow
&\circ&&&&\cr
&&&&&&&&&&&&&&\cr
&&&&&&&&&&&&&&\cr
&&&&&&&&&&&&&&\cr
C_n& & \circ &\ld& \circ & \ld & \cdots & \ld& \circ & \Leftarrow
&\circ&&&&\cr
&&&&&&&&&&&&&&\cr
&&&&&&&&\circ&&&&&&\cr
&&&&&&&&|&&\cr
D_n&& \circ &\ld& \circ & \ld & \cdots & \ld& \circ & \ld
&\circ&&&&\cr
&&&&&&&&&&&&&&\cr
&&&&&&\circ&&&&&&&&&&\cr
&&&&&&| &&&&&&&&&\cr
E_6&& \circ & \ld & \circ & \ld & \circ & \ld & \circ & \ld
& \circ&&&&\cr
&&&&&&&&&&&&&&\cr
&&&&&&\circ&&&&&&&&&&\cr
&&&&&&| &&&&&&&&&\cr
E_7&& \circ & \ld & \circ & \ld & \circ & \ld & \circ & \ld
& \circ&\ld&\circ&&\cr
&&&&&&&&&&&&&&\cr
&&&&&&\circ&&&&&&&&&&\cr
&&&&&&| &&&&&&&&&\cr
E_8&& \circ & \ld & \circ & \ld & \circ & \ld & \circ & \ld & \circ&\ld&\circ&\ld&\circ\cr
&&&&&&&&&&&&&&\cr
&&&&&&&&&&&&&&\cr
&&&&&&&&&&&&&&\cr
F_4 && \circ& \ld &\circ &\Rightarrow & \circ& \ld &\circ&&&&&&\cr
&&&&&&&&&&&&&&\cr
&&&&&&&&&&&&&&\cr
&&&&&&&&&&&&&&\cr
G_2 && \circ & \Rrightarrow & \circ &&&&&&&&&&\cr}$$
\centerline{Figure 2}

\vfill\eject
\noindent \bf PART 2. REPRESENTATION THEORY. \rm
\vskip .1in
\it In this part we develop the representation theory of an arbitrary
compact simple matrix group analogously to that of $SU(2)$. Passing to
the Lie algebra, we classify representations by their highest weight
and give a description of the representation in terms of lowering and
raising operators. WE then prove the Weyl character formula by studing
the supercharge or Dirac operator. This same procedure will be
followed for affine Kac--Moody algebras. The most important example to
follow is $SU(N)$, but we use the language of root systems so that the
proofs apply in general. 
\vskip .1in
\noindent \bf 10. ROOT AND WEIGHT LATTICES. \rm Let $G$ be
a compact simple matrix group with lie algebra $\g$. Let $T$ be a
maximal torus in $G$ with Lie algebra $\h$ and Weyl group $W=N(T)/T$. 
Let $\Lambda\subset \h$ be the kernel of the map
$\h\rightarrow T$ $X\mapsto e^{2\pi T}$. Thus $\h/\Lambda \cong T$ so that
$\Lambda$ is a lattice in $\h$, called the unit or integer
lattice. We claim that $T^\vee={\rm Hom}({\Bbb
T}, T)\cong \Lambda$. In fact 
any homomorphism $f:{\Bbb T}\rightarrow T$ has the form $f(e^{it}) =
\exp(tX)$ for a unique $X\in \h$ (simply take the infinitesimal
representations). Setting $t=2\pi$, we get $\exp(2\pi X)=1$, so that $X\in
\Lambda$. 
\vskip .1in
\bf \noindent The weight lattice. \rm The weight lattice $P(G)$ of $T$ is
the group 
$\hat{T}={\rm Hom}(T,{\Bbb T})$. Looking at the corresponding
infinitesimal homomorphism, we see that any $\chi\in\hat{T}$ 
has the form $\chi(\exp(X))=e^{i\lambda(x)}$
for a unique $\lambda\in \h^*$. Let $P(G)$ be the subgroup of $\h^*$
consisting of weights. Clearly $\chi(\exp(2\pi X))=1$ if $X\in
\Lambda$, so that $\lambda\in \h^*$ defines a character or weight iff
$\lambda(X)\in {\Bbb Z}$ for all $X\in {\Bbb Z}$. Thus $P(G)$ forms a
lattice, the weight lattice, and $P(G)$ and $\Lambda$ are dual
lattices. We write $P(G)=\Lambda^*$ and $\Lambda=P(G)^*$. 
If we restrict the faithful representation of $G$ to $V$ and
decompose it in characters, we get a finite set of homomorphisms
$\chi_i:T\rightarrow {\Bbb T}$. Let $\Gamma_0$ be the subgroup of
$\hat{T}$ generated by the $\chi_i$'s. The next result shows that
$\Gamma_0=\hat{T}$. 
\vskip .1in
\noindent \bf Lemma. \it If $\Gamma_0$ is a subgroup of $\Gamma={\rm
Hom}({\Bbb T}^n,{\Bbb T})$ distinguishing the points of ${\Bbb
T}^n$, then $\Gamma_0=\Gamma$.
\vskip .1in
\noindent \bf Proof. \rm $\Gamma/\Gamma_0$ is a finitely generated
Abelian group, so there admits a non--trivial homormorphism $\theta$
into ${\Bbb T}$. Now $\Gamma={\Bbb Z}^n$ with ${\Bbb Z}$--basis
$e_1,\dots,e_n$. Let $t_i=\theta(e_i)\in {\Bbb T}$. Thus $t=(t_i)\in
{\Bbb T}^n$. By definition $\theta(\sum m_i e_i)=t^m=e_m(t)$ for $m\in
{\Bbb Z}^n$. But then $e_m(t)=1$ for all $m\in \Gamma_0$, a
contradiction. 
\vskip .1in
\bf \noindent The root lattice. \rm Let $\Phi$ be the set of non--zero
weights appearing in the complexified adjoint representation on
$\g_{\Bbb C}$. Since the adjoint representation is real, if $\alpha\in
\Phi$, then $-\alpha\in \Phi$. Let $Q$ be the sublattice of $\h^*$
spanned by $\Phi$. Thus $Q\subseteq P_G$.
\vskip .1in
\noindent \bf The centre of $G$. \rm The centre $Z(G)$ is a closed
subgroup of $G$, so a Lie group. Its Lie algebra is just the centre of
$\g$, so trivial. Hence $Z(G)$ is finite. 
\vskip .1in
\noindent \bf Lemma. \it $Z(G)\cong Q^*/P^*=(P/Q)^*$.
\vskip .1in
\noindent \bf Proof. \rm Note that $Z(G)\subset T$. Thus we may write
any $z\in Z(G)$ as $z=e^{2\pi X}$ for $x\in \h/P^*$. Now
$e^{2\pi X}\in Z(G)$ iff ${\rm ad} (e^{2\pi X})=1$ iff
${\rm Ad}(e^{2\pi X})$ fixes $\g_\alpha$ for all $\alpha\in \Phi$ iff
$e^{2\pi i \alpha(X)} =1$ for all $\alpha$ iff $\alpha(X)\in {\Bbb Z}$
for all $\alpha\in Q$ iff $X\in Q^*$. Hence $Z(G)\cong Q^*/P^*$, as
required. 
				
\vskip .1in
\noindent \bf The generalised weight lattice. \rm For each $\alpha$,
we have a copy of $\sl_2$ corresponding to $\alpha$, namely
$H_\alpha,E_\alpha,F_\alpha$, where
$H_\alpha=-2iT_\alpha/\|\alpha\|^2$ and
$\mu(T_\alpha)=(\alpha,\mu)$. Let $V$ be the defining representation
of $G$ (so that $G\subset U(V)$). By the $SU(2)$ theory
$\lambda(H_\alpha)\in {\Bbb Z}$ for each weight $\lambda$ of
$V$. Hence $2(\lambda,\alpha)/\|\alpha\|^2\in {\Bbb Z}$ for each
root $\alpha$, i.e.~$(\lambda,\alpha^\vee)\in {\Bbb Z}$. This defines
a lattice $P(\g)=\{\lambda\in \h^*:(\lambda,\alpha^\vee)\in {\Bbb Z}\}$,
containing $P(G)$. It is called the generalised weight lattice.

\vskip .1in 
\noindent \bf \noindent \bf The dual root system. \rm The inner
product on the real inner product space $\h$ allows us to identify
$\h$ and $\h^*$. Recall that if
$\alpha\in \h^*$, then $\alpha^\vee\in \h$ is defined by
$\alpha^\vee(\lambda) =2(\alpha,\lambda)/(\alpha,\alpha)$. 
\vskip .1in
\noindent \bf Proposition. \it (1) If $\Phi$ is a root
system, then $\Phi^\vee$ is also a root system with $W(\Phi^\vee)
=W(\Phi)$. 

\noindent (2) If  $\alpha\in \Phi^\vee$, then $\pm \alpha$ is a
non--negative integer combination of the $\alpha_i^\vee$'s,
$\alpha=\sum n_i\alpha_i$.

\vskip .1in
\noindent bf Proof. \rm (1) Recall that $\alpha^\vee
2\alpha/(\alpha,\alpha)$. Thus $\alpha^{\vee\vee}=\alpha$. Also
$t\alpha\in \Phi$ iff $t^{-1}\alpha^\vee\in \Phi^\vee$; thus
$\Phi^\vee$ is reduced. Since $\Phi$ spans $V$, so too does
$\Phi^\vee$. Note also that
$\sigma^{\beta^\vee}(\alpha^\vee)=\sigma_\beta(\alpha^\vee)
=\sigma_\beta(\alpha)^\vee)$. Finally since $\Phi$ is a root system,
$(\alpha,\beta^\vee)\in {\Bbb Z}$ for all $\alpha,\beta\in
\Phi$. Since $\alpha^{\vee\vee}=\alpha$, the same condition holds for
$\Phi^\vee$. Hence $\Phi^\vee$ is an abstract root system. The
Weyl groups coincide because $\sigma_{\alpha^\vee}=\sigma_\alpha$. 

\noindent (2) Since $\alpha_i^\vee$ is proportional to $\alpha_i$, the
$\alpha_i^\vee$'s are a basis of $V$. So it suffices to show that
$\Phi^\vee \subset \bigoplus {\Bbb Z}\alpha_i^\vee=Q_1$ say. Clearly
$Q_1$ is invaraint under the $\sigma_{\alpha_i}$'s, so $W$. If
$\alpha\in \Phi^+$, then $\alpha=\sigma\alpha_i$ some $i$ (since
$\h_\alpha^+$ is the wall of some chamber). So $\alpha^\vee
=\sigma\alpha_i^\vee\in Q_1$, as required.

\vskip .1in
\noindent \bf The coroot lattice. \rm This is the lattice in $\h$ spanned by
the coroots $\alpha^\vee$ over ${\Bbb Z}$. The relationship between
the lattices is best illustrated in a 
picture.

$$\matrix{ Q^* & \qquad & P(\g) \cr
           \cup    & & \cup \cr
           \Lambda & & P(G)\cr
           \cup    & & \cup \cr
           Q^\vee & & Q\cr}$$
\vskip .1in
\noindent \bf Toplogical remark. \rm It turns out that the fundamental
group of $G$ is just $P(\g)/P(G)$. Thus $G$ is simply connected iff
$P(\g)=P(G)$, i.e.~every generalised weight is a weight. This result is
due to Hermann Weyl; its proof requires extra analytic or topological
tools. In this chapter we shall content ourselves with constructing
matrix group $\tilde{G}$ and a homomorphism $f:\tilde{G}\rightarrow G$
which is a covering map (${\rm ker}(f)$ is finite and central) such
that $P(\tilde{G})=P(\g)$. (Note that since $f$ is a covering map, the
Lie algebra of $\tilde{G}$ is just $\g$, just like the 
double cover ${\rm Spin}(V)\rightarrow SO(V)$.) 
\vskip .1in
\noindent \bf 11. POINCARE--BIRKHOFF--WITT THEOREM. \rm Let $\g$ be a Lie
algebra, possibly infinite dimensional, with basis
$X_1,X_2,\dots$. Then
the universal eneveloping algebra of $\g$ has basis $X_{i_1} X_{i_2}
\dots X_{i_k}$ with $i_1\le i_2\le \cdots \le i_k$.
\vskip .1in
\noindent \bf Proof (Jacobson). \rm We define $U(\g)$ to be the
quotient of the tensor algebra $T(\g)$ by the two--sided ideal
generated by $X\otimes Y - Y\otimes X -[X,Y]$ with $X,Y\in \g$. Let
$X_1,X_2,\dots$ be a basis of $\g$ (possibly
infinite--dimensional). Let $S(\g)$ be the symmetric algebra of $\g$
with basis $X_{i_1} \dots X_{i_n}$ with $i_1\le i_2 \le \le i_n$. We
claim that there is a unique linear `symbol' map
$\sigma:T(\g)\rightarrow S(\g)$ such that 
$$\sigma(a_1\otimes a_2 \otimes \cdots \otimes( a\otimes b-b\otimes a
-[a,b])\otimes \cdots \otimes a_m)=0\eqno{(1)}$$ 
and
$$\sigma(X_{i_1}\otimes \cdots \otimes X_{i_n})=
X_{i_1}\cdots X_{i_n}\eqno{(2)}$$ 
if $i_1\le i_2\le \cdots \le
i_n$. Suppose that 
such a map has been constructed for the linear span of all monomials
of degree $\le n-1$. We proceed by induction on the length of the
permutation required to put a monomial of degree $n$ in `correct'
order. If this length is $0$, we simply use (2) to define
$\sigma$. Otherwise if the order is wrong we can find a transposition
of two adjacent terms which decreases the length of the permutation. We
then define
$$\sigma(X_{i_1}\otimes \cdots
\otimes X_{i_n})
=\sigma(X_{i_1}\otimes \cdots\otimes X_{i_{k+1}}\otimes X_{i_k}\otimes \cdots 
\otimes X_{i_n})+
\sigma(X_{i_1}\otimes \cdots\otimes [X_{i_k}, X_{i_{k+1}}]\otimes \cdots 
\otimes X_{i_n}).$$
This has to be the case and therefore proves uniqueness of
$\sigma$. We must show that $\sigma$ is independent of the choice  of
$k$. If we used another transposition disjoint from $(k,k+1)$. Then we
could apply the same process to the right hand side of the above
equation with $(\ell,\ell+1)$ in place of $(k,k+1)$. We would clearly
get the same answer if we did the $(\ell,\ell+1)$ transposition first
followed by $(k,k+1)$, since they are disjoint. If the other
transpoition was not disjoint, we may suppose that the two
transpositions are $(k,k+1)$ and $(k-1,k)$. Thus
$i_{k-1}<i_k<i_{k+1}$.  Write $Y_j=X_{i_j}$. Then
if we use $(k,k+1)$ to define $\sigma$, we get
$$\eqalign{\sigma(\cdots Y_{k-1}&\otimes Y_{k} \otimes Y_{k+1} \cdots )\cr
&=\sigma( \cdots Y_{k-1}\otimes Y_{k+1} \otimes Y_k \cdots) 
+\sigma(\cdots Y_{k-1}\otimes [Y_k,Y_{k+1}]\cdots)\cr
&=\sigma(\cdots Y_{k+1} \otimes Y_{k-1} \otimes Y_k \cdots)
+ \sigma(\cdots [Y_{k-1}, Y_{k}]\otimes Y_k \cdots)
+\sigma(\cdots Y_{k-1}\otimes [Y_k,Y_{k+1}]\cdots)\cr
&=\sigma(\cdots Y_{k+1} \otimes Y_k \otimes Y_{k-1}\cdots )
+ \sigma(\cdots Y_{k+1} \otimes [Y_{k-1},Y_k]\cdots)\cr
&\quad + \sigma(\cdots [Y_{k-1}, Y_{k}]\otimes Y_k \cdots)
+\sigma(\cdots Y_{k-1}\otimes [Y_k,Y_{k+1}]\cdots).\cr}$$
If we use $(k-1,k)$ to define $\sigma$, we get
$$\eqalign{\sigma&(\cdots Y_{k-1}\otimes Y_{k} \otimes Y_{k+1} \cdots )\cr
&=\sigma(\cdots Y_{k+1} \otimes Y_k \otimes Y_{k-1}\cdots )
+ \sigma(\cdots Y_{k} \otimes [Y_{k-1},Y_{k+1}]\cdots)\cr
& \quad+ \sigma(\cdots [Y_{k}, Y_{k+1}]\otimes Y_{k-1} \cdots)
+\sigma(\cdots [Y_{k-1},Y_k]\otimes Y_{k+1}\cdots).\cr}$$
The difference between these expression is zero by (1) and the Jacobi
identity for $Y_{k-1}$, $Y_k$ and $Y_{k+1}$. The fact that the
extension is well--defined means that it has all the stated
properties.

Property (2) implies the map $\sigma$ induces a map $\sigma
:U(\g)\rightarrow S(\g)$. Property (1) implies that $\sigma $ carries
the span of the monomials of 
degree $\le N$ onto $\bigoplus_{j\le N} S^j(\g)$. Since the monomials
of degree $\le N$ are spanned by the monomials $X_1^{\alpha_1} \cdots 
X_n^{\alpha_n}$ with $|\alpha|\le N$ which map into identical
linearly independent symbols under $\sigma$, the theorem follows.

\vskip .1in 
\noindent \bf 12. HIGHEST WEIGHT VECTORS. \rm Let $\pi:G\rightarrow
U(V)$ be a continuous homomorphism. We call $(\pi,V)$ a unitary
representation of $G$ or $G$--module. We may decompose $V$ as a
$T$--module, $V=\bigoplus_{\mu\in \widehat{T}} V_\mu$ where $e^X\in T$
acts on $V_\mu$ as the scalar $e^{i\mu(X)}$. We call $\mu$ a weight of
$V$ and $V_\mu$ the corresponding weight space. Because of the
identification of $\widehat{T}$ and the weight lattice $P(G)$, we have
$\mu\in P(G)$. Notice that the root system is invariant under the Weyl
group. Hence $Q$ and $Q^\vee$ are invariant under the Weyl group. It
follows that $P(\g)$ and $Q^*$ are also invaraint under the Weyl
group. Since ${\rm Ad}(g) \exp(2\pi X)=\exp(2\pi \,{\rm Ad}(g)\cdot
X)$, for $X\in 
\h$ and $g\in N(T)$, it follows that the integer lattice is invariant
under $W$. Hence the weight lattice is invariant under $W$. This also
follows by applying the following lemma to a faithful representation. 
\vskip .1in
\noindent \bf Lemma. \it The weights of $V$ are invariant under the
Weyl group. In fact if $\sigma=gT$ for $g\in N(T)$, then
$g V_\lambda=V \sigma \lambda$. 

\vskip .1in
\noindent \bf Proof. \rm We have
$$\pi(e^X) \pi(g) v=
\pi(g) \pi(e^{g^{-1} Xg}) v= \pi(g)
e^{i\lambda{g^{-1}Xg}}v=e^{i\sigma\cdot\lambda(X)}\pi(g)v,$$
so that $\pi(g)v$ lies in $V_{\sigma\lambda}$. 

\vskip .1in
If we pass to the infinitesimal representation of $\g$ on $V$,
$\pi(\g)$ and $\pi(G)$ have the same commutant since $G$ is
connected. As with $SU(2)$, we shall temporarily drop the
skew--adjointness assumption on the matrices $\pi(X)$. Thus we shall
study finite--dimensional representations $\pi:\g\rightarrow {\rm
End}(V)$ of the Lie algebra such that the matrices $\pi(X)$
($X\in \h$) are simultaneously diagonalisable. We may therefore
decompose $V$ as a direct sum $\bigoplus V_\mu$ where
$\pi(X)v=i\mu(X)v$ for $X\in \h$ with $\mu\in {\rm Hom}(\h,{\Bbb
C})$. We call $\mu$ a weight of $V$ and $V_\mu$ the corresponding
weight space. On the other hand for each root $\alpha$ we can
construct a copy of $\sl_2$ iinside $\g_{\Bbb C}$, namely $E_\alpha,
F_\alpha, H_\alpha$ with $H_\alpha=-2iT_\alpha/\|\alpha\|^2$ and
$(T_\alpha,X)=\alpha(X)$ for $X\in \h$. We may consider $V$ as an
$\sl_2$--module: plainly 
$$\pi(H_\alpha)v= (\mu,\alpha^\vee) v$$
for $v\in V_\mu$. This immediately gives the following result.
\vskip .1in
\noindent \bf Lemma. \it (a) Any weight of $V$ lies in $P(\g)$, so is
a generalised weight.

\noindent (b) $E_\alpha V_\mu \subseteq V_{\mu+\alpha}$ and $F_\alpha
V_\mu \subseteq V_{\mu-\alpha}$, so that $E_\alpha$ and $F_\alpha$ act
as raising and lowering operators if $\alpha>0$.

\vskip .1in
\noindent \bf Proof. \rm The first statement follows because
$\pi(H_\alpha)$ can only have integer eigenvalues by the $SU(2)$
theory. The second statement follows because
$[X,E_\alpha]=\alpha(X)E_\alpha$ for $X\in \h$.
\vskip .1in

The appearance of lowering and raising operators leads us to define a
partial order on $\h^*$. We say that $\lambda\ge \mu$ iff $\lambda-\mu=\sum
t_i\alpha_i$ with $t_i\ge 0$. (In all cases the $t_i$'s will be
integers.) Clearly every finite--dimensional representation has a
highest weight space $V_\mu$, not necessarily unique. Any vector in
it, besides being an eigenvector for $X\in \h$, will be annihilated by
any raising operator $E_\alpha$ for $\alpha>0$. 
\vskip .1in
\noindent \bf Theorem. \it (1) Every irreducible representation $V$ has a
unique highest weight vector. 

\noindent (2) Every vector is in linear span of the vectors obtained
by successively applying lowering operators to the highest weight
vector. 

\noindent (3) The corresponding weight space is one--dimensional.

\noindent (4) The highest weight $\lambda$ satisifies
$(\lambda,\alpha)\ge 0$ for all $\alpha>0$; we say that $\lambda$ is
{\rm dominant}. 

\noindent (5) Two irreducible representations with the same highest
weight are isomorphic. 
\vskip .1in
\noindent \bf Proof. \rm  There certainly is at least one highest
weight vector, since ${\rm dim}(V)<\infty$.
Let $v\in V_\lambda$ be a non--zero highest
weight vector. Thus $\pi(E_\alpha)v=0$ for $\alpha>0$.
Every elemnt in the enveloping algebra ${\cal
U}(\g)$ can be written as a sum of monomials $LDR$ where $L$ is a
product of lowering operators $\pi(F_\alpha)$ ($\alpha>0$), $D$ is a
diagonal operator (a monomial in $\pi(X)$'s for $X\in \h$) and $R$ is
a product of raising operators $\pi(E_\alpha)$ ($\alpha>0$). Since
$Rv=0$ and $Dv$ is proprtional to $v$, it follows that $LDRv$ is
either $0$ or proportional to $Lv$, so (2) holds. Since applying a
lowering operator always lowers the weight, (2) implies uniqueness in
(1) as well as (3). (4) follows from
$\lambda(H_\alpha)=(\lambda,\alpha^vee)$, $\pi(E_\alpha)v=0$
($\alpha>0$) and the $SU(2)$ theory. To prove (5), we may suppose we
have two irreducible representations $V_1$ and $V_2$ with highest
weight $\lambda$. Let $v_1$ and $v_2$ be corresponding non--zero
highest weight vectors. Then $v=v_1\oplus v_2\in V_1\oplus V_2$ is a
highest vector of weight $\lambda$. Let $E$ be the $\g$--submodule
generated by $v$. Since the raising operators annihilate $v$, $E$ is
spanned by vectors obtained by applying lowering operators to $v$. 
Hence $E_\lambda ={\Bbb C}v$. Let
$f$ be the restriction of the propjection $V_1\oplus V_2\rightarrow
V_1$ to $E$. Since $f(v)=v_1$ and $V_1$ is irreducible, $f(E)=V_1$. On
the other hand ${\rm ker}(f)\subset V_2\cap E$ by definition of
$f$. Now $V_2\cap E$ is a submodule of $V_2$ and $E$. It does not
contain $v_2$, because if it did $v_1=v-v_2$ would lie in $E$ which
would contradict $E_\lambda={\Bbb C}v$. Thus $V_2\cap E\ne V_2$, so
that $V_2\cap E=(0)$, by irreducibility of $V_2$. It follows that $f$
is an isomorphism of $E$ onto $V_1$, so that $E\cong V_1$ as
$\g$--modules. Similarly $E\cong V_2$ and hence $V_1\cong V_2$.
\vskip .1in
Our next goal will be to prove a converse of this theorem, namely to show
that every dominant generalised weight is the highest weight of a
finite--dimensional irreducible representation of $\g$ (see the next
two sections). 
\vskip .1in
\noindent \bf The fundamental weights. \rm We know that if $\alpha_i$
are simple roots, then $\alpha^\vee_i$ are simple coroots. Thus they
form a ${\Bbb Z}$--basis for the coroot lattice $Q^\vee$. Since
$Q^\vee=P(\g)^*$, we get a dual basis $\lambda_i\in P_\g$ defined by
$\lambda_i(\alpha_j^\vee)=\delta_{ij}$. Thus
$(\lambda_i,\alpha_j)=\delta_{ij} (\alpha_i,\alpha_i)/2$. The
generalised weights $\lambda_i$ are called the fundamental weights. 
\vskip .1in
\noindent \bf Simple reflections and positive roots. \it If $\alpha$
is a simple root, then $\sigma_\alpha$ permutes all the positive roots
not equal to $\alpha$.
\vskip .1in
\noindent \bf Proof. \rm Suppose $\beta$ is a positive root with
$\beta\ne \alpha$. Then $\beta=\sum_{\gamma\in \Delta} n_\gamma \gamma$. Since $\beta\ne
\alpha$, it is not proportional to $\alpha$, so $n_\gamma>0$ for some
$\gamma \ne \alpha$. But $s_\alpha \beta =\beta - t\alpha$, so the
coefficient of $\gamma$ in $s_\alpha\beta$ is also $n_\gamma$. Thus
$s_\alpha\beta$ must be positive. 
\vskip .1in
\noindent \bf Corollary~1. \rm If $\rho$ is half the sum of the positive
roots and $\alpha$ is simple, then $\rho-\sigma_\alpha\rho=\alpha$.
\vskip .1in
\noindent \bf Proof. \rm Immediate since $\sigma_\alpha 
\alpha=-\alpha$. 
\vskip .1in
\noindent \bf Corollary~2. \it If $\rho={1\over 2} \sum_{\alpha>0} \alpha$,
then $\rho=\sum \lambda_i$. Hence $\rho\in P(\g)$ and, if $\sigma\rho=\rho$ for
$\sigma \in W$, then $\sigma=1$.
\vskip .1in
\noindent \bf Proof. \rm By Corollary~1,
$(\rho,\alpha_i^\vee)=1$. The result follows because $(\lambda_i)$ is
the dual basis to $(\alpha_j^\vee)$. If $\sigma\rho=\rho$, then,
because $\rho$ is in the positive Weyl chamber,
$\sigma$ is in the subgroup generated by simple reflections fixing
$\rho$, of which there are none. So $\sigma=1$. Note that this result
is obvious for $SU(N)$. Alternatively, we will see in section~14 that
each $\lambda_i$ is a highest weight vector of a representation whose
weights are invariant under $W$. Thus $\sigma\lambda_i\le
\lambda_i$. Hence $\sigma\lambda_i =\lambda_i$ for all $i$. Hence
$\sigma=1$.  
\vskip .1in
\noindent \bf Corollary~3. \it $\rho-\sigma\rho=\sum_{\alpha\in
\Phi_{\sigma^{-1}}} \alpha$. 
\vskip .1in

\noindent \bf Proof. \rm Clearly $\rho-\sigma\rho=\sum \beta$, where
the sum is over all $\beta>0$ such that $\beta=-\sigma\alpha$ with
$\alpha>0$. Hence $\rho-\sigma\rho =-\sum_{\alpha\in \Phi_\sigma}
\sigma\alpha$. We get the result by changing $\sigma$ to $\sigma^{-1}$
and applying $\sigma$.

\vskip .1in

\noindent \bf 13. EIGENVALUES OF THE CASIMIR OPERATOR. \rm Let $\g$ be
a Lie algebra and 
$(x,y)$ an invariant real inner product on $\g$. Let $(X_i)$ be an
orthonormal basis. Then $Z=-\sum X_i\otimes X_i$ is an invariant element
in $\g\otimes \g$; for clearly $Z$ is independent of the choice of
orthonormal basis and ${\rm ad}(g)\cdot X_i$ is also an orthonormal
basis. Now let $V$ be any representation of $\g$. Then $\g\otimes \g
\rightarrow {\rm End}(V)$, $X\otimes Y\mapsto XY$ is
$\g$--equivariant. Hence the image of $Z$ commutes with $\g$. This
image is called the Casimir operator $\Omega=-\sum \pi(X_i)^2$. Note
that we take the minus sign because, if $\pi(X_i)^*=-\pi(X_i)$, then
$\Omega =\sum \pi(X_i)^*\pi(X_i)$ is a positive operator.								
\vskip .1in

Recall that if
$(x,y)$ is an invariant real inner product on $\g$, $(X_i)$ is an
orthonormal basis of $\g$ and $\pi:\g\rightarrow {\rm End}(V)$ a
representation, then the Casimir operator $C=-\sum \pi(X_i)^2$
commutes with $\g$. As for $\sl(2)$, we can express $C$ in terms of the
elements $H_i$ and the 
lowering and raising operators $E_\alpha$, $F_\alpha$
($\alpha>0$). Let $(T_i)$ be any orthonormal basis of $\h$ and let
$H_\alpha=-2iT_\alpha/\|\alpha\|^2$ with corresponding elements
$E_\alpha$, $E_{-\alpha}=E_\alpha^*$. Since $H_\alpha=4/\|\alpha\|^2$
we have $\|E_{\pm\alpha}\|^2=2/\|\alpha\|^2$. Thus $(T_i)$ and
$\|\alpha\|E_{\pm\alpha}/\sqrt{2}$ is an orthonormal basis of
$\g_{\Bbb C}$. If we take any orthonormal basis $(X_j)$ of $\g_{\Bbb C}$, we
still have $\Omega=\sum X_j^*X_j$. (This is independent of the choice
of orthonormal basis). 

\vskip .1in
\noindent \bf Lemma. \it $\Omega=-\sum T_i^2 -i\sum_{\alpha>0}
T_\alpha + \sum_{\alpha>0}\|\alpha\|^2
E_{\alpha}^* E_\alpha$. It acts on the representation $V_\lambda$ as 
the scalar $\|\lambda\|^2 +2(\lambda,\rho) =\|\lambda +\rho\|^2 -\|\rho\|^2$, where
$\rho={1\over 2}\sum_{\alpha>0} \alpha$. 
\vskip .1in
\noindent \bf Proof. \rm We have 
$$\Omega=-\sum T_i^2 +\sum_{\alpha\in \Phi} {\|\alpha\|^2\over 2}
E_\alpha^* E_\alpha .$$
On the other hand $E_\alpha F_\alpha =H_\alpha +F_\alpha E_\alpha
=-2iT_\alpha/\|\alpha\|^2 +F_\alpha E_\alpha$. Hence
$$\Omega=-\sum T_i^2 -i\sum_{\alpha>0}
T_\alpha + \sum_{\alpha>0}\|\alpha\|^2
E_{\alpha}^* E_\alpha.$$ 
Applying $\omega$ to the highest weight vector $v_\lambda\in
V_\lambda$ (which is annihilated by $E_\alpha$ for $\alpha>0$), we get
$$\Omega v_\lambda = [(\lambda,T_i)^2 +\sum_{\alpha>0}
(\alpha,\lambda)]v_\lambda=
[\|\lambda\|^2 +2(\lambda,\rho)]v_\lambda.$$
\vskip .1in
\noindent \bf Freudenthal's Lemma. \it If $\mu$ is a weight of
$V_\lambda$ and $\nu$ 
is a weight of $V_\rho$, then 
$|\mu+\nu|^2\le |\lambda +\rho|^2$ with equality iff
$\mu=\sigma\lambda$ and $\nu=\sigma \rho$ for $\sigma\in W$,
necessarily unique.
\vskip .1in
\noindent \bf Proof. \rm Take $\sigma\in W$ such that
$\sigma^{-1}(\mu+\nu)\ge 0$. Since $\mu$ and $\nu$ are weights of $V_\lambda$
and $V_\rho$ respectively, we have $\mu_1=\sigma^{-1}\mu \le \lambda$ and
$\nu_1=\sigma^{-1} 
\nu \le \rho$, so that  $\lambda +\rho - \mu_1 -\nu_1$ is a sum of
positive roots. But then
$$0=\|\lambda+\rho\|^2-\|\mu_1+\nu_1\|^2=(\lambda +\rho -\mu_1-\nu_1,
\lambda+ \rho +\mu_1 + nu_1)\le (\lambda +\rho -\mu_1-\nu_1,\rho).$$
Thus $(\lambda-\mu_1,\rho)=0=(\rho-\nu_1,\rho)$ and hence
$\lambda=\mu_1 =\sigma^{-1} \mu$, $\rho=\nu_1=\sigma^{-1}\nu$, as
required. Uniqueness follows because $\sigma\rho=\rho$ imples
$\sigma=1$. 
\vskip .1in 
\noindent \bf Corollary. \it If $\mu$ is a weight of $V_\lambda$, then 
$|\mu+\rho|^2\le |\lambda +\rho|^2$ with equality iff
$\mu=\lambda$.
\vskip .1in
\noindent \bf Proof. \rm In this case $\nu=\rho$. On the other hand if
$\sigma\rho=\rho$, we must have $\sigma=1$. Hence
$\mu=\sigma\lambda=\lambda$. 

\vskip .1in

\noindent \bf 14. LIE ALGEBRAIC CONSTRUCTION OF IRREDUCIBLE
REPRESENTATIONS. 
\vskip .1in
\bf \noindent Generation by simple root vectors. \rm Let
$\alpha_1,\dots,\alpha_m$ be the simple positive roots and set
$E_i=E_{\alpha_i}$, $F_i=F_{\alpha_i}$ and $H_i=H_{\alpha_i}$. 
\vskip .1in
\noindent \bf Lemma. \it $\g$ is generated by the $E_i$'s and $F_i$'s
as a Lie algebra.
\vskip .1in
\noindent \bf Proof. \rm Let $\g_c$ be the complex Lie algebra
generated by all $E_i$, $F_i$ and hence $H_i=[E_i,F_i]$. Clearly
$\g_0$ is *--invariant and hence the complexification of its skew
adjoint part $\g_0$. Since the
$\alpha_i$'s are a basis of $\h^*$, the $H_i$'s are a basis of
$\h_{\Bbb C}$. Hence $\h\subset \g_0$. Let $g_i\in G$ be the Weyl group
element in the copy of $SU(2)$ corresponding to $\alpha_i$. Thus
$g_i\in \exp \g_0$ and ${\rm Ad}(g_i)$ is an automorphism of
$\g_0$. On the other hand $g_i$ is a representative of the simple
reflection $\sigma_i$ in $N(T)$ and the $\sigma_i$'s generate the Weyl
group. Hence $\g_0$ and $\g_c$ are invariant under $W$. But if $g\in
N(T)$ corresponds to $\sigma\in W$, ${\rm Ad}(g)\cdot \g_\alpha =
\g_{\sigma\alpha}$. Now by definition every root space $\g_{\pm
\alpha_i}$ lies in $\g_c$. Since every positive root is in the $W$--orbit
of a simple root, each root space $\g_\alpha$ lies in $\g_c$. Hence
$\g_c=\g_{\Bbb C}$ and $\g_0=\g$.

\vskip .1in
\noindent \bf Corollary. \it The $E_i$'s generate
$\bigoplus_{\alpha>0} \g_\alpha$ and the $F_i$'s generate
$\bigoplus_{\alpha<0} \g_\alpha$. 
\vskip .1in
\noindent \bf Proof. \rm Let $\g_+$ and $\g_-$ be the Lie algebras
generated by the $E_i$'s and $F_i$'s respectively. The relation
$[E_i,F_j]=\delta_{ij} H_i$
shows that $\g_+\oplus \h \oplus \g_-$ is a Lie subalgebra of
$\g_{\Bbb C}$. Since it contains $E_i,F_i,H_i$, it must be the whole
of $\g_{\Bbb C}$ so the result follows.

\vskip .1in
\bf \noindent Lemma (Serre relations). \it  The generators
$E_i,F_i,H_i$ satisfy the following relations:

\item{S1.} $[H_i,H_j]=0$.

\item{S2.} $[E_i,F_j]=\delta_{ij}H_i$.

\item{S3.} $[H_i,E_j]=n(i,j) E_j$ and $[H_i,F_j]=-n(i,j) F_j$ where
$n(i,j)= 2(\alpha_i,\alpha_j)/\|\alpha_i\|^2$. 

\item{S$^+_{ij}$.} ${\rm ad}E_i^{-n(i,j)+1}E_j=0$ for $i\ne j$.

\item{S$^-_{ij}$.} ${\rm ad}F_i^{-n(i,j)+1}F_j=0$ for $i\ne j$.
\vskip .1in
\noindent \bf Proof. \rm We already know S1, S2 and S3.
To prove the S$^-_{ij}$ ($i\ne j$), note that ${\rm ad}(E_i)\cdot
F_j=0$, ${\rm ad}(H_i)\cdot F_j=-n(i,j)F_j$. Thus the result follows
from $SU(2)_i$--theory, because $F_j$ is a highest weight vector. (In
particular $n(i,j)\le 0$.) S$^+_{ij}$ follows by taking adjoints.
\vskip .1in
\noindent \bf Remark. \rm The affine Kac--Moody algebra is given by
similar relations but this time indexed by the extended Cartan
matrix. We have to add an extra triple of generators
$E_0,F_0,H_0$. The theory of this chapter then proceeds almost without
change. 
\vskip .1in
\noindent \bf Verma module construction (induced modules). \rm Let $\g$ be a
Lie algebra, possibly infinite--dimensional, and $\g_1$, $\g_2$
subalgebras such that $\g=\g_1 \oplus \g_2$. If $W$ is any
finite--dimensional $\g_1$--module, then the Verma module is just
the induced module ${\cal U}(\g)\otimes_{{\cal U}(\g_2)} W$. Since
${\cal U}(\g)={\cal U}(\g_1)\otimes {\cal U}(\g_2)$ as vector space by
the Poincar\'e--Birkhoff--Witt theorem, the Verma module is also
isomorphic to ${\cal U}(\g_1)\otimes W$ as a vector space. This
description makes the action of $\g_1$ clear, but the action of $\g_2$
is harder to describe. We therefore give a more down to earth
computational recipe specialised to the case where $W$ is
one--dimensional. (The same arguments apply in general.) Let $f:\g_2
\rightarrow {\Bbb C}$ be a one--dimensional representation of
$\g_2$. Set $V={\cal U}(\g_1)$. We have to make $V$ into an ${\cal
U}(\g)$--module or equivalently a $\g$--module. 

Choose a basis $(b_i)$ of $\g_1$ consisting of monomials. By the
Poincar\'e--Birkhoff--Witt theorem, any element in ${\cal U}(\g)$ can
be written uniquely as $\sum b_i a_i$ with $a_i\in {\cal
U}(\g_2)$. Now take $x\in {\cal U}(\g)$. Then we have 
$xb_i=\sum b_j a_{ij}(x)$ with $a_{ij}(x)\in {\cal U}(\g_2)$. Let
$A(x)=a_{ij}(x))$ an infinite matrix with entries in ${\cal U}(g_2)$
with finitely many entries in any row. By uniqueness $A(xy)=A(x)A(y)$ for
$x,y\in {\cal U}(\g)$. In particular if $X,Y\in \g$, we have
$[A(X),A(Y)]=A([X,Y])$. Now define $x (b_i\otimes w) =\sum b_j \otimes
a_{ij}(x)w$. This is a representation because $A(xy)=A(x)A(y)$. By
definition the Verma module is a cyclic representation of $\g$
generated by a vector $w$ such that $Xw=f(X)w$ for $X\in
\g_2$. Conversely if $V^\prime$ is any other such cyclic
representation there is clearly a unique homorphism of $V$ onto
$V^\prime$ taking the cyclic vector $w$ onto the cylic vector
$w^\prime$. The homomorphism is given by $b\otimes w\mapsto bw^\prime$. 

\vskip .1in

Let $\g$ be the Lie algebra of a compact
semisimple Lie group $G$. For each $\alpha>0$ let
$E_\alpha,F_\alpha,H_\alpha$ be the basis of the 
Lie algebra $s\ell(2)_\alpha$
corresponding to the simple root $\alpha$. Let $\lambda\ge 0$ in
$P(\g)$ be a generalised highest weight. Let $\g_2=\h_{\Bbb C}\oplus  
\bigoplus_{\alpha>0} \g_\alpha$ and $\g_1=\bigoplus_{\alpha<0}
\g_\alpha$. These are Lie subalgebras of $\g_{\Bbb C}$ with
$\g_{\Bbb C}=\g_1\oplus \g_2$. Consider the 1--dimensional
representation sending $E_\alpha$ to $0$ and $H\in \h$ to
$i\lambda(H)$. Let $M(\lambda)$ be the corresponding Verma module.
Thus if $v=v_\lambda$ is the highest weight
vector of $M(\lambda)$, we have $E_\alpha v_\lambda=0$ and $ Hv_\lambda =
\lambda(H) v_\lambda$ where $\lambda(H_\alpha)\in {\Bbb Z}_+$ for all
$\alpha>0$. If $\alpha_1,\dots,\alpha_k$ is a numbering of the
positive roots, then a basis of $M(\lambda)$ is given by
$F_{\alpha_1}^{m_1} \cdots f_{\alpha_k}^{m_k} v_\lambda$. 
We know that $M(\lambda)$ has a unique maximal submodule $N$ such that
$L(\lambda)= M(\lambda)/N$ is irreducible as a $\g$--module. In fact,
since $\h$ is diagonalisable, every submodule is the sum of its weight
spaces. Hence if we take $N$ to be the algebraic sum of all proper
submodules, we muat have $v\notin N$, so that $N$ is the unique
maximal proper submodule.
By the $s\ell(2)$ theory, if $\ell_i=(\lambda,\alpha_i^\vee)$, then
$w_i=F_i^{\ell_i +1}v_\lambda$ is a singular 
vector i.e.~$E_iw=0$ and $w$ is an eigenvector for $\h$. It therefore
generates a proper submodule (all weights are strictly less than
$\lambda$). Hence $w_i\in N$ for all $i$. Let $N_0$ be the submodule
generated by the $w_i$'s. 
\vskip .1in
\noindent \bf Theorem (Harish--Chandra). \it $L(\lambda)$ is the
quotient of $M(\lambda)$ by the submodule generated by $F_i^{\ell_i
+1}v_\lambda$ and is finite--dimensional.
\vskip .1in
\noindent \bf Proof. \rm We have to show that $N=N_0$ and $L(\lambda)$ is
finite dimensional. Set $L=M(\lambda)/N_0$. Thus $L$ is a cyclic
module for $\g$ generated 
by $v=v_\lambda$ satisfying $X v=\lambda(X) v$ for $X\in \h$,
$E_\alpha v=0$ and 
$F_i^{\ell_i+1}v=0$. The identity 
$$[a^n,b]=\sum_{r=1}^n {n\choose r} [({\rm ad}\,a)^r b] a^{n-r}\eqno{(*)}$$
implies that the action on $L$ is
locally nilpotent, i.e.~some power of each $E_i$ or $F_i$ kills any
vector. For the $E_i$'s this follows because the $E_i$'s lower energy.
For the $F_i$'s it follows because $L$ is spanned by vectors
$F_{i_1}\cdots F_{i_k}v$ where $i_1,\dots,i_k$ are arbitrary (recall
that the $F_i$'s generate the $F_\alpha$ subalgebra). Starting from
the relation $F_i^{\ell_i+1}v=0$, 
successive application of $(*)$ and the Serre relations show that each
$F_i$ is nilpotent on any such monomial vector. This local nilpotence
shows that any vector in $L$ lies in a finite dimensional $s\ell(2)_i$
module for each $i$. 

We claim that the weights of $L$ are invariant under the Weyl group
$W$. In fact suppose $w\in L$ has weight $\mu$, $H_i w=m_i w$ with
$m_i=\mu(H_i)=(\mu,\alpha_i^\vee)$. Then 
$w$ lies in a sum of $s\ell(2)_i$ modules. If $m_i\ge 0$, set
$u=f_i^{m_i} w$ and if $m_i<0$, set $u=e_i^{-m_i}$. Thus $u\ne 0$ by
the $s\ell(2)$ theory and $u$ has weight $\lambda -m_i \alpha_i
=\sigma_i \lambda$. Thus the set of weights is invariant under each
simple reflection $\sigma_i$ and hence the whole of $W$.

To see that $L$ is finite--dimensional, we take the unique
$\sigma\in W$ such that $\sigma
C=-C$ (note that $-C$ is also a Weyl chamber and the Weyl group acts
simply transitively on these). By uniqueness $\sigma^2=1$. {\it We claim
that $\sigma \lambda$ is the lowest weight of $L$.} Since there are
plainly only finitely many weights $\mu$ such that $\sigma\lambda \le \mu \le
\lambda$, each of finite multiplicity, this proves that $L$ is
finite--dimensional. To prove the claim, note that $\alpha_i \mapsto
-\sigma\alpha_i$ 
must be a permutation of the simple roots because the walls of $C$ and $-C$
correspond to the same simple roots. (In particular $\sigma
\Phi^+=-\Phi^+$, so that $\sigma$ takes the positive roots onto
the negative roots.) Let $\mu\le \lambda$ be a weight of $L$. But then
$\sigma\mu$ is also a weight, so that
$\sigma\mu=\lambda -\sum n_i, \alpha_i$ with $n_i\ge 0$. Applying
$\sigma$, we get
$\mu=\sigma\lambda +\sum n{i^\prime} \alpha_i$, where
$\sigma\alpha_i =-\alpha_{i^\prime}$. Thus $\mu\ge \sigma\lambda$,
showing that $\sigma\lambda$ is the lowest weight.

\vskip .1in
Irreducibility of $L$ is now a consequence of the following result. 
\vskip .1in
\noindent \bf Casimir lemma. \it Let $V$ be a finite--dimensional cyclic
representation of $\g$ generated by a highest weight vector $v$ of
weight $\lambda$. Then $V$ is irreducible. 
\vskip .1in
\noindent \bf Proof. \rm Note that $V$ must be completely reducible
for each $\sl(2)_i$, so that as above the weights are integrable and
invariant under the Weyl group. Suppose that $V$ is not
irreducible. Then $V$ 
must contain a singular vector $w$ of weight $\mu$ strictly lower than
$\lambda$: thus $e_iv=0$ and $h_iw= m_iw$ where $m_i=\mu(h_i)\le
\ell_i$. But then $w$ is a highest weight vector generating an irreducible
representation of each $s\ell(2)_i$. On the other
hand let $\Omega$ be the Casimir operator of $\g$. Then $\Omega
v=(|\lambda +\rho|^2 -|\rho|^2)v$, so by cyclicity $\Omega=
(|\lambda +\rho|^2 -|\rho|^2)I$. Since $\Omega w=(|\mu +\rho|^2
-|\rho|^2) w$, we must have $|\lambda +\rho|^2 =|\mu +\rho|^2$. 
Choose $\sigma\in W$ such that $\sigma (\mu +\rho)\ge 0$. Then
$\sigma \mu <\lambda $ and $\sigma \rho <\rho$. But then by
Freudenthal's lemma $\lambda=\mu$, a contradiction. Hence $V$ is
irreducible. 

\vskip .1in
Finally we use the Casimir to argue that $V_\lambda=M(\lambda)$ admits
an invariant inner product. 
\vskip .1in
\noindent \bf Theorem. \it The representation $V_\lambda$ is
unitary. 
\vskip .1in
\noindent \bf Proof (Garland). \rm The representation
$\overline{V_\lambda}^*$ is irreducible with highest weight
$\lambda$. So there is an isomorphism $T:V_\lambda \rightarrow
\overline{V_\lambda}^*$, unique up to a scalar by Schur's lemma. Hence
we get an essentially unique $\g$--invariant sesquilinear form on
$V_\lambda$, $\phi(v,w)$. But $\phi^\prime(v,w)=\overline{\phi(w,v)}$
is another such form, so that $\overline{\phi(v,w)}=c \phi(w,v)$ for
some constant $c\in {\Bbb C}$. Clearly $|c|^2=1$. Multiplying $\phi$
by $a$ with $\overline{a}/a=c$, we get
$\overline{\phi(v,w)}=\phi(w,v)$ with $\phi$ non--degenerate. We claim
that $\pm\phi$ is positive definite. Clearly all the weight spaces are
orthogonal. Since the $\lambda$ weight space is 1--dimensional, we mus
have $\phi(v_\lambda,v_\lambda)\ne 0$. Since it is real, we may
rescale $\phi$ so that $\phi(v_\lambda,v_\lambda)=1$. To prove that
$\phi$ is positive definite, it suffices to show that $\phi(v,v)\ge 0$
for any weight vector. We prove this by downwards induction on the
weights under the usual ordering on weights.
In fact if $v$ is a weight vector of weight $\mu\le \lambda$, then
$$\phi(\Omega v,v)=(|\lambda+\rho|^2 -|\rho|^2) \phi(v,v)\eqno{(1)}$$
while
$$\phi(\Omega v,v ) = \sum \phi(T_iv,T_i,v)  -i \sum_{\alpha>0}
\phi(T_\alpha v,v)  + \sum_{\alpha>0} \|\alpha\|^2
\phi(E_\alpha v, E_\alpha v)=(|\mu +\rho|^2-|\rho|^2)\|v\|^2 +\sum_{\alpha>0}
\phi(E_\alpha v,E_\alpha v).\eqno{(2)}$$
Comparing (1) and (2) we get
$$(|\lambda +\rho|^2 -|\mu+\rho|^2)\phi(v,v) =\sum \|\alpha\|^2
\phi(E_\alpha v, E_\alpha v).$$
Assume $(v,v)=\phi(v,v)>0$. since $\mu<\lambda$, we have $\|\mu
+\rho|^2<|\lambda +\rho|^2$. Since the right hand side is
non--negative (by induction), we deduce that $\phi(v,v)\ge 0$. The
induction argument shows that $\phi(v,w)$ is positive
semi--definite. Since $\phi$ is non--degenerate, $\phi$ must be
positive definite, i.e.~a complex inner product.
\vskip .1in
\noindent \bf Proposition. \it Any finite dimensional representation $W$
of $\g$ is completely reducible. If $V$ is a non--trivial irreducible
representation $C\ne 0$ on $V$.
\vskip .1in
\noindent \bf Proof. \rm We may assume $C$ has only one eigenvalue on
$W$. Note that $W$ is completely reducible for 
each $\sl(2)_i$. This implies that $\h$ is diagonalisable and the
weights of $W$ are in $P(\g)$. Let $W_1$ be sum of all the irreducible
submodules of $W$. If $W_1\ne W$, find an irreducible subspace
$\overline{V}$ of $W/W_1$. Suppose $\overline{V}=V/W_1$ of weight $\mu$. Let
$\overline{v}$ be a highest weight vector in $\overline{V}$. Since
$\h$ is diagonalisable, we can lift $\overline{v}$ to an eigenvector $v$
of $\h$ of the same weight. We claim that $E_i v=0$ for all $i$.
If this is so, we may apply the Casimir Lemma to the cyclic module
generated by $v$. It must be irreducible. This contradicts the
maximality of $W_1$. Hence $W$ is completely reducible.

To prove the claim, observe that if $E_iv\ne 0$, then $\mu+\alpha_i$
is a weight of $W_1$, so is less than some highest weight
$\lambda$. Thus $\lambda\ge \mu \ge 0$ with $\|\lambda+\rho\|^2=\|\mu
+\rho\|^2$. By Freudenthal's lemma, $\lambda=\mu$, a contradiction. 
The last assertion is obvious from the formula
$\Omega|_{V_\lambda}=(\lambda+2\rho,\lambda)I$.

\vskip .1in
\bf \noindent 15. PROJECTIVE REPRESENTATIONS AND COVERING GROUPS. \rm
For each 
weight $\lambda\in P(\g)$ with $\lambda\ge 0$, we have constructed a
representation of $\g$ in ${\rm End}(V_\lambda)$ with
$\pi(X)^*=-\pi(X)$. Since $[\g,\g]=0$, each $\pi(X)$ has trace zero,
i.e.~$\pi(X)\in su(V_\lambda)$. Now consider the representation
$\pi_g(X)= \pi({\rm Ad}(g)\cdot X)$. This is also an irreducible
representation with highest weight $\lambda$, so there exists $U_g\in
SU(V_\lambda)$ unique up to an element of $Z_\lambda$, the finite
centre of $SU(V_\lambda)$, such that $\pi(gXg^{-1})=U_g \pi(X) U_g^*$. 
Let $G_1 =\{(u,g):\pi(gXg^{-1})=u \pi(X) u^{-1} \}$, a closed
subgroup of the compact matrix group $SU(V_\lambda)\times G$, so a
matrix group itself. Let $\widetilde{G}_\lambda = G_1^0$, the connected
component of the identity in $G_1$. This is also a matrix group. 
There is a natural homomorphsism $f:G_1\rightarrow G$ with kernel $G_1\cap
Z$. Since $G$ is connected and $f\exp(\g_1)=\exp(\g)$, we see that
$f(\widetilde{G}_\lambda)=f(G_1^0)=G$. Thus we have an exact sequence
$1\rightarrow C_\lambda \rightarrow \widetilde{G}_\lambda \rightarrow
G\rightarrow 1$ with $C_\lambda =Z_\lambda \cap
\widetilde{G}_\lambda$, a cyclic central subgroup. The Lie algebras of
$\tilde{G}_\lambda$ and $G$ can naturally be identified. Moreover
$P(G_\lambda)=\langle P(G),\lambda\rangle$. 

In general if $\pi$ is an irreducible projective unitary
representation of $G$, the same argument shows that we can find
$\tilde{G}$ a connected compact matrix group with a covering
homomorphism $f:\tilde{G}\rightarrow G$ and a representation
$\tilde\pi$ such that $\tilde\pi =\pi \circ f$ as projective
representations. 

Let $\lambda_i$ be a choice of cosets for $P(\g)/P(G)$ and let $G_i$
be the central extension (by $C_i$) constructed above with projection
maps $f_i:G_i\rightarrow G$. Let $\tilde G$ be the connected component
of the identity of the closed subgroup of
elements $(x_i)\in\prod G_i$ such that $f_i(x_i)=f_j(x_j)$. Then
there is a natural projection $f:\tilde{G}\rightarrow G$ with kernel
$Z=\prod Z_i \cap \tilde{G}$, finite Abelian. By construction
$\lambda_i\in P(\tilde G)$, so that $P(\tilde G)=P(\g)$. Hence
$Z(\tilde G) = P(\tilde G)/Q$. The group $\tilde G$ has the property
that every generalised weight is now a weight.
\vskip .1in
\noindent \bf Remark. \rm Weyl's theorem implies that $\tilde G$ is
simply connected.
\vskip .2in
\bf \noindent 16. THE DIRAC OPERATOR AND SUPERSYMMETRY RELATIONS. \rm
Our aim now is to prove Weyl's character formula for the character of
an irreducible representation with highest weight $\lambda$. The proof
is based on supersymmetry and the coset construction of
Goddard--Kent--Olive. The supercharge or Dirac operator and its
properties lie at the heart of the method which is manifestly
unitary. Later we will 
explain why it also gives a geometric construction of all irreducible
representations as twisted harmonic spinors. The Dirac operator or
supercharge operator also exists for affine algebras, as shown by
Kazama and Suzuki. The same technique can therefore be used for to
prove the Weyl--Kac character formula. 
\vskip .1in
\noindent \bf The spin representation. \rm Let $V$ be a real inner
product space with orthonormal basis $(e_i)$. 
For $T\in so(V)$, define $s(T)={1\over 4} \sum c(T\cdot e_i) c(e_i)$.
Recall that $[s(T),c(X)]=c(TX)$.
\vskip .1in
Let $G$ be a
compact simple matrix group with Lie algebra $\g$. Let $T$ be a
maximal torus in $G$ with Lie algebra $\h$. We assume that the
invariant inner product on $\g$ has been normalised so that the
highest root $\theta\in \h^*$ has $\|\theta\|^2=2$ in the induced
norm. Let $\m=\h^\perp\subset \g$. Thus $\m$ is even dimensional and
invariant under ${\rm Ad}(T)$ and ${\rm ad}(\h)$. Since it is an inner
product space, we may consider the real Clifford algebra ${\rm
Cliff}\, \m$. It has a unique irreducible representation $W_\m$ which
is ${\Bbb Z}_2$--graded: $W_\m=W^+_\m \oplus W^-_\m$. 
Let $V$ be an irreducible representation of $G$. We
define the {\it supercharge} or {\it Dirac operator} on $V\otimes W_\m$ by 
$$Q=\sum (\pi(X_i)+{1\over 3}s(\overline{X}_i))c(X_i).$$
By definition $Q$ takes $V\otimes W_\m^\pm $ into $V\otimes W_m^\mp$
and commutes with $T$. As we shall explain below, this operator is
really the Dirac operator for a very special connection restricted to
an isotypic subspace. However to understand why we take this
particular formula we will have to take a supersymmetric path rather
than a geometric one. 
\vskip .1in

\noindent \bf A. Computations for ${\rm Cliff}(\g)$. \rm Take
generators $c(X)$ ($X\in\g$) for ${\rm Cliff}(\g)$. Thus
$\{c(X),c(Y)\}= 2(X,Y)$. Let $Q_0=\sum s(X_i) c(X_i)$ so that $s(X)
={1\over 4} \sum 
c([X,X_i])c(X_i)$. Evidently $Q_0$ is independent of the choice of
orthonormal basis, so that $s(g) Q_0 s(g)^{-1} =Q$ since
$s(g)X_is(g)^{-1}$ is another orthonormal basis of
$\g$. Differentiating this relation, we get 
$$[Q_0,s(X)]=0.\eqno{(1)}$$
This is one of the supersymmetry relations; we get the other as
follows:
$$\{Q_0,c(X)\}
= \sum s(X_i) \{c(X_i),c(X)\} -[s(X_i),c(X)]c(X_i)
=\sum s(X_i) 2(X_i,X) +c([X,X_i])c(X_i)
=6s(X).$$
Thus 
$$\{Q_0,c(X)\}=6s(X).\eqno{(2)}$$

\noindent \bf B. Computations in ${\rm End}(V)\otimes {\rm
Cliff}(\g)$. \rm Let $V$ any unitary $G$--module and set
$$\widetilde{Q}=\sum (\pi(X_i) +{1\over 3} s(X_i)) c(X_i).$$
Then (1) implies that $[\widetilde{Q}, \pi(X)+s(X)]=0$,  while (2)
implies
$$\{\widetilde{Q},c(X)\}={1\over 3} \{Q,c(X)\} +\sum \pi (X_i)
\{c(X_i),c(X)\} = 2(\pi(X) + s(X)).$$
\vskip .1in
\noindent \bf C. Coset construction of the Dirac operator. \rm If
$(X_i)$ is an orthonormal basis of $\m$ and $U$ is an irreducible
representation of $G$, then the Dirac (or supercharge) operator on
$U\otimes W_\m$ is given by
$$Q=\sum c(X_i) (\pi(X_i) + {1\over 3} s(\overline{X_i})).$$
We can consider the same operator acting on $U\otimes W_\m\otimes
W_\h=U\otimes W_\g$ and denote it by the same symbol (strictly
speaking it should be $Q\otimes I$). The supercharge operator $Q_\g$
also acts on this space as does the supercharge operator for $\h$
acting on $U^\prime\otimes W_\h$, where $U^\prime =U\otimes W_\h$. 
These operators are given by
$Q_\g=\sum_a (\pi(X_a)+{1\over 3} S(X_a)) c_a$ and $Q_\h=\sum
(\pi(X_A) +s_\m(X_A))c(X_A)$ since $s_\h(X)=0$ for $X\in\h$. 
\vskip .1in
\noindent \bf Theorem (coset construction). \it $Q=Q_\g - Q_h$.
\vskip .1in
\noindent \bf Proof. \rm Let $(X_a)$ be an orthonormal basis of $\g$,
made up of an orthonormal basis $(X_i)$ of $\m$ and $(X_A)$ of
$\h$. Thus we have structure constants $f_{abc}$ given by
$[X_a,X_b]=\sum_c f_{abc} X_c$. The invariance of the inner product
and the orthonormality of $(X_a)$ imply that $f_{abc}$ is {\it totally
antisymmetric} in its three arguments. In particular we have
$s(X_a) ={1\over 4} \sum f_{abc} c(X_a) c(X_b)$. Thus
$$\eqalign{ Q_\g - Q_\h& = \sum (\pi(X_a) +{1\over 3} s_\g(X_a))c(X_a)
-\sum (\pi(X_A) + s_\m(X_A)) c(X_A)\cr
& =\sum \pi(X_a) c(X_a) +{1\over 12} \sum_{a,b,c} c(X_a) c(X_b) c(X_c) 
-\sum \pi(X_A) c(X_A) -{1\over 4} \sum_{A,i,j} f_{Aij} c(X_A) c(X_i) c(_j)\cr
&=\sum \pi(X_i) c(X_i) +{1\over 12} \sum_{i,j,k} c(X_i)c(X_j)c(X_k)\cr
&=\sum (\pi(X_i) +{1\over 3} s_\m(\overline{X}_i))c(X_i),\cr}$$ 
since there are three ways that terms $c(X_A)c(X_i)c(X_j)$ can appear
and $f_{abc}=0$ if two or more coefficients corresponds to basis
elements in $\h$.

\vskip .2in
\bf\noindent 17. THE SQUARE OF THE DIRAC OPERATOR. \rm 
\vskip .1in
\noindent \bf A. Lichnerowicz's lemma. \rm Let $c_i=c(X_i)$, so that
$\{c_i,c_j\}=2\delta_{ij}$, and let $R_{ijk\ell} = ([X_i,X_j],
[X_k,X_\ell])$. Apart from the obvious symmetry
properties
$$R_{ijk\ell} = R_{k\ell ij} = -R_{jik\ell} = - R_{ij\ell k},$$
the tensor $R_{ijk\ell}$ satisfies Bianchi's
first identity, namely
$$R_{ijk\ell} + R_{ik\ell j} + R_{i\ell jk} =0.$$
This follows immediately from the Jacobi identity, since
$R_{ijk\ell}=(X_i,[X_j,[X_k,X_\ell]])$. 
\vskip .1in
\noindent \bf Lemma (Lichnerowicz). \it $\sum_{i,j,k,l}
R_{ijkl}c_ic_jc_kc_l = -2\sum_{i,j} R_{ijij}$.
\vskip .05in
\noindent \bf Proof. \rm By Bianchi's identity we have
$$\eqalign{\sum_{i,j,k,\ell} R_{ijk\ell} c_ic_jc_kc_\ell& = 
-\sum_{i,j,k,\ell} R_{ijk\ell} c_i(c_kc_\ell c_j + c_\ell c_jc_k)\cr
& = 2\sum_{i,j,k,\ell} R_{ijk\ell} c_i(\delta_{jk}c_\ell +\delta_{k\ell}c_j
     -2\delta_{j\ell} c_k - c_jc_kc_\ell).\cr}$$
The symmetry properties allow one to delete the $\delta_{k\ell}$ term and,
upon rearrangement, we get
$$\eqalign{\sum_{i,j,k,\ell} R_{ijk\ell}c_ic_jc_kc_\ell &= 
{2\over 3} \sum_{i,j,k,\ell} R_{ijk\ell}(\delta_{jk}c_ic_\ell - 2\delta_{j\ell}
c_ic_k)\cr 
& = 2\sum_{ij\ell} R_{ijj\ell} c_ic_\ell
 = 2\sum_{ijl} R_{ijj\ell} (c_ic_\ell+c_\ell c_i)/2
 = 2\sum_{ij} R_{ijji}\cr}$$
as required.
\vskip .1in
\noindent \bf B. Computations in ${\rm Cliff}(\g)$. \rm If $(X_i)$ is
an orthonormal basis of $\g$, 
let $Q=\sum s(X_i) c(X_i)$ in ${\rm 
Cliff}(\g)$. 
\vskip .1in
\noindent \bf Lemma~1. \it $Q_0^2$ is a scalar operator.
\vskip .1in
\noindent \bf Proof. \rm The supersymmetry relations
$\{Q_0,c(X)\}=6 s(X)$ and $[Q_0,s(X)]=0$ imply
$$[Q_0^2, c(X)]=\{Q_0,c(X)\}Q_0 -Q_0\{Q_0,c(X)\}=6 s(X)Q_0 - Q_0s(X) =0.$$
So $Q_0^2$ is central in ${\rm Cliff}(\g)$. Since $Q_0^2$ lies in
${\rm Cliff}^+(\g)$, it must be a scalar operator.
\vskip .1in
\noindent \bf Corollary. \rm $Q_0^2=3 \sum s(X_i)^2$, so $\sum s(X_i)^2$
is a scalar.
\vskip .1in
\noindent \bf Proof. \rm We have
$$\{Q_0,Q_0\} = \{Q_0,\sum s(X_i)c(X_i)\} = \sum [Q_0,s(X_i)]c(X_i) + 3
\sum s(X_i)^2= 6 \sum s(X_i)^2,$$
so $\sum s(X_i)^2$ is a scalar, as claimed. 
\vskip .1in
\noindent \bf Lemma~2. \it $\sum s(X_i)^2 = -{1\over 8} \sum
\|[X_i,X_j]\|^2 =-{1\over 8} {\rm dim}(\g)\cdot h_\g$, where the dual
Coxeter number $h_\g$ equals the value of the Casimir in the adjoint
representation. 
\vskip .05in
\noindent \bf Remark. \rm Thus $h_\g=-\sum {\rm Tr}({\rm
ad}(X_i)^2)=\|\theta+\rho\|^2 -\|\rho\|^2$, where $\theta$ is the
highest root, i.e~the highest weight in the adjoint representation.
\vskip .05in
\noindent \bf Proof. \rm We have
$$s(X_a)={1\over 4} \sum ([X_a,X_j],X_i)c_i c_j
=-{1\over 4} \sum (X_a,[X_i,X_j])c_ic_j,$$
so that 
$$\eqalign{\sum s(X_i)^2&={1\over 16} \sum
([X_i,X_j],[X_k,X_\ell])c_jc_jc_kc_\ell\cr
&=-{1\over 8}\sum \|[X_i,X_j]\|^2\cr
&=-{1\over 8}\sum ([X_i,X_j],[X_i,X_j])\cr
&={1\over 8} \sum ({\rm ad}(X_i)^2 \cdot X_j,X_j)\cr
&={1\over 8} {\rm Tr} (\sum {\rm ad}(X_i)^2)\cr
&={1\over 8} {\rm Tr} (-h_\g I)\cr
&={1\over 8} {\rm dim}(\g)\cdot h_\g.\cr}$$
\vskip .1in 
\noindent \bf Corollary. \it $Q_0^2=-{3\over 8} h_\g \cdot {\rm
dim}(\g)$.
\vskip .1in
\bf \noindent Proof. \rm We have $Q_0^2=3\sum s(X_i)^2 = -{3\over 8}
h_\g\cdot {\rm dim}(\g)$.  
\vskip .1in
\noindent \bf Lemma~3. \it The weights of $W_\m^{\pm}$ as a representation
of $\h$ are exactly $\rho-\sum_{\alpha\in S} \alpha$, where $S$
is an arbitrary subset of positive roots with $|S|$ even for $W_\m^+$
and odd for $W_\m^-$. Thus ${\rm ch}_s (W_\m) = {\rm ch}(W^+_\m) -{rm
ch}(W^-_\m) = e^\rho \prod_{\alpha>0} (1-e^{-alpha})$. 
\vskip .1in
\noindent \bf Proof. \rm This follows immediately from the formulas for
the weights of the spin representations. Recall that every matrix
$g\in SO(N)$ is conjugate to a matrix with $2\times 2$ 
blocks $\pmatrix{\cos\theta_j & \sin\theta_j \cr -\sin\theta_j
&\cos\theta_j\cr}$ down the diagonal where $j=1,\dots,[N/2]$. Note that
there is an additional $1$ on the diagonal if $N$ is odd. Thus the
complex eigenvalues of $g$ are $e^{\pm i \theta_j}$. In any irreducible
projective representation $\pi(g)$ of $SO(N)$, the eigenvalues of
a generic block diagonal element $\pi(g)$ are called the {\it weights}
of the representation.
\vskip .1in
\noindent \bf Lemma~3$^\prime$ (weights of spin representation). \it (1)
If ${\rm dim}(V)$ is even, the weights of the spin representation $W^\pm$ are
$\exp{i {1\over 2} \sum \pm \theta_k}$ where the number of plus signs is even
for $W^+$ and odd for $W^-$.

\noindent (2) If ${\rm dim}(V)$ is odd, the weights of the spin
representation $W$ are $\exp{ {1\over2} \sum \pm i\theta_k}$. 
\vskip .1in
\bf \noindent Proof. \rm Let $v_1,\dots ,v_m$ be an orthonormal basis
for $V_J$ and 
set $v_I=v_{i_1}\wedge \cdots \wedge v_{i_k}$ for
$i_1<\cdots<i_k$. Then
$c(v_k)=e(v_k) + e(v_k)^*$ and $c(iv_k) = i(e(v_k)-e(v_k)^*)$. The
generators of Lie algebra of the maximal torus are given by
$T_j={1\over 4} (c(v_j)c(iv_j) - c(iv_j)c(v_j))={1\over 2} c(v_j)
c(iv_j)$. Thus $T_j v_I = i/2 \, v_I$ if $j\in I$ and $T_j v_I =-i/2\,
v_I$ if $j\notin I$. The corresponding self--adjoint operators 
$S_j$ satisfy $T_j=iS_j$ so that $S_jv_I=\pm {1\over 2} v_I$. Since
$W^+=\Lambda^{\rm even}(V)$ is spanned by $v_I$'s with $I$ even and 
$W^-=\Lambda^{\rm odd}(V)$ is spanned by $v_I$'s with $I$ odd, the
result follows. 

\noindent (2) Note that the maximal torus of $SO(V_0)$ coincides with
the maximal torus of $SO(V)$. On the other hand the spin
representation of $SO(V)$ equals $W^+\oplus W^-$, where $W^\pm$ are
the spin representations of $SO(V_0)$. So the result follows
immediately from (1).
\vskip .1in
\bf \noindent Corollary of Lemma~3. \it As a representation of $\g$, the
irreducible representation 
$W$ of ${\rm Cliff}(\g)$ contains a highest weight
vector of weight $\rho$; in fact it contains exactly $N=2^{M}$
such vectors where $M=[m/2]$ with $m={\rm dim}\, \h$. 
\vskip .1in
\noindent \bf Proof. \rm We have $W_\g=W_\m\otimes W_\h$ with $\h$
acting trivially On $W_\h$ and with highest weight $\rho$ on
$W_\m$; for as we have seen the weights of $W_\m$ are exactly
$\rho -\sum_{\alpha \in S} \alpha$, where $S$ is a subset of the
positive roots. Clearly this means that there are exactly $N={\rm
dim}(W_\h)$ vectors of weight $\rho$ in $W$. 
\vskip .1in
\noindent\bf Corollary~2 (Kostant, Atiyah--Schmid, Wallach, \dots). \it
$W=V_\rho\otimes {\Bbb C}^N$. In particular we get the ``strange''
isomorphism $W_\m\cong V_\rho$ as $\h$ or $T$ modules, so the weights
and weight multiplicities of $W_\m$ are invariant under $W$. 

\vskip .1in
\noindent \bf Remark. \rm Invariance of the weights and weight
multiplicities of $W_\m$ under the Weyl group also
follows because $N(T)$ acts on $\m$ and hence $W_\m$ (by
quantisation); it is not hard to show that $N(T)$ permutes the weight
spaces according to the corresponding Weyl group elements.
\vskip .1in
\noindent \bf Proof. \rm Clearly it suffices to show that all
highest weight vectors in $W$ have weight $\rho$. We have already seen
that the Casimir of $\g$ has value $3\|\rho\|^2$ on $W$, so any other
highest weight $\lambda=\rho-\sum_{\alpha\in S} \alpha=\rho-\nu$ must satisfy
$\|lambda+\rho\|^2 =4\|\rho^2\|$. But
$$\|\lambda+\rho\|^2 = 4\|\rho\|^2 -(\nu,rho) - 3(\lambda,\nu),$$
with $(\nu,\rho)\ge 0$ and $(\lambda,\nu)\ge 0$. Thus
$\|lambda+\rho\|^2 =4\|\rho^2\|$ only if $(\nu,\rho)=0$ and
$(\lambda,\nu)= 0$, so that $\nu=0$ and $\lambda=0$.
\vskip .1in
\noindent \bf Lemma~4. \it $\sum s(X_i)^2=3\|\rho\|^2$.
\vskip .1in
\noindent \bf Proof. \rm Let $W$ be an irreducible representation of
${\rm Cliff}(\g)$. Now $W=W_\m\otimes W_\h$ and $\h$ acts trivially on
$W_\h$. On the other hand $W_\m =\Lambda^*(\m_{\Bbb C}^+)$, so that
$\rho$ is a highest weight of $W$ as a $\g$--module. In fact
although $\h$ is canonically quantised, $\pi(T)={1\over 4} \sum
c([T,X_i])c(X_i)$ does not agree with the canonical quantisation. A highest
weight vector is given by $v=\bigwedge_{\alpha>0} X_\alpha$ and $T
X_\alpha =\alpha(T) X_\alpha$, so that $\pi(T)v=2\rho(T)v$. On the
other hand $\pi(T)\Omega=0$, so that $\Omega$ is a lowest weight
vector. Since the canonical quantisation has the form $\pi(T) +
\lambda(T)$ and the highest and lowest weight of $\pi$ have the form
$\pm \mu(T)$ for some $\mu>0$, we must have
$\lambda(T)=\mu(T)=\rho(T)$. The Casimir $\sum s(X_i)^2$ therefore acts on the
corresponding submodule as $3\|\rho\|^2$. Since it acts as a scalar on
$W$, the lemma follows.
\vskip .1in
\noindent \bf Corollary (Freudental--de Vries ``strange
formula''). \it $\|\rho\|^2= h_\g\cdot {\rm dim}\, \g/12$.
\vskip .1in
\bf \noindent Proof. \rm Immediate from lemmas~2 and 3. 
\vskip .1in
\noindent \bf C. Computations in ${\rm End}(V)\otimes {\rm
Cliff}(\g)$. \rm Let $\pi:\g\rightarrow so(V)$ be a representation of
$\g$ and set $Q_\g=\sum (\pi(X_i) +{1\over 3}
s(X_i))c(X_i)$. 
\vskip .1in
\noindent \bf Theorem. \it $Q_g^2=\sum \pi(X_i)^2 +{1\over 3}\sum
s(X_i)^2$. 
\vskip .1in
\noindent \bf Proof. \rm Let $Q_0=\sum s(X_i) c(X_i)$. Thus
$[Q_0,s(X)]=0$ and $\{Q_0,c(X)\}=6s(X)$. Moreover
$[Q_\g,\pi(X)+s(X)]=0$ and $\{Q,c(X)\}=2\pi(X)+
+2 s(X)$. Hence
$$\{Q_0,Q_0\}=
\{Q_\g, \sum (\pi(X_i)+s(X_i))c(X_i) + ({1\over 3}
-1)Q_0\}.\eqno{(1)}$$
The first term gives 
$$\{Q_\g, \sum (\pi(X_i)+s(X_i))c(X_i)\}
=\sum (\pi(X_i)+s(X_i)) (2\pi(X_i) + {2\over 3} s(X_i))\eqno{(2)}$$
while the second term gives
$$\{Q_\g, Q_0\}=\sum \pi(X_i)\{Q_0,c(X_i)\} + \varepsilon
\{Q_0,Q_0\} =6\sum \pi(X_i) s(X_i) + 2 \sum
s(X_i)^2.\eqno{(3)}$$
The result follows by substituting in (1) from (2) and (3). 
\vskip .1in
\noindent \bf D. Coset computation of $Q^2$. \rm  By the coset construction 
$Q=Q_\g-Q_\h$. Thus $Q_\g=Q+Q_\h$, where $Q$ and $Q_\h$ anticommute,
i.e.~$\{Q,Q_\h\}= 0$. Hence 
$$\{Q_\g,Q_\g\}=\{Q+Q_\h,Q+Q_\h\} = \{Q,Q\} + \{Q_\h,Q_\h\}.$$
Thus $Q_\g^2 = Q^2 +Q_\h^2$, so that $Q^2=Q_\g^2-Q_\h^2$.

\vskip .1in

\bf \noindent 18. WEYL'S CHARACTER AND DENOMINATOR FORMULAS. \rm 

\vskip .1in

\noindent \bf Lemma. \it If $\sigma\in W$, then
$\varepsilon(\sigma)={\rm det}_\h(\sigma)=(-1)^{n(\sigma)}$ where
$n(\sigma)=|\{\alpha>0: \sigma \alpha<0\}|$. 
\vskip .1in
\noindent \bf Proof. \rm Since $G$ is connected, ${\rm Ad}(G)\subset
SO(\g)$. Hence $N(T)\subset SO(\g)$. But $N(T)$ normalises $\h$ and
therefore acts on $\m$. If $g\in N(T)$ it follows, that ${\rm det}_\m
(g) ={\rm det}_\h(g)$. Since ${\rm det}_\g(t)={\rm det}_\h(t)=1$, we
have ${\rm det}_\m(t)=1$ for $t\in T$. Thus if $\sigma\in W=N(T)/T$,
${\rm det}_\m(\sigma)$ is well--defined and equals ${\rm
det}_\h(\sigma)$. Let $X_\alpha,Y_\alpha$ be a basis for $\m$. Clearly
$\sigma X_\alpha \wedge \sigma Y_\alpha$ is
well--defined (independently of the choice of representative in
$N(T)$). If $\sigma\alpha>0$, it equals $X_{\sigma\alpha}\wedge
Y_{sigma\alpha}$, while if $\sigma\alpha<0$ it equals
$-X_{-\sigma\alpha} \wedge Y_{-\sigma\alpha}$. Thus ${\rm
det}_\h(\sigma) =(-1)^N$ where $N$ is the number of $\alpha$ such that
$\sigma\alpha<0$, as required.

\vskip .1in
\noindent \bf Corollary. \it $\pi(\sigma)\in {\rm End}(W_\m)$ is even
or odd according as $\varepsilon(\sigma)=1$ or $-1$. If $\xi\in W_\mu$
has weight $\mu$, then $\pi(\sigma)\xi$ has weight $\sigma\mu$.
\vskip .1in
\noindent \bf Proof. \rm We already know that if $T\in O(\m)$,
then $\pi(T)$ is even or odd according to the sign of ${\rm det}(T)$,
so the first result follows. Thus $\pi(\sigma) \pi(t)
\pi(\sigma)^{-1} =\pm \pi(\sigma t)^{-1}$ where we regard $\pi$ as a
representation of $T$ (recall that $\rho$ is a weight). Since $T$ is
connected, only the plus sign is possible, so the result follows.
\vskip .1in
\noindent \bf Lemma (Euler--Poincar\'e Principle). \it 
Let $W=W^+\oplus
W_-$ be a vector space ${\Bbb Z}_2$--graded and let $A$ and $B$ be
even and odd commuting operators on $W$ with $B$ diagonalisable. Then
${\rm Tr}_s\, A = {\rm Tr}_s A|_{{\rm ker}(B)}$.
\vskip .1in
\noindent \bf Proof. \rm Since $B$ is diagonalisable, so is
$B^2$. Moreover ${\rm ker}(B)={\rm ker}(B^2)$. Let
$W^\pm_\lambda=\{\xi\in W^\pm : B^2 \xi =\lambda \xi\}$. If
$\lambda\ne 0$, then $B$ gives an isomorphism between $W^+\lambda$ and
$W^-_\lambda$. Since $A$ commutes with $B$, $A$ leaves $W^\pm_\lambda$
and the isomorphism given by $B$ intertwines the two actions of
$A$. Hence ${\rm Tr}_{W_\lambda^+} (A) ={\rm Tr}_{W_\lambda^-} (A)$
for $\lambda\ne 0$. Hence
$${\rm Tr}_s\, A ={\rm Tr}_{W_+}(A) -{\rm Tr}_{W^-}(A)=\sum_\lambda
{\rm Tr}_{W^+_\lambda}(A) -{\rm Tr}_{W^-_\lambda}(A)
={\rm Tr}_{W_0^+}(A) -{\rm Tr}_{W_0^-}(A)=
 {\rm Tr}_s A|_{{\rm ker}(B)},$$
as required. 
\vskip .1in
\noindent \bf Theorem (Weyl's character formula). \it ${\rm
ch}(V_\lambda) = \Pi ^{-1} \sum_{\sigma\in W} \varepsilon(\sigma)
e^{\sigma(\lambda +\rho)-\rho}$, where $\Pi=\prod_{\alpha>0}
(1-e^{-\alpha})$.

\vskip .1in
\noindent \bf Proof. \rm Consider the $\h$ or $T$ module
$V_\lambda\otimes W_\m$. Evidently
$${\rm ch}_s V_\lambda\otimes W_\m= {\rm ch}(V) \cdot {\rm ch}_s(W_\m)
={\rm ch}(V_\lambda)\cdot e^\rho \prod_{\alpha>0} (1-e^{-\alpha}).\eqno{(1)}$$
On the other hand, by the Euler--Poincar\' e principle,
$${\rm ch}_s V_\lambda \otimes W_\m = {\rm ch}_s {\rm ker}(Q)\cap
(V_\lambda\otimes W_\m).\eqno{(2)}$$
Since $Q$ commutes with $\h$, we can assume $\xi\in 
{\rm ker}(Q)\cap
(V_\lambda\otimes W_\m)$ is an $\h$--eigenvector of weight $\mu +\nu$
say, where $\mu<\lambda$ and $\nu<\rho$. Since $Q$ is skew--adjoint,
${\rm ker}(Q)={\rm ker}(Q)$. But
$Q^2\xi=(\|\lambda+\rho\|^2-\|\mu+\nu\|^2)\xi$. Therefore we must have
$\mu=\sigma \lambda$ and $\nu=\sigma\rho$ for some unique $\sigma\in
W$. Note that if $g\in N(T)$ corresponds to $\sigma\in W$, then
$g v_\lambda\otimes \pi(\sigma)v_\rho$ has weight $\sigma\lambda +
\sigma \rho$. We claim, up to scalar multiples this is the only vector
with this weight. In fact suppose that $\mu_1+ \nu_1 =
\mu + \nu$
with $\lambda \ge \mu_1$ and $\rho \ge 
\nu_1$. Since $|\mu_1 +\nu_1|^2 = |
\mu+\nu|^2 =|\lambda+\rho|^2$, the previous argument
implies that $ \mu_1=\tau \lambda$ and
$ \nu_1=\tau \rho$ for some $\tau\in  W$.
But then $\gamma=\tau^{-1}\sigma$ fixes $\lambda
+\rho$. Since $\lambda\ge \gamma  \lambda$ and
$\rho\ge \gamma\rho$, we get
$\gamma\rho=\rho$ so that $\gamma={\rm id}$. Thus $\mu=\mu_1$ and
$\nu=\nu_1$, so that the vector must lie in the tensor product of the
$\mu$--weight space of $V_\lambda$ and the $\nu$--weight space of
$W_\m$. But each of hese weight spaces is obtained by applying
$\sigma$ to the $\lambda$ and $\rho$ weight spaces; they therefore have
multiplicity one. Thus the kernel of $Q$ (or equivalently $Q^2$) is
indexed by elements of $\sigma$ and has a basis consisting of vectors 
$g_\sigma v_\lambda \otimes \pi(g_\sigma)v_\rho$. The vector
$\pi(g_\sigma)v_\rho$ lies in $W^\pm_\m$ according as
$\varepsilon(\sigma)=\pm 1$. Thus 
$${\rm ch}_s({\rm ker}(Q)\cap(V_\lambda\otimes W_\m))= \sum_{\sigma\in
W} \varepsilon(\sigma) e^{\sigma(\lambda+\rho)}.\eqno{(3)}$$
The character formula follows from (1), (2) and (3). 

\vskip .1in
\noindent \bf Corollary (Weyl's denominator formula). \it
$\sum_{\sigma\in W} \varepsilon(\sigma)
e^{\sigma\rho-\rho}=\prod_{\alpha>0} 
(1-e^{-\alpha})$.
\vskip .1in
\noindent \bf Proof. \rm This follows by setting $\lambda=0$, since
the character of the trivial representation is identically $1$.

\vskip .1in
\noindent \bf Remark. \rm Using the denominator formula, we can write
${\rm ch}(V_\lambda)=A(\lambda+\rho)/A(\rho)$ where $A(\mu)= \sum
\varepsilon(\sigma) e^{\sigma\mu}$.

\vskip .1in
\noindent \bf Corollary (Weyl's dimension formula). \it  
${\rm dim}(V_\lambda) =\prod_{\alpha>0}
(\lambda+\rho,\alpha)/\prod_{\alpha >0} (\rho,\alpha)$.
\vskip .1in
\noindent \bf Proof. \rm Let $X,Y\in \h$ be the elements such that
$\mu(X)=(\mu,\rho)$ and $\mu(Y)=(\mu,\rho+\lambda)$ for $\mu\in
\h^*$. Then by Weyl's denominator formula
$$\sum \varepsilon(\sigma)e^{i(\rho+\lambda)(Xt)}=\sum
\varepsilon(\sigma) e^{i\rho(Yt)} =\prod_{\alpha>0} (e^{i\alpha(Y)t}
-e^{-i\alpha(Y)t}),$$
and
$$\sum \varepsilon(\sigma) e^{i\rho(Xt)}=\prod_{\alpha>0} (e^{i\alpha(X)t}
-e^{-i\alpha(X)t}).$$
Dividing these we get
$${\rm Tr}_{V_\lambda} (e^{Xt}) = \prod_{\alpha>0} {\sin
(\lambda+\rho,\alpha)t \over \sin(\rho,\alpha)t}.$$
The result follows by letting $t\rightarrow 0$.

\vskip .1in
\bf \noindent 19. REMARKS ON CONNECTIONS AND DIRAC OPERATORS. \rm Let $G$ be a
with simple compact Lie group with Lie 
algebra $\g$ with invariant inner product $(X,Y)$. Let $H$ be a closed
subgroup with Lie algebra $\h$ and let $\m=\h^\perp$, so that
$\g=\h\oplus \m$. Set $M=G/H$. Then ${\rm Vect}(M)=
(C^\infty(G)\otimes \m)^H$, so that $X(gh)={\rm ad}(h)^{-1} X(g)$ for
$X\in {\rm Vect}(M)$. If $P=P_\m$ is the orthogonal projection onto
$\m$ and $Y\in \m$, then $\widetilde Y(g)=P(g^{-1} Y g)$ defines a
vector field with $\widetilde Y(1) =Y$; translating on the left, we
can produce a vector field equal to $Y$ at a given point. If $X(g)$ is
a vector field and $f\in C^\infty(M)=C^\infty(G)^H$, we define
$$Xf(g) = {d\over dt} f(g\exp(X(g)t))|_{t=0}.$$
It is immediate from the definitions that $Xf(gh)=Xf(g)$, so that
$Xf\in C^\infty(M)$. The canonical connexion is defined by
$$\nabla_X\xi(g) = {d\over dt} \xi(g\exp(X(g)t))|_{t=0},$$
for $\xi\in {\rm Vect}(M)=(C^\infty(G)\otimes \m)^H$. Clearly
$\nabla_X(f\xi) =(Xf)\xi + f\nabla_X \xi$. Any other $G$--invariant
connexion is 
given by
$$\nabla^\prime_X \xi = \nabla_X\xi +{\rm id} \otimes \alpha_X \xi,$$
where $\alpha:\m\rightarrow {\rm End}(\m)$ is an $H$--invariant linear
mapping; note that if $X\in (C^\infty(G)\otimes \m)^H$, then
$({\rm id}\otimes \alpha)(X)\in (C^\infty(G)\otimes {\rm
End}(\m))^H$. In our case $\m$ has an ${\rm Ad}\,H$--invariant inner
product $(X,Y)$. This induces a $G$--invariant hermitian structure on
tangent vectors, $(X,Y)(g) =(X(g),Y(g))\in C^\infty(G/H)$. A connexion
$\Delta^\prime_X$ is a metric connexion (or compatible with the metric)
$X(\xi,\eta) =(\nabla_X^\prime,\eta) +
(\xi,\nabla_X^\prime\eta)$. It is obvious that $\nabla_X$ is
compatible with the metric; and $\Delta_X +\alpha_X$ is compatible
with the metric iff $\alpha(\m)\subset so(\m)$. We call $\alpha$ the
connexion $1$--form. We shall take $\alpha_X(Y)=-\varepsilon [X,Y]_\m
=\varepsilon P_\m[X,Y]$ with $\varepsilon\in {\Bbb R}$. 

Let $W$ be an irreducible ${\rm Cliff}_{\Bbb C}\m$--module, ${\Bbb
Z}_2$--graded if ${\rm dim}\m$ is even. Thus $W$ is a complex inner
product space. Consider $(C^\infty(G)\otimes W\otimes V)^H$, where $V$
is any unitary $H$--module. Given $X\in (C^\infty(G)\otimes \m)^H$, a
tangent vector, and $\xi\in (C^\infty(G)\otimes W\otimes V)^H$, define
$$c(X)\xi=c(X(g))\otimes {\rm id})xi(g).$$
Then clearly $c(X)^*=c(X)$, $c(X)c(Y)+c(Y)c(X)=2(X,Y)I$ and 
$[\widetilde\nabla_X,c(Y)]=c(\nabla_XY)$, where the spin connection
$\widetilde\nabla_X$ is given by $\nabla_X + s(\alpha_X)$. (Recall
that if $T\in so(V)$, $s(T)= {1\over 4} \sum c(T\cdot v_i) c(v_i)$,
where $(v_i)$ is an orthonormal basis of the real inner product space
$V$.) 

The twisted Dirac operator $D_V$ is defined by
$D_V=\sum c(X_i) \nabla_{X_i}$, where $(X_i)$ is locally an
orthonormal basis of vector fields near $x$. $D_V$ is clearly
independent of the local choice of orthonormal basis, so globally
defined. By definition $D_V$ is $G$--invariant. We want to find a
simpler expression for $D_V$ in terms of an orthonormal basis of
$\m$. 
\vskip .1in
\noindent \bf Proposition. \it If $(X_i)$ is an orthonormal basis of
$\m$, then on $(C^\infty(G)\otimes W \otimes V)^H$, the Dirac operator
is given by
$$D_V= \sum ({\rm id} \otimes c(X_i) \otimes {\rm id})(r(X_i) \otimes
{\rm id}\otimes {\rm id} -\varepsilon {\rm id} \otimes
s(\overline{Y_i}) \otimes{\rm id}).$$
\vskip .1in
\noindent \bf Proof. \rm Both $D_V$ and the right hand side
$D^\prime_V$ are evidently $G$--invariant and act on the correct
spaces. Take a section $\xi\in (C^\infty(G)\otimes W\otimes V)^H$. By
invariance, it is enough to show that $(D_V\xi)(1)=
(D^\prime_V\xi)(1)$. As above we have  $\widetilde X_i$, vector fields
on $G/H$, orthonormal at $g=1$. We then have
$$\eqalign{(D_V\xi)(1) &= \sum (c(\widetilde X_i)\nabla_{\widetilde X_i}\xi)(1)
=\sum c(\widetilde X_i(1)) {d\over dt} \xi(\exp \widetilde
X_i(1)t)|_{t=0} \cr
&=\sum c(X_i) [(r(X_i)\xi)(1)-\varepsilon s(\overline{X_i}) \xi(1)]
=(D^\prime_V\xi)(1),\cr}$$
so the result follows.
\vskip .2in

\bf \noindent 20. REMARKS ON DIRAC INDUCTION AND BOTT'S PRINCIPLE. \rm
Consider ${\rm ker}(D_V^\pm )\subset (C^\infty(G)\otimes W^\pm\otimes
V)$. This is a closed $G$--invariant subspace (in the $C^\infty$
topology). Since $D_V$ is an elliptic operator,  we know it is
finite--dimensional; in any event it is the closure of the sum of its
irreducible subspaces. Let $G$ be a compact matrix group and let $A$
be the *--algebra 
generated by the matrix coefficients of a finite--dimensional
faithfull representation. Thus $A\subset C^\infty(G)\subset
C(G)\subset L^2(G)$. By the Stone--Weierstrass theorem, $A$ is
uniformly dense in $C(G)$ and hence dense in $L^2(G)$. Let $V\subset
L^2(G)$ be a finite--dimensional left invariant subspace. If $f\in
C^\infty(G)$ and $\xi\in V$, we have $f\star \xi\in V$. On the other
hand $f\star\xi\in C^\infty(G)$, so that $V\subset C^\infty(G)$. Now
let $W$ be any finite--dimensional irreducible representation of
$G$. Then the map of taking matrix coefficients defines a $G\times
G$--equivariant embedding $W\otimes W^*\rightarrow C^\infty(G)$.  
The algebra $A$ is a algebraic direct sum of representations $V\otimes W$ of
$G\times G$. If $W\otimes W^*$ does not appear in this list, the
corresponding elements would have to be orthogonal to $A$, a
contradiction. Hence $W\otimes W^*\subset A$. On the other hand
${\rm Hom}_G(V,C(G))=V^*$ under the map $f\mapsto f^*$ with $f^*(v)=
f(v)(1)$. So $A=\oplus V_i \otimes V_i^*$ (algebraic direct sum) as a
$G\times G$--module. Thus $L^2(G)=\oplus V_i\otimes V_i^*$ (Hilbert
space direct sum). The multiplicity space of $U$ in $C^\infty(G)\otimes V$
is ${\rm Hom}_G(U,C^\infty(G)\otimes V)= Hom(U,V)$ and the
multiplicity of $U$ in $(C^\infty(G)\otimes V)^H$ is
${\rm Hom}_G(U,(C^\infty(G)\otimes V)^H)={\rm Hom}_H(U,V)$. 
Moroever $(L^2(G)\otimes V)^H=\oplus V_i\otimes {\rm Hom}_H(V_i,V)$ as
a $G$--module.

Thus every representation $U$ of $G$ appears in $(C^\infty(G)\otimes
W\otimes V)^H$ with finite multiplicity. The multiplicity space is
given by ${\rm Hom}_H(U,W\otimes V)$. The operator $D_V$ commutes with
$G$ and therefore carries each multiplicity space onto itself. In fact
on the space ${\rm Hom}_G(U,C^\infty(G)\otimes W\otimes V)$, 
$$D_V= \sum ({\rm id} \otimes c(X_i) \otimes {\rm id})(r(X_i) \otimes
{\rm id}\otimes {\rm id} -\varepsilon {\rm id} \otimes
s(\overline{X_i}) \otimes{\rm id})
=\sum -c(X_i)( \pi(X_i)+\varepsilon 
s(\overline{X_i}) ).$$
It follows by direct computation (see below) or ellipticity of $D_V$
that ${\rm ker}(D_V)$ is finite--dimensional. Note that $D_V$ is
skew--adjoint but breaks up as two operators
$D_V^\pm:(C^\infty(G)\otimes W^\pm \otimes V)^H
\rightarrow (C^\infty(G)\otimes W^\mp \otimes V)^H$ with
finite--dimensional kernels. The index is the formal difference ${\rm
ind}(D_V^+)= [{\rm
ker}(D_V^+)]-[{\rm ker}(D_V^-)]$ in the representation ring $R(G)$. 
Let $(C^\infty(G)\otimes W^\pm \otimes V)^H=\oplus V_i\otimes M_i^\pm$,
where $M_i^\pm$ is the multiplicity space of $V_i$. Then
$D_V^\pm:M_i^\pm \rightarrow M_i^\mp$. Since ${\rm ker}(D_V^\pm)$ is
finite--dimensional, $D_V^2$ is almost everywhere an isomorphism of
$M_i^\pm$ onto $M_i^\pm$. Thus
$${\rm ind}(D_V^+)=\sum ({\rm dim}(M_i+) -{\rm dim}(M_i^-))\cdot
[V_i].$$
This proves Bott's principle: the index depends only on the underlying
bundles and not on the elliptic operator between them. one can also
compute ${\rm ind}(D_V^+)$ using supersymmetry and show that if $V$ is
positive (in a certain sense), then ${\rm ker}(D_V^-)=(0)$. This shows
that every irreducible representation of $G$ can be
realised on a space of twisted harmonic spinors ${\rm ker}(D_V^+)$ and
gives a uniform geometric construction of all irreducible
representations.

\vfill\eject

  \bf \centerline{CHAPTER~III. REPRESENTATIONS OF AFFINE KAC--MOODY
ALGEBRAS}
\vskip .1in
\noindent \bf Background. \it If $G$ is a compact simply connected
group, the corresponding loop 
group is $LGC^\infty(S^1,G)$ under pointwise multiplication. This group
acts by multiplication on $C^\infty(S^1,V)$ whenever
$V$ is a finite--dimensional representation of $V$. Let ${\rm Diff}(S^1)$ be a
the group of orientation--preserving diffeomorphisms of $S^1$ and
${\rm Rot}\, S^1$ the rotation subgroup. These groups act by
automorphisms on $LG$. They also acts on $C^\infty(S^1,V)$ compatibly
with $LG$. The semidirect product acts unitarily on $H=L^2(S^1,V)$ (after
correcting the action of ${\rm Diff}(S^1)$ by a Radon--Nikodym
cocycle). Although $H$ is already a complex Hilbert space, we regard
it as a real Hilbert space, taking a new complex structure given by
the Hilbert transform $J=i(I-2P)$, where $P$ is the Hardy space
projection onto the space with vanishing negative Fourier
coefficients. This complex structure defines an irreducible
representation of the corresponding real Clifford algebra with
generators $c(f)$ ($f\in H$). Since each element $g$ in the semidirect
product $LG\rtimes 
{\rm Diff}(S^1)$ commutes with $J$ modulo Hilbert--Schmidt operators,
each such $g$ is implemented on fermionic Fock space by a unitary
$U_g$. This gives a projective representation $LG\rtimes {\rm
Diff}(S^1)$ which lifts to an ordinary representation of ${\rm
Rot}(S^1)$ of positive energy (there are only non--negative
eigenspaces each of finite multiplicity). As in the
finite--dimensional case, there is a 
corresponding infinitesimal projective representation of the Lie
algebras of $LG$ and ${\rm Diff}(S^1)$ by quadratics in fermions. In
general, if $C^\infty(S^1,\g)={\rm Lie}(LG)$, one can look directly for
positive energy projective representations of $C^\infty(S^1,\g)\rtimes {\rm
Rot}(S^1)$. We do this below for the Lie subalgebra $L\g\subset
C^\infty(S^1,\g)$ consisting of (trigonometric) polynomial maps. In
view of the fermionic construction, it is not so surprising that every
positive energy projective representation of $L\g\rtimes {\rm
Rot}(S^1)$ extends the Witt algebra. This is the Lie algebra of polynomial
vector fields on $S^1$, a subalgebra of ${\rm Lie}\, {\rm Diff}(S^1)$.
(We will establish an infinitesimal version of this fermionic
construction in the course of this chapter.)
\vskip .1in

\noindent \bf 1. LOOP ALGEBRAS AND THE WITT ALGEBRA. \rm Let $\g$ be a
simple compact Lie algebra with complexification ${\g}_{\Bbb C}$. The
(trigonometric) polynomial loops $S^1\rightarrow \g$ are spanned by
$X\sin m\theta$, $X\cos m \theta$ ($X\in \g$, $m\ge 0$) and form a
real Lie algebra under pointwise Lie bracket,
$[X(\theta),Y(\theta)]$. The complexification has a slightly easier
spanning set $X_n =Xe^{in\theta}$. It is called the {\it loop algebra}
and denotes by $L\g$. Note that
$$[X_m,Y_n]=[X,Y]_{m+n}.\eqno{(1)}$$
Now there is a natural involution $X\mapsto X^*$ on $\g_{\Bbb C}$ as
well as the usual involution (conjugation) on complex functions. This
leads to the (conjugate--linear) involution on $L\g$ given by
$(X_n)^*=(X^*)_{-n}$. The Witt algebra corresponds to the
complexification of the real Lie algebra of (trigonmetric) polynomial
vector fields $a(\theta) \, d/d\theta$ on $S^1$. It has basis $d_n=
ie^{in\theta} \, d/d\theta$. We can use Leibniz' rule to compute the
Lie brackets:
$$[d_m,d_n] =(m-n)d_{m+n}. \eqno{(2)}$$
Note that $[d_m,X_n]=ie^{im\theta}{d\over d\theta} (e^{in\theta}) X=-n
e^{i(m+n)\theta} X = -nX_{m+n}$, so that
$$[d_m,X_n]=-n X_{n+m}.\eqno{(3)}$$
Note that $d\equiv d_0=i\, d/d\theta$ is the vector field corresponding to
rotations and $[d,X_m]=-m X_m$. Note that the rotation group ${\rm
Rot}\, S^1\cong {\Bbb T}$ acts on
$L\g$ by $(r_\alpha f)(\theta)=f(\theta-\alpha)$ and
$r_\alpha=e^{i\alpha d}$ by Taylor's theorem. 
Thus the Witt algebra ${\rm Vect}(S^1)$ acts by Lie algebra
derivations on $L\g$. We extend the involution to ${\rm Vect}(S^1)$ by
declaring that $d_n^*=d_{-n}$. This picks out the usual real structure
on ${\rm Vect}(S^1)$. 

\vskip .1in
\noindent \bf 2. POSITIVE ENERGY REPRESENTATIONS AND KAC
MOODY--ALGEBRAS. \rm We shall be interested in {\bf projective}, {\bf
unitary}, {\bf positive energy} representations of $L\g \rtimes {\Bbb
C}d$ or $L\g\rtimes {\rm Rot}\, S^1$. Thus if $\a=L\g \rtimes {\Bbb
C}d$ is the semidirect product, we look for inner product spaces $H$
(not complete!!) such that:
\vskip .1in
\noindent (1) \bf Projective: \rm $\a$ acts  projectively by operators
$\pi(A)$ ($A\in 
\a$), i.e. $A\mapsto \pi(A)$ is linear and
$[pi(A),\pi(B)]-\pi([A,B]) $ lies in ${\Bbb C}I$ for $A,B\in \a$.

\noindent (2) \bf Unitary: \rm $\pi(A)^*=\pi(A^*)$. 

\noindent (3) \bf Positive energy: \rm $H$ admits an orthogonal
decomposition $H=\bigoplus_{k\ge 0} H(k)$ such that $D=\pi(d)$ acts on
$H(k)$ as multiplication by $k$, $H(0)\ne 0$ and ${\rm dim}\, H(k)<\infty$
\vskip .1in
Note that we can take $R_\theta=\pi(r_\theta)=e^{i\theta D}$, so that
$R_\theta$ acts on $H(k)$ as multiplication by
$e^{ik\theta}$. A more general version of the positive energy
condition requires only that $H(k)=(0)$ for $k<<0$. Tensoring $H$ be a
representation of ${\rm Rot}\, S^1$, we 
may always convert this more general positive energy representation
into the normalised form given in (3). The subspace $H(k)$ are called
the energy subspaces with energy $k$; the operator $D$ has many names,
including the energy operator or hamiltonian operator.

Since the representation is projective, $[\pi(A),\pi(B)]-\pi([A,B])=b(A,B) I$
where $b(A,B)\in {\Bbb C}$. We call $b$ a 2--cocycle --- in fancy
language it gives a class in $H^2(\a,{\Bbb C})$. The definition
immediately implies the antisymmetry condition
 $$b(A,B)=-b(B,A)$$
because Lie brackets are antisymmetric; and the Jacobi identity
immediately implies that
$$b([A,B],C) + b([B,C],A) + b([C,A],B)=0$$
for all $A,B,C\in \a$. On the other hand we are free to adjust the
operators $\pi(A)$ by adding on scalars. Thus to preserve linearity,
we change $\pi(A)$ to $\pi(A)+f(A)I$ where $f:\a\rightarrow {\Bbb C}$
is linear. This changes $b(A,B)$ to $b(A,B) -f([A,B])$. We shall now
make such adjustments so that $b$ has a canonical form. We normaliseS
the inner product on $\g$ (and hence $h$ and $\h^*$) so that $\|\theta\|^2=2$,
where $\theta$ is the highest root. 
\vskip .1in
\noindent \bf Theorem. \it Representatives of $X(n)=\pi(X_n)$ and $D$
can be chosen so that
$$[X(m),Y(n)]=[X,Y](m+n) -m\ell (X,Y) \delta_{m+n,0}, \qquad
[D,X(n)]=-n X(n)$$
and $D=D^*$, $X(n)^*=-X(-n)$ for $X\in \g$. Here $\ell$ is a
non--negative integer called the {\rm level} of $H$.
\vskip .1in
\noindent \bf Remark. \rm We can extend the inner product
$(\cdot,\cdot)$ on $\g$ to a complex inner product on $\g_{\Bbb
C}$. If we set $<X,Y>=(X,Y^*)$ for $X,Y\in {\Bbb C}$, then we obtain
an invariant complex symmetric bilinear form on $\g_{\Bbb C}$. The
commutation relation then becomes 
$[X(m),Y(n)]=[X,Y](m+n) +m\ell <X,Y> \delta_{m+n,0}$ for $X,Y\in
\g_{\Bbb C}$. When $\g=\su(n)$, this bilinear form is given by
$<X,Y>={\rm Tr}(XY)$ and is in general a positive multiple of the
Killing form. 
\vskip .1in
\noindent \bf Proof. \rm We start by adjusting the operator $X=X(0)$
for $X\in \g$ (the zero modes). Set $b(X(0),Y(0))=ic(X,Y)$. Thus
$c([X,Y],Z) +c([Y,Z],X) + c([Z,X],Y)=0$ and $c(X,Y)=-c(Y,X)$. Since
$X^*=-X$ and $\pi(X)^*=-\pi(X)$, it follows that $c(X,Y)$ is real.
We show that $c(X,Y)=f([X,Y])$ for $f\in \g^*$. Let $(\cdot,\cdot)$ be
an invariant inner product on $\g$ and write $b(X,Y)=(\delta(X),Y)$
for some linear operator $\delta\in {\rm End}(\g)$. The antisymmetry
of $b$ implies that $\delta^*=-\delta$, i.e.~$\delta$ is
skew--adjoint, and the cocycle relation translates into
$$\delta([X,Y])=[\delta(X),Y] + [X,\delta(Y)].$$
Thus $\delta$ is a skew--adjoint derivation of $\g$. We saw in
Chapter~2 that any such derivation is inner, i.e.~$\delta(X)=[A,X]$
for some $A\in \g$. Hence $b(X,Y)=([A,X],Y)=-([X,A],Y) =
(A,[X,Y])=f([X,Y])$ with $f(X)=(A,X)$. Thus we may adjust the $X$'s by
purely imaginary scalars so that $[\pi(X),\pi(Y)]=\pi([X,Y])$ with
$\pi(X)^*=-\pi(X)$. 

Now consider the operators $[D,\pi(X)]=ig(X)I$. Since $D=D^*$ and
$\pi(X)^*=-\pi(X)$, $g(X)$ is real. It is also linear and the cocycle
relation implies that $g([X,Y])=0$. Since $[\g,\g]=\g$, it follows
that $g\equiv 0$ and hence that $[D,\pi(X)]=0$.

This completes the zero mode adjustments. If $n>0$, we can choose
$X(n)$ so that $[D,X(n)]=-n X(n)$, for the left hand side is
independent of any choices. Taking $X(-n)=-X(n)^*$ for $X\in \g$, we
still have $[D,X(-n)]=nX(-n)$. These equations imply that the operator
$X(n)$ takes the energy space $H(k)$ into $H(k-n)$. Thus the
operators $X(-n)$  and $X(n)$ ($n>0$) raise  and lower energy.
  We thus have
$[\pi(X),\pi(Y)]=\pi([X,Y]) + ib(x,y)I$ for some 
2--cocycle $b(x,y)\in {\Bbb C}$. We now compute $b(x,y)$.
\vskip .05in
\noindent (1) We have already imposed the condition $[D,X(n)]=-nX(n)$
for $X\in \g$ which  uniquely specifies the choice of $X(n)$'s.
\vskip .05in
\noindent (2) $[X(n),Y(m)]=[X,Y](n+m)$ if $n+m\ne 0$. For these
elements lower energy by $n+m$ and scalars preserve energy.
\vskip .05in
\noindent (3) $[X(n),Y(-n)] = [X,Y](0) + \delta_n (X,Y)\cdot
I$. For if 
$[X(n),Y(-n)]-[X,Y](0)=\lambda(X,Y)$, then taking Lie brackets with
$Z(0)$, we find $\lambda([Z,X],Y) + \lambda(X,[Z,Y])=0$.  But any
$\g$--invariant bilinear form on $\g$ is a multiple of $(X,Y)$.
\vskip .05in
\noindent (4) $\delta_{n}=n\delta_1$ for $n>0$. It suffices to show
that $\delta_{n+1} =\delta_n + \delta_1$. But
$$[A(-1),[X(n+1),Y(-n)]]= [A(-1),[X,Y](1)]=([A,[X,Y]])(0)
-\delta_1(A,[X,Y]).$$ 
On the other hand
$$\eqalign{ [A(-1),[X(n+1),Y(-n)]]& =
-[X(n+1),[Y(-n),A(-1)]]-[Y(-n),[A(-1),X(n+1)]] \cr
&=-[X(n+1),[Y,A](-n-1)] -[Y(-n),[A,X](n)]\cr
&=-[X,[Y,A]](0) -[Y,[A,X]](0) -
\delta_{n+1} (X,[Y,A]) - \delta_n (Y,[A,X]).\cr}$$
Since $(\cdot,\cdot)$ is $\g$--invariant, and $[\g,\g]=\g$, the result
follows. 
\vskip .05in
\noindent (5) $-\delta_1$ is a non--negative integer $\ell$. Suppose
that $H(0)$ has a summand $V_\lambda$  with highest weight $\lambda$.  
Let $v_\lambda$ be a highest weight vector in $V_\lambda$, so that
$H_i(0)v_f =\lambda(H_i)v_\lambda$ and $E_i(0)v_\lambda=0$
for $i>0$. Now consider $E=E_{-\theta}(1)$,
$F=E_\theta(-1)$, $H=[E,F]= iT_\theta(0)+ \delta_1$, where
$(T_\theta,T)=\theta(T)$. Then
$[H,E]=2E$, $[H,F]=-2F$, $H=H^*$ and $E^*=F$. Moreover $Ev_\lambda=0$ and
$Hv_\lambda =(-(\theta,\lambda)+\delta_1)v_\lambda$. So by the usual
$\sl_2$ lemma, $-\delta_1-(\theta,\lambda)$ is a non--negative
integer. Hence $-\delta_1$ must be a non--negative integer $\ell$.
\vskip .1in
\noindent \bf Corollary of proof. \it Each energy space $H(k)$ is a
$\g$--module. For level $\ell\ge 0$, any
highest weight $\lambda$ appearing in $H(0)$ must satisfy
$(\lambda,\theta)\le \ell$.
\vskip .1in
In the light of this theorem, we define the {\it affine Kac--Moody
algebra} $\widehat{g}$ by
$$\widehat{g} =L\g \oplus {\Bbb C}d\oplus {\Bbb C}c,$$
with $[X(m),Y(n)]=[X,Y](m+n) +m <X,Y> \delta_{m+n,0} c$ and
$[d,X(n)]=-nX(n)$ (note the minus sign!). All other brackets
zero. This contains $L\g \oplus {\Bbb C}c={\cal L}\g$ as an
ideal. Note that ${\Bbb C}c$ is a central subalgebra and ${\cal
L}\g/{\Bbb C}c=L\g$. Thus ${\cal L}\g$ is a central extension of $L\g$
by ${\Bbb C}$ and ${\rm Rot}\,S^1$ acts on ${\cal L}\g$. 
${\cal L}\g$ is called an {\it affine Lie algebra}. 
\vskip .1in
\noindent \bf Invariant symmetric bilinear form. \rm We define a
bilinear form on $\widehat{\g}$ by 
$$(X_1(m)+\delta_1 d +\gamma_1 c, X_2(n)+\delta_2 d + \gamma_2
c)= (X_1,X_2) \delta_{m+n,0} +\delta_1\gamma_2 +\delta_2\gamma_1
=-<X_1,X_2> \delta_{m+n,0} +\delta_1\gamma_2 +\delta_2\gamma_1.$$
It is straightforward to check that this form is 
${\rm ad}$--invariant. It is important that it is not positive
definite on the ''real part''; for example it is clearly
Lorentzian on $i{\Bbb R}d +i{\Bbb R}c$.

\vskip .1in
\noindent \bf 3. COMPLETE REDUCIBILITY. \rm
\vskip .1in
\bf \noindent Positive energy theorem. \it (a) Let $H$ be a positive
energy representation of $\widehat{g}$ and let $V=H(0)$ be the lowest
energy subspace. If $\xi\in H(0)$ is cyclic, then some vector in
$H(0)$ generates a lowest energy subspace.

\noindent (b) Any positive energy representation is an orthogonal
direct sum of irreducible positive energy representations.
\vskip .1in
\noindent \bf Proof. \rm (a) Let $V$ be the subspace of lowest
energy. Let $K$ be any invariant subspace of $H$. So
$K=\bigoplus K(n)$ and $H=K\oplus K^\perp$. Note also that $K^\perp$
is invariant. Now $H(0)=K(0)\oplus K^\perp(0)$. Moreover the
$\widehat{g}$--module generated by $K(0)$ or $K^\perp(0)$ is contained
in $K$ or $K^\perp$. But the $\widehat{g}$--module generated by
$K(0)\oplus K^\perp(0)=H(0)$ equals $H$. Hence $K$ must be the
$\widehat{g}$--module generated by $K(0)=K\cap H(0)$. Thus there is a
1--1 correspondence between submodules and certain invariant subspaces of
$H(0)$. Taking $K$ so that $K\cal H(0)$ has minimal dimension, we see
that $K$ must be irreducible. Any non--zero vector in $K\cap H(0)$
must be cyclic, by irreducibility.

\noindent (b) Take the cyclic module generated by
a vector of lowest energy. This contains an irreducible submodule
generated by another vector of lowest energy $H_1$ say. Now repeat
this process for $H_1^\perp$, to get $H_2$, $H_3$, etc. The positive
energy assumption shows that $H=\bigoplus H_i$.
\vskip .1in
\noindent \bf 4. CLASSIFICATION OF POSITIVE ENERGY
REPRESENTATIONS. \rm 
\vskip .1in
\noindent \bf Theorem (Uniqueness). \it Let $(\pi,H)$ be an
irreducible positive energy representation of $\widehat{g}$ of level
$\ell$.
\vskip .05in
\item{(1)} $H(0)$ is irreducible as an $SU(N)$--module.

\item{(2)} If $H(0)=V_\lambda$, then $(\theta,\lambda)\le \ell$.

\item{(3)} (Uniqueness) If $H$ and $H^\prime$ are irreducible positive
energy representations of level $\ell$ of the above form with
$H(0)\cong H^\prime(0)$ as $\g$--modules, then $H$ and $H^\prime$ are
unitarily equivalent as representations of $\widehat{g}$.
\vskip .1in

\noindent \bf Proof. \rm (1) Let $V$ be an irreducible
$SU(N)$--submodule of $H(0)$. By irreducibility the
$\widehat{g}$--module generated by $V$ is the 
whole of $H^0$. Since $D$ fixes $V$, it follows that the ${\cal L}\g$--module
generated by $V$ is the whole of $H$. The commutation rules 
show that any monomial in the $X(n)$'s can be written as a sum of
monomials of the form $P_- P_0 P_+$, where $P_-$ is a monomial in the
$X(n)$'s for $n<0$ (energy raising operators), $P_0$ is a monomial in
the $X(0)$'s (constant energy operators) and $P_+$ is 
a monomial in the $X(n)$'s with $n>0$ (energy lowering operators).
Hence $H$ is spanned by 
products $P_- v$ ($v\in V$). Since the lowest energy subspace of this
$L^0\g$--module is $V$, 
$H(0)=V$, so that $H(0)$ is irreducible as a $\g$--module. 

\noindent (2) We have already proved this in section 2 by introducing
$E=E_{-\theta}(1)$, $F=E_\theta(-1)$ and $H=[E,F]=\ell +iT_\theta(0)$,
where  $(T_\theta,T)=\theta(T)$. [Thus, if $v\in H(0)$ has highest weight
$\lambda$, then $T(0)v=i\lambda(T)v$ and $Ev=0$. Thus $Hv=(\ell +i
T_\theta(0)) v=(\ell-(\theta,\lambda))v$. Since $E,F,H$ give a copy of
$\sl_2$, we get $\ell\ge (\theta,\lambda)$.]

\noindent (3) Any monomial $A$ in operators from $\g$ is a sum of
monomials $RDL$ with $R$ a monomial in energy raising operators, $D$ a
monomial in constant energy operators and $L$ a monomial in energy
lowering operators. Observe that if $v,w\in H(0)$, the inner products 
$(A_1v,A_2w)$ are uniquely determined by $v,w$ and the monomials $A_i$:
for $A_2^*A_1$ is a sum of terms $RDL$ and $(RDLv,w)=(DLv,R^*w)$ with
$R^*$ an energy lowering operator. Hence, if $H^\prime$ is another
irreducible positive energy representation with $H^\prime(0)\cong
H(0)$ by a unitary isomorphism 
$v\mapsto v^\prime$, $U(Av)=Av^\prime$ defines a unitary map of $H$
onto $H^\prime$ intertwining $\widehat{\g}$. 
\vskip .1in

\noindent \bf 5. SUGAWARA'S FORMULA FOR $L_0$. \rm 
\vskip .1in

\bf \noindent Sugawara's formula for $L_0$. \it Let $H$ be a cyclic positive
energy representation at level $\ell$ and let $(X_i)$ be
an orthonormal basis of $\g$. Let $L_0$ be the
operator defined on $H^0$ by
$$L_0={1\over N+\ell} \left(-\sum_{i} {1\over 2} X_i(0)X_i(0)
-\sum_{n>0}\sum_i  X_i(-n)X_i(n)\right).$$
Then $L_0 = D + C/2(N+\ell)$ where
the Casimir $\Delta=-\sum_i X_i(0)X_i(0)$ acts on $H(0)$ as
multiplication by $C$ and on $\g$ as $2g$ (where $g$ is the dual
Coxeter number, equal to $N$ for $\su(N)$). 
\vskip .1in

\noindent \bf Proof. \rm Since $\sum_i X_i(a) X_i(b)$ is
independent of the orthonormal basis $(X_i)$, it commutes with $G$ and
hence each $X(0)$ for $X\in \g$. Thus $\sum_i[X,X_i](a) X_i(b) + X_i(a)
[X,X_i](b) =0$ for all $a,b$. If
$A=\sum_i{1\over 2}  X_i(0)X_i(0) +\sum_{n>0} X_i(-n) X_i(n)$, then
using the above relation we get
$$\eqalign{[X(1),A]
&=-\ell X(1) +\sum_i {1\over 2}([X,X_i](1) X_i(0) + X_i(0)
[X,X_i](1)) \cr
& \quad +\sum_n[X,X_i](-n+1) X_i(n) + X_i(-n) [X,X_i](n+1)\cr
& =- \ell X(1) +{1\over 2} \sum _i [[X,X_i](1),X_i(0)]= \ell X(1)
+{1\over 2} \sum_i [[X,X_i],X_i](1),\cr}$$ 
since $([X,X_i],X_i)=0$ by invariance of $(\cdot,\cdot)$. 
Hence $[X(1),A]=-(g+\ell) X(1)$, since $-\sum_i {\rm ad}(X_i)^2 =2g$.
Similarly $[X(-1),A]=(g+\ell) X(-1)$. [Note that if $H$ were a unitary
representation, so that $X(n)^*=-X(-n)$, then $A^*=A$ and taking adjoints
we get $[X(-1),A]=(g+\ell)X(-1)$. A similar argument could be applied
in general using the pairing between $H$ and its algebraic dual.] Thus
$(g+\ell)D+A$ commutes with all $X(\pm 1)$'s. Since
$[\g,\g]=\g$, these generate $L^0\g$, and hence $(g+\ell)D +A$
commutes with $L\g$. Since $(g+\ell)D+A=-{\Delta\over 2}\cdot I$ on
$H(0)$ and the cyclic subspace generated by $H(0)$ is the whole of
$H$, we get $(g+\ell)D+A=-{\Delta\over 2}\cdot I$ on $H$ as required.
\vskip .1in
\bf \noindent Corollary. \it Let $H$ be a positive energy
representation of $\widehat{\g}$. 
\vskip .05in
\item{(a)} If $H$ is irreducible as an $\widehat{\g}$---module, then 
it is irreducible as an  ${\cal L}\g$--module.
\item{(b)} If $H_1$ and $H_2$ are irreducible $\widehat{\g}$--modules 
which are isomorphic as ${\cal L}\g$--modules, then one is obtained from the
other by tensoring with a character of ${\Bbb T}$.
\vskip .1in
\noindent \bf Proof. \rm (a) The Sugawara formula show that given
$k\ge 0$, there is a finite linear combination $T$ of operators
$X(a)Y(b)$ such that $T\xi =D\xi$ for all $\xi \in H(0)\oplus \cdots \oplus 
H(k)$. Hence the submodule generated by any such $\xi$ also contains
the submodules generated by any of the components $\xi_j\in H(j)$ ($j\le
k$). However it is clear that the $\widehat{g}$--module generated by
any $\xi_j$ is the same as the ${\cal L}\g$ module generated by
$\xi_j$. By irreducibility, it follows that the ${\cal L}\g$--module
generated by $\xi$ is the whole of $H$. 

\noindent (b) Let $T:H_1\rightarrow H_2$ be a unitary
intertwiner for $\widehat{g}$. Then $V_t^*TU_t$ is also a unitary intertwiner, 
so must be of the form $\lambda(t) T$ for $\lambda(t)\in {\Bbb T}$ by
Schur's lemma. Since $TU_tT^*=\lambda(t) V_t$, $\lambda(t)$ must be a
character of ${\Bbb T}$. 

\vskip .1in
\noindent \bf Remark. \rm The previous corollary is important because
it shows that positive energy representations are classified by up to
tensoring with a character of ${\rm Rot}\, S^1$. This will appear as
an important feature in our discussion of roots and weights for
$\widehat{g}$ below. 

\vskip .1in
\noindent \bf 6. SUGAWARA'S CONSTRUCTION OF THE VIRASORO ALGEBRA. \rm 
\vskip .1in
\noindent \bf Theorem. \it Let $H$ be an irreducible  positive
energy prepresentation at level $\ell$. If $L-0$ is defined
as above and in addition we set
$$L_m=-{1\over 2(\ell +g)} \sum_i \sum_{a+b=m} X_i(a) X_i(b),$$
then
$$[L_m,X(n)]=-n X(n+m), \qquad [L_m,L_n]=(m-n)L_{m+n} + {c\over 12}
(m^3-m) \delta_{m+n,0},$$
where $c={\rm dim}\, \g \cdot \ell/(\ell+g)$ and $g$ is the dual
Coxeter number.
\vskip .1in
\noindent \bf *Proof. \rm Since $\sum_i X_i(a) X_i(b)$ is
independent of the orthonormal basis $(X_i)$, it commutes with $G$ and
hence each $X(0)$ for $X\in \g$. Thus $\sum_i[X,X_i](a) X_i(b) + X_i(a)
[X,X_i](b) =0$ for all $a,b$. If
$B={1\over 2} \sum_{a+b=m} X_i(a) x_i(b)$, 
then using the above relation we get $[X(0),B]=0$.
Similarly writing 
$$B=\sum_i ({1\over 2}\sum _{a+b=m,a=b} +\sum_{a+b=m,b>a})
X_i(a)X_i(b),$$
we get
$$[X(1),B]=\ell X(1)+ \sum_i ({1\over 2}\sum _{a+b=m,a=b} +\sum_{a+b=m,b>a})
[X,X_i](a+1)X_i(b) +X_i(a)[X,X_i](b+1).$$
In this sum, consider terms $P(\alpha)Q(\beta)$ with
$\alpha+\beta=m+1$: if $\alpha<\beta+1$, the first term gives
a contribution $\sum_i [X,X_i](\alpha)X_i(\beta)$, while if $\alpha=\beta+1$,
it gives 
${1\over 2}\sum_i [X,X_i](\alpha)X_i(\beta)$; if $\beta=\alpha+1$, the second
terms gives a contribution
${1\over 2} \sum _i X_i(\alpha)[X,X_i](\beta)$, while if
$\beta>\alpha+1$, it gives $\sum_i X_i(\alpha)[X,X_i](\beta)$. Adding these
contributions (when they occur), we get a total of $0$ if
$\beta>\alpha+1$; $\sum_i {1\over 2}  X_i(\alpha) [X,X_i](\beta)
+[X,X_i](\alpha)X_i(\beta)$ if $\beta=\alpha+1$;
$\sum_i [X,X_i](\alpha)X_i(\beta)$ if $\beta=\alpha$; ${1\over 2}
\sum_i [X,X_i](\alpha)X_i(\beta)$ if $\beta=\alpha-1$; and $0$ if
$\beta<\alpha-1$. If $m+1=2k$, then we must have $\alpha=\beta=k$ and
$$\eqalign{[X(1),B]+\ell X(1)&=\sum_i [X,X_i](k)X_i(k)\cr
&={1\over 2} \sum_i ([X,X_i](k)X_i(k)
-X_i(k) [X,X_i](k))\cr
&={1\over 2} [[X,X_i],X_i](m+1) -\ell \delta_{k,0}
([X,X_i],X_i) \cr
&=gX(m+1).\cr}$$
If $m=2k$, then
$$\eqalign{[X(1),B]+\ell X(1)& =\sum_i {1\over 2} X_i(k)[X,X_i](k+1) +
[X,X_i](k) X_i(k+1) +{1\over 2} [X,X_i](k+1) X_i(k)\cr
&={1\over 2}\sum_i [[X,X_i](k+1),X_i(k)]\cr
&=-gX(m+1).\cr}$$
Thus in both cases $[X(1),B]=-(\ell +g)X(m+1)$. Hence
$[L_m,X(1)]=-X(m+1)$. Since $L_m^*=L_{-m}$ and $X(1)^*=-X(-1)$,
taking adjoints we get $[L_m,X(-1)]=X(m-1)$. Since
$[\g,\g]=\g$, $L^0\g$ is generated by the $X(1)$'s and $Y(-1)$'s. The
relation $[L_m,X(n)]=-nX(m+n)$ then follows easily by induction and
the Jacobi identity. Hence
$$[[L_m,L_n],X(p)]=[L_m,[L_n,X(p)]]-[L_n,[L_m,X(p)]]
=-p(m-n)X(p+m+n)=(m-n)[L_{m+n},X(p)].$$
Since the $X(p)$'s act irreducibly, Schur's lemma implies that
$[L_m,L_n]=(m-n)L_{m+n}+\lambda(m,n)I$ for some scalar
$\lambda(m,n)$. Note that by definition $L_m$ carries $H(k)$ into
$H(k-m)$. Thus both $[L_m,L_n]$ and $L_{m+n}$ carry $H(k)$ into
$H(k-m-n)$. Thus $\lambda(m,n)=0$ if $m+n\ne 0$. Clearly
$\lambda(m,-m)=-\lambda(-m,m)$ by antisymmetry of the Lie bracket. 
Take $\xi\in H(0)$ and $m>0$. Then
$([L_m,L_{-m}]\xi,\xi)=2m (L_0\xi,\xi)
+\lambda(m,-m)(\xi,\xi)$. 
On the other hand $L_0\xi=\mu\xi$, where $\mu=\Delta/2(\ell+g)$, and
$([L_m,L_{-m}]\xi)=(L_mL_{-m}\xi,\xi)$. For $m>0$, we have
$$L_{-m}\xi={1\over 2(\ell +g)} \sum_i \sum_{a+b=m; a,b\ge 0}
X_i(-a)X_i(-b)\xi.$$
Since $[L_m,X(p)]=-pX(p+m)$ and $L_m\xi=0$, we get
$$L_m L_{-m}\xi ={1\over 2(\ell +g)} \sum_i\sum_{a+b=m;a,b\ge 0} aX_i(-a+m)
X_i(-b)\xi +b X_i(-a)X_i(-b+m)\xi.$$
Now $S_j=\sum-(X_i(j)X_i(-j)\xi,\xi)=0$ if $j<0$ and $=\Delta\|\xi\|^2$ if
$j=0$. Since $X_i(j)X_i(-j)=X_i(-j)X_i(j) +j\ell\cdot I$, 
$\sum-(X_i(j)X_i(-j)\xi,\xi)= j\ell \cdot{\rm dim}\, \g \, \|\xi\|^2$
for $j>0$. Hence
$$\eqalign{(L_m L_{-m}\xi,\xi) &=
{1\over 2(\ell+g)} [ \sum_{b=0}^m (m-b) S_b +\sum_{a=0}^m (m-a) S_a]\cr
&={1\over (\ell+g)} [m\Delta\|\xi\|^2 + \sum_{b=0}^m b(m-b) \ell \cdot
{\rm dim}\, \g \, \|\xi\|^2]\cr
&={\|\xi\|^2\over \ell +g} (m\Delta + (m^3-m)/6 \cdot \ell {\rm dim}\,
\g),\cr}$$
since $\sum_{a=0}^m a(m-a)=(m^3-m)/6$. 
Thus $\lambda(m,n)=\delta_{m+n,0} c(m^3-m)/12$ with $c=\ell {\rm
dim}\, \g/(\ell+g)$, as required.

\vskip .1in
\noindent \rm This result is an example of quantisation. The Witt
algebra acts by derivations on $L\g$ preserving the
central extension. Thus by Schur's lemma there is at most one
covariant projective representation of it compatible with the
group. Thus if $\pi(d_n)=L_n$, we require
$[L_n,X(m)]=-mX(m+n)$. By uniqueness, we must have
$[L_m,L_n]=(m-n)L_{n+m} =\lambda(m,n)I$. Here $\lambda(n,m)$ is a
2--cocycle. As we show below, by appropriate adjustment of the $L_n$'s
by scalars, $\lambda$ can always be normalised so that
$\lambda(m,n)={c\over 12} \delta_{m+n,0}$, where $c$ is called the
{\it central charge}. Note that the Sugawara construction was automtically
normalised with central charge ${\rm dim}\, \g \cdot \ell/(\ell+N)$.
The representation was also unitary, in that $L_n^*=L_{-n}$ for all
$n$. This central extension of the Witt algebra is usually called the
Virasoro algebra and is classified by its central charge. 
\vskip .1in
\noindent \bf Virasoro cocycle lemma. \it If $[L_0,L_n]=-nL_n$ for all
$n$, then 
$[L_m,L_n]=(m-n)L_{m+n} +(am^3 +bm)\delta_{m+n,0}I$. If we choose
$L_0$ so that $[L_1,L_{-1}]=L_0$, then $a+b=0$. 
\vskip .1in
\noindent \bf Proof. \rm Note that $-n^{-1} [L_0,L_n]$ is independent
of adding scalars onto $L_0$ or $L_n$, so we may always choose $L_n$
so that $[L_0,L_n]=-nL_n$. By the Jacobi idenity for $L_0$, $L_m$ and
$L_n$, we have
$[L_0,[L_m,L_n]]=-(n+m)[L_m,L_n]$. On the other hand
$[L_0,L_{m+n}]=-(m+n) L_{m+n}$. Since $[L_m,L_n]=(m-n)L_{m+n}
+\lambda(m,n)I$, we must have $\lambda(m,n)=0$ if $m+n\ne 0$. Thus
$$[L_m,L_n]=(m-n) +A(m) \delta_{m+n,0}.$$
Clearly $A(m)=-A(-m)$ nad $A(0)=0$. Writing out the Jacobi
identity for $L_k$, $L_n$ and $L_m$ with $k+n+m=0$, we get
$$(n-m)A(k) + (m-k)A(n) +(k-n)A(m)=0.$$
Setting $k=1$ and $m=-n-1$, we get
$$(n-1)A(n+1) =(n+2)A(n)-(2n+1) A(1).$$
This recurrence relation allows $A(n)$ to be determined from $A(1)$
and $A(2)$. Since $A(n)=n$ and $A(n)=n^3$ give solutions, we see that
$A(m)=am^3+bm$ for some constants $a$ and $b$. Clearly we are free to
choose $L_0=[L_1,L_{-1}]$ (since we have made no adjustment to $L_0$
so far). But then $A(1)=0$ and hence $a+b=0$. 
\vskip .1in

\noindent \bf 7. WEIGHTS, ROOTS AND THE QUANTUM CASIMIR OPERATOR. \rm 
\vskip .1in

\noindent \bf Weights. \rm It is immediately verified that
$\overline{\h}=\h\oplus i{\Bbb R}d \oplus i{\Bbb R}c$ is a maximal
Abelian subalgebra of $\widehat{g}$. If $H$ is a positive energy
representation (in the generalised sense), it first has an energy
decomposition $\bigoplus H(k)$ where $H(k)$ is the $k$--eigenspace of
$d$. If $H$ has level $\ell$, then $c=\ell I$ on $H$. Each $H(k)$
breaks up as a sum of $\h$--modules with weights $\mu\in P(\g)$. Thus
the weights of $H$ are triples $\overline{\mu}=(\mu,k,\ell)\in
\overline{h}^*$. The dimension of the corresponding weight space in
$H$ is called the multiplicity of the weight.

\vskip .1in
\noindent \bf The Lorentzian inner product on weights. \rm We
introduce a real symmetric bilinear form on 
$\overline{\h}^*$ via $(\mu_1,k_1,\ell_1)\cdot
(\mu_2,k_2,\ell_2)=(\mu_1,\mu_2) +k_1\ell_2 +k_2\ell_1$. This bilinear
form is thus obtained by taking the direct sum of the eculidean space
$\h^*$ with the Lorentzian lattice ${\Bbb R}^{1,1}$, with indefinite symmetric
form $(x_1,y_1)\cdot (x_2,y_2)=(x_1y_2+y_1x_2)$. 

\vskip .1in

\noindent \bf Roots and multiplicities. \rm The Lie algebra
$\overline{\h}$ acts on $\widehat{g}$ through the adjoint
representation preserving the Lorentzian form introduced before
(recall $(X_1(m)+\gamma_1 c+\delta_1 d, X_2(n) +\gamma_2 c+\delta_2 d)
=\delta_{m+n,0} (X_1,X_2) +\gamma_1\delta_2 +\delta_2\gamma_1$). The
inner product is non--degenerate on $\overline{\h}$ and the orthogonal
complement splits as a direct sum of non--zero eigenspaces of
$\overline{\h}$, each of finite multiplicity. Indeed
$$\overline{\h}^\perp = \bigoplus_{n\ne 0} \h(n) \oplus
\bigoplus_{\alpha\in \Phi^+} \g_\alpha(0) \oplus \bigoplus{\alpha\in
\Phi,n\ne 0} \g_\alpha(n).$$
These give weights of $\overline{\h}$ which we call affine roots. 
Since $c$ is central, they all have the form $(*,*,0)$. We can list
all the roots: $(\alpha,0,0)$ with $\alpha\in \Phi$; $(\alpha,n,0)$
with $\alpha\in \Phi$ and $n\ne 0$; and $(0,n,0)$. We denote the set
of affine roots by $\overline{\Phi}$. We define the positive roots by 
$\overline{\Phi}^+$ to be $(0,n,0)$ or $(\alpha,n,0)$ with $n<0$ or,
if $n\ne 0$, $(\alpha,0,0)$ with $\alpha\in \Phi^+$. If
$\overline{\alpha}$ is an affine root we denote by
$m_{\overline{\alpha}}$ the multiplicity of the corresponding root
space. Thus $(0,n,0)$ has multiplicity $m={\rm dim}\, \h$ while all
other roots have multiplicity one. (These conventions are adopted so
that we can use highest weight theory painlessly.) The roots are of
two types those of form $\overline{\alpha}=(\alpha,n,0)$ with
$\alpha\in \Phi$ and those of the form
$\overline{\alpha}=(0,n,0)$. The former satisfy
$(\overline{\alpha},\overline{\alpha})>0$ and are called {\it
space--like}; the latter satsfy
$(\overline{\alpha},\overline{\alpha})=0$ and are called {\it
time--like}. 

\vskip .1in
\noindent \bf The simple roots. \rm Let $\alpha_1,\dots,\alpha_m$ be the
simple roots of $\g$ (with respect to a standard Weyl chamber). Let
$\theta$ be the highest root. Set $\overline{\alpha}_i=(\alpha_i,0,0)$
for $i=1,\dots,m$ 
and $\overline{\alpha}_0=(-\theta,-1,0)$. The $\overline{\alpha}_i$
are all roots called the simple roots. 
\vskip .1in
\noindent \bf Lemma. \it A space--like root is positive iff it is a
non--negative integer combination of simple roots. 
\vskip .1in
\noindent \bf Proof. \rm  Let $\overline{\alpha}=(\alpha,n,0)$ be a
space--like root. If $n=0$, the result is known from the
finite--dimensional case. If $n\ne 0$ and $\overline{\alpha}
=\sum_{i=0}^m n_i\overline{\alpha_i}$, then $n_0=-n$ and
$\alpha=\sum_{i=1}^m n_i\alpha_i +n\theta$. If $n\ge 1$, then
$n\theta -\alpha=(n-1)\theta +(\theta-\alpha) =-\sum n_i\alpha_i$ with
$n_i\le 0$ for $i\ge 1$. In this case $n_0<0$. If $n\le -1$,
then $n\theta -\alpha =-\sum n_i\alpha_i$ with $n_i\ge 0$ for $i\ge
1$. In this case $n_0>0$. This prove the lemma. 

\vskip .1in

\noindent \bf Kac--Moody--Serre generators and relations. \rm Let
$\alpha_1,\dots,\alpha_m$ be the simple positive roots of $\g$ and set
$E_i=E_{\alpha_i}(0)$, $F_i=F_{\alpha_i}(0)$ and
$H_i=H_{\alpha_i}(0)$.  In addition let $E_0=E_{-\theta}(1)$,
$F_0=E_\theta(-1)$ and $H_0=[E_0,F_0]=\ell +iT_\theta(0)$, 
where  $(T_\theta,T)=\theta(T)$. Thus $E_0^*=F_0$, $H_0^*=H_0$,
$[H_0,E_0]=2E_0$, $[H_0,F_0]=-2F_0$ and $[E_0,F_0]=H_0$. (Similar
relations of course hold for each other $E_i,F_i,H_i$.) Note that
$E_0$ corresponds to the positive root $(-\theta,-1,0)$ and $F_0$ to
the negative root $(\theta,1,0)$. Since every affine root or its
negative is a non--negative combination of simple roots, it follows
that $[E_0,F_i]=0$ for $i\ne 0$ and $[F_0,E_i]=0$ for $i\ne 0$. Thus
$[E_i,F_j]=\delta_{ij} H_i$ as in the finite--dimensional case. 
\vskip .1in
\noindent \bf Lemma. \it ${\cal L}\g$ is generated as a complex Lie algebra by
the $E_i$'s and $F_i$'s where $i=0,\dots,m$.
\vskip .1in
\noindent \bf Proof. \rm Let ${\cal L}\g_c$ be the complex Lie algebra
generated by all $E_i$, $F_i$ and hence $H_i=[E_i,F_i]$. Clearly
${\cal L}\g_c$ is *--invariant and hence the complexification of its skew
adjoint part ${\cal L}\g_0$. We already know that $\g_{\Bbb C}$ is
generated by the $E_i$'s and $F_i$'s with $i\ge 1$. We prove by
induction that $\g_{\Bbb C}(n)\subset {\cal L}\g_c$. The Lie algebra
$\g=\g(0)$ acts by the adjoint representation on $\g_{\Bbb C}(\pm
1)$. Since this action is isomorphic to the adjoint representation of
$\g$, it is irreducible. Hence the $\g(0)$--module generated by
$E_0=E_{-\theta}(1)$ is the whole of $\g_{\Bbb C}(1)$ and hence
$\g_{\Bbb C}(1)\subset {\cal L}\g_c$. Now suppose that $\g_{\Bbb
C}(n)\subset {\cal L}\g_c$ for $n>0$. Then
$[X(1),Y(n)]=[X,Y](n+1)$. Since $[\g,\g]=\g$, we see that $\g_{\Bbb
C}(n+1)\subset {\cal L}\g_c$. Taking adjoints, it follows that ${\cal
L}\g_c={\cal L}\g$, as required.

\vskip .1in
\noindent \bf Corollary. \it The $E_i$'s generate
${\cal
L}^+\g=\bigoplus_{\overline{\alpha}>0}\widehat{\g}_{\overline\alpha}$
and the $F_i$'s generate ${\cal L}^-\g=\bigoplus_{\alpha<0}
\widehat{\g}_\alpha$.  
\vskip .1in
\noindent \bf Proof. \rm Let ${\cal L}\g_+$ and ${\cal L}\g_-$ be the
Lie algebras 
generated by the $E_i$'s and $F_i$'s respectively. The relation
$[E_i,F_j]=\delta_{ij} H_i$
shows that $\g_+\oplus \h \oplus \g_-$ is a Lie subalgebra of
$\g_{\Bbb C}$. Since it contains $E_i,F_i,H_i$, it must be the whole
of ${\cal L}\g$ so the result follows.

\vskip .1in
\bf \noindent Lemma (Serre relations). \it  The generators
$E_i,F_i,H_i$ satisfy the following relations:

\item{S1.} $[H_i,H_j]=0$.

\item{S2.} $[E_i,F_j]=\delta_{ij}H_i$.

\item{S3.} $[H_i,E_j]=n(i,j) E_j$ and $[H_i,F_j]=-n(i,j) F_j$ where
$n(i,j)= 2(\overline{\alpha}_i,\overline{\alpha}_j)/
(\overline{\alpha}_i,\overline{\alpha}_i)$.

\item{S$^+_{ij}$.} ${\rm ad}E_i^{-n(i,j)+1}E_j=0$ for $i\ne j$.

\item{S$^-_{ij}$.} ${\rm ad}F_i^{-n(i,j)+1}F_j=0$ for $i\ne j$.
\vskip .1in
\noindent \bf Proof. \rm We already know S1 and S2. S3 follows from
the definition of the Lorentzian inner product and the fact that $E_0$
has weight $(-\theta,-1,0)$ and $F_0$ has weight $(\theta,1,0)$. 
To prove the S$^-_{ij}$ ($i\ne j$), note that ${\rm ad}(E_i)\cdot
F_j=0$, ${\rm ad}(H_i)\cdot F_j=-n(i,j)F_j$. Thus the result follows
from $SU(2)_i$--theory, because $F_j$ is a highest weight vector. (In
particular $n(i,j)\le 0$.) S$^+_{ij}$ follows by taking adjoints.
\vskip .1in
\noindent \bf Definition. \rm $n_{ij}=n(i,j)$ ($i,j\ge 0$) is called the
extended Cartan matrix of $\g$. It is the matrix obtained by taking
$\alpha_0=-\theta$ together with the simple roots of $\g$ and
therefore coincides with our previous definition of the extended
Cartan matrix.

\vskip .1in

If $H$ is irreducible, we 
know that $H(0)$ is an irreducible $\g$--module. Let $\lambda$ be its
highest weight. The corresponding highest weight vector $v_\lambda$
has weight $\overline{\lambda}=(\lambda,0,\ell)$. The space $H$ is
spanned by all vectors 
obtained by applying lowering operators $F_i$ to $v_\lambda$. But if
$w\in H$ has weight $\overline{\mu}$, $F_iw$ has weight
$\overline{\mu} -\overline{\alpha}_i$. Thus all the weights of $H$
have the form $\overline{\mu}=\overline{\lambda}-\sum_{i=0}^m
n_i\overline{\alpha}_i$ with $n_i\ge 0$, in analogy to the
finite--dimensional case. Again we denote this relation by
$\overline{\lambda}\ge \overline{\mu}$.

\vskip .1in
\noindent \bf The fundamental weights. \rm If $\theta=\sum m_i
\alpha_i^\vee$, we define the fundamental weights by
$\overline{\lambda}_0=(0,0,1)$ (the ``vacuum'' weight) and
$\overline{\lambda}_i=(\lambda_i,0,m_i)$ for $i\ge 1$. Note that
these satisfy $(\overline{\lambda}_i,\overline{\alpha}_j^\vee)
=\delta_{ij}$ with respect to the Lorentzian form where for a
space--like vector $W$ (i.e.~one with $(W,W)>0$) we set $W^\vee =2
(W,W)^{-1} W$ as in the euclidean case. If we include $(0,0,-1)$ with
the simple roots and $\delta=(0,1,0)$ with the fundamental weights, we get
dual bases of $\overline{\h}$. (The inclusion of $\delta=(0,1,0)$ amongst the
fundamental weights again reflects the freedom to tensor positive
energy representations by characters of ${\rm Rot}\, S^1$.) 
\vskip .1in
\noindent \bf Dominant weights. \rm A weight
$\overline{\lambda}=(\lambda,k,\ell)$ is 
said to be dominant if $\lambda$ is dominant and $(\lambda,\theta)\le
\ell$. This is the permissibility condition for the highest weight of
$H(0)$. Plainly $\overline{\lambda}$ is dominant iff 
$(\overline{\lambda},\overline{\alpha}_i)\ge 0$ for all $i$. 
We will usually normalise the weight to have $k=0$,
using the freedom to tensor by a character of ${\rm Rot}\, S^1$. 
We then have the following analogue of the finite--dimensional
result. 

\vskip .1in
\noindent \bf Lemma. \it $\overline{\lambda}=(\lambda,0,\ell)$ is a
dominant weight iff $\overline{\lambda}=\sum_{i=0}^m n_i
\overline{\lambda}_i$ with $n_i\ge 0$.

\vskip .1in
\noindent \bf Proof. \rm This is immediate because
$(\overline{\lambda}_i,\overline{\alpha}_j^\vee)= \delta_{ij}$.

\vskip .1in
\noindent \bf The Weyl weight $\overline{\rho}$. \rm We define
$\overline{\rho}= (\rho,0,g)$. We already know that $\rho$ is a
weight, in fact $\rho=\sum_{i=1}^m \lambda_i$. An analogous statement
holds for $\overline{\rho}$, called the Weyl weight.

\vskip .1in
\noindent \bf Lemma. \it $\overline{\rho}=\sum_{i=0}^m
\overline{\lambda}_i$. 
\vskip .1in
\noindent \bf Proof. \rm We have
$\overline{\lambda_0}=(0,0,1)$ and $\overline{\lambda}_i=(\lambda_i,
0, m_i)$ for $i\ge 1$. Thus we have to prove that $g=1+\sum_{i\ge 1}
m_i$. Recall that $2g=(\theta,\theta) + 2 (\rho,\theta)$. We have
$\theta=\sum m_i \alpha_i^\vee$ (since $\theta^\vee=\theta$), so that
$$2g=\|\theta\|^2 + 2 \sum_{i=1}^m (\lambda_i,\theta)
=2 + 2\sum_{i=1}^m m_i,$$
as required. 

\vskip .1in

\noindent \bf The quantum Casimir operator. \rm We define the quantum
Casimir operator in any positive energy representation to be the
operator $\Omega=L_0 -d$.  
\vskip .1in
\noindent \bf Theorem. \it The quantum Casimir operator commutes
with ${\cal L}\g$. If $H$ is a positive energy
representation generated by a highest weight vector $v$ of
weight $\overline{\mu}= (\mu,k,\ell)$, then
$\Omega$ acts as the scalar
$[(\overline{\mu}+\overline{\rho},\overline{\mu}+\overline{\rho}
-(\overline{\rho},\overline{\rho})]/2(g+\ell)$ on $H$.

\vskip .1in
\noindent \bf Proof. \rm By construction $L_0-d$ commutes with
$\g$. Since $v$ is cyclic for ${\cal L}\g$, it suffices to show that
$\Omega$ acts as the given scalar on $v$.
Since $v$ is a highest weight vector
$X(n)v=0$ for $n>0$, so that 
$$\Omega v= [2(g+\ell)]^{-1}(-\sum X_i(0)^2v,v)]=[2(g+\ell)]^{-1}
(\|\mu+\rho\|^2 -\|\rho\|^2)v.$$
Since $dv=kv$, we are reduced to showing that
$$(\overline{\mu} +\overline{\rho},\overline{\mu}+\overline{\rho})
-(\overline{\rho},\overline{\rho}) = (\mu+\rho,\mu+\rho) -(\rho,\rho) 
-2(\ell+g)k,$$
which is immediate from the definition of the Lorentzian form.

\vskip .1in

\noindent \bf 8. THE AFFINE WEYL GROUP. \rm 
\vskip .1in
\noindent \bf Hyperbolic realisation. \rm Given a space--like vector
$\overline{\alpha}=(\alpha,k,0)$ with $\alpha\ne 0$, we define the
hyperbolic reflection by 
$$\sigma_{\overline{\alpha}} \overline{\mu} =
\overline{\mu} -2 {(\overline{\mu},\overline{\alpha})\over
(\overline{\alpha},\overline{\alpha})}\overline{\alpha}.$$
Let
$\widehat{W}$ be the group 
of transformations on $\overline{\h}^*$ generated by hyperbolic
reflections in the space--like roots. We call $\widehat{W}$ the affine
Weyl group.
\vskip .1in
\noindent \bf Proposition~1. \it The affine Weyl group permutes the
weights and weight multiplicities in an irreducible positive energy
presentation. 
\vskip .1in
\noindent \bf Proof. \rm If $\alpha\in \Phi$, set $E=E_\alpha(n)$,
$F=E_{-\alpha}(n)$ and $H=[E,F]=-2i\|\alpha\|^{-2} T_\alpha(0) +2n
\|\alpha\|^{-2}$. Then $E^*=F$, $H^*=H$, $[H,E]=2E$ and
$[H,F]=-2F$. Thus we have a copy of $\sl_2$. Now suppose that $w\in H$
has weight $\overline{\mu}=(\mu,k,\ell)$. Then $Hw=Mw$ where 
$M=\|\alpha\|^{-2} (2(\alpha,\mu) + 2n\ell)$. If $M\ge 0$, set $uF^M
w$. If $M<0$, set $u=E^{-M}w$. By the $\sl_2$ theory, we know that
$u\ne 0$. The weight of $u$ is $(\mu-\alpha M, k+nM,\ell)
=\sigma_{(\alpha,-n,0)} (\overline{\mu})$. Thus the weights are
invariant under the generators of $\widehat{W}$. Since the inner
product is invariant, orthogonal $w$'s of weight $\overline{\mu}$ give
rise to orthogonal $u$'s. This proves that the multiplicity of
$\sigma(\overline{\mu})$ is greater than or equal to the multiplicity
of $\overline{\mu}$. Applying $\sigma^{-1}$, we get the reverse
inequality, so the affine Weyl group preserves multiplicities. 
\vskip .1in
We compute $T=\sigma_{(\alpha,0,0)}^{-1}\sigma_{(\alpha,x,0)}$. We
have $T=T_{x\alpha^\vee}$ where $T_\beta(\mu,k,\ell)= (\mu + \ell
\beta, k-(\mu,\beta) -\ell \|\beta\|^2/2, \ell)$. It is easy to verify
that $T_\alpha \circ T_\beta =T_{\alpha + \beta}$ 
and $\sigma T_\alpha \sigma^{-1} =T_{\sigma\alpha}$ for $\alpha\in
h^*$ and $\sigma\in W$. For $\alpha\in Q^*$, the $T_{x\alpha}$'s
generate a copy of $Q^\vee$. The group $W$ normalises this translation
group, so the two groups together generate a group isomorphic to
$Q^\vee \rtimes W$. 

\vskip .1in
\noindent \bf Proposition~2. \it $\widehat{W}\cong Q^\vee \rtimes W$. 
\vskip .1in
\noindent \bf Proof. \rm We have just seen that $\widehat{W}$ lies in
the group generated by $W$ and the $T_{x\alpha}$'s ($\alpha\in
Q^\vee$). Clearly $W\subset \widehat{W}$ since
$\sigma_\alpha=\sigma_{(\alpha,0,0)}$.  In addition we can get
$T_{\alpha}\in \widehat{W}$ by taking $n=1$ in the above
discussion. Thus $\widehat{W}$ is generated by
$W$ and the $T_{x\alpha}$'s ($\alpha\in
Q^\vee$), as required.

\vskip .1in

\noindent \bf Euclidean realisation. \rm  For fixed level, we consider
the first component of the hyperbolic realisation. An affine root 
$(\alpha,n,0)$
gives rise to a transformation $X \mapsto \sigma_\alpha(X) + n\ell
\alpha^\vee$. Conjugating by the homothety $R(X)=\ell X$, we get
transformations 
$X\mapsto \sigma_\alpha(X)$ which give the usual euclidean generators
of the affine Weyl group.Let $\h_{\alpha,n}=\{X\in \h: 
\alpha(X)=n\}$, an affine hyperplane. Set $\h^\prime=\h\backslash
\bigcup \h_{\alpha,n}$. The connected components of $\h^\prime$ are
called Weyl alcoves $A$. The affine Weyl group $\hat{W}=\hat{W}(\Phi)$
acting in $\h$
is the group generated by reflections $\sigma_{\alpha,n}$ in the
hyperplane $\h_{\alpha,n}$.  Evidently
$$\sigma_{\alpha,n}(x)=\sigma_\alpha(x) +n\alpha^\vee.\eqno{(*)}$$
Let $\alpha^\vee=2 \alpha/(\alpha,\alpha)$, identifying $\h$ and
$\h^*$ using the normalised inner product. Let $Q$ be the lattice generated by
$\Phi$ in $\h^*$ and $Q^\vee$ the lattice generated by $\Phi^\vee$ in
$\h$. The Weyl group $W$ acts in both $\h$ and $\h^*$; it is generated
by $\sigma_{\alpha_i}=\sigma_{\alpha_i^\vee}$ where $\alpha_i$ are
walls of a Weyl chamber in $\h$. Under the identification $\h=\h^*$
induced by the inner product, the $\alpha_i^\vee$ are walls of the
corresponding chamber in $\h^*$. The Weyl group acts simply
transitively on the chambers. Every element of $\Phi$ is a positive or
negative integer combination of the $\alpha_i$'s. The following result
is an immediate corollary of Proposition~2; it can also easily be proved
directly.
\vskip .1in
\noindent \bf Proposition~2*. \it $\hat{W}=Q^\vee \rtimes W$, where
$Q^\vee$ acts on $\h$ by translation.
\vskip .1in
\noindent \bf Proof. \rm Formula $(*)$ shows that $\hat{W}\subseteq
Q^\vee \rtimes W$. Clearly $W\subset \hat{W}$, since
$\sigma_\alpha=\sigma_{\alpha,0}$. Moreover
$\sigma_{\alpha,1}\sigma-\alpha(x)=x+\alpha^\vee$. Hence
$Q^\vee\subset \hat{W}$.
\vskip .1in
\noindent \bf Corollary~1. \it Restriction to the first component with
$\ell=1$ gives an isomorphism between the hyperbolic and euclidean
affine Weyl groups. 
\vskip .1in
\noindent \bf Corollary~2. \it $Q^*\rtimes W$
permutes the $\h_{\alpha,n}$ 
and hence the Weyl alcoves. Hence the same is true of
$\widehat{W}=Q^\vee \rtimes W\subseteq Q^*\rtimes W$.
\vskip .1in
\noindent \bf Proof. \rm Clearly
$$\sigma(\h_{\alpha,n})=\h_{\sigma_\alpha,n},\eqno{(1)}$$ 
so $W$ permutes the
hyperplanes. Since $(Q^*,Q)\subset {\Bbb Z}$, $Q^*$ also
permutes the hyperplanes;
in fact if $X\in Q^*$ we have
$$X+ \h_{\alpha, n}=\h_{\alpha,n+\alpha(X)}.\eqno{(2)}$$
Hence $Q^*\rtimes W$ permutes the hyperplanes. 
\vskip .1in

\noindent \bf The highest root. \rm Fix a Weyl chamber $C$ and
and let $\alpha_1,\dots,\alpha_\ell$ be the corresponding simple
roots. Let $\theta$ be the highest weight 
of the adjoint representation on $\g$. Since $\theta$ is the highest
root and
$\sigma_{\alpha_i}\theta=\theta-(\theta,\alpha^\vee_i)\alpha_i$ is
also a root, we must have $(\theta,\alpha_i)\ge 0$. 
Since $\theta$ is a positive root, we may
write $\theta=\sum_{i=1}^n d_i\alpha_i$ with $d_i\ge 0$. Since
$\alpha_i$ is also a weight of $\g$, $\theta-\alpha_i\ge 0$. Thus
$d_i\ge 1$ for all $i$. Take $X_i\in \h$ with
$\alpha_i(X_j)=\delta_{ij}$, the dual basis of $Q^*$. Recall that we
have normalised the 
inner product on $\g$ (and hence $\h$ and $\h^*$) so that
$(\theta,\theta)=2$. As usual we identify $\h$ and $\h^*$ using
this normalised inner product. 
\vskip .1in

\noindent \bf Proposition 3. \it Let $A$ be the unique Weyl alcove
contained in $C$ with $0\in\overline{A}$. Then $A=\{X\in \h:
(X,\alpha_i)>0, (X,\theta)<1\}$ with closure
$\overline{A}=\{X\in \h:(X,\alpha_i)\ge 0, (X,\theta)\le 1\}$.
Moreover $A$ is a simplex with vertices $m_i^{-1} \lambda_i$, where
$\lambda_i$ are the fundamental weights with level $m_i$.
\vskip .1in
\noindent \bf Proof. \rm Clearly $A=\{X\in \h: \alpha(X)\in (0,1)
(\alpha\in \Phi^+)\}$. Let $A^\prime=\{X\in \h: \alpha_i(X)>0,
\theta(X)<1\}$. Plainly $A\subseteq A^\prime$. Conversely if $X\in A$,
then $\alpha(X)>0$ for all $\alpha\in \Phi^+$; and, since
$\theta-\alpha$ is a non--negative combination of simple roots for any
root $\alpha$, $\alpha(X)\le \theta(X)<1$ for $\alpha\in \Phi^+$. Thus
$X\in A$, so that $A^\prime\subseteq A$. Hence $A=A^\prime$. Since
$\theta=\sum d_i \alpha_i$ with $d_i\ge 0$, this means that
$A$ is a simplex with vertices at $0$ and the points
$Y_i$ given by $(Y_i,\alpha_j)=\delta_{ij}$. Since
$(\lambda_i,\alpha_j) =m_i\delta_{ij}$ (where $\theta^\vee =\sum m_i
\alpha_i^\vee$), we see that $Y_i=m_i^{-1}\lambda_i$ as claimed.

\vskip .1in
\noindent \bf Proposition~4. \it  The affine Weyl group permutes the
Weyl alcoves transitively and is generated by the reflections in the
walls of $A$ (``simple reflections''). In particular every weight is
in the affine Weyl group orbit of a dominant weight.
\vskip .1in
\noindent \bf Proof. \rm Let $\widehat{W}_0$ be the subgroup of
$\widehat{W}$ generated by the simple reflections. We first prove that
$\widehat{W}_0$ permutes the Weyl alcoves transitively. Let $A$ be the
standard alcove and let $A^\prime$ be another alcove. Take $X\in A$,
$X^\prime \in 
A^\prime$. Choose $\sigma\in \widehat{W}_0$ minimising $\|X-\sigma
X^\prime \|$ and set $Y=\sigma X^\prime$. If the line segment $[X,Y]$
crosses a hyperplane it must also cross a wall $\h_{i}$ of $A$;
but then $\|X-Y\|> 
\|X-\sigma_{\sigma_i}(Y)\|$, contradicting minimality since the simple
reflection $\sigma_i$ is in $\widehat{W}_0$. Thus $Y$ must
lie in $A$ and hence $A^\prime=\sigma A$ with $\sigma\in
\widehat{W}_0$. 

Since $\widehat{W}_0$ permutes the Weyl alcoves transitively, facts
about $A$ can be transported to $A^\prime \sigma A$. In particular any
alcove $A^\prime$ has well--defined walls. Since $\sigma A$ is an
alcove for any $\sigma\in Q^*\rtimes W$, it follows that every
hyperplane $\h_{\alpha,n}$ is the wall of some alcove; for 
(1) and (2) every hyperplane is the image of a wall of $A$ under
$Q^*\rtimes W$. Let $A^\prime$ be an alcove having $\h_{\alpha,n}$ as
a wall and take $\sigma\in \widehat{W}_0$ such that $A^\prime=\sigma
A$. Then $\sigma^{-1}\h_{\alpha,n}$ is a wall $\h_i$ of $A$. It follows that
$\sigma^{-1}\sigma_{\alpha,n}\sigma=\sigma_i$. Hence
$\sigma_{\alpha,n}=\sigma\sigma_i\sigma^{-1}$ lies in
$\widehat{W}_0$ and so $\widehat{W}=\widehat{W}_0$.

\vskip .1in

\bf \noindent Remark. \rm It is also true that the affine Weyl group
permutes the Weyl alcoves simply transitively and that the Weyl alcove
is a fundamental domain. Although we shall not need these results, we
note that they are easy to prove directly for the affine Weyl group of
$SU(N)$. Recall that the integer lattice
$\Lambda={\Bbb Z}^N$ acts by translation on ${\Bbb R}^n$. The symmetric
group $S_N$ acts on ${\Bbb R}^N$ by permuting the coordinates and
normalises $\Lambda$, so we get an action of the semidirect product
$\Lambda\rtimes S_N$. The subgroup $\Lambda_0=\{(N+\ell)(m_i):\sum
m_i=0\}\subset \Lambda$ is 
invariant under $S_N$, so we can consider the semidirect product
$W=\Lambda_0\rtimes S_N$. This is essentially the affine Weyl group of
$SU(N)$.

\vskip .1in
\noindent \bf Lemma~1. \it If $\sigma_i$ is the hyperbolic reflection
corresponding to a simple root $\overline{\alpha}_i$, then $\sigma_i$
permutes $\overline{\Phi}^+\backslash \{\overline{\alpha}_i\}$ (and
$\sigma_i\overline{\alpha}_i =-\overline{\alpha}_i$). Moreover the
$\sigma_i$ preserves the multiplicity of a root.
Hence each $\sigma\in \widehat{W}$ permutes the roots, preserving
their multiplicities.

\vskip .1in
\noindent \bf Proof. \rm We know that $\overline{\alpha}$ is a
positive root iff $\overline{\alpha} =\sum_{i=0}^m n_i
\overline{\alpha}_i$ with $n_i\ge 0$. If $\overline{\alpha} \ne
\overline{\alpha}_i$, it cannot be a multiple of $\overline{\alpha}_i$
so $n_j>0$ for some $j\ne i$. But $\sigma_i \overline{\alpha}
=\overline{\alpha} -t \overline{\alpha}_i$, so the coefficient of
$\overline{\alpha}_j$ in $\sigma_i\overline{\alpha}$ is also
$n_j$. Thus $\sigma_i\overline{\alpha}$ must be positive. Note that
$\sigma_i(0,n,0) =(0,n,0)$, so the $\sigma_i$'s preserve root
multiplicities.

\vskip .1in
\noindent \bf Lemma~2. \rm $\sigma_i
\overline{\rho}=\overline{\rho}-\overline{\alpha}_i$. 

\vskip .1in
\noindent \bf Proof. \rm We have $\sigma_i
\overline{\lambda}_j=\overline{\lambda}_j 
-\delta_{ij} \overline{\alpha}_i$ because
$(\overline{\lambda}_i,\overline{\alpha}_j^\vee) =\delta_{ij}$. 
This implies the result because $\overline{\rho}=\sum_{i=0}^m
\overline{\lambda}_i$. 

\vskip .1in

\noindent \bf 9. CONSTRUCTION OF IRREDUCIBLE REPRESENTATIONS. \rm  
We now prove an analogue of the Harish--Chandra theorem for an affine
Kac--Moody algebra: the proof is almost identical to the
finite--dimensional case. Let $\overline{\lambda}$ be a dominant weight.
For each simple root $\overline{\alpha}_i$,  let
$E_i,F_i,H_i$ be the basis of the 
Lie algebra $s\ell(2)_i$
corresponding to the simple root $\overline{\alpha}_i$. 
Let $\g_2=\overline{\h}_{\Bbb C}\oplus  
\bigoplus_{\alpha>0} \g_\alpha$ and $\g_1=\bigoplus_{\alpha<0}
\g_\alpha$. These are Lie subalgebras of $\widehat{\g}$ with
$\widehat{\g}=\g_1\oplus \g_2$. We know that $\g_2$ is generated by
the $E_i$'s and $\overline{\h}$ and $\g_1$ is generated by the
$F_i$'s. Consider the 1--dimensional
representation sending $E_i$ to $0$ and $H\in \overline{\h}$ to
$i\overline{\lambda}(H)$. Let $M(\overline{\lambda})$ be the
corresponding Verma module. 
Thus if $v=v_{\overline{\lambda}}$ is the highest weight
vector of $M(\overline{\lambda})$, we have $E_i
v_{\overline{\lambda}}=0$ and $ Hv_{\overline{\lambda}} = 
\overline{\lambda}(H) v_{\overline{\lambda}}$ where
$\overline{\lambda}(H_i)\in {\Bbb Z}_+$ for all 
$i$. We know that $M(\overline{\lambda})$ has a unique maximal
submodule $N$ such that 
$L(\overline{\lambda})= M(\overline{\lambda})/N$ is irreducible as a
$\widehat{\g}$--module. In fact, 
since $\overline{\h}$ is diagonalisable, every submodule is the sum of
its weight 
spaces. Hence if we take $N$ to be the algebraic sum of all proper
submodules, we muat have $v\notin N$, so that $N$ is the unique
maximal proper submodule.
By the $s\ell(2)$ theory, if
$\ell_i=(\overline{\lambda},\overline{\alpha}_i^\vee)$, then 
$w_i=F_i^{\ell_i +1}v_{\overline{\lambda}}$ is a singular 
vector i.e.~$E_iw=0$ and $w$ is an eigenvector for $\overline{\h}$. It
therefore 
generates a proper submodule (all weights are strictly less than
$\overline{\lambda}$). Hence $w_i\in N$ for all $i$. Let $N_0$ be the
submodule generated by the $w_i$'s. 
\vskip .1in
\noindent \bf Theorem (Harish--Chandra--Kac). \it
$L(\overline{\lambda})$ is the 
quotient of $M(\overline{\lambda})$ by the submodule generated by $F_i^{\ell_i
+1}v_\lambda$.
\vskip .1in
\noindent \bf Proof. \rm We have to show that $N=N_0$.
Set $L=M(\overline{\lambda})/N_0$. Thus $L$ is a cyclic
module for $\widehat{\g}$ generated 
by $v=v_{\overline{\lambda}}$ satisfying $X v=i \overline{\lambda}(X)
v$ for $X\in 
\overline{\h}$, 
$E_i v=0$ and 
$F_i^{\ell_i+1}v=0$. The identity 
$$[a^n,b]=\sum_{r=1}^n {n\choose r} [({\rm ad}\,a)^r b] a^{n-r}\eqno{(*)}$$
implies that the action on $L$ is
locally nilpotent, i.e.~some power of each $E_i$ or $F_i$ kills any
vector. For the $E_i$'s this follows because the $E_i$'s lower energy.
For the $F_i$'s it follows because $L$ is spanned by vectors
$F_{i_1}\cdots F_{i_k}v$ where $i_1,\dots,i_k$ are arbitrary (recall
that the $F_i$'s generate the $F_\alpha$ subalgebra). Starting from
the relation $F_i^{\ell_i+1}v=0$, 
successive application of $(*)$ and the Serre relations show that each
$F_i$ is nilpotent on any such monomial vector. This local nilpotence
shows that any vector in $L$ lies in a finite dimensional $s\ell(2)_i$
module for each $i$. 

We claim that the weights of $L$ are invariant under the affine Weyl group
$\widehat{W}$. In fact suppose $w\in L$ has weight $\overline{\mu}$. Then
$H_i w=m_i w$ with 
$m_i=\mu(H_i)=(\mu,\alpha_i^\vee)$. Then 
$w$ lies in a sum of $s\ell(2)_i$ modules. If $m_i\ge 0$, set
$u=F_i^{m_i} w$ and if $m_i<0$, set $u=E_i^{-m_i}$. Thus $u\ne 0$ by
the $s\ell(2)$ theory and $u$ has weight $\overline{\lambda} -m_i
\overline{\alpha}_i 
=\sigma_i \overline{\mu}$. Thus the set of weights is invariant under each
simple reflection $\sigma_i$ and hence the whole of $\widehat{W}$. As
a consequence of this reasoning we have the following result.
\vskip .1in
\noindent \bf Lemma. \it If
$\overline{\lambda}$ is dominant, then $\sigma\overline{\lambda} \le
\overline{\lambda}$ for all $\sigma\in\widehat{W}$.
\vskip .1in
\noindent \bf Corollary. \it $\tau\overline{\rho}=\overline{\rho}$
for $\tau\in 
\widehat{W}$ iff $\tau=1$. 

\vskip .1in
\noindent \bf Proof. \rm Let $\overline{\lambda}_i$ be the
fundamental weights of 
$\widehat\g$. We can apply the lemma to these.  Since $\overline{\rho}
=\sum_{i\ge 0} \overline{\lambda}_i$, 
the equality $\tau \overline{\rho} =\overline{\rho}$ and the inequalities
$\overline{\lambda}_i \ge \tau \overline{\lambda}_i$ force
$\tau\overline{\lambda}_i =\overline{\lambda}_i$ for $i\ge 0$. Hence
$\tau=1$.
\vskip .1in
Now suppose that $L$ is not irreducible. Then $V$ 
must contain a singular vector $w$ of weight $\overline{\mu}$ strictly
lower than 
$\overline{\lambda}$: thus $E_iv=0$ and $H_iw= m_iw$ where
$m_i=\overline{\mu} (h_i)\le
\ell_i$. But then $w$ is a highest weight vector generating an irreducible
representation of each $s\ell(2)_i$. On the other
hand let $\Omega$ be the quantum Casimir operator of
$\widehat{\g}$ and set $C=2(\ell+g)\Omega$. Then $C
v=((\overline{\lambda} +\overline{\rho},\overline{\lambda}
+\overline{\rho}) -(\overline{\rho},\overline{\rho}))v$, 
so by cyclicity $\Omega=
((\overline{\lambda} +\overline{\rho},\overline{\lambda}
+\overline{\rho}) -(\overline{\rho},\overline{\rho}))I$.
Since $\Omega w=
((\overline{\mu} +\overline{\rho},\overline{\mu}
+\overline{\rho}) -(\overline{\rho},\overline{\rho}))w$, 
we must have 
$((\overline{\lambda} +\overline{\rho},\overline{\lambda}
+\overline{\rho}) -(\overline{\rho},\overline{\rho}))
=
((\overline{\mu} +\overline{\rho},\overline{\mu}
+\overline{\rho}) -(\overline{\rho},\overline{\rho}))$. As in the
finite--dimensional case, the proof is completed by the contradiction
implied by Freudenthal's lemma. 
\vskip .1in
\noindent \bf Freudenthal's Lemma. \it Let
$\overline{\lambda}=(lambda,0,\ell)$ be a dominant and let
$\overline{\mu}=(\mu,k,\ell)$ be another weight such that
$\overline{\lambda}-\overline{\mu}=\sum_{i=0}^m  n_i
\overline{\alpha_i}$ with $n_i\ge 0$. Then
$(\overline{\lambda}+\overline{\rho},\overline{\lambda}+\overline{\rho})
\ge 
(\overline{\mu}+\overline{\rho},\overline{\mu
}+\overline{\rho})$ with equality iff
$\overline{\lambda}=\overline{\mu}$.
\vskip .1in 
\noindent \bf Proof. \rm Take $\tau\in \widehat W$ such that
$\tau(\overline{\mu}
+\overline{\rho})$ is dominant. Thus $\tau \overline{\mu}\le
\overline{\lambda}$ and $\tau \overline{\rho}\le \overline{\rho}$ by
the lemma. But then
$$0=(\overline {\lambda} +\overline{\rho} +\tau \overline{\mu} +\tau
\overline{\nu}, \overline{\lambda} - \tau \overline{\mu}  +\overline{\rho}-\tau
\overline{\rho}) \ge (\overline{\rho}, \overline{\lambda} -\tau\overline{\mu}) 
+(\overline{\rho}, \overline{\rho}-\tau \overline{\rho}) \ge 0.$$
Hence $(\overline{\rho}, \overline{\lambda} -\tau\overline{\mu})=0$ and
$(\overline{\rho}, \overline{ \rho}-\tau \overline{\rho})=0$, so that
$\overline{\lambda} =\tau\overline{\mu}$ and $\overline{\rho}=\tau \overline{
\rho}$. By the corollary above, $\tau=1$ and hence
$\overline{\lambda}=\overline{\mu}$, as required.

\vskip .1in

\noindent \bf 11. GARLAND'S `NO--GHOST' THEOREM ON UNITARITY. \rm Let
$H=L(\overline{\lambda})$ be the irreducible representation just
constructed as a 
quotient of the Verma module $M(\overline{\lambda})$. Consider the algebraic
dual $L^*(\overline{\lambda})=\bigoplus H(k)^*$; if we take complex 
multiplication on the dual to be given by $z\cdot
\xi=\overline{z}\xi$, it is easy to verify that the canonical action
of $\widehat{g}$ is positive energy of level $\ell$ with highest
weight $\overline{\lambda}$. By Schur's lemma there is a unique isomorphism
of $L(\overline{\lambda})$ onto this module, we get a $\widehat{g}$--equivariant
linear map from $L(\overline{\lambda})$ onto its conjugate dual. Since
any such 
map gives and is equivalent to an invariant sesquilinear form on
$L(\overline{\lambda})$, we deduce that there is an essentially unique
invariant 
sesquilinear form $(v,w)$ on $L(\overline{\lambda})$. Since the form is
invariant, its kernel is $\widehat{g}$--invariant and hence trivial by
irreducibility. Thus the form is non--degenerate on
$L(\overline{\lambda})$. Since $D$ is self--adjoint, the 
energy spaces $H(k)$ must be orthogonal and the form non--degenerate
on each of these. In particular the form is non--degenerate on $H(0)$,
an irreducible $\g$--module. By the finite--dimensional theory and
Schur's lemma, it is proportional to an inner product on $H(0)$. We
may therefore assume its restriction to $H(0)$ is positive definite. 

\vskip .1in
\noindent \bf Theorem. \it If
$L(\overline{\lambda})$ is an irreducible positive energy presentation
at level $\ell$
with $(\lambda,\theta) \le \ell$, then the canonical invariant
sesquilinear form on $L(\overline{\lambda})$ is positive definite.
\vskip .1in
\noindent \bf Proof. \rm Let $(v,w)$ be the invariant sequilinear form
on $H=L(\overline{\lambda})$. By irreducibility $(\cdot,\cdot)$ is
non--degenerate on $H$. Clearly the spaces $H(k)$ are orthogonal with
respect to $(\cdot,\cdot)$ because $D=\pi(d)$ is self--adjoint.
The form must also be non--degenerate on each $H(k)$. Each $H(k)$ is a
finite--dimensional $\g$--module and therefore completely reducible. In
particular the action of $\h$ is diagonalisable on $H(k)$, so that
$H(k)$ breaks up as a sum of weight spaces for $\h$. Since $\h$ acts
as skew--adjoint operators with respect to $(\cdot,\cdot)$, these
eigenspaces must be mutually orthogonal. To prove that $(\cdot,\cdot)$
is positive definite, it therefore suffices to show that $(v,v)\ge 0$
for any vector in $H(k)$ that is a highest weight vector for $\g$. We
prove this by induction on 
$H(k)$. For $k=0$, this follows from the no--ghost theorem for $\g$
proved in Chapter~2. We therefore assume that $(\cdot,\cdot)$ is
positive definite on $\bigoplus_{j\le k}$ and show that $(v,v)\ge 0$
for $v\in H(k+1)$ of weight $\mu$. Now
$$2(\ell+g)((L_0-D)v,v)=(\ell+g)[(\overline{\lambda}
+\overline{\rho},\overline{\lambda}+\overline{\rho})
-(\overline{\rho},\overline{\rho})](v,v).\eqno{(1)}$$ 
But we also have
$$\eqalign{2(\ell+g)((L_0-D)v,v)&={1\over 2} \sum_i (X_i(0)v,X_i(0)v)
-2(\ell+g)(Dv,v) +\sum_{i,n>0} (X_i(n)v,X_i(n)v)\cr
&\ge (\ell+g)(\|\mu +\rho\|^2 -2k)(v,v)\cr
&=(\ell+g)[(\overline{\mu}+\overline{\rho},\overline{\mu}+\overline{\rho})
-(\overline{\rho},\overline{\rho})](v,v),\cr}$$ 
since $(X_i(n)v,X_i(n)v)\ge 0$ for $n>0$ by the induction hypothesis.
Combining this equation with (1) we get
$$[(\overline{\lambda}+\overline{\rho},\overline{\lambda}+\overline{\rho})
   -(\overline{\mu}+\overline{\rho},
\overline{\mu}+\overline{\rho})](v,v)\ge 0.$$ 
By Freudenthal's lemma, $
(\overline{\lambda}+\overline{\rho},\overline{\lambda}+\overline{\rho})
   -(\overline{\mu}+\overline{\rho},
\overline{\mu}+\overline{\rho})>0$, so we obtain $(v,v)\ge 0$ as required.
\vskip .1in
\noindent \bf 11. THE CHARACTER OF A POSITIVE ENERGY
REPRESENTATION. \rm Our aim now is to determine the character of a
unitary irreducible positive energy representation
$H=L(\overline{\lambda})$. If $H$ is a positive energy representation of
$\widehat{\g}$ or $\widehat{\h}$ and $H=\bigoplus H(n)$, we define the
character ${\rm ch}\, L(\overline{\lambda})$ to be the formal power
series
$\sum_{n\ge 0} q^n {\rm Tr}_{H(n)}(z)$ for $z=e^T$ with $T\in
\h$. Although defined as a formal power series in $q$, the character
converges absolutely for $|q|<1$ and $T\in \h_{\Bbb C}$. We can write
the character as ${\rm Tr}(q^d z)$. It turns out that $q^d$ is a
trace--class operator for $|q|<1$ (for $0<q<1$, this means that the
positive operator $q^d$ is diagonalisable with summable
eigenvalues). Actually to make the characters invariant under the
modular group, it is more natural to take the normalised characters
${\rm Tr}(q^{L_0-c/24} z)$. 

\vskip .1in

\noindent \bf 12. BOSONS AND FERMIONS ON THE CIRCLE. \rm 
\vskip .1in
\noindent \bf BOSONS. \rm Consider the operators the Lie algebra
generated by the operators $X(n)$ and $d$, where $X\in \h$. If we take
an orthonormal basis $(X_i)$ of $\h$ and set $X_i(n)=iA_i(n)$, then we
obtain the boson algebra, also often called the
oscillator or Heisenberg algebra. If $M={\rm dim}\, \h$, then we have
$M$ commuting bosonic fields, which physicists usually suggestive
write as $X_i(z)=\sum X_i(n) z_{-n-1}$. The boson algebra is the
infinite--dimensional Lie 
algebra with basis $(A_i(n))$ ($n\in {\Bbb Z}$)
satisfying the commutation relations 
$[A_i(m),A_j(n)]=m\delta_{ij}\delta_{m+n,0}$. It has a
conjugate--linear involution 
given by $A_i(n)^*=A_i(-n)$ and a derivation $d$ given
by $[d,A_i(n)]=-nA_i(n)$ with $d^*=d$. Note that the zero modes
$A_i(0)$ are central in the 
semidirect product. Just as with affine algebras, we may consider
positive energy unitary representations. As before every positive
energy representation is a direct sum of irreducible positive energy
representations. 
\vskip .1in
\noindent \bf Theorem. \it The boson algebra has a unique irreducible
positive energy representation, up to tensoring by characters of $d$
and the $A_i(0)$'s.
\vskip .1in
\noindent \bf Proof. \rm Let $H$ be a positive energy irreducible
representation. Then the $A_i(0)$'s  leave the lowest energy subspace $H(0)$
invariant and form a commuting self--adjoint set. Thus we can find
$\Omega\in H(0)$ such 
that $A_i(0)\Omega =\mu_i\Omega$. The vector $\Omega$ is annihilated by the
creation operators $A_i(n)$ for $n>0$, so by the commutation relations (or
the Poincar\'e--Birkhoff--Witt theorem) the submodule generated by $\Omega$
is spanned by vectors $\prod A_i(-n)^{k_{i,n}}\Omega$ where
$k_{i,n}\ge 0$, $n\ge 1$. (Note that the operators $A_i(-n)$ commute
for $n<0$.)  It is invariant under 
$d$ and the $A_i(0)$'s and therefore coincides with $H$ by irreducibility.
Clearly $H(0)={\Bbb C}\Omega$. Note that each pair of elements
$A=\sqrt{n} A_i(n)$ and $A^*=\sqrt{n}A_i(-n)$ ($n>0$) satisfies the
Heisenberg commutation relations
$AA^*-A^*A=I$. This is a copy of the usual Heisenberg Lie algebra. $A$
is called an annilation operator and $A^*$ a creation operator. The
following is a standard computation in quantum mechanics.

\vskip .1in

\noindent \bf Lemma. \it Let
$D=A^*A$ and $\xi_n=A^{*n}\xi_0$, where $\xi_0$ is a vector
satisfying $A\xi_0=0$. 

\noindent (a) $AA^{*n} -A^{*n}A= n A^{*(n-1)}$.

\noindent (b) $D\xi_n =n \xi_n$. 

\noindent (c) $(\xi_n,\xi_m)=\delta_{nm} n! (\xi,\xi)$. 
\vskip .05in
\bf \noindent Proof. \rm (a) follows by
induction from $[A,A^*]=I$, since
$AA^{*n}-A^{*n} A= (AA^{*n-1} -A^{*n-1}A) A^* +A^{*n-1} (AA^*-A^*A) =
n A^{*n-1}$.  Next by (a),  $D\xi_n=A^*AA^{*n} \xi_0 = (A^{*n+1} A +
n A^{*n})\xi_0 = n\xi_n$, since $A\xi_0=0$. So (b) follows. Since the
$\xi_n$'s correspond to different 
eigenvalues of the self--adjoint operator $D$, they must be pairwise
orthogonal.  Moreover we have
$$(\xi_n,\xi_n)=(A^{*n} \xi_0,A^{*n} \xi_0) =(AA^{*n} \xi_0, A^{*n-1} \xi_0)
=(A^{*n} Axi_0 + n A^{*n-1}\xi_0,A^{*n-1} \xi_0)=n(A^{*n-1} \xi_0,A^{*n-1}
\xi_0).$$
Thus $\|\xi_n\|^2 = n \|\xi_{n-1}\|^2$, so
the result follows by induction.
\vskip .1in
\noindent \bf Prophetic remark. \rm Note that $[D,A]=-A$ and
$[D,A^*]=A^*$. Mathematically this may be viewed as 
part of the metaplectic action of the symplectic
Lie algebra, the bosonic version of spin quantisation for
fermions. Below we will use bilinears in bosons to give the
Fubini--Veneziano construction of an $L_0$ operator implementing $d$.
In fact we will give a construction of the entire Virasoro algebra
with central charge $1$. Like the Sugawara construction, this is
another example of quantisation in infinite dimensions. 
\vskip .1in
The computation above implies that all the vectors
$\prod A_i(-n)^{k_{i,n}}\Omega$ are orthogonal with
$$\|\prod A_i(n)^{k_{i,n}}\Omega\|^2 = \prod  k_{i,n}! \,
n^{k_{i,n}}.$$
Conversely the standard Verma module
construction produces a similar basis generated from a vector $\Omega$
with $A_i(0)\Omega =\mu_i\Omega$,
$d\Omega=0$. As before, there is unique invariant sesquilinear form on
it. By the previous computation it is positive definite and coincides
with the formula given above. It may be realised on polynomials
${\Bbb C}[z_{i,n}]$ by $A_i(n)=n\partial_{z_{i,n}}$, $A_i(-n)=z_{i,n}$,
$A_i(0)=\mu I$ and $dz_{i,n}=nz_{i,n}$. To prove the representation is
irreducible we proceed as in Chapter~1, acting by annihilation
operators $\partial_{z_{i,n}}$ until we get the vacuum and then acting by
creation operators until we get all vectors. 
\vskip .1in
\noindent \bf Corollary. \it The character of the irreducible
representation $H$ of $\widehat{\h}$ with $X_j(0)$ acting as $i\mu_j$ and
$d$ as $h$ on $H(0)$ is $e^{i\mu(X)} q^h \prod_{n\ge 1} (1-q^n)^{-1}$. 

\vskip .1in
As promised we now produce the quantised action of the Virasoro
algebra, which is a simpler case of the Sugawara construction (for an
Abelian Lie algebra rather than a simple Lie algebra). 
\vskip .1in
\noindent \bf Proposition (Fubini--Venziano construction). \it Let
$$L_0^\h = \sum_i {1\over 2} A_i(0)^2 +\sum_{j>0} A_i(-j)A_i(j)
=-\sum_i {1\over 2} X_i(0)^2 +\sum_{j>0} X_i(-j)X_i(j)$$ and $L_m={1\over
2}\sum_{i,j} A_i(-j) A_i(j+m)$ for $m\ne 0$. Then $[L_m,A_i(n)]=-n A_i(n+m)$
and $[L_m,L_n]=(m-n)L_{m+n} +M{m^3-m\over 12} \delta_{m+n,0}$. Thus
the central charge is $M$, the number of bosons.
\vskip .1in
\noindent \bf Proof. \rm To check that $[L_0,A_k(n)]=-n A_k(n)$, we compute
$$[L_0,A_k(n)]=\sum_i \sum_{j>0} [A_i(-j),A_k(n)]A_i(j) +A_i(-j)
[A_i(j),A_k(n)] =-nA_k(n).$$ 
A similar computation shows that $[L_m,A_k(n)]=-n AS_k(n+m)$. It follows
that $[L_m,L_n]-(m-n)L_{m+n}$ commutes with all $A_k(n)$'s and hence
equals a scalar. Since this operator lowers energy by $m+n$, this
scalar is zero unless $m+n=0$. Thus if $m>0$, we have
$[L_{m},L_{-m}]=2mL_0 +\lambda(m)I$. Now $L_{-m}\Omega={1\over
2}\sum_{j=0}^m 
A_i(-j) A_i(j-m) \Omega$, so that
$$\eqalign{[L_m,L_{-m}]\Omega&=L_mL_{-m}\Omega ={1\over 2}
\sum_i\sum_{j=0}^m L_m 
A_i(-j) A_i(j-m)\Omega\cr
&={1\over 2}\sum_i \sum_{j=0}^m j
A_i(-j+m)A_i(j-m)\Omega + 
(m-j) A_i(-j) A_i(j)\Omega.\cr}$$
Thus 
$$([L_m,L_{-m}]\Omega,\Omega)=m \mu^2/2 +{1\over 2} \sum_{j=0}^m 
j(a_{-j}\Omega,a_{-j}\Omega)=m\mu^2/2 +{1\over 2}\sum_{j=0}^m j(m-j) =m\mu^2
+M(m^3-m)/12.$$ 
Hence $\lambda(m)=M(m^3-m)/12$ as required.
\vskip .1in
\noindent \bf FERMIONS. \rm The fermion algebra (or Clifford algebra)
is the Lie superalgebra with basis $\psi_i(n)$ ($n\in {\Bbb Z}$, $1\le
i\le N$)
satisfying the anticommutation relations
$\{\psi_i(m),\psi_j(n)\}=\delta_{m+n,0}\delta_{ij}$ and the adjoint condition
$\psi_i(n)^*=\psi_i(-n)$. It has a derivation $d$ given by
$[d,\psi_-i(n)]=-n\psi_i(n)$. Note that this is really a collection of
$N$ independent fermi fields on the circle. The zero modes form a
subalgebra which can be identified with the real Clifford algebra on
an $N$--dimensional real inner product space (let $c_i=\sqrt{2}
\psi_i(0)$). As we have seen in Chapter~I, these finite--dimensional
algebras behave differently for $N$ even or odd, so for this reason our
account is not entirely parallel to the bosonic case. 

We start by
defining a cyclic representation of the fermion algebra on the
exterior algebra on the unit vectors $v_{i,j}$ ($j\ge 0$). $\psi_i(n)$
acts as $e(v_{i,n})^*$ for $n>0$ and $e(v_{i,-n})$ for $n<0$;
$\psi_i(0)$ acts as ${1\over \sqrt{2}} (e(v_{i,0})+e(v_{i,0})^*)$. The
operator $d$ acts as the (even) derivation $dv_{i,n}=nv_{i,n}$. It is
clear that we have defined a positive energy representation $H$ with
lowest energy space $H(0)$, the exterior algebra on the
$v_{i,0}$'s. This is not an irreducible representation of the zero
mode algebra, but the usual creation--annihilation argument shows that
any irreducible submodule $W\subset H(0)$ generates an irreducible
representation of the fermion algbebra. When $N$ is odd, this
representation will not be graded, because the grading operator does
not lie inside the Clifford algebra. 

However when $N$ is even, the
grading operator is proportional to $\psi_1(0)\cdots \psi_N(0)$ and
$W$ is automatically graded. The irreducible module generated by $W$ is
clearly isomorphic to graded tensor product of $W$ and the exterior
algebra on the generators $v_{i,n}$ with $n\ge 1$. Conversely any
positive energy irreducible representation must have this form: for
$H(0)$ must be irreducible as a module over the zero modes and if
$w_j$ is an orthonormal basis of $W=H(0)$, then the vectors in the
Verma module obtained by applying products of distinct
$\sqrt{2}\psi_i(-n)$'s to different $w_j$'s are all orthonormal.

The fermionic version of the Fubini--Venziano construction is defined
on the exterior algebra via 
$$L_0={N\over 16} +\sum_{i;j>0} j\psi_i(-j)\psi_i(j),\qquad
L_k={1\over 2}\sum_{i,j} j\psi_i(-j)\psi_i(j+k) \quad (k\ne 0).$$
As before we prove that $[L_0,\psi_i(n)]=-n\psi_i(n)$ and
$[L_k,\psi_i(n)]=-(n+{k\over 2}) \psi_i(n+k)$. Again
$[L_m,L_n]-(m-n)L_{m+n}$ commutes with $\psi_i(n)$ and lowers energy
by $m+n$. If $m+n>0$, it therefore acts trivially on the vacuum vector
$\Omega=1$ and hence everywhere. Thus $[L_m,L_n]=(m-n)L_{m+n}$ for
$m+n>0$. Taking adjoints the same is true for $m+n<0$. We check
directly that $[L_m,L_{-m}]\Omega={N\over 24} (m^3-m) +{mN\over 16}$,
by a computation similar to the bosonic one. Hence 
$$[L_m,L_n]=(m-n) L_{m+n} + {N\over 24} (m^3-m) \delta_{m+n,0},$$
so that we get a quantised representation of the Virasoro algebra,
with central charge $c=N/2$. 

\vskip .1in
\noindent \bf 13. THE KAZAMA--SUZUKI SUPERCHARGE OPERATOR. \rm 
\vskip .1in
\noindent \bf A. GKO supercharge operator in ${\rm Cliff}\, L\g$. \rm
Let $(X_a)$ be an orthonormal 
basis of $\g$ 
with $[X_a,X_b]=\sum f_{abc}X_c$. Let $(\psi_a(n))$ ($n\in
{\Bbb Z}$) be Ramond
fermions on $L\g$, so that
$\{\psi_a(m),\psi_b(n)\}=\delta^{ab}\delta_{m+n,0}$.
Let ${\cal F}_\g$ be the fermionic Fock space giving a cyclic
representation representation of the $\psi_a(n)$'s. Let $S_a(n)=-{1\over 2}
\sum f_{abc}\psi_b(m)\psi_c(n-m)$ on ${\cal F}_\g$. Then the
anticommutation relations for $\psi$ immediately imply that
$$[S_a(n),\psi_b(m)]= \sum f_{abc}\psi_c(n+m),$$
and $[d,S_a(n)]=-nS_a(n)$ so that
by uniqueness $[S_a(m),S_b(n)]=\sum f_{abc} S_c(n+m) + \delta_{n+m,0} C(m)$,
where $C(m)$ is a constant. Taking vacuum expectations, i.e.~computing
$([S_a(m),S_b(-m)]\Omega,\Omega)$, we find $C(m)=
-mg$, with $g$ the dual Coxeter number. Let
$Q_0={1\over 3}\sum_{a,m} \psi_a(m) S_a({-m})$ (note that
$[\psi_a(m),S_a(n)]=0$ so that
no normal ordering is required). Clearly $S_a(m)^*=-S_a(m)$ and
$Q_0^*=-Q_0$. Then 
$$\eqalign{3\{Q_0,\psi_b(n)\}
&= \{\psi_a(m),\psi_b(n)\} S_a({-m}) + \psi_a(m)
[S_a(-m), \psi_b(n)]\cr
&= S_b(n) + \sum f_{abc}\psi_a(m)\psi_c(n-m)\cr
&=S_b(n) -\sum f_{bac}\psi_a(m)\psi_c({n-m})\cr
&=3S_b(n).\cr}$$
Thus $\{Q_0,\psi_b(n)\}= S_b(n)$.

It is not possible to compute $[Q_0, S_b(n)]$ by this method,
so we proceed indirectly. Let $\phi_b(n) = [Q_0, S_b(n)]$.
Then 
$$\eqalign{\{\phi_b(n),\psi_a(m)\} 
&=\{[Q_0,S_b(n)],\psi_a(m)\}\cr
& =[\{Q_0,\psi_a(m)\}, S_b(n)] +\{Q_0,[S_b(n),\psi_a(m)]\}\cr
&=[S_a(m),S_b(n)] +\sum \{Q_0,f_{bac}\psi_c(n+m)\}\cr
&=\sum f_{abc}S_c({m+n}) +\sum f_{bac}S_c({m+n})
-mg\delta_{ab}\delta_{m+n,0} \cr 
&=-mg\delta_{ab} \delta_{m+n}.\cr}$$
If $n\ne 0$, let $\xi=-(gn)^{-1}\phi_b(n)+\psi_b(n)$. Then $[d,\xi]=-n\xi$ and
$[\xi, \psi_a(m)]=0$. So by cyclciity $\xi$ must be a scalar and
hence zero. Thus $[Q_0,S_n^b]=+ng\psi_b(n)$; this relation also holds
for $n=0$ because the construction is manifestly $\g$--invariant. In
summary we have obtained the ``supersymmetry relations'' 
$$\{Q_0,\psi_b(n)\}= S_b(n),\qquad [Q_0,S_b(n)]=ng\psi_b(n);\eqno{(*)}$$
the supercharge operator thus interchanges fermions $\psi_b(n)$
and bosons $S_b(n)$. 

\vskip .1in
\noindent \bf B. GKO supercharge operator in ${\rm End}(H)\otimes {\rm
Cliff}\, L\g$. \rm Let $L(\overline{\lambda})$ be a level $\ell$
positive energy 
representation of ${\cal L}\g$ with corresponding generators
$T_a(m)$. These satisfy the commutation relations
$$[T_a(m),T_b(n)]= \sum f_{abc} T_c({m+n}) -\ell m \delta_{ab}\delta_{m+n,0}.$$
We extend the operator $Q_0$ to $L(\overline{\lambda})\otimes {\cal F}_\g$ as
$Q_0\equiv I\otimes Q_0$. Now consider the operator $Q_1=\sum
T_a(m)\psi_a(-m)$ on 
$L(\overline{\lambda})\otimes {\cal F}\g$. Then
$$\{Q_1,\psi_b(n)\} =\sum T_a(m)\{\psi_a(-m),\psi_b(n)\} = \sum T_a(m)
\delta_{ab} 
\delta_{mn} =T_b(n),$$
so that $\{Q_1,\psi_b(n)\}=T_b(n)$. Moreover
$$[Q_1, S_b(n)]= \sum T_a(m) [\psi_a({-m}), S_b(n)]=-\sum f_{bac}
T_a(m)\psi_c({n-m})\eqno{(1)}$$ 
and $[Q_1, T_b(n)] = \sum [T_a(m),T_b(n)]\psi_a({-m})$, so that
$$[Q_1,T_b(n)] =\sum f_{abc}T_c({m+n})
\psi_a({-m}) -\ell n\psi_b (n)
=\sum f_{bac} T_a(m)\psi_c({n-m})-\ell n\psi_b(n).\eqno{(2)}$$
Adding (1) and (2), we get $[Q_1,T_b(n) +S_b(n)]=n\ell \psi_b(n)$. Thus we
have 
$$\{Q_1,\psi_b(n)\}=T_b(n), \qquad [Q_1,T_b(n) +S_b(n)]=n\ell
\psi_n^b.\eqno{(**)}$$ 
Let $Q=Q_0+Q_1$ and $X_b(n)=T_b(n) +S_b(n)$. The supercharge operator is
therefore given by the formula
$$Q= \sum (T_a (m) + {1\over 3} S_a(m)) \psi_a({-m});$$
the factor  $1/3$ is very important here and is often given incorrectly in
much of the literature. Combining $(*)$ and $(**)$,
we get the supersymmetry relations
$$\{Q,\psi_a(n)\}= X_a(n),\qquad \{Q, X_a(n)\}=n(\ell+g) \psi_a(n).$$
\vskip .1in
\noindent \bf C. Kazama--Suzuki supercharge operator. \rm We keep the
above notation, but write $Q=Q^\g$ to show the dependence of the
construction on
$\g$. Thus
$$\{Q^\g,\psi_a(n)\}= X_a(n),\qquad \{Q^\g, X_a(n)\}=n(\ell+g)
\psi_a(n).\eqno{(3)}$$ 
Let $\h$ be a maximal torus in $\g$ with orthogonal complement
$\m$. We may choose the orthonormal basis of $\g$ to be made up of
orthonormal bases $(X_A)$ for $\h$ and $(X_i)$ for $\m$. We then have
fermions $\psi_A(n)$ and $\psi_i(n)$. We take the submodule of ${\cal F}_\g$
given by the tensor product ${\cal F}_\m\otimes {\cal F}_\h$, where the
$\psi_A(n)$'s act only on the first factor, irreducibly, and the
$\psi_i(n)$'s act only on the second factor. Since $S_A(n)$
commutes with all $\psi_B(n)$'s, they act only on the first factor;
indeed we have $S_A(n)=-{1\over 2} \sum f_{Abc} \psi_b({m}) \psi_c({n-m})$ and
$f_{ABc}\equiv 0$, so that $S_A(n)=-\sum {1\over 2} f_{Aij} \psi_j({m})
\psi_j({n-m})$.  Consider the representation of ${\cal L}\h$ on 
$K(\lambda)=L(\lambda)\otimes {\cal F}_\m$ given by $Y_A(n) = T_A(n) +
S_A(n)=X_A(n)$. Thus $[Y_A(m),Y_B(n)] =-m(\ell+g)
\delta_{AB}\delta_{m+n,0}$, since 
$\h$ is Abelian. The supercharge operator $Q^\h$ corresponding to $\h$
on $K(\lambda) \otimes {\cal F}_\h$ is just
$$Q^\h= \sum Y_A(m) \psi_A(-m).$$
It satisfies
$$\{Q^\h,\psi_A(n)\}=Y_A(n),\qquad \{Q^\h, Y_A(n)\}=n(\ell+g)\psi_A(n).
\eqno{(4)}$$
The Kazama--Suzuki supercharge operator is defined by $Q=Q^\g-Q^\h$
(the supersymmetric coset construction). Comparing (3) and (4), we see
that 
$$\{Q,\psi_A(n)\}=0,\qquad [Q,X_A(n)]=0.$$
The first equation tells us that $Q$ really acts on
$L(\overline{\lambda})\otimes 
{\cal F}_\m$, while the second implies that $Q$ commutes with the
natural action of ${\cal L}\h$ there. To see why $Q$ explicitly why $Q$ acts
on $L(\overline{\lambda})\otimes {\cal F}_\m$, recall that
$Q^\g=\sum  (T_a(m) + {1\over 3} S_a(m)) \psi_a(-m)$
and $Q^\h=\sum  Y_A(m) \psi_A(-m) = \sum (T_A(m) + S_A(m))\psi_A(-m)$.
Thus we get
$$Q=\sum (T_i(m) +{1\over 3} S_i(m))\psi_i(-m),$$
where 
$$S_i(n) = -{1\over 2} \sum f_{ijk} \psi_i(m) \psi_j({n-m}).$$
Thus 
$$Q=\sum T_i(m) \psi_i({-m}) -{1\over 6} \sum f_{ijk}\psi_i(m) \psi_j(n)
\psi_k(-m-n),$$
which evidently acts on $L(\overline{\lambda})\otimes {\cal F}_\m$.

\vskip .1in

\noindent \bf 14. THE SQUARE OF THE SUPERCHARGE OPERATOR. \rm
\vskip .1in
\noindent \bf A. Computations in ${\rm Cliff}\, L\g$. \rm We already
have proved the supersymmetry relations
$\{Q_0,\psi_b(n)\}=S_b(n)$ and $[Q_0,S_b(n)]=ng\psi_b(n)$. 
\vskip .1in
\noindent \bf Proposition. \it ${-1\over g} Q_0^2 =L_0^\psi -{\rm
dim}\g/48$.
\vskip .1in
\noindent \bf Proof. \rm Note that ${-1\over g} Q_0^2 -L_0^\psi$
commutes with $\psi_n^b$. Recall that $Q_0^*=-Q_0$ and that
$$Q_0=-{1\over 12}\sum  f_{abc} \psi_a(m)\psi_b(n)\psi_c({-m-n}).$$
The
corresponding finite--dimensional supercharge operator is given by
$G_0=-{1\over 12} \sum f_{abc} \psi_a(0)\psi_b(0)\psi_c(0)$. Thus
$Q_0\xi=G_0\xi$ for any $\xi\in H(0)$. Hence
$Q_0^2\Omega=Q_0G_0\Omega=G_0^2\Omega=-{g\, {\rm dim}\,\g\over 24}
\Omega$. Thus $(-{1\over g} Q_0^2 -L_0^\psi)\Omega= ({1\over
16}-{1\over 48}){\rm dim}\, \g \, \Omega={{\rm dim}\, \g\over 48}\cdot
\Omega$. Since $\Omega$ is cyclic and ${-1\over g} Q_0^2 -L_0^\psi$
central, the result follows.
\vskip .1in
\bf \noindent Remark. \rm If we have a Sugawara operator $L_0$ with
corresponding central charge $c$, we define ${\cal L}_0=L_0 -{c\over
24}$. If we redefine ${\cal L}_n=L_n$ for $n\ne 0$, then
$[{\cal L}_m,{\cal L}_n]=(m-n){\cal L}_{m+n} +{c\over 12} m^3
\delta_{m+n,0} $. We call ${\cal L}_0$ the normalised Sugawara
operator. It is the operator needed to make various characters modular
invariant. 

\vskip .1in
\noindent \bf B. Computation in ${\rm End}(H)\otimes {\rm Cliff}\,
L\g$. \rm Let $Q=\sum T_a(m) \psi_a({-m}) +{1\over 3} S_a(m)
\psi_a({-m})$. Set $X_a(m)=T_a(m)+S_a(m)$. Then we already have proved the
supersymmetry relations $\{Q,\psi_a(n)\}=X_a(n)$ and
$[Q,X_a(n)]=n(\ell+g) \psi_a(n)$.
\vskip .1in
\noindent \bf Proposition. \it $-(\ell +g)^{-1} Q^2 ={\cal L}_0^\g
+{\cal L}_0^\psi$, where ${\cal L}_0=L_0-c/24$ with $c_\g=\ell {\rm
dim}\, \g/ (\ell +g)$ and $c_\psi={\rm dim}\,\g/16$. 
\vskip .1in
\noindent \bf Proof. \rm If we apply $\{Q,\cdot\}$ to the
supersymmetry relations, we get $[Q^2,\psi_a(n)]=n(\ell+g) \psi_a(n)$ and 
$[Q^2,X_a(n)]=n(\ell+g)X_a(n)$. Thus 
$-(\ell+g)^{-1}Q^2-L_0^\g -L_0^\psi$ commutes with the $\psi_a(n)$'s
and $X_b(m)$'s. We claim that these operators act cyclically with
cyclic vector $\xi\otimes \Omega$, where $\xi\in H(0)$ is a highest
weight vector. Since the operator $S_a(n)$ are combinations of biinears
in $\psi_b(m)$'s, the cyclic module generated by $\xi\otimes \Omega$
must also be invariant under the $S_a(n)$'s. But then it must also be
invariant under $T_a(n)=X_a(n)-S_a(n)$. Since $\xi\otimes \Omega$ is
obviously cyclic for the comuuting actions of $\psi_a(n)$'s and
$T_b(m)$'s, our claim follows. Thus it will suffice to show that
$-(\ell+g)^{-1} Q^2 -{\cal L}_0^\g - {\cal L}_0^\psi$ annihilates
$\xi\otimes \Omega$. But if $G$ is the finite--dimensional version of
$Q$, we have as before that 
$$-Q^2(\xi\otimes
\Omega)=-G^2(\xi\otimes\Omega)= ({g\, {\rm dim}\, \g\over 24} +{1\over
2}(\|\lambda+\rho\|^2 -\|\rho\|^2)(\xi\otimes \Omega).$$
But $({\cal L}_0^\g+{\cal L}_0^\psi)(\xi\otimes \Omega)$
$$=({{\rm dim}\,
\g\over  16} + 
{\|\lambda+\rho\|^2 -\|\rho\|^2\over 2(\ell+g)} - 
{\ell\, {\rm dim}\, \g\over 24 (\ell+g)} - {{\rm dim}\,\g\over 48}
)(\xi\otimes\Omega)
=({g\, {\rm dim}\, \g\over 24} +{1\over
2}(\|\lambda+\rho\|^2 -\|\rho\|^2))(\xi\otimes \Omega),$$
as required.
\vskip .1in
\noindent \bf C. Coset construction of $Q^2$. \rm By the coset construction 
$Q=Q_\g-Q_\h$. Thus $Q_\g=Q+Q_\h$, where $Q$ and $Q_\h$ anticommute,
i.e.~$\{Q,Q_\h\}= 0$. Hence 
$$\{Q_\g,Q_\g\}=\{Q+Q_\h,Q+Q_\h\} = \{Q,Q\} + \{Q_\h,Q_\h\}.$$
Thus $Q_\g^2 = Q^2 +Q_\h^2$, so that $Q^2=Q_\g^2-Q_\h^2$. Hence we
have:
\vskip .1in
\noindent \bf Theorem. \it $-{1\over \ell+g} Q^2 = {\cal L}_0^\g +{\cal
L}_0^{\psi,\m} -{\cal L}_0^{\h}$ on $H\otimes {\cal F}_\m$. 
\vskip .1in
\noindent \bf Proof. \rm If we use the formula in B for $Q^\h$, we get
$-{1\over \ell+g} Q_\h^2= {\cal L}_0^\h +{\cal L}_0^{\psi,\h}$. By
definition ${\cal L}_0^{\psi,\g} ={\cal L}_0^{\psi,\m} +{\cal
L}_0^{\psi,\h}$, so the result follows. (Note that when the tensor
product of $H$ and ${\cal F}_\m$ is
restricted to $\widehat{\h}$, it splits up as a direct sum of positive
energy irreducible representations $H_j$. The formula for $Q_\h^2$ is
valid on each tensor product $H_j\otimes {\cal F}_\h$.)
\vskip .1in
\noindent \bf Corollary. \it Let ${{\cal F}}_\m$ be the
irreducible representation of ${\rm Cliff}\, L\m$. Then on 
$H\otimes {{\cal F}}_\m$, we have
$$-{1\over \ell+g} Q^2= {\|\lambda+\rho\|^2\over 2(\ell+g)} -(L_0^\h
-d^\g -d^\psi),$$
where $d^\g$ and $d^\psi$ are the energy operators on $H$ and
${{\cal F}}_\m$.
\vskip .1in
\noindent \bf Proof. \rm Using Freudenthal's strange formula, we get
$$\eqalign{-{1\over \ell +g} Q^2 & = L_0^\g + L_0^{\psi,\m} - L_0^\h -
(c_\g +c_{\psi,\m} -c_{\h})/24 \cr
&=(L_0^\g-d^\g) +(L_0^{\psi,\m} -d^{\psi}) -(L_0^\h-d^\g-d^\psi) -
(c_\g +c_{\psi,\m} -c_{\h})/24 \cr
&={\|\lambda+\rho\|^2-\|\rho\|^2\over 2(\ell+g)} +{{\rm dim}\, \m\over 16} 
-{\ell \, {\rm dim}\, \g\over 24(\ell+g)} -{{\rm dim}\, \m\over 24}
+{{\rm dim}\, \h\over 24} -(L_0^\h-d^\g-d^\psi)\cr
& ={\|\lambda+\rho\|^2\over 2(\ell+g)} -{g\, {\rm dim}\, \g\over
24(\ell+g)} +{{\rm dim}\, \m\over 16} 
-{\ell \, {\rm dim}\, \g\over 24(\ell+g)} -{{\rm dim}\, \m\over 24}
+{{\rm dim}\, \h\over 24} -(L_0^\h-d^\g-d^\psi)\cr
&={\|\lambda+\rho\|^2\over 2(\ell+g)} -(L_0^\h
-d^\g -d^\psi),\cr}$$
as required.
\vskip .1in

\noindent \bf 15. KAC'S CHARACTER AND DENOMINATOR FORMULAS. \rm
\vskip .1in

\noindent \bf Proposition. \it (1) Let $\overline{\mu}
+\overline{\nu}$ be a weight of $\widehat\h$ appearing in
$L(\overline{\lambda}) \otimes {\cal F}_\m$. Then
$(\overline{\lambda} +\overline{\rho},\overline{\lambda}
+\overline{rho}) \ge (\overline{\mu}
+\overline{\rho},\overline{\mu}+\overline{\rho})$. 

\noindent (2) Equality occurs in (1) iff $\overline{\mu}
=\sigma\overline{\lambda}$ and $\overline{\nu} =\sigma\overline {\rho}$ for
some $\sigma\in \widehat W$. In this case
$\overline{\mu}+\overline{\nu}=\sigma(\overline{\lambda} +\overline{\rho})$.

\noindent (3) $\sigma\mapsto \sigma(\overline{\lambda} +\overline{\rho})$ is
a bijection from $\widehat W$ onto the solutions of (2).

\noindent (4) $\sigma(\overline{\lambda} +\overline{\rho})$
appears in $L(\overline{\lambda}) \otimes {\cal F}_\m$ with
multiplicity one and corresponds to a tensor $\xi\otimes \eta$, where
$\xi$ has weight $\overline{\mu}$ and $\eta$ has weight
$\overline{\nu}$.

\vskip .1in
\noindent \bf Proof. \rm (1) Let $\overline{\mu}
+\overline{\nu}$ be any weight in the tensor product with
corresponding vector $\xi$. The corollary in
the last 
section gives the following formula for the square of the supercharge
operator 
$$-(\ell+g)^{-1} Q^2= {\|\lambda+\rho\|^2\over 2(\ell+g)} -(L_0^\h
-d^\g -d^\psi).$$
Since $Q*=-Q$, the operator $-Q^2$ is positive.  Thus 
$$\eqalign{(2\ell +2g)^{-1}(\|\lambda +\rho\|^2) \|\xi\|^2 
-((L_0^\h -d^\psi-d^\g)\xi,\xi)& \ge ((L_0^\h 
-d^\g -d^\psi)\xi,\xi) \cr 
&\ge ((\Omega^\h -d^\g-d^\psi)\xi,\xi)\cr
& =
(2\ell+2g)^{-1}|\overline{\mu}
+\overline{\nu}|^2 \|\xi\|^2,\cr}$$
where $\Omega^\h$ is the zero--mode Casimir contribution to $L_0^\h$. 
Clearly $\widehat\xi$ is a highest weight vector for $\widehat{\h}$ if
and only if $(L_0^\h \xi,\xi) = (\Omega^\h\xi,\xi)$ if and only if
$L_0^\h \xi = \Omega^\h\xi$. Hence
$(\overline{\lambda}+\overline{\rho},\overline{\lambda}+\overline{\rho})\ge
\ge (\overline{\mu}+\overline{\nu},\overline{\mu}+\overline{\nu})$;
equality is possible only if $\overline{\mu}+\overline{\nu}$ is a highest
weight of $\widehat{\h}$.

\noindent (2) Suppose that $|\overline{\lambda} +\overline{\rho}|^2 \ge
|\overline{\mu} + \overline{\nu}|^2$ with
$\overline{\lambda} \ge \overline{\mu}$ and $\overline {\rho} \ge \overline
{\nu}$. Take $\tau\in \widehat W$ such that $\tau(\overline{\mu}
+\overline{\nu})\ge 0$. Since $\tau \overline{\mu}$ is a weight of
$L(\overline{\lambda})$ and $\tau \overline{\nu}$ is a weight of ${\cal
F}$, we have $\overline{\lambda} -\tau\overline{\mu}\ge 0$ and 
$\overline{\rho} - \tau\overline{\nu} \ge 0$. But then
$$0=(\overline {\lambda} +\overline {\rho} +\tau \overline{\mu} +\tau
\overline{\nu}, \overline {\lambda} - \tau \overline{\mu}
+\overline{\rho}-\tau 
\overline{\nu}) \ge (\overline{\rho}, \overline{\lambda} -\tau\overline{\mu}) 
+(\overline{ \rho}, \overline {\rho}-\tau \overline {\nu}) \ge 0.$$
Hence $(\overline{\rho}, \overline{\lambda} -\tau\overline{\mu})=0$ and
$(\overline{ \rho}, \overline {\rho}-\tau \overline {\nu})=0$, so that
$\overline{\lambda} =\tau\overline{\mu}$ and $\overline {\rho}=\tau
\overline {\nu}$. 

\noindent (3) Let $\overline{\lambda}_i$ be the fundamental weights of
$\widehat\g$.  Since $\overline{\rho} =\sum_{i\ge 0} \overline{\lambda}_i$,
the equality $\tau \overline {\rho} =\overline {\rho}$ and the ineqalities
$\overline{\lambda}_i \ge \tau \overline{\lambda}_i$ force
$\tau\overline{\lambda}_i =\overline{\lambda}_i$ for $i\ge 0$. Hence
$\tau=1$.

\noindent (4) Suppose that $\overline {\mu}_1+\overline {\nu}_1 =\overline
{\mu} +\overline{ \nu}$
with $\overline{ \lambda} \ge \overline {\mu}_1$ and $\overline{\rho}
\ge \overline 
{\nu}_1$. Since $|\overline{\mu}_1 +\overline{\nu}_1|^2 = |\overline
{\mu}+\overline{\nu}|^2 =|\overline{\lambda}+\overline{\rho}|^2$, (2) and
(3)imply that $\overline {\mu}_1=\sigma \overline{ \lambda}$ and
$\overline {\nu}_1=\sigma \overline{\rho}$ for some $\sigma\in \widehat W$.

But then $\gamma=\tau^{-1}\sigma$ fixes $\overline{\lambda}
+\overline{\rho}$. Since $\overline{\lambda}\ge \gamma \overline {\lambda}$ and
$\overline{\rho}\ge \gamma\overline{\rho}$, we get
$\gamma\overline{\rho}=\overline{\rho}$ so that $\gamma=1$.
\vskip .1in
\bf \noindent Lemma. \it The formal power series
$D=e^{\overline{\rho}}\prod_{\alpha\in \overline{\Phi}^+}
(1-e^{-\alpha})^{m_\alpha}$ satisfies $D\circ \sigma
=\varepsilon(\sigma) D$ for $\sigma\in \widehat{W}$, where
$\varepsilon:\widehat{W}\rightarrow \{\pm 1\}$ is the sign character
obtained by the sign character of $W$ and the projection
$\widehat{W}\rightarrow W$. 
\vskip .1in
\noindent \bf Proof. \rm It clearly suffices to show that 
$D\circ \sigma_i
=\varepsilon(\sigma_i) D$ for each simple reflection $\sigma_i$. We
know that $\sigma_i\overline{\rho} =\overline{\rho}
-\overline{\alpha}_i$. Moreover $\sigma_i$ permutes
$\overline{\Phi}^+\backslash\{\overline{\alpha}_i\}$ and satisfies
$\sigma_i\overline{\alpha}_i = -{\alpha}_i$. Hence
$D\circ \sigma_i =-D$, as required.

\vskip .1in
\bf \noindent Theorem (Kac Character Formula). \it 
$${\rm ch} L(\overline\lambda) = \sum_{\sigma\in \widehat W}
 \varepsilon (\sigma)
e^{\sigma(\overline{\lambda}+\overline{\rho})-\overline{\rho}}/D,$$
where the denominator $D$ is given by
$$D=\prod_{n\ge 1} (1-q^n)^m\cdot \prod_{\alpha\in \Phi^+}
(1-e^{-\alpha}) \cdot \prod_{n\ge 1,\alpha\in \Phi}(1-e^\alpha q^n)
=\prod_{\overline{\alpha}\in\overline{\Phi}^+ } 
(1-e^{-\overline{\alpha}})^{{\rm mult} \, \alpha}.$$

\vskip .1in
\noindent \bf Proof. \rm Let ${\rm ch}_s W$ be the supercharacter of a
${\Bbb Z}_2$--graded module. Clearly 
$${\rm ch}_s
L(\overline{\lambda})\otimes {\cal F}_\m ={\rm ch}\, L(\overline{ \lambda})
\cdot {\rm ch}_s \, {\cal F}_\m.\eqno{(1)}$$
Moreover, because ${\cal F}_\m$ can be constructed as the tensor
product of an irreducible representation of ${\rm Cliff}(\m)$ on
$\Lambda^*\m_+$ and the irreducible representation of the non--zero
modes $\psi_i(n)$ on the exterior algebra with generators $v_{i,n}$,
we easily check that
$${\rm ch}_s\, {\cal F}_\m = e^{\overline{\rho}} \prod_{\alpha\in
\Phi^+} 
(1-e^{-\alpha}) \cdot \prod_{n\ge 1, \alpha\in \Phi} (1-e^\alpha
q^n).\eqno{(2)}$$
On the other hand ${\rm ch}_s \, L(\overline{\lambda})\otimes {\cal F}_\m$
can be computed by decomposing the tensor product as a direct sum of
$\widehat\h$--modules, according to the previous proposition. We find
$${\rm ch}_s \, L(\overline{\lambda})\otimes {\cal F}_\m
=\sum_{\sigma\in\widehat W}  \varepsilon^\prime(\sigma) {\rm ch}
V(\sigma(\overline{\lambda}+\overline{\rho})),$$
where $V(\overline{\mu})$ is the irreducible representation of
$\widehat \h$ with 
highest weight $\overline{\mu}$
and $\varepsilon^\prime(\sigma)=\pm 1$ according to whether the weight
$\sigma\overline{\rho}$ appears in the even or odd part of ${\cal
F}_\m$. But ${\rm ch}\, V(\mu)= e^\mu \cdot \prod (1-q^n)^m$, where 
$m={\rm dim}\,\h$ is the rank of $\g$. Hence
$${\rm ch}_s \, L(\overline{\lambda})\otimes {\cal F}_\m
=\sum_{\sigma\in\widehat W} \varepsilon^\prime(\sigma)
e^{\sigma(\overline{\lambda}+\overline{\rho})} \cdot \prod
(1-q^n)^m.\eqno{(3)}$$ 
The character formula follows by combining (1), (2) and (3).
It still remains to show that
$\varepsilon^\prime(\sigma)=\varepsilon(\sigma)$. 
Specialising to $\overline{\lambda}=0$, we find
$$e^{\overline{\rho}} \prod_{\overline{\alpha}>0} (
1-e^{-\overline{\alpha}}) =\sum_{\sigma\in\widehat{W}}
\varepsilon^\prime(\sigma) e^{\sigma\overline{\rho}}.$$
We know that all the exponents $\sigma\overline{\rho}$ are distinct
and we know that the left hand side is antisymmetric. Since the
coefficient of $e^{\overline{\rho}}$ is $1$ and $\varepsilon^\prime$
must coincide with $\varepsilon$, as required.

\vskip .1in
\noindent \bf Theorem (Macdonald's identity/Kac denominator
formula). \it  
$$\sum_{\sigma \in \widehat W} \varepsilon(\sigma)
e^{\sigma(\overline{\rho}) -\overline{\rho}}=
=\prod_{\overline{\alpha}\in \overline{\Phi} }
(1-e^{-\overline{\alpha}})^{{\rm mult} \, \alpha}.$$ 
\vskip .1in
\noindent \bf Proof. \rm Since ${\rm ch} \, L(0) =1$, the result
follows immediately from the character formula with $\lambda=0$.

\vskip .1in
\bf \noindent Comments on the denominator formula. \rm Macdonald's
identities for the classical groups were first proved by Dyson 
independently of root systems 
and then for general root systems by Macdonald, without
representations. In the case of $\sl_2$, we see that the action of
$\widehat{W}= {\Bbb Z}\rtimes {\Bbb Z}_2$ is given by
$\sigma(j,k,\ell) =(-j,k,\ell)$ for $\sigma\in {\Bbb Z}_2$ and
$T_n(j,k\ell)= (j+2n\ell,k+jn+\ell m^2,\ell)$ for $n\in {\Bbb Z}$. 
From the denominator formula, we retrieve Jacobi's
celebrated triple product identity:
$$\sum_{k\in {\Bbb Z}} (-1)^k q^{(k-1)k/2} t^k = \prod(1-q^{m-1} t)(1-q^m
t^{-1}) (1-q^m).$$
We can get another a formula for the characters in terms of theta
functions. Let 
$$\Theta_{n,m}(q,z)=\sum_{k\in {n\over 2m} +{\Bbb Z}} q^{mk^2} e^{2\pi
i km z}.$$
Then if $0\le j\le \ell/2$ is a half--integer, the normalised character of the
irreducible positive energy representation of level $\ell$ with spin
$j$ is given by
$${\rm ch}\, L(\ell,j)=q^{-c/24}{\rm Tr}(q^{L_0} z^H)=
{\Theta_{2j+1,\ell+2}(q,z) - 
\Theta_{-2j-1,\ell+2}\over
\Theta_{1,2}(q,z) - \Theta_{-1,2}(q,z)},$$
where $s=(j+{1\over 2})^2/(\ell+2)-{1\over 8}$.
A similar formula holds for any simple algebra: the sum over the
affine Weyl group is first performed as a sum over the coroot lattice
followed by an antisymmetrisation over the finite Weyl group. Each sum
over the Weyl lattice results in a theta function and we thus get a
product formula for an alternating sum of theta functions. 

The formula
for the level one vacuum character of affine $SU(2)$ can be further
simplified. (This 
simplification is related to the boson--fermion correspondence and the
Frenkel--Kac--Segal construction; an analogue holds for each simply laced
simple Lie algebra of types A, D or E.)
\vskip .1in
\noindent \bf Product formula. \it $\Theta_{n,m}(q,z)
\Theta_{n^\prime,m^\prime}(q,z) =\sum_{j\in {\Bbb Z}/(m+m^\prime){\Bbb Z}}
F_j(q) \Theta_{n+n^\prime +2mj, m+m^\prime} (q,z)$, where
$F_j(q)=\sum_{k\in {\Bbb Z}+x} q^{mm^\prime (m+m^\prime)k^2}$ and
$x=(m^\prime n -m n^\prime + 2j m m^\prime) /2mm^\prime (m+m^\prime)$.
\vskip .1in
\noindent \bf Proof. \rm  We have
$$\Theta_{n,m}\Theta_{n^\prime,m^\prime} = \sum_{k,k^\prime} q^{mk^2
m^\prime k^{\prime 2}}{\bf e}_{mk +m^\prime k^\prime}(z),$$
where $k\in {n\over 2m} +{\Bbb Z}$ and $k^\prime \in {n^\prime \over
2m^\prime} +{\Bbb Z}$. Set $k=j+{n\over 2m}$ and $k^\prime =j^\prime+
{n^\prime \over 2m^\prime}$. Define
$s=(k-k^\prime)/(m+m^\prime)$ and $s^\prime(mk+m^\prime k^\prime)
/(m+m^\prime)$.  Write $j-j^\prime=(m+m^\prime)a +b$ with $a\in {\Bbb
Z}$ and $0\le b <m+m^\prime$. Then
$$s\in {nm^\prime - n^\prime m + 2 mm^\prime b\over
2mm^\prime(m+m^\prime)} + {\Bbb Z}, \qquad s^\prime\in {n+ n^\prime +
2mj\over 2(m+m^\prime)} +{\Bbb Z}.$$
This gives a bijection between pairs $(k,k^\prime)$ and triples
$(s,s^\prime,b)$. Since
$mk^2+m^\prime k^{\prime 2}=mm^\prime(m+m^\prime)s^2
+(m+m^\prime)s^{\prime 2}$, we get
$$\Theta_{n,m}\Theta_{n^\prime,m^\prime}=\sum_b(\sum_s
q^{mm^\prime(m+m^\prime)s^2}) (\sum_{s^{\prime}}
q^{(m+m^\prime)s^{\prime 2}}{\bf e}_{m+m^\prime}(s^\prime z),$$
as required.
\vskip .1in 
\bf \noindent Corollary. \it ${\rm ch}\,
L(1,0)=\Theta_{0,1}/\eta(q)$,  where $\eta(q) =q^{1/24}\prod_{n\ge 1}
(1-q^n)$ is Dedekind's eta function. 
\vskip .1in
\noindent \bf Proof. \rm By Jacobi's triple product identity (the
Weyl--Kac denominator formula for $SU(2)$), we have
$$\sum_{k\in {\Bbb Z}} q^{k^2} t^k = \prod(1+q^{2m-1} t)(1+ q^{2m-1}
t^{-1}) (1-q^{2m}).$$
If we specialise $(q,t)$ to $(q^{3/2},-q^{-1/2})$, we obtain Euler's
pentagonal identity:
$$\varphi(q)\equiv \prod_{m\ge 1} (1-q^m) =\sum_{m\in {\Bbb Z}} (-1)^k
q^{3k^2-k}/2=\sum_{m\in {\Bbb Z}} (-1)^k
q^{3k^2+k}/2.$$
To prove the corollary, we must show that
$$\Theta_{0,1}(\Theta_{1,2}-\Theta_{-1,2})=q^{1\over 24} \varphi(q)
(\Theta_{1,3} -\Theta_{-1,3}).$$
By the product formula, the left hand
side is
$$(\Theta_{1,3}-\Theta_{-1,3})(\sum_{k\in {-1\over 12} +{\Bbb Z}}
q^{k^2} -\sum_{k^\prime\in {5\over 12} +{\Bbb Z}} q^{k^{\prime 2}})=
(\Theta_{1,3}-\Theta_{-1,3})q^{1\over 24}
\sum_{k\in {\Bbb Z}}(-1)^k q^{(3k^2+k)/2} 
=(\Theta_{1,3}-\Theta_{-1,3})q^{1\over 24}\varphi(q),$$
as required. 
\vfill\eject
\def\d{{\hbox{\goth d}}}
\bf \centerline{CHAPTER~IV.~REPRESENTATIONS OF THE VIRASORO
ALGEBRA}
\vskip .1in
\it In this chapter we discuss positive energy unitary
representations of the Virasoro algebra. We show that they are
classified by a lowest energy $h$ and a central charge $c$. Of
particular interest are the representations with $0<c<1$. A series of
these can be constructed by the coset construction of
Goddard--Kent--Olive. For these representations $c=1-6/m(m+1)$ with
$m\ge 3$. The $h$ values are then $h_{p,q}=[(p(m+1)-qm)^2-1]/4m(m+1)$
where $1\le q\le p\le m-1$. The coset construction gives
representations of the Virasoro algebra on multiplicity spaces and the
characters of these multiplicity spaces can be computed explicitly in
terms of theta functions and Dedekind's eta function. Using this
information, we present the short proof of Kac's determinant formula
due to Kac and Wakimoto. This is then used to prove the easy part of
the unitarity criterion of Friedan--Qiu--Shenker: the only irreducible
unitary representations with $c=1-6/m(m+1)$ have $h=h_{p,q}$ as
above. We then use this unitarity criterion to give our own direct
proof that the multiplicity spaces are irreducible. (Such a
`multiplicity one' theorem seems to hold more generally for
W--algebras.) This gives a very
short proof of the Feigin--Fuchs character formula for these values of
$c$ and $h$. Our method uses unitarity properties rather than a
detailed knowledge of null vectors and Verma module
resolutions. (These resolutions are not easy; indeed the lengthy
derivation of Feigin and Fuchs needs to be supplemented 
with arguments of Astashkevich on the Jantzen filtration.) Finally we
prove the hard part of the unitarity criterion of
Friedan--Qiu--Shenker: the only values of $c\in (0,1)$ yielding 
unitary representations are $c=1-6/m(m+1)$. Our techniques
extend easily to treat the case $c=1$. 

\vskip .1in
\noindent \bf 1. POSITIVE ENERGY REPRESENTATIONS OF THE VIRASORO
ALGEBRA. \rm We shall be interested in {\bf projective}, {\bf
unitary}, {\bf positive energy} representations of the Witt 
algebra $\d={\rm Vect}\, S^1$. Recall that $\d$ is the
complexification of the real Lie algebra of (trigonometric) polynomial
vector fields $a(\theta) \, d/d\theta$ on $S^1$. It has basis $d_n=
ie^{in\theta} \, d/d\theta$. We can use Leibniz' rule to compute the
Lie brackets: $[d_m,d_n] =(m-n)d_{m+n}$. We set $d_n^*=d_{-n}$. This
extends to a conjugate--linear involution on $\d$. Thus
we are looking for inner 
product spaces $H$ (not complete!!) such that:
\vskip .1in
\noindent (1) \bf Projective: \rm $\d$ acts  projectively by operators
$\pi(A)$ ($A\in 
\d$), i.e. $A\mapsto \pi(A)$ is linear and
$[\pi(A),\pi(B)]-\pi([A,B]) $ lies in ${\Bbb C}I$ for $A,B\in \d$.

\noindent (2) \bf Unitary: \rm $\pi(A)^*=\pi(A^*)$. 

\noindent (3) \bf Positive energy: \rm $H$ admits an orthogonal
decomposition $H=\bigoplus_{k\ge 0} H(k)$ such that some (necessarily
unique) representative $D$ for $\pi(d_0)$ acts on
$H(k)$ as multiplication by $k$, $H(0)\ne 0$ and ${\rm dim}\, H(k)<\infty$
\vskip .1in
The subspaces $H(k)$ are called
the energy subspaces with energy $k$; the operator $D$ has many names,
including the energy operator or hamiltonian operator.
Since the representation is projective, $[\pi(A),\pi(B)]-\pi([A,B])=b(A,B) I$
where $b(A,B)\in {\Bbb C}$. We call $b$ a 2--cocycle --- in fancy
language it gives a class in $H^2(\d,{\Bbb C})$. The definition
immediately implies the antisymmetry condition
 $$b(A,B)=-b(B,A)$$
because Lie brackets are antisymmetric; and the Jacobi identity
immediately implies that
$$b([A,B],C) + b([B,C],A) + b([C,A],B)=0$$
for all $A,B,C\in \d$. On the other hand we are free to adjust the
operators $\pi(A)$ by adding on scalars. Thus to preserve linearity,
we change $\pi(A)$ to $\pi(A)+f(A)I$ where $f:\d\rightarrow {\Bbb C}$
is linear. This changes $b(A,B)$ to $b(A,B) -f([A,B])$. We shall now
make such adjustments so that $b$ has a canonical form. 
We start by choosing the canonical representative $D$ for $\pi(d_0)$
as above. By uniqueness, we must have
$[L_m,L_n]=(m-n)L_{n+m} +\lambda(m,n)I$. Here $\lambda(n,m)$ is a
2--cocycle. As we now show, by appropriate adjustment of the $L_n$'s
by scalars, that $\lambda$ can always be normalised so that
$\lambda(m,n)={c\over 12} (m^3-m) \delta_{m+n,0}$, where $c$ is called the
{\it central charge}. The corresponding central extension of the Witt
algebra is usually called the Virasoro algebra. 
\vskip .1in
\vskip .1in
\noindent \bf Virasoro cocycle lemma. \it Representatives $L_n$ of
$\pi(d_n)$ may be chosen uniquely so that $[D,L_n]=-nL_n$ for all
$n$. In this case
$[L_m,L_n]=(m-n)L_{m+n} +(am^3 +bm)\delta_{m+n,0}I$. If we choose
$L_0$ so that $[L_1,L_{-1}]=L_0$, then $a+b=0$ and
$$[L_m,L_n]=(m-n)L_{m+n} +{c\over 12} (m^3-m)
\delta_{m+n,0}.\eqno{(1)}$$
If instead we take ${\cal L}_0=L_0-c/24$ and ${\cal L}_n=L_n$ ($n\ne 0$), then
$$[{{\cal L}}_m,{{\cal L}}_n]=(m-n){{\cal L}}_{m+n} +{c\over 12} m^3
\delta_{m+n,0}.\eqno{(2)}$$
\vskip .1in

\noindent \bf Proof. \rm Note that $-n^{-1} [L_0,L_n]$ is independent
of adding scalars onto $L_0$ or $L_n$, so we may always choose $L_n$
so that $[L_0,L_n]=-nL_n$. Thus $[D,L_n]=-nL_n$, so that $L_n$ lowers
energy by $n$, i.e.~takes $H(k)$ into $H(k-n)$. But then $[L_m,L_n]$
and $L_{m+n}$ lower energy by $n+m$. Since $[L_m,L_n]-(m-n)L_{m+n}$ is
a scalar, it must be $0$ if $n+m\ne 0$. Thus 
$$[L_m,L_n]=(m-n) +A(m) \delta_{m+n,0}.$$
Clearly $A(m)=-A(-m)$ and $A(0)=0$. Writing out the Jacobi
identity for $L_k$, $L_n$ and $L_m$ with $k+n+m=0$, we get
$$(n-m)A(k) + (m-k)A(n) +(k-n)A(m)=0.$$
Setting $k=1$ and $m=-n-1$, we get
$$(n-1)A(n+1) =(n+2)A(n)-(2n+1) A(1).$$
This recurrence relation allows $A(n)$ to be determined from $A(1)$
and $A(2)$. Since $A(n)=n$ and $A(n)=n^3$ give solutions, we see that
$A(m)=am^3+bm$ for some constants $a$ and $b$. Clearly we are free to
choose $L_0=[L_1,L_{-1}]$ (since we have made no adjustment to $L_0$
so far). But then $A(1)=0$ and hence $a+b=0$. This gives (1) and (2)
follows by an easy manipulation. 
\vskip .1in

\bf \noindent Complete reducibility theorem. \it (a) If $H$ is a positive
energy unitary representation of $\d$, then any non--zero vector in
the lowest energy subspace $H(0)$ generates an irreducible submodule.

\noindent (b) Any positive energy representation is an orthogonal
direct sum of irreducible positive energy representations.
\vskip .1in
\noindent \bf Proof. \rm (a) Take $v\ne 0$ in $H(0)$ and let $K$ be
the $\d$--invariant subspace it generates. Clearly since $L_n v=0$ for
$n>0$ and $L_0v=hv$, we see that $K$ is spanned by all products $Rv$
with $R=L_{-n_k} \cdots L_{-n_1}$ with $n_k\ge \cdots \ge n_1\ge 1$. But
then $K(0)={\Bbb C}v$. We claim that $K$ is irreducible. If not let
$K^\prime$ be a submodule and let $P$ be the orthogonal projection
onto $K^\prime$. By unitarity $P$ commutes with $\d$
and hence $D$. Thus $P$ leaves $K(0)={\Bbb C}v$ invariant, so that
$Pv=0$ or $v$. But $P(Rv)=RPv$. Hence $K^\prime=(0)$ or $K$, so that
$K$ is irreducible.

\noindent (b) Take the irreducible module generated by
a vector of lowest energy $H_1$. Now repeat
this process for $H_1^\perp$, to get $H_2$, $H_3$, etc. The positive
energy assumption shows that $H=\bigoplus Hi$.
\vskip .1in

\noindent \bf Uniqueness Theorem. \it If $H$ and $H^\prime$ are
irreducible positive 
energy representations of $\d$ with central charge $c$ and
$h=h^\prime$, then $H$ and $H^\prime$ are
unitarily equivalent as representations of $\d$. 
\vskip .1in

\noindent \bf Proof. \rm Any monomial $A$ in operators from $\d$ is a sum of
monomials $RDL$ with $R$ a monomial in energy raising operators, $D$ a
monomial in constant energy operators and $L$ a monomial in energy
lowering operators. Observe that if $v,w\in H(0)$, the inner products 
$(A_1v,A_2w)$ are uniquely determined by $v,w$ and the monomials $A_i$:
for $A_2^*A_1$ is a sum of terms $RDL$ and $(RDLv,w)=(DLv,R^*w)$ with
$R^*$ an energy lowering operator. Hence, if $H^\prime$ is another
irreducible positive energy representation with $h=h^\prime$, with
$H(0)={\Bbb C}v$ and $H^\prime(0)={\Bbb C}v^\prime$ for unit vectors
$v,v^\prime$, then $U(Av)=Av^\prime$ defines a unitary map of $H$
onto $H^\prime$ intertwining $\widehat{\d}$. 
\vskip .1in
\noindent \bf 2. THE GODDARD--KENT--OLIVE CONSTRUCTION. \rm We have
seen a variety of construction of positive energy representations of
the Virasoro algebra in the chapter on affine Kac--Moody algebras: the
Segal--Sugawara construction; the Fubini--Veneziano construction using
bosons; and the Fubini--Veneziano construction using fermions. 
We now describe a further ``coset'' construction of
Goddard--Kent--Olive on multiplicity spaces. 

\vskip .1in
\noindent \bf Lemma. \it Let $\h$ be a Lie algebra acting unitarily on
the inner product space $H$. Suppose that $H$ is a direct sum of
irreducible submodules and that there are only finitely many
isomorphism types of irreducible summands $H_i$. Let $K_i={\rm
Hom}_\h(H_i,H)$. Then $K_i$ is naturally an inner product space and
the map $\bigoplus K_i\otimes H_i\rightarrow H$, $\sum \xi_i\otimes
\eta_i\mapsto \sum \xi_i\eta_i$ is a unitary map of $\h$--modules.
The operators $A$ on $H$ which commute with $\h$ act naturally on each
$K_i$ by $A_i$. This action is a *--homomorphism. Under the
unitary isomorphism above, $A$ corresponds to $\oplus 
A_i\otimes I$ and $X\in \h$ to $\oplus I\otimes \pi_i(X)$.
\vskip .1in
\noindent \bf Proof. \rm If $S,T\in {\rm Hom}(H_i,H)$, then $T^*S\in
K_i={\rm End}_\h(H_i) ={\Bbb C}$ by Schur's lemma. Thus the canonical
inner product on $K_i$ is defined by $(S,T)=T^*S$. It is then easy to
check the assertions about the map $\oplus K_i\otimes H_i\rightarrow
H$, since by assumption this map is surjective. The action of $A$ on
$K_i={\rm Hom}_\h(H_i,H)$ is defined by $A\xi$. 
\vskip .1in
\noindent \bf Proposition (coset construction). \it Let $\g$ be a Lie
algebra with subalgebra $\h$. Let $\d$ be a Lie algebra of derivations
action on $\g$ by $D,X\mapsto [D,X]$ ($D\in \d, X\in \g$) such that
$\d$ leaves $\h$ invariant. Suppose that $\g$ acts irreducible on the
inner product space $H$ and that $(H,\h)$ satisfy the hypotheses of
the previous lemma. Suppose in addition that $H$ and the $H_i$'s admit
projective unitary actions of $\d$ compatible with the action of $\g$
and $\h$. If $D\in \d$ acts by $\pi(D)$ on $H$ and $\pi_i(D)$ on
$H_i$, then $\pi(D)=\sum I\otimes \pi_i(D) + \sigma_i(D)\otimes I$,
where $\sigma_i$ is a projective unitary representation of $\d$ on
$K_i$; the cocycle 
of $\d$ on $K_i$ is the difference of the cocycles of $\d$ on $H$ and
on $H_i$.
\vskip .1in
\noindent \bf Proof. \rm Let $\sigma(D)=\pi(D)-\sum I\otimes
\pi_i(D)$. By construction this operator commutes with $\h$. Therefore
by the previous lemma $\sigma(D)=\sum \sigma_i(D)\otimes I$. 
Now suppose that $c(D_1,D_2)I=[\pi(D_1),\pi(D_2)]-\pi([D_1,D_2])$ 
and $c_i(D_1,D_2)I=[\pi_i(D_1),\pi_i(D_2)]-\pi_i([D_1,D_2])$ for
$D_1,D_2\in \d$. 
If $D_1,D_2\in \d$, then
$$\eqalign{c(D_1,D_2)I&= \sum [\pi(D_1),\pi(D_2)]-\pi([D_1,D_2])\cr
& =\sum
([\sigma_i(D_1),\sigma_i(D_2)]-\sigma_i([D_1,D_2]) )\otimes I
+\sum I \otimes ([\pi_i(D_1),\pi_i(D_2)]-\pi_i([D_1,D_2]))\cr
&=\sum
([\sigma_i(D_1),\sigma_i(D_2)]-\sigma_i([D_1,D_2]) )\otimes I
+\sum I\otimes c_i(D_1,D_2).\cr}$$
Looking at this equation on $K_i\otimes H_i$, we get
$$[\sigma_i(D_1),\sigma_i(D_2)]-\sigma_i([D_1,D_2]) )=
c(D_1,D_2) -c_i(D_1,D_2),$$ 
as required. 
\vskip .1in
We shall apply the coset construction in the following setting. Let
$H_0$ be the vacuum representation of $\widehat{\sl(2)}$ at level one
and let $H_{j,\ell}$ be any irreducible representation of $\widehat{\sl(2)}$
at level $\ell$ with lowest energy space of spin $j\in {1\over 2}{\Bbb
Z}$. Thus $H_0\otimes H_{j,\ell}$ gives a positive energy 
representation of $\widehat{\sl(2)}$ of level $\ell+1$. Thus we may
write $H_0\otimes H_{j,\ell} = \bigoplus M_k \otimes H_{k,\ell+1}$
with the $M_k$ multiplicity spaces. In this case $\g={\cal
L}\sl_2\oplus {\cal L}\sl_2$ and $\h={\cal L}\sl_2$, embedded
diagonally via $X\mapsto X\otimes I +I\otimes X$. The Witt algebra
acts on $\g$ and $\h$ and is implemented in both cases by the Sugawara
constructions. For $\g$, it has central charge $3\ell/(\ell+2) +1$
while for $\h$ it has central charge $3(\ell+1)/(\ell+3)$. Let
$m=\ell+2$. Subtracting the central charges, we see that there are
canonical projective representations of the Virasoro algebra on the
multiplicity spaces $M_k$ with central charge 
$$c=1 -3[\ell (\ell+3)
-(\ell+1)(\ell+2)]/(\ell+2)(\ell+3)= 1-6/m(m+1).$$

\vskip .1in

\noindent \bf 3. CHARACTER OF THE MULTIPLICITY SPACE. \rm We recall the
formula for the characters of the positive energy representations of
$LSU(2)$ (in normalised form). Let 
$$\Theta_{n,m}(q,\zeta)=\sum_{k\in {n\over 2m} +{\Bbb Z}} q^{mk^2}
\zeta^{2mk}.$$ 
Then if $0\le j\le \ell/2$ is a half--integer, the character of the
irreducible positive energy representation of level $\ell$ with spin
$j$ is given by
$${\rm ch}\, L(\ell,j)= {\Theta_{2j+1,\ell+2}(q,z) -
\Theta_{-2j-1,\ell+2}(q,z)\over
\Theta_{1,2}(q,z) - \Theta_{-1,2}(q,z)}.$$
Note that the character of a representation of $\widehat{\sl_2}$ is
${\rm Tr}(q^{{\cal L}_0}z)$ where ${\cal L}_0=L_0-c/24$ and $z$
corresponds to the element $\pmatrix{z & 0\cr 0 & z^{-1}\cr}$ in
$SU(2)$ or $SL(2)$. Similarly the character of a positive energy
representation 
\vskip .1in
\noindent \bf Product formula. \it $\Theta_{n,m}(q,z)
\Theta_{n^\prime,m^\prime}(q,z) =\sum_{j\in {\Bbb Z}/(m+m^\prime){\Bbb Z}}
F_j(q) \Theta_{n+n^\prime +2mj, m+m^\prime} (q,z)$, where
$F_j(q)=\sum_{k\in {\Bbb Z}+x} q^{mm^\prime (m+m^\prime)k^2}$ and
$x=(m^\prime n -m n^\prime + 2j m m^\prime) /2mm^\prime (m+m^\prime)$.
\vskip .1in
\noindent \bf Proof. \rm  We have
$$\Theta_{n,m}\Theta_{n^\prime,m^\prime} = \sum_{k,k^\prime} q^{mk^2
m^\prime k^{\prime 2}}z^{2(mk +m^\prime k^\prime)},$$
where $k\in {n\over 2m} +{\Bbb Z}$ and $k^\prime \in {n^\prime \over
2m^\prime} +{\Bbb Z}$. Set $k=j+{n\over 2m}$ and $k^\prime =j^\prime+
{n^\prime \over 2m^\prime}$. Define
$s=(k-k^\prime)/(m+m^\prime)$ and $s^\prime(mk+m^\prime k^\prime)
/(m+m^\prime)$.  Write $j-j^\prime=(m+m^\prime)a +b$ with $a\in {\Bbb
Z}$ and $0\le b <m+m^\prime$. Then
$$s\in {nm^\prime - n^\prime m + 2 mm^\prime b\over
2mm^\prime(m+m^\prime)} + {\Bbb Z}, \qquad s^\prime\in {n+ n^\prime +
2mj\over 2(m+m^\prime)} +{\Bbb Z}.$$
This gives a bijection between pairs $(k,k^\prime)$ and triples
$(s,s^\prime,b)$. Since
$mk^2+m^\prime k^{\prime 2}=mm^\prime(m+m^\prime)s^2
+(m+m^\prime)s^{\prime 2}$, we get
$$\Theta_{n,m}\Theta_{n^\prime,m^\prime}=\sum_b(\sum_s
q^{mm^\prime(m+m^\prime)s^2}) (\sum_{s^{\prime}}
q^{(m+m^\prime)s^{\prime 2}}z^{2(m+m^\prime)s^\prime},$$
as required.
\vskip .1in 
\bf \noindent Corollary. \it ${\rm ch}\,
L(1,0)=\Theta_{0,1}/\eta(q)$.
\vskip .1in
\noindent \bf Proof. \rm By Jacobi's triple product identity (the
Weyl--Kac denominator formula for $SU(2)$), we have
$$\sum_{k\in {\Bbb Z}} q^{k^2} t^k = \prod(1+q^{2m-1} t)(1+ q^{2m-1}
t^{-1}) (1-q^{2m}).$$
If we specialise $(q,t)$ to $(q^{3/2},-q^{-1/2})$, we obtain Euler's
pentagonal identity:
$$\varphi(q)\equiv \prod_{m\ge 1} (1-q^m) =\sum_{m\in {\Bbb Z}} (-1)^k
q^{3k^2-k}/2=\sum_{m\in {\Bbb Z}} (-1)^k
q^{3k^2+k}/2.$$
To prove the corollary, we must show that
$$\Theta_{0,1}(\Theta_{1,2}-\Theta_{-1,2})=q^{1\over 24} \varphi(q)
(\Theta_{1,3} -\Theta_{-1,3}).$$
By the product formula, the left hand
side is
$$(\Theta_{1,3}-\Theta_{-1,3})(\sum_{k\in {-1\over 12} +{\Bbb Z}}
q^{k^2} -\sum_{k^\prime\in {5\over 12} +{\Bbb Z}} q^{k^{\prime 2}})=
(\Theta_{1,3}-\Theta_{-1,3})q^{1\over 24}
\sum_{k\in {\Bbb Z}}(-1)^k q^{(3k^2+k)/2} 
=(\Theta_{1,3}-\Theta_{-1,3})\eta(q),$$
as required. 
\vskip .1in
\noindent \bf Theorem. \it ${\rm ch}\, L(0,1) \cdot {\rm ch}\, L(j,\ell)
=\sum_{0\le k\le (\ell +1)/2} \psi_k(q) \cdot {\rm
ch}\, L(k,\ell+1)$, where $k-j\in {\Bbb Z}$,
$$\psi_k(q)=\eta(q)^{-1} (\Theta_{a_+,b}(q,1)-\Theta_{a_-,b}(q,1)),$$
$a_\pm=r(m+1)\mp sm$, $b=m(m+1)$, $r=2j+1$ and $s=2k+1$.
\vskip .1in
\noindent \bf Proof. \rm By the product formula
$$\Theta_{0,1}(q,z)\Theta_{s,m}(q,z)
=\sum_{r\equiv s \, (2), \, |r|\le m} \Theta_{r,m+1}(q,z) f_{rs}(q),$$
where 
$$f_{rs}(q)=\sum_{r^\prime\in r +2m{\Bbb Z}}
q^{(r^\prime(m+1)-qm)^2/4m(m+1)}.$$
Note that $f_{rs}(q)=f_{-r,-s}(q)$. Thus
$$\eqalign{\Theta_{0,1}(q,z) (\Theta_{s,m}(q,z) -\Theta_{-s,m}(q,z)
&=\sum_{r\equiv s\, (2),\, |r|\le m} (\Theta_{r,m+1}(q,z)
-\Theta_{-r,m+1}(q,z)) f_{rs}(q)\cr
&=\sum_{r\equiv s\, (2), \, 1\le r\le m} (\Theta_{r,m+1}(q,z)
-\Theta_{-r,m+1}(q,z)) (f_{rs}(q)-f_{-r,s}(q)),\cr}$$
since there is no contribution for $r=0$. Thus the character of the
multiplicity space for $L(k,\ell+1)$ is
$${(f_{r,s}(q)-f_{-r,s}(q))\over \eta(q)}=
{(\Theta_{r(m+1)-sm,m(m+1)}(q,1)
-\Theta_{r(m+1)+sm,m(m+1)}(q,1))\over\eta(q)},$$
as required.

\vskip .1in
\noindent \bf 4. THE KAC DETERMINANT FORMULA \rm Note that as in Chapter~II, section~14,
a Verma module $V(c,h)$ can be constructed which is a representation of the Virasoro algebra
with central charge $c$, generated by a cyclic vector $\xi_0$ such that 
$L_n\xi=0$ for $n>0$ and $L_0\xi =h\xi$. It has the universal property that, for any representation 
generated by a cyclic vector
satisfying similar relations, there is  unique equivariant map sending 
$\xi_0$ to the cyclic vector. A basis of the Verma module
is given by monomials in the raising operators $L_{-k}^{n_k} \cdots L_{-2}^{n_2} L_{-1}^{n_1}\xi_0$, 
where $n_i\ge 0$. Clearly the Verma module is a positive energy representation. 

We now assume that $c$ and $h$ are real.
Let $f:V(c,h)\rightarrow {\Bbb C}$ be the linear map 
picking out the coefficient of $\xi_0$ and extend the involution $L_n^*=L_{-n}$ to a complex involution
on the universal enveloping algebra, so that $(AB)^*=B^*A^*$. We can then define a hermitian form
on $V(c,h)$ by $(A\xi_0,B\xi_0)=f(B^*A\xi_0)$. By definition it satisfies the invariance condition
$(L_n\xi,\eta)=(\xi,L_{-n}\eta)$. Since $L_0^*=L_0$, the eigenspaces of $L_0$ are orthogonal. Moreover
$(\cdot,\cdot)$ is the unique invariant hermitian form on $V(c,h)$ with $(\xi_0,\xi_0)=1$: for the orthogonality 
conditions force $(A\xi_0,B\xi_0)=(B^*A\xi_0,\xi_0)=f(B^*A)$. In particular if $L(c,h)$ is a unitary irreducible 
representation and $T:V(c,h)\rightarrow L(c,h)$ is the canonical map, then the invariant hermitian form on $V(c,h)$ 
is just the pull back of the inner product on $L(c,h)$. Let $K=\{\xi\in V(c,h)|(\xi,V(c,h))=0\}$. 
Then $K$ is invariant under the Virasoro algebra and the hermitian form passes to a non--degenerate invariant
hermitian form on $L=V(h,c)/K$. Since $K$ is invariant under $L_0$, $L$ is itself a positive energy representation. 
We claim that $L$ is irreducible. In fact let $L^\prime$ be a submodule of $ L$ and let 
$v_0$ be the image of $\xi_0$ in $L$. 
$L^\prime$ can be written as the direct sum of eigenspaces $L_0$. Choose $v\ne 0$ in $L^\prime$, 
the image of $A\xi_0$. Thus $(A\xi_0,B\xi_0)\ne 0$ for some monomial $B$, by nondegeneracy. 
Hence $(B^*Av,v_0)\ne 0$. But then $L^\prime(0)\ne 0$, so that $v_0$ lies in $L^\prime$ and thus $L^\prime=L$.

It follows that the Verma
module $V(c,h)$ is irreducible iff $(\cdot,\cdot)$ is
non--degenerate. In particular this happens iff $(\cdot,\cdot)$ is
non--degenerate on every energy subspace $V(N)$. Let $M_N(c,h)$ be the
$P(N)\times P(N)$ matrix 
$$(L_{-i_p} \cdots L_{-i_1}\xi_0, L_{-j_q} \cdots L
_{-j_1}\xi_0)$$
where $1\le i_1 \le \cdots \le i_p$ and $1 \le j_1 \le \cdots \le
j_q$ with $N=\sum i_s =\sum j_t$. The Kac determinant ${\rm
det}_N(c,h)$ is the determinant of this matrix. Note that $L(c,h)$ is
unitary iff $M_N(c,h)$ is positive semi--definite for all $N$. In this
case, it is is necessary that ${\rm det}_N(c,h)\ge 0$ for all
$N$. Using raising and lowering operators to compute the entries of
$M_N(c,h)$, we see that they are all polynomials in $c$ and $h$ if
$c,h\in {\Bbb R}$. 
\vskip .1in
\noindent \bf Examples. \rm (0) ${\rm det}_0(x,h)=\|\xi_0\|^2=1$.

\noindent (1) ${\rm det}_1(c,h)=(L_{-1}\xi_0,L_{-1}\xi_0)=2h$.

\noindent (2) ${\rm det}_2(c,h)=\pmatrix{4h+c/2 & 6h \cr
                                          6h  & 8h^2 +4h\cr}=2h(16h^2
+2hc -10h +c).$
\vskip .1in 
\noindent \bf Lemma. \it If $L(c,h)$ is unitary, then $h\ge 0$ and
$c\ge 0$.
\vskip .1in
\noindent \bf Proof. \rm The computation of ${\rm det}_1=2h$ shows
that $h\ge 0$. Now for $n>0$ we compute 
$$\|L_{-n}\xi_0\|^2=([L_n,L_{-n}]\xi_0,\xi_0)=2nh +c(n^3-n)/12.$$
For this to be positive for all values of $n$, we must have $c\ge 0$. 

\vskip .1in

\noindent \bf Proposition. \it For fixed $c$, ${\rm det}_N$ is a
polynomial in $h$ of degree $\sum_{1\le rs\le N} p(N-rs)$. (The
coefficient of the highest power of $h$ is independent of $c$.)
\vskip .1in
\noindent \bf Proof. \rm  We prove the result by ``degenerating to
bosons''. For $c$ fixed, let $h=t^{-2}$ and
 $a_0=h^{-1}L_0$, $a_n=(2h)^{-1/2} L_n$. Thus
$a_0\xi_0=\xi_0$, $[a_m,a_{-m}]=m a_0 +t^2 c(m^3-m)/12$ and
$[a_m,a_n]=(m-n)ta_{m+n}$ if $m\ne -n,0$, $[a_0,a_m]=-m t^2 a_m$.
Moreover $a_n^*=a_{-n}$. If we look at monomials $(a_{-i_p}\cdots
a_{-i_1}\xi_0, a_{-j_q} \cdots a_{-j_1}\xi_0)$, these are polynomials
in $t$. We extend these to $t=0$; this obviously gives the leading
order terms in $h$ in the original problem. In the limit $t=0$, we get
the system of oscillators $a_0=I$, $[a_m,a_{-m}]=mI$. For this bosonic
system it is immediate that
$x=(a_{-p}^{m_p} \cdots a_{-1}^{m_1}\xi_0,a_{-q}^{n_q} \cdots
a_{-1}^{n_1} \xi_0)$ is zero unless $m_s=n_s$ for all
$s$, in which case $x=\prod m_s! s^{m_s}$. This is independent of
$c$. If we substitute these terms into the determinant for ${\rm
det}_N$, we see that the off--diagonal terms vanish when $t=0$, so the
determinant is given by the product of the diagonal entries, all
non--zero. Thus $\lim_{h\rightarrow \infty} h^{-M} {\rm det}_N\ne 0$
is indepedent of $c$, where $M$ is the sum of all $\sum j_k$'s with
$\sum j_k k =N$. Let $m(r,s)$ be the number of partitions of $N$ in
which $r$ appears exactly $s$ times. Cearly $M=\sum_{1\le rs\le N}
s\cdot m(r,s)$. Now the number of partitions of $N$ in which $r$
appears $\ge s$ times is $P(N-rs)$. Thus
$m(r,s)=P(N-rs) -P(N-r(s+1))$ (where $P(0)=1$ and $P(-k)=0$ for
$k>0$). Thus 
$$M=\sum s\cdot m(r,s)=\sum \sum_{s=1}^{[{n\over r}]}
s\cdot(P(n-rs)-P(n-r(s+1)) =\sum_{1\le rs\le N} P(N-rs).$$
Since $h^M$ is evidently the highest power of $h$ with a non--zero
coefficient in ${\rm det}_N$, the result follows.
\vskip .1in
\noindent \bf Definitions. \rm Let 
$$h_{p,q}(c) ={1\over 48} [(13-c)(p^2 + q^2) +\sqrt{(c-1)(c-25)}
(p^2-q^2) -24pq -2 +2c].$$
Set $\varphi_{p,p}(c,h) =h-h_{p,p}(c) = h +(p^2-1)(c-1)/24$ and
$\varphi_{p,q}(c,h)=(h-h_{p,q}(c))(h-h_{q,p}(c))$
$$=(h-(p-q)^2/4)^2 +{h\over 24} (p^2+q^2 -2)(c-1) +{1\over 576}
(p^2-1)(q^2-1)(c-1)^2 +{1\over 48} (c-1) (p-q)^2 (pq+1).\eqno{(*)}$$
If we parametrise $c$ as $c=1-6/m(m+1)$, then
$$h_{p,q}(c)={((m+1)p-mq)^2 -1\over 4m(m+1)}.$$
\vskip .1in
\noindent \bf Kac determinant formula. \it ${\rm det}_N(c,h)=C_N
\prod_{1\le rs\le N} (h-h_{rs}(c))^{P(N-rs)}$, where $C_N>0$ is
independent of $c$ and $h$.
\vskip .1in

\noindent \bf Lemma~1. \it If $t\mapsto
A(t)$ is a polynomial mapping into $N\times N$ matrices and ${\rm
dim}\, {\rm ker}\, A(t_0)=k$, then $(t-t_0)^k$ divides ${\rm
det}\,A(t)$. 
\vskip .1in
\noindent \bf Proof. \rm Take a basis $v_i$ such that $A(t_0)v_i=0$ for
$i=1,\dots,k$. Thus the first $k$ columns of $A(t)$ are divisible by
$t-t_0$ and hence $(t-t_0)^k$ divides ${\rm det}\, A(t)$. 
\vskip .1in
\noindent \bf Lemma~2. \it Fix $c$ and regard ${\rm det}_N(c,h)$ as a
polynomial in $h$. If ${\rm det}_k$ vanishes at $h=h_0$, then 
$(h-h_0)^{P(N-k)}$ divides ${\rm det}_N$.  
\vskip .1in
\noindent \bf Proof. \rm We may take $k$ minimal subject to ${\rm
det}_k(c,h_0)=0$. Thus $V(c,h_0)$ has a singular vector $v$ at energy
level $k$. By the Poincar\'e--Birkhoff--Witt theorem, the vectors
$L_{-i_t} \cdots L_{-i_1} v$ are all linearly independent for $i_t\ge
\cdots\ge i_1\ge 1$. So at level $N$, this submodule has dimension
$P(N-k)$. On the other hand this submodule is contained in ther kernel
of $(\cdot,\cdot)$. Thus $M_N$ has a kernel of dimension at least
$P(N-k)$ at $h_0$. The assertion therefore follows from Lemma~1.

\vskip .1in
\noindent \bf Lemma~3. \it ${\rm det}_N$ vanishes at $h_{r,s}(c)$ for
$1\le rs\le N$. 
\vskip .1in
\noindent \bf Proof. \rm By the GKO construction, 
$${\rm ch}\,
L(c^{(m)}, h^{(m)}_{r,s})\le {q^h\over \varphi(q)} (1-q^{rs}
-q^{r^\prime s^\prime} +\cdots),$$
where $r^\prime=m-r$, $s^\prime=m+1-s$ and the inequality is to be
understood in terms of coefficients of $q^i$. It follows that the
kernel of $(\cdot,\cdot)$ in $V(c,h)$ has a non--zero component at
each energy level $N\ge \min(rs,r^\prime s^\prime)$. Thus ${\rm
det}_N$ vanishes at $h_{rs}^{(m)}$ for $m$ sufficiently large. But
then ${\rm det}_N$ vanishes at infinitely many points of the
curve $\varphi_{rs}(c,h)=0$ (namely $(c^{(m)},h^{(m)}_{r,s})$). Since
$\varphi_{r,s}(c,h)$ is irreducible in ${\Bbb C}[c,h]$, we see that
$\varphi_{r,s}$ divides ${\rm det}_N$ for $N\ge rs$. Thus
 ${\rm det}_N$ vanishes at $h_{r,s}(c)$ for
$1\le rs\le N$.

\vskip .1in
\noindent \bf Proof of determinant formula (Kac--Wakimoto). \rm By
Lemmas~2 and 3, ${\rm det}_N$ is divisible by 
$$\prod_{1\le pq\le N} (h-h_{p,q}(c))^{P(N-pq)},$$
since the $h_{p,q}(c)$'s are distinct for
generic $c$. Since both sides have the same degree in $h$ and the
highest order term in $h$ is indepedent of $c$ (by the Proposition),
the result follows.
\vskip .1in

\noindent \bf 5. THE FRIEDAN--QIU--SHENKER UNITARITY CRITERION FOR
h. \rm We prove the 
easy part of the FQS criterion for 
unitarity. (The harder part of their criterion gives the restrictions
on the values of $c$. It depends on a detailed knowledge of the
representations $L(1,m^2/4)$; in this sense, their proof is analogous
to the proof of Jones' index theorem that uses detailed knowledge of
the limiting $SU(2)$ subfactor.)
\vskip .1in
\noindent \bf FQS Theorem. \it Let $L(c,h)$ be a unitary representation
of the Virasoro algebra with $c=1- 6/m(m+1)$ for $m\ge 3$. Then
$h=h_{p,q}$ with $1\le q\le p \le m-1$ and $h_{p,q}=[(p(m+1) -
qm)^2-1]/4m(m+1)$.
\vskip .1in
\noindent \bf Proof. \rm Let 
$$h_{p,q}(c) ={1\over 48} [(13-c)(p^2 + q^2) +\sqrt{(c-1)(c-25)}
(p^2-q^2) -24pq -2 +2c].$$
Set $\varphi_{p,p}(c,h) =h-h_{p,p}(c) = h +(p^2-1)(c-1)/24$ and
$\varphi_{p,q}(c,h)=(h-h_{p,q}(c))(h-h_{q,p}(c))=$
$$(h-(p-q)^2/4)^2 +{h\over 24} (p^2+q^2 -2)(c-1) +{1\over 576}
(p^2-1)(q^2-1)(c-1)^2 +{1\over 48} (c-1) (p-q)^2 (pq+1).\eqno{(*)}$$
If we parametrise $c$ as $c=1-6/m(m+1)$, then
$$h_{p,q}(c)={((m+1)p-mq)^2 -1\over 4m(m+1)}.$$
It will sometimes be more convenient to use the variable
$x=m+{1/2}$. Thus $c= 1-6/(x^2 -{1\over 4})$ and
$h_{p,q}(c)=[(x(p-q) +{1\over 2} (p+q))^2-1]/(4x^2 -1)$. Let $C_{p,q}$
be the real curve $\{(c,h):\varphi_{p,q}(c,h)=0\}$. By the symmetry of
$(*)$ in $p$ and $q$, $C_{p,q}=C_{q,p}$. The form of $(*)$ also shows
that $\phi_{p,q}(c,h)>0$ for $c>1$ and $h>0$. Hence, for $h\ge 0$, the curve
$C_{p,q}$ lies in the region $c\le 1$. We shall only consider it in
the domain $h\ge 0$ and $c\ge 0$. The curve $C_{p,q}$ is parametrised
by $x\in {\Bbb R}$ with $x=\pm \infty$ giving its intersection with
the line $c=1$ which it touches at $h=(p-q)^2/4$, giving two branches
$C_{p,q}^\pm$ according to the sign of $x$. Clearly $C_{p,q}^\pm
=C_{q,p}^\mp$ and $C_{p,q}^+$ is the upper branch if $p>q$. The curve
$C_{p,q}$ arises at level $N=pq$. Thus there are infinitely many
curves through the point $(1,M^2/4)$. The figure below shows the possible 
curves schemtically (with a rescaling in the $h$ direction).
\vskip .1in
\epsfxsize=3in
\centerline{\epsfbox{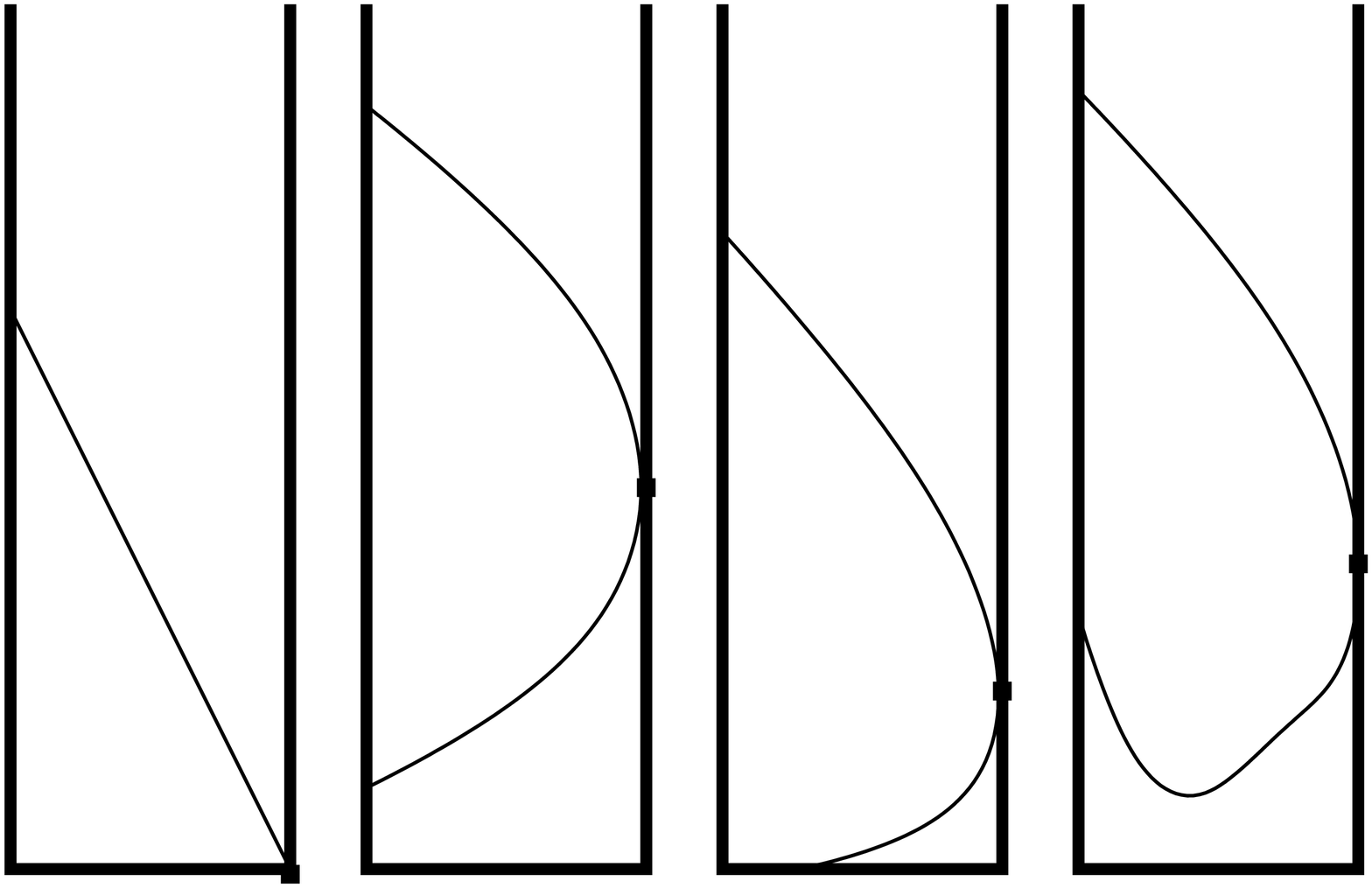}}
\vskip .1in
\noindent In the degenerate case $p=q$, the curves $C_{p,p}$ become straight
lines joining the point $(1,0)$ to $(0,(p^2-1)/24)$, the first case in
the figure above. The remaining three cases give all possible shapes
for the curve $C_{p,q}$. The curve $C^+_{p,q}$ cuts the $h$--axis
at $h=[(3p-2q)^2-1]/24$ and the $c$-axis at
$0$ if $|3p-2q|=1$, between $0$ and $1$ if $3p=2q$ and not all otherwise.

We next make the important observation that at through the special points
$(1,M^2/4)$ there at most one curve of a given level and the upper branch 
of a higher level curve lies above that of curves of lower level and the lower
branch lies below that of the lower level curves the special point. 
In fact if $C_{p,q}$ occurs at a higher
level than $C_{r,s}$ with $(p-q)^2/4=(r-s)^2/4$ and $p\ge q$, $r\ge
s$, then $p> r$, $q> s$, so that $p-q=r-s>0$ and $p+q>r+s$. Since we
evidently have
$$((p-q)x + {1\over 2}(p+q))^2 > ((r-s) x + {1\over 2}(r+s))^2$$
for $x\ge 0$ and the opposite inequality for $x<0$, sufficiently
large, we see that near $(1,(p-q)^2/4)$, the curve $C_{r,s}$ lies
between the upper and lower branches of the curve $C_{p,q}$. Note also
that it is not possible for $r-s=p-q$ with $(r,s)$ and $(p,q)$ having
the same level. For if $pq=rs$, $p,-q$ and $r,-s$ are roots of the
same quadratic. Since $p,r\ge 1$ and $-q,-s\le -1$, it follows that
$p=r$ and $q=s$. Hence at any fixed level, there is at most one curve
through $(1,M^2/4)$. 

Let $X$ be the strip $[0,1]\times [0,\infty)$ in the $(c,h)$ plane and
let $U_{p,q}$ be the open subset where $\varphi_{p,q}(c,h)<0$. Its
closure $\overline{U_{p,q}}$ is compact and connected with boundary
made up of the segment of the curve $C_{p,q}$ in $X$ and parts of the
lines $c=0$ and possibly $h=0$. In fact if $p=q$, $\overline{U_{p,q}}$
is a rightangle triangle with vertices $(1,0)$, $(0,0)$ and
$(0,(p^2-1)/24)$. While if $p>q$, as may be assumed without loss of
generality, $C^+_{p,q}$ increases as $c$ decreases, cutting the
$h$--axis at $(0,[(3p-2q)^2-1]/24)$; similarly the lower branch 
$C_{p,q}^-$ cuts the $h$--axis at $(0,[(3q-2p)^2-1]/24)$, which lies
outside $X$ if $|3q-2p|<1$, in which case it cuts the $c$--axis at
$(1-6/(x^2-{1\over 4}),0)$ with $x={1\over 2}(p+q-2)/(p-q)$. 

Let $X_N=\bigcup_{pq\le N} \overline{U_{p,q}}$. We will show that thia
is the region bounded a explicit sequence of pieces of the curves $C_{p,q}$
as depicted schematically in the diagram below, with the $h$ direction
rescaled to make the special points equally spaced on the line $c=1$.  
\vskip .1in
\epsfxsize=3in
\centerline{\epsfbox{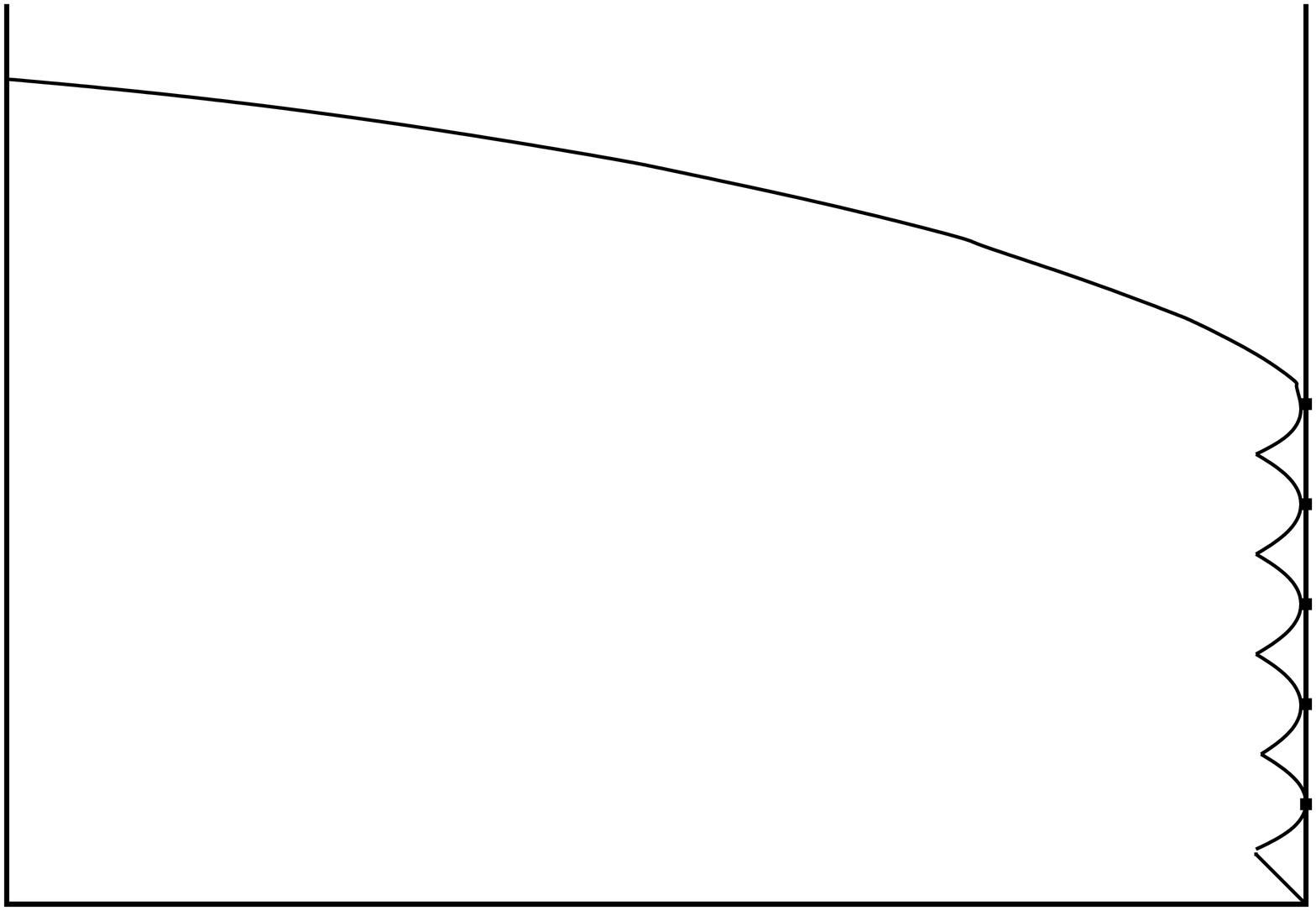}}
\vskip.1in

Then for
$M=0,1,2,\dots,N-1$ there are curves $C_{pq}$ with $p>q$, $p-q=M$ and $pq\le N$.
When $M=N-1$, $p=N$ and $q=1$, so the upper branch of the curve cuts the $h$--axis at 
$c=[(3N-2)^2-1]/24$. The lower branch is decreasing and cuts the $h$--axis between
$0$ and $(N-1)^2/4$. At the other extreme $m=0$, the curve of highest level through
$(1,0)$ is $C_{p,p}$ where $p=[\sqrt{N}]$, the largest integer 
less than or equal to $\sqrt{N}$.

More generally if we take a highest level curve 
$C_{p,q}$ through a special point $(1,M^2/4)$ with $M=p-q>0$, then 
$pq\le M < (p+1)(q+1)$. The highest level curve through the special point with
$(1,(M+1)^2/4$ is
$C_{p+1,q}$ if $(p+1)q\le N$ and $C_{p,q-1}$ otherwise. The highest level curve
$(1,(M-1)^2/4$ is $C_{p,q+1}$ if $p(q+1)\le N$ and $C_{p-1,q}$ otherwise.
Note that if $p(q+1)>N$ then necessarily $(p+1)q>N$ since $p\ge q$, 
so that for are are at most three
possibilities for the highest level curves through the two adjacent 
special points. In particular if a highest curve through a special point has
level $N$, then the highest level curves through the adjacent 
points have level $<N$.

We now work out the points of intersection of these curves: 
in fact not only do we calculate them but we prove that the key observation of
Friedan--Qiu--Shenker 
that they are these points are characterized as being the 
first intersections of the
highest level curve of level $N$ through $(M^2/4,1)$ with any other curves 
of level $\le N$. 
\vskip .1in
\noindent \bf Lemma~1. \it The curve $C^+_{N,1}$ does not 
intersect any of the other curves at level $N$.
\vskip .05in
\noindent \bf Proof. \rm It suffices to shows that for each value of $x$ 
the corresponding point on  $C^+_{N,1}$ lies above any point on $C^+_{p,q}$ for
$p\ge q$ and $pq\le N$. But this follows immediately because
$$x(N-1) +N/2 \ge x(p-q) +(p+q)/2.$$ 
Indeed if $q=1$, then $p< N$ and the result is obvious. If $q\ge 2$, then
$N-1> p-q$. Moreover $p\le N/q\le N/2$. Since $q\le p$, this implies 
$p+q\le 2p\le N$.

\vskip .1in
\noindent \bf Lemma~2. \it If $|p-q|\ne |r-s|$ and $(p,q)$ is not
proportional to $(r,s)$, the curves 
$C^+_{p,q}$ and $C^+_{r,s}$
intersect transversely in the $(c,h)$--plane in the distinct points 
$$(1-6/(x^2-{1\over 4}),[(x(p-q) +{1\over 2}(p+q))^2-1]/(4x^2 -1))$$ 
with $x={1\over 2}(r+s - p -q)/(p-q + s -r)$ and
$x={1\over 2} (r+s+p+q)/(q-p + s -r)$. 
\vskip .05in
\noindent \bf Proof. \rm  The points of intersection are given by the
solutions of
$$x(p-q)+{1\over 2}(p+q) = \pm [x(r-s)+{1\over 2}(r+s)].$$
Transversality occurs if and only if the derivatives at $x$ are equal, i.e.
$$(p-q)[x(p-q)+{1\over 2}(p+q)] = (r-s)[x(r-s)+{1\over 2}(r-s)].$$
Since
$$|x(p-q)+{1\over 2}(p+q)| = |x(r-s)+{1\over 2}(r+s)|$$
while $|p-q|\ne|r-s|$, this can only happen if both sides above
vanish and the two points of intersection coincide. In this case
$$x={1\over 2}(p+q)/(q-p)={1\over 2}(r+s)/(s-r)$$
so that $p/q=r/s$. 
\vskip .1in

\noindent \bf Lemma~3. \it If $C^+_{p,q}$ is a highest level curve at level $N$
with 
$p\ge q$, $pq\le N$ and $(p+1)(q+1)>N$
then the first intersection with a curve of level
$\le N$ is 
with $C^-_{p, q-1}$ if $(p+1)q> N$ and with $C^-_{p+1,q}$ otherwise. 
The intersection
with $C^-_{p+k, q+k-1}$ takes place at $x=p+q +k -1/2$ where $k=0$ or $1$. 
The first intersection of
$C^-_{p,q}$ with a curve of level $\le N$ with $C^+_{p-1, q}$ if $p(q+1)>N$ and
with $C^+_{p+1,q}$ otherwise. 
The intersection
with $C^+_{p+k-1, q+k-1}$ takes place at $x=p+q +k -1/2$ where $k=0$ or $1$.
The interections are transverse.

\vskip .05in
\noindent \bf Proof. \rm Immediate from Lemma~2. 
\vskip .1in
 
In the prvious diagram we have marked the parts of the 
highest degree curves lying between
these intersections. Let $X^\prime_N$ be the area bounded by these curves
and the $c$ and $h$ axes. We now verify the assertion made above about the
boundary of the closed region $X_N$.
\vskip .1in
\noindent \bf Theorem. \it $X_N=X^\prime_N$.
\vskip .05in

\noindent \bf Proof. \rm 
We assume the result by induction on $N$, the result being obvious for $N=1$. 
\vskip .1in
\epsfxsize=3in
\centerline{\epsfbox{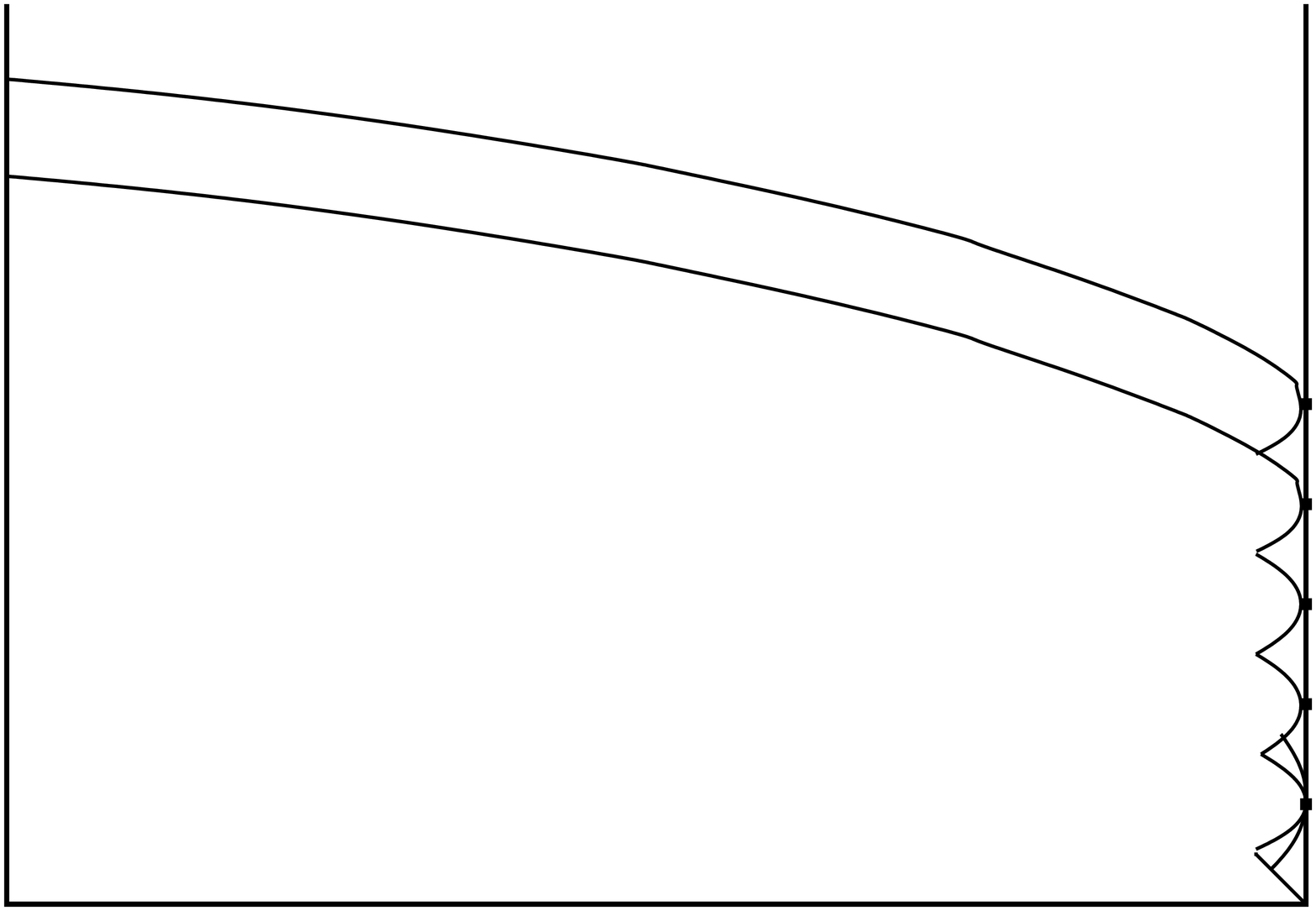}}
\vskip .1in
At level $N$, the regions bounded by the
curves $C_{pq}$ with $pq=N$ are added. 
Now when one of these new curves occurs at level $N$ through $m^2/4$, 
there are no new curves through $(m\pm 1)^2/4$. 
We claim that the new areas added are just the areas between the new curve, 
the old curve $C_{p-1,q-1}$ and the places where the new curve cuts the 
highest level curves to either side. For a curve $C_{p,q}$ with $q>1$,
After the new curve cuts the boundary of the old region it cannot intersect 
the curved part of the old boundary again by Lemma~3. 
The same arguments apply to the curve $C_{N,1}$ 
using Lemma~1 and Lemma~3.
\vskip .1in
\bf \noindent Remark. \rm As observed by Friedan--Qiu--Shenker, 
an immediate consequence of the 
theorem is that the new regions added to $X_{N-1}$ 
to produce $X_N$ are 
exactly the regions bounded by the curves $C_{p,q}$ of level $N$, the previous
highest level curve $C_{p-1,q-1}$ through that point and 
the two adjacent highest level curves $C_{p-1,q}$ and 
$C_{p,q+1}$ through adjacent special points; 
for the special point $(1,(N-1)^2/4$, 
the new region is the one between $C_{N,1}$ and $C_{N-1,1}$. These two types
of region are indicated in the diagram above: the first is like a bow--tie; 
the second like a curved strip.     
\vskip .1in
\bf \noindent Proof of FQS Theorem. \rm 
We first prove that for fixed $c\in (0,1)$ with $c=1-6/m(m+1)$, the only
values of $h$ for which $L(c,h)$ can be unitary are the
$h_{p,q}^{(m)}$ with $p,q\ge 1$ and $p+q\le m$. At level $N$ define the
$N$th excluded region $R_N$ as being the open set where
${\rm det}_n$ is negative for some $n\le N$. We claim that 
$$\bigcup R_N = \{(c^{(m)},h^{(m)}_{p,q}): m>p+q-1\}.$$
Note first that every point $(c,h)$ with $0<c<1$ and $h>0$ lies in the
interior of one of the curves $C_{p,q}$ for $pq$ sufficiently
large. Indeed the regions $U_{p,1}$ sweep out the region since the points of intersection
of $C_{p,1}$ with the $h$--axis tend to infinity.

Consider the closure $F_N$ of all the interiors of the curves
$C_{p,q}$ for $pq\le N$. We shall prove that the only parts of $F_N$
that might not lie in $R_N$ are the parts of the curves $C_{p,q}$,
$pq=N$ with $|m+1/2|>|p+q-{1/2}|$. In fact the functions ${\rm det}_N$
do not vanish for $c>1$ and therefore have the same sign. Since the
Segal--Sugawara construction gives positive values for all $N$, it
follows that ${\rm det}_N$ is positive for all $c>1$, $h>0$ and $N\ge
0$. Hence if $p\ne q$, near $(1,(p-q)^2/4)$, ${\rm det}_N$ is positive
for $c>1$. Thus on one side of the curve $C_{p,q}$ near 
$(1,(p-q)^2/4)$, ${\rm det}_N$ is positive. Note that near
$(1,(p-q)^2/4)$, the function ${\rm det}_N$
changes sign as $C_{p,q}$ is crossed, because at the $N$th stage
$h_{p,q}$ is a simple zero of ${\rm det}_N$. We can therefore exclude
all parts of the interior of $C_{p,q}$ which do not meet $F_{N-1}$ or
the closure of the interior of any other $C_{rs}$ with $rs=N$. There
are two such regions. Since the first curves of level $\le N$ met by
$C_{p,q}$ are $C_{q-1,p}$ and $C_{q,p-1}$, both of these regions are
bounded on one of their sides by a segment of $C_{p,q}$ with
$|m+{1\over 2}|>p+q-{1\over 2}$ (starting at $(1,(p-q)^2/4)$) and on
their other sides by segments of curves $C_{rs}$ with $rs<N$. It
follows that the new parts of the closure of the interiors of the
curves $C_{p,q}$ where ${\rm det}_N$ vanishes are exactly the boundary
parts with $|m+{1\over 2}|>p+q-{1\over 2}$. Hence the only possible
points of unitarity in $h>0$ and $c\in (0,1)$ have $h=h_{p,q}^{(m)}$
with $p+q<m+1$. 

Finally we specialise to the case when $m\ge 3$ is an integer. 
We have shown that if $L(c,h)$ is unitary then
$h=h_{p,q}^{(m)}=[(p(m+1)-qm)^2-1]/4m(m+1)$ with $p,q\ge 1$ and
$p+q\le m$. We want to show that $h=h_{rs}^{(m)}$ with $1\le s\le r\le
m-1$. Note that since $p,q\ge 1$ and $p+q\le m$, we have $p,q\le
m-1$. If $p\ge q$, we take $r=p$, $s=q$. Otherwise $q\ge p+1\ge
2$. Let $p^\prime=m-p$ and $q^\prime=m+1-q$. Then
$h_{pq}^{(m)}=h_{p^\prime,q^\prime}^{(m)}$ and $1\le q^\prime\le
p^\prime \le m-1$. so in this case we may take $r=p^\prime$,
$s=q^\prime$. This completes the proof of the FQS proposition. 
\vskip .1in
\noindent \bf Corollary. \it For $c=1-6/m(m+1)$ with $m\ge 3$, the
values of $h$ from which $L(c,h)$ is unitary are given by $h=h_{r,s}=
[(r(m+1) - sm)^2 -1]/4m(m+1)$ with $1\le s\le r\le m-1$.
\vskip .1in

\noindent \bf 6. THE MULTIPLICITY ONE THEOREM. \it Each of the
multiplicity spaces $K_h$ appearing in the $SU(2)\times SU(2)/SU(2)$
decomposition gives an irreducible representation of the Virasoro
algebra, so that $K_h=L(c,h)$. 
\vskip .1in
\noindent \bf Lemma. \it Let $h=h_{rs}$  and $M=rs+(m-r)(m+1+s)=m(m+1)-(m+1)r +ms$. Then
$${\rm ch}\,L(c,h)={\rm ch}\, K_{h} \, {\rm mod}\, q^{h+M}.$$ 
\vskip .1in
\noindent \bf Proof. \rm By the Kac determinant formula we have 
$${\rm det}_N(c^{(m)},h)=\prod_{1\le pq\le N}
(1-h_{pq}^{(m)})^{P(N-pq)}.$$
Note that $h_{rs}^{(m)}=h_{r^\prime s^\prime}^{(m)}$. This is the only
possible such coincidence, because if $((m+1)p-qm)^2=((m+1)r-ms)^2$
with $1\le p,q\le m$ and $1\le s\le r \le m-1$, we would have
$(m+1)(p\pm r) =m (q\pm s)$. Since $m$ and $m+1$ are coprime, we would
have $p\pm r=am$ and $q\pm s=a(m+1)$ for some integer $a$. Hence
either $p=r$ and $q=s$ or $p=m-r$ and $q=m+1-s$. 
Similar reasoning shows that $rs\ne r^\prime s^\prime$. Indeed if
$rs=r^\prime s^\prime$, we get 
$rs=(m-r)(m+1-s)$. Thus $(m-r)(m+1)=sm$. Since $m$ and $m+1$ are
coprime, we would have $r$ would be divisible by $m$ and $s$ by
$m+1$. This contradicts $1\le r,s\le m-1$. Thus we may assume that
$rs<r^\prime s^\prime$. It follows that $h=h_{rs}^{(m)}$ is first a
zero of ${\rm det}_N$ when $N=rs$. It has multiplicity one. It has
multiplicity $P(N-rs)$ for $rs\le N< r^\prime s^\prime$ and multiplicity
$P(N-rs)+P(N-r^\prime s^\prime$ when $r^\prime s^\prime \le N<M$. By
Lemma~1 in section~4, this gives an upper bound for the dimension of
the kernel of the form $(\cdot,\cdot)$ and hence a lower for the character:
$${\rm ch}\, L(c,h_{rs}^{(m)}) \ge {q^{h_{rs}}\over
\varphi(q)}(1-q^{rs} - q^{r^\prime s^\prime}) \,{\rm mod}\, q^{h_{rs}+M}.$$
Here $\varphi(q) =\prod_{n\ge 1}(1-q^n)$ and the inequality means that
the coefficient of $q^{i+h}$ on the left hand side is greater than or
equal to the coefficient on the right hand side for $i\ge M$.  On the other hand
the right hand side agrees with ${\rm ch}\, K_h$ mod $q^{h_{rs}+M}$.
Since ${\rm ch}\,L(c,h)\le {\rm ch}\,K_h$, the result follows.

\vskip .1in
\noindent \bf Proof of Theorem. \rm By the preceding lemma, we have
${\rm ch}\, L(c,h_{rs}^{(m)})$ and ${\rm ch}\, K_h$ agree for energy levels $<
M=rs-(m+1)r+ms$, where
$K_h$ is the $SU(2)\times SU(2)/SU(2)$ multiplicity space. Now suppose that the
representation on the multiplicity space $K_h$ is not irreducible. The
character computation we have made so far shows that there can be no
singular vectors with level $< M$ in $K_h$. For any such would
already appear in the cyclic module generated by the lowest energy
vector $v_h$ by the equality above. This module is irreducible (by
unitarity), so has no singular vectors apart from its lowest energy
vector $v_h$. On the other hand any lowest energy singular vector in
$K_h$ would have energy $h^\prime=[((m+1)p-mq)^2-1]/4m(m+1)$ with
$1\le q\le p\le m-1$ by the FQS criterion. We will check that $h^\prime<M+h$, i.e.
$$[((m+1)p-mq)^2-1]/4m(m+1)<  rs -(m+1)r +ms+ [((m+1)r-ms)^2-1]/4m(m+1),
\eqno{(1)}$$
so that such a vector
would have to lie within the energy range discussed above. It follows
that $K_h$ is irreducible. 

To prove (1), note the left hand side is maximised by taking $p=m-1$ and $q=1$, 
so we must show that
$${(m^2-m -1)^2-1\over 4m(m+1)} < m(m+1) -(m+1)r +ms  +{((m+1)r-ms)^2-1\over 4 m(m+1)}.$$
The quadratic expression in $(r,s)$ is minimised on the triangle $1\le s\le r\le m-1$ 
at its vertices. For fixing $s$, the derivative in $r$ of the right hand side is 
$-m-1 + ((m+1)r-ms)/2m <0 $; fixing $r$ the derivative in
$s$ is $+m +(ms-(m+1)r)/2m>0$. Thus in the interior of the triangle we can always decrease the right hand
side by moving towrds an edge parallel  the axes. On interior points of edges parallel to the axes we can decrease
the right hand side by moving towards a vertex. On the segment $(r,s)=(t,t)$ with $1\le t\le m-1$,
the derivative in $t$ of the right hand side is $-1+ t/4m(m+1)<0$. Thus the minimum occurs at a vertex. 
At the vertex $(r,s)=(m-1,1)$, the two extreme terms agree and the middle terms are positive. When
$(r,s)=(1,1)$, the right hand side becomes $m^2+m-1$ and the inequality is immediate. When $(r,s)=(m-1,m-1)$,
the inequality becomes
$${(m^2-m -1)^2-1\over 4m(m+1)} < (m-1)(m+1)  +{(m-1)^2-1\over 4 m(m+1)}.$$
The left hand side equals $(m-1)(m-2)/4$ which is less than the first term on the right hand side.
\vskip .1in
\noindent \bf 7. THE FEIGIN--FUCHS CHARACTER FORMULA FOR THE DISCRETE
SERIES. \it The unitary representation $L(c,h)$ with $c=1-6/m(m+1)$ 
and $h= [(p(m+1)-qm)^2-1]/4(m+1)m$ has normalised character
$${\rm ch}\,
L(c,h)=\eta(q)^{-1}(\Theta_{a_+,b}(q,1)-\Theta_{a_-,b}(q,1)),$$
where $\eta(q)=q^{1\over 24}\prod_{n\ge 1} (1-q^n)$, $a_\pm=p(m+1)\mp
qm$ and $b=m(m+1)$. 
\vskip .1in 
\noindent \bf Proof. \rm We have just shown that the multiplicity
space with lowest energy $h=h_{p,q}$ is irreducible. Thus the character
formula for $L(c,h)$ is given by the character of the multiplicity
space. 
\vskip .1in
\noindent \bf 8. THE FRIEDAN--QIU--SHENKER UNITARITY CRITERION FOR c. \rm 
Our aim now is to prove the complete version of the
Friedan--Qiu--Shenker unitarity theorem.
\vskip .1in
\noindent \bf Theorem (Friedan--Qiu--Shenker). \it If $0<c<1$, then a
representation with 
central charge $c$ is unitary iff $c=1-6/m(m+1)$ with $m\ge 3$.
\vskip .1in
\bf Remark. \rm From our previous work, the permitted values of $h$
for a particular 
$m\ge 3$  are $h_{p,q}(m)$ with $1\le q\le p\le m-1$. 
\vskip .1in
We already know that any point not on a curve $C_{p,q}$ cannot be
unitary. The points $h_{p,q}(m)$ are precisely the intersections of
the different curves $C_{p,q}$. Each intersection $P$ may be described by
first choosing a curve
$C_{p,q}$ with $P\in C_{p,q}$ and $pq$ minimal and then choosinf
$C_{p^\prime,q^\prime}$ with $P\in C_{p^\prime,q^\prime}$ and $p^\prime q^\prime$ minimal. It is easy to check that these
intersections are obtained by taking $p^\prime=p+k$,
$q^\prime=q+k$ for $k\ge 1$ and $m=p+q+k-1$. We now rule out all the
points between these intersection points. Note that every curve
$C_{p,q}$ touches $c=1$ at $h=(p-q)^2/4$. Thus if
$C_{p^\prime,q^\prime}$ intersects $C_{p,q}$,
then the part of $C_{p,q}$ on
the $c<1$ side of $C_{p^\prime,q^\prime}$ is the part with $h$
decreasing if $h^\prime <h$ and $h$ increasing if $h^\prime>h$.

Take a point $P_0$ on $C_{p,q}$ corresponding
to $p^\prime,q^\prime$ as above. Let $N=pq<N^\prime =p^\prime
q^\prime$. Starting from the asymptote at $c=1$ through
$(p^\prime-q^\prime)^2/4$, we may follows the $C_{p^\prime,q^\prime}$
curve as it travels to $P_0\in C_{p,q}$. We get a stright line
parametrisation of $C_{p^\prime,q^\prime}$ by taking $y=m$ as
coordinate. Along the way to $P_0$, the curve will cross other curves
$C_{p^{\prime\prime}, q^{\prime\prime}}$ transversely and simply at
points $P_1,\dots, P_k$. At level $N^\prime$, the dimension of the
null space is $1$ on $C_{p^\prime,q^\prime}$ away from the
intersection points. Near the asymptote on the $c>1$ side of
$C_{p^\prime,q^\prime}$, the matrix of inner products $A(c,h)$ is
positive--definite. We shall find a open neighbourhood $U$ of the part of
the curve $C_{p^\prime,q^\prime}$ above $P_0$ (containing
$P_1,\dots,P_k$ and part of the asymptote to $c=1$) and a rank one
spectral projection $P(c,h)$ of $A(c,h)$ in this strip depending
continuously on $(c,h)$ such that $P(c,h)A(c,h)=\lambda(c,h)P(c,h)$
with $\lambda(c,h)=0$ only on $C_{p^\prime,q^\prime}$. It will follow
that $\lambda(c,h)<0$ on the $c<1$ side of $C_{p^\prime,q^\prime}$. In
particular $A(c,h)$ will have a negative eigenvalue on $C_{p,q}$ in
the segment between $P_0$ and the next intersection. This clearly
will prove the unitarity theorem.
\vskip .1in
Using $m$ as parameter, we can replace the curve
$C_{p^\prime,q^\prime}$ by the $y$--axis. The following result (with
$M=1$) shows the existence and uniqueness of the spectral projection
$P(z)$ for $z$ in an open neighbourhood of the $y$--axis with $z\ne
P_i$. (This open neighbourhood should of course contain $P_0,\dots,
P_k$.)

\vskip .1in
\noindent \bf Lemma~1. \it Let $A(z)$ be a continuous self--adjoint
matrix--valued function on a topological space $Z$ such that ${\rm
ker}(A(z))$ has 
constant rank $1$ (or rank $M$ more generally) for $z\in Z_0$, a closed subset of $Z$, and is
invertible otherwise. For each $z\in Z_0$, there is an open
neigbourhood $U$ of $z$ such that 
if $P(z)$ is the orthogonal projection onto ($M$th) the
lowest eigenspace(s) of $A(z)^2$, then $z\mapsto P(z)$ is continuous on
$U$. If $Z$ is an open 
subset of ${\Bbb R}^n$ and $A(z)$ is a smooth (or analytic) function
of $z$, then $z\mapsto P(z)$ is also smooth (or analytic) on $U$. 
\vskip .1in
\noindent \bf Proof. \rm Take $z\in Z_0$. By minimax there is a
neighbourhood $U$ of $z$ such that the lowest eigenvalue of $A(t)^2$
is less than $r/2$ and the next eigenvalue is greater than
are greater than $2r>0$. Let $\chi$ be a continuous bump
function supported in $(-r,r)$ with $\chi(0)=1$. Then
$P(z)=\chi(A(z)^2)$ for $t\in U$, since the only the lowest eigenvalue
of $A(z)$ occurs in $(-r,r)$. Since $\chi$ can uniformly 
approximated by polynomials on any compact interval, it follows that
$z\mapsto A(z)$ is continuous on $U$. The
second assertion follows immediately from the contour integral
expression for the spectral projection $P(z)$:
$$P(z)={1\over 2\pi i} \int_{|w|=\varepsilon} (wI-A(z))^{-1} \, dw.$$
\vskip .1in
\noindent \bf Corollary. \it There is an open subset $U$ of $Z$
containing $Z_0$ on which $P(z)$ can be defined (uniquely).
\vskip .1in
\noindent \bf Proof. \rm Take an open neighbourhood $U_z$ for each
point $z\in Z_0$ and set $U=\bigcup_{z\in Z_0} U_z$. By uniqueness the
different $P(z)$'s must agree on intersections of these opens.
\vskip .1in
We next need to use information from the Kac character formula to
continue $P(z)$ across the points $P_i$. We simply have to define
$P(z)$ in an open neighbourhood of each point $P_i$. Since the
intersection at $P_i$ is transverse, we may assume that the transverse
curve is the $x$--axis and $P=P_i$ corresponds to the point $(0,0)$. 
\vskip .1in
\noindent \bf Proposition. \it In an open neighbourhood of $P=P_i$ there
are 

\item{(a)} a unique continuous determination of a rank one projection $P(z)$
such that $P(z)$ is a spectral projection of $A(z)$ coinciding with
the projection onto the kernel of $A(z)$ for $(0,y)$ with $y\ne 0$.

\item{(b)} a unique continuous determination of a rank $m$ projection $Q(x)$
on $y=0$ such that $Q(x)$ is the spectral projection onto the kernel
of $A(x,0)$.

\noindent Moreover $P$ and $Q$ are orthogonal on $y=0$ and $P(0)+Q(0)$
is the projection onto the kernel of $A(0,0)$. 
\vskip .1in
\noindent \bf Proof. \rm Note that $Q(z)$ could be constructed on
$y=0$, $x\ne 0$ using the method of the previous lemma. We need a
variant of this construction. Let $C^\prime=C_{p^\prime,q^\prime}$ and
let the transverse curve be
$C^{\prime\prime}=C_{p^{\prime\prime},q^{\prime\prime}}$ with
$N^{\prime\prime}=p^{\prime\prime}q^{\prime\prime} <N^\prime$. The
kernel of $M_{N^{\prime\prime}}$ is rank one on
$C^{\prime\prime}$. We choose parameters such that $P=(0,0)$,
$C^\prime$ is the $y$--axis and $C^{\prime\prime}$ the $x$--axis. 
As in the previous lemma, let $R(z)$ be a smooth or analytic
determination of a spectral subspace of $M_{N^{\prime\prime}}(z)$
giving the kernel on $y=0$. Let $u(z)=R(z)u_0/\|R(z)u_0\|$ be a smooth
or analytic choice of eigenvector near $z=0$. Define 
vectors $u_j(z)=L_{-j_1} \cdots L_{-j_r} u(z)$ at level $N^\prime$ for
$\sum j_i=N^\prime-N^{\prime\prime}$. There are
$P(N^\prime-N^{\prime\prime})$ such vectors and they form a basis of
$A(z)$ for $z=(x,0)$ with $x\ne 0$. The Gram--schmidt
orthonormalisation process shows that the orthogonal projection $Q(z)$
onto the subspace spanned by the $u_j(z)$'s is smooth or
analytic. Thus $Q(z)$ is a projection of rank
$P(N^\prime-N^{\prime\prime})$ defined in a neighbourhood of $0$. Note
that $Q(z)$ is a spectral projection of $A(z)$ for $y=0$, $x\ne 0$,
but not necessarily otherwise. Let $B(z)=Q(z)A(z)Q(z)$ considered as a
self--adjoint operator on ${\rm im}(Q(z))$. Since ${\rm
im}(Q(x,0))\subseteq {\rm ker}(A(x,0))$, 
we must have $B(x,0)=0$. Hence $B(x,y)=y
B_0(x,y)$ where $B_0$ is smooth or
analytic. The next lemma shows that $B_0(0)$ is invertible.

\vskip .1in
\noindent \bf Lemma~2. \it ${\rm det}B_0(0)\ne 0$.
\vskip .1in
\noindent \bf Proof. \rm We claim that ${\rm det}
M_{N^{\prime\prime}}(c,h+pq)\ne 0$ where $N^{\prime\prime}=N^\prime
-N$. In fact if the determinant vanished, $(c,N^{\prime\prime} +pq)$ would
have to lie on some $C_{rs}$ with $rs\le N^{\prime\prime}=pq
-p^\prime q^\prime$. Thus
$$(m+1)p+mq=\pm [(m+1)r-ms].\eqno{(1)}$$
By assumption $p^\prime=q-1+k$, $q^\prime=p+k$ for
some $k\ge 1$. Hence 
$$rs\le p^{\prime} q^\prime -pq =m(m+1) -(m+1)p -mq.\eqno{(2)}$$
Combining (1) and (2) yields $rs\pm(m+1)r-ms) \le m(m+1)$, so that
$(r\pm m)(s\mp (m+1)) \le 0$. Since $1\le r,s\le m$, it follows that
$r=m$ or $s=m+1$ and equality holds in (2). Reducing modulo $m$ or
$m+1$, we deduce that $p=m$ or $q=m+1$, neither of which is compatible
with $m=p+q+k-1$. Thus the claim holds.

Let $u(x)$ be the null vector at level $N^{\prime\prime}$ and as above
extend $u$ to $u(z)$. Since the
submodule generated by $u(z)$ is isomorphic to the Verma module
$M(c,h+pq)$, the corresponding matrix of inner products is
$\psi(z)\cdot M_{N^\prime -N}$ where
$\psi(z)=(u(z),u(z))$. But $\psi(x,0)=0$, so that $y|\psi(x,y)$. From
the Kac determinant formula, ${\rm det}\,
M_{N^\prime\prime}(x,y)=yf(x,y)$ with $f(0)\ne 0$. Hence $\psi(z)=yg(z)$
with $g(0)\ne 0$. But then ${\rm det} \,
B(z)=(u(z),u(z))^{P(N^\prime-N^{\prime\prime})} h(z)$ where
$h(0)\ne 0$. Thus $y^{P(N^\prime-N^{\prime\prime})}$ is the highest
power of $y$ dividing ${\rm det}\, B(z)$. Since $B(z)=yB_0(z)$, we
must have ${\rm det}\, B_0(0)\ne 0$, as required.

\vskip .1in
By continuity we deduce that $B_0(z)$ is invertible (possibly by
shrinking the neighbourhood of $0$). With respect to the orthogonal
decomposition corresponding to $I=Q(z)\oplus (I-Q(z))$, we may write
$A(z)=\pmatrix{B(z) & C(z) \cr C(z)^* & D(z)\cr}$. As above we have
$C(z)=yC_0(z)$, so that 
$A(z)=\pmatrix{yB_0(z) & yC_0(z) \cr y C_0(z)^* & D(z)\cr}$. 
We have already seen that $B_0(z)$ is invertible. If we try to solve
$A(z)v=0$ with $v=\pmatrix{a \cr b\cr}$, we find $a= -B_0^{-1} C_0 b$
and $(D-yC^*B^{-1} C)b=0$. Looking at the kernel at $(0,0)$, we
see that $D(0,0)$ has one--dimensional kernel. Likewise $D(0,y)$ has
one--dimensional kernel for $y\ne 0$. On the other hand $D(x,y)$ must
have zero kernel for $x\ne 0$. Thus we may define $P(z)$ near $z=0$ as
the spectral projection of $F(z)=D-yC^*B^{-1}B$ corresponding to the lowest
eigenvalue, just as in Lemma|1. By definition $F(z)$ and hence $P(z)$
is orthogonal to $Q(z)$ on the $x$--axis. By construction $P(z)$ and
$Q(z)$ have all the required properties.

\vskip .1in
\noindent \bf Remark. \rm Recall that if $P$ and $Q$ are orthogonal
projections with $\|P-Q\|<1/2$, then $T=PQ+(I-P)(I-Q)$ satisfies
check that
$\|T-I\|< 1$ and $PT=TQ$. Thus $T$ is an
invertible operator conjugating $P$ into $Q$. Since $TT^*$ commutes
with $P$, $U=(TT^*)^{-1/2}T$ gives a unitary such that $UQU^*=P$. This
means that on a sufficiently small neighbourhood of $0$
the projections $(I-Q(z))$ can be identified using a unitary gauge
change. This is not true for self--adjoint maps $A(z)$!
\vskip .1in
\noindent \bf Corollary. \it The rank one projection--valued function
$P(z)$ can be defined on a neighbourhood of the
$C_{p^\prime,q^\prime}$ containing the points $P_0,\dots, P_k$. It is
continuous and satisfies $P(z)A(z)=A(z)P(z)=\lambda(z)P(z)$ with
$\lambda(z)=0$ iff $z\in C_{p^\prime,q^\prime}$. $\lambda(z)<0$ on the
$c<1$ side of $C_{p^\prime,q^\prime}$. 
\vskip .1in
\noindent \bf Proof. \rm These first part follows taking the open to
be a (finite) union of neighbourhoods of the $P_i$'s and fintely many
other points on $C_{p^\prime,q^\prime}$. For $c$ near $1$, the matrix
$M_{N^\prime}(c,h)$ is positive--definite on the $c>1$ side of
$C_{p^\prime,q^\prime}$, invertible off $C_{p^\prime,q^\prime}$ and
has one--dimensional kernel on $C_{p^\prime,q^\prime}$. On the other
hand ${\rm det}\, M_{N^\prime}<0$ on the $c<1$ side of
$C_{p^\prime,q^\prime}$. By minimax, at most one eiegenvalue of $A(z)$
can change sign crossing $C_{p^\prime,q^\prime}$. The determinant
condition therefore implies that it is the lowest eigenvalue of $A$
that changes sign. (This evidently corresponds to the lowest
eigenvalue of $A^2$.) Thus $\lambda(z)<0$ on the $c<1$ side of
$C_{p^\prime,q^\prime}$. Since $\lambda(z)$ is real and non--zero on
the $c<1$ side of $C_{p^\prime,q^\prime}$, the last assertion
follows. 

\vskip .1in
This last corollary completes the proof of the unitarity theorem. 
\vfill\eject
\end
In between for $0<M<N-1$, the curve of highest level through $M^2/4$ is $C_{M+a,a}$ with
$a>0$ the largest integer such that $a(M+a)\le N$. Thus 
$a=\psi(M)=[F(M)]$,
with 
$$F(M)=\sqrt{N+M^2/4}-m/2]={N\over \sqrt{N+M^2/4} +M/2}.$$
Thus $\psi(M)$ decreases as $M$ increases. 

We claim that $\psi(M) -\psi(M+1)$ is either $0$ 
or $1$. Indeed 
$$F(M)-F(M+1)={1\over 2} -{M/2 + 1/4\over \sqrt{N+ M^2/4} + \sqrt{N+(M+1)^2/4}},$$
so that $0<F(M)-F(M+1) <1/2$. It follows that $0\le [F(M)]-[F(M+1)]\le 1$. 
Similarly
$$F(M)-F(M+2) = 1-{M + 1\over \sqrt{N+ M^2/4} + \sqrt{N+(M+2)^2/4}},$$
so that $0<F(M)-F(M+2) <1$. It follows that $0\le [F(M)]-[F(M+1)]\le 1$.
Thus if $0<M<N-1$, only one of $\psi(M-1)-\psi(M)$ and $\psi(M)-\psi(M+1)$ can be equal to
$1$. There are three cases:
\vskip .05in
\item{(a)} $\psi(M-1)=psi(M)=\psi(M+1)$
\item{(b)} $\psi(M-1)=\psi(M) +1 =\psi(M+1) +1$
\item{(c)} $\psi(M-1)-\psi(M) =\psi(M+1) +1$
\vskip .05in
We now compute the intersection of the curves $C^+_{M+\psi(M),M}^+$ and $C^-_{M+1+\psi(M+1),\psi(M+1)}$.
First suppose $a=\psi(M)=\psi(M+1)$. Then
$$[xM +{1\over 2}(M+2a)]^2 =[-x(M+1) +{1/2}(M+1+2a)]^2$$
so that
$$x=M+2a +{1\over 2}.$$
If on the other hand $a=\psi(M)=\psi(M+1)+1$, then
$$[xM +{1\over 2}(M+2a)]^2 =[-x(M+1) +{1/2}(M-1+2a)]^2$$
and
$$x=M+2a -{1\over 2}.$$
Similarly, replacing $M$ by $M-1$, the intersection 
of the curves $C^-_{M+\psi(M),M}^+$ and $C^+_{M-1+\psi(M+1),\psi(M+1)}$
occurs at
$$x=M+2a -{1\over 2}$$
if $\psi(M)=\psi(M-1)=a$ and at
$$x= M+ 2a -{3\over 2}$$
if $\psi(M)+1=\psi(M)=a$.
These calculations also apply in the extreme cases when $M=0$ and $M=1$.

\noindent \bf Corollary~1. \it 
If $p\ge q>1$, the first intersection of the curve $C_{p,q}^+$ with 
with a curve of level $\le N=pq$ is with $C_{p,q-1}^-$ for $x=p+q-1/2$. 
The first intersection of
$C_{p,q}^-$ with a curve of level $\le N=pq$ is with 
$C_{p-1,q}^+$ for $x=p+q-1/2$. The intersections are transverse.
\vskip .05in
\noindent \bf Proof. \rm Immediate from the lemma.   
\vskip .1in

\vskip .05in
[A point of intersection would satisfy
$((m+1)p-m)^2=((m+1)r-sm)^2$ with $rs\le p$. So either
$(m+1)(p-r) +(s-1)m=0$ or $(m+1)(p+r) =(s+1)m$. Thus the intersection
occur when $m=(\pm r-p)/(p\mp r-1-s)$. Clearly $-1<m<0$, so that
$p(m+1) -m=(p-1)(1+m) +1>1$. Thus $h=[(p(m+1) -m)^2-1]/4m(m+1) <0$, so
the points of intersection lie outside the domain.]
\noindent \bf Lemma~A. \it If $|p-q|\ne |r-s|$ and $(p,q)$ is not
proportional to $(r,s)$, the curves 
$C^+_{p,q}$ and $C^+_{r,s}$
intersect transversely in the distinct points 
$$(1-6/(x^2-{1\over 4}),[(x(p-q) +{1\over 2}(p+q))^2-1]/(4x^2 -1))$$ 
with $x={1\over 2}(r+s - p -q)/(p-q + s -r)$ and
$x={1\over 2} (r+s+p+q)/(q-p + s -r)$. The points coincide iff the pairs
$(p,q)$ and $(r,s)$ are proportional; in this case the
regions $\overline{U_{p,q}}$ and $\overline{U_{r,s}}$ intersect in $X$ only
at this common point, where they are tangent.
\vskip .05in
\noident \bf Proof. \rm  The points of intersection are given by the
solutions of
$$x(p-q)+{1\over 2}(p+q) = \pm [x(r-s)+{1\over 2}(r+s)].$$
Transversality occurs if and only if the derivatives at $x$ are equal, i.e.
$$$(p-q)[x(p-q)+{1\over 2}(p+q)] = (r-s)[x(r-s)+{1\over 2}(r-s)].$
Since
$$|x(p-q)+{1\over 2}(p+q)| = |x(r-s)+{1\over 2}(r+s)|$$
while $|p-q|\ne|r-s|$, this can only happen if both sides above
vanish and the two points of intersection coincide. In this case
$$x={1\over 2}(p+q)/(q-p)={1\over 2}(r+s)/(s-r)$$
so that $p/q=r/s$. 
\vskip .1in
\noindent \bf Lemma~B. \it The highest value of $c<1$
for which $C^+_{pq}$ intersects a curve $C^+_{rs}$ in $X$ with $rs\le pq$ and
$|p-q|\ne |r-s|$ 
corresponds to $m=p+q-1$ and 
occurs with $(r,s)=(q-1,p)$ at $h=h_{p,q}(m)$ provided $q>1$. The
intersection is transverse. 
If $q=1$ it has no intersection in $X$ with any of the curves $C_{rs}^+$ 
with $rs\le pq$.
\vskip .05in
\noindent Proof. \rm The value of $x>0$ in Lemma~A is maximised by
taking
$r=q-1$ and $s=p$ when $x=p+q -{1\over 2}$ and hence $m=p+q-1$. When
$q=1$, neither value of $x$ is positive.
\vskip .1in
\noindent \bf Lemma~C. \it If $rs\le pq$, as soon as the curve $C_{p,q}^\pm$
enters the region $\overline{U_{r,s}}$ it remains there if $|p-q|\ne
|r-s|$. If $|p-q|=|r-s|=M$ say, it intersects the region 
only at the point $(1,M^2/4)$ and $\overline{U_{r,s}}$ is wholly
contained in $\overline{U_{p,q}}$.
\vskip .05in
\noindent \bf Proof. \rm If $|p-q|=|r-s|$, then the second assertion
has already been noted when $M>0$. When $M=0$, both curves are
straight lines through $(1,0)$ and the result is immediate.
Now suppose that $|p-q|\ne |r-s|$. We first treat the case 
when $p\ne q$ and $r\ne s$. 
\vskip .1in
\noindent \bf Lemma. \it (a) $C_{p,1}$ does not intersect any curves of
level $\le p$ in $0\le c<1$ and $h\ge 0$.

\noindent (b) If $q>1$, $C_{p,q}$ intesects curves of level $\le pq$
in points $(c,h)$ in $0\le c<1$ and $h\ge 0$. The largest value of $c$
corresponds to $m=p+q-1$ and is with the curves $C_{q-1,p}$ and
$C_{p-1,q}$. The points of intersection occur at $h_{p,q}(c)$ and
$h_{q,p}(c)$. 
\vskip .1in
\noindent \bf Proof. \rm (a) A point of intersection would satisfy
$((m+1)p-m)^2=((m+1)r-sm)^2$ with $rs\le p$. So either
$(m+1)(p-r) +(s-1)m=0$ or $(m+1)(p+r) =(s+1)m$. Thus the intersection
occur when $m=(\pm r-p)/(p\mp r-1-s)$. Clearly $-1<m<0$, so that
$p(m+1) -m=(p-1)(1+m) +1>1$. Thus $h=[(p(m+1) -m)^2-1]/4m(m+1) <0$, so
the points of intersection lie outside the domain.

\noindent (b) Two curves $C_{p,q}$ and $C_{p^\prime,q^\prime}$
intersect where $((m+1)p-mq)^2=((m+1)p^\prime -mq^\prime)^2$ with
$p^\prime q^\prime \le pq$ and $p,q,p^\prime,q^\prime\ge 1$.
Thus either
$$m+1/2={ (p+p^\prime + q+ q^\prime)\over 2(q+q^\prime -p
-p^\prime)}\eqno{(1)}$$
or
$$m+{1\over 2} ={(p-p^\prime + q-q^\prime)\over 2(p^\prime + q-p
-q^\prime)}.\eqno{(2)}$$ 
Note that we may interchange $p^\prime$ and $q^\prime$ if desired. The
assertion of (b) is equivalent to the assertion that
any solution of (1) or (2) satisfies 
$$|m+{1\over 2}| \le p+q -{1/2}.\eqno{(3)}$$
Suppsoe first that we have a solution of (1). Interchanging
$p,p^\prime$ and $q,q^\prime$ if necessary, we may assume that
$q+q^\prime \ge p+p^\prime +1$. Inequality (3) is equivalent to
$$q+q^\prime \le (p+q) (q+q^\prime - p -p^\prime).\eqno{(4)}$$
If $q^\prime \le p$, (4) holds since
$q+q^\prime \le p+q\le (p+q)(q+q^\prime-p-p^\prime)$. So we may assume
that $q^\prime >p$. Since $pq\ge p^\prime q^\prime$, we must have
$p^\prime<q$. Thus $q+q^\prime-p-p^\prime\ge 2$ in this case, so that
(4) would follow if $q+q^\prime\le 2(p+q)$, i.e.~$q^\prime\le
2p+q$. If $q^\prime>2p+q$, then $q^\prime>q,2p$, so that
$p^\prime<p$. But also $p^\prime<q$. Hence
$$(p+p^\prime)(p+q) \le (p+q-1)(q^\prime + q)$$
so (4) is true.

Now suppose that (2) holds. Then 
$$|m+{1\over 2}|\le {1\over 2}
|p+q-p^\prime -q^\prime|,\eqno{(5)}$$ 
since the denominator in (2) is a non--zero
integer. If $p+q\ge p^\prime + q^\prime$, the right hand side is less
than $p+q-{1\over 2}$, so (3) holds. Thus we may assume that
$p+q\le p^\prime+q^\prime+1$. This implies that $p^\prime>p$ or
$q^\prime >q$. Without loss of generality we may assume $p^\prime >p$.
Since $p^\prime q^\prime \le pq$, we must have $q^\prime < q$. Thus
$|m+{1\over 2}|= {1\over 2} |\alpha-\beta|/|\alpha+\beta|$ where
$\alpha=p^\prime -p>0$ and $\beta=q-q^\prime>0$. Thus
$|m+{1\over 2}|\le {1\over 2} \le p+q-{1\over 2}$, as required.
\vskip .1in